\documentclass[10pt,a4paper]{article}
 \usepackage{graphics}
 \usepackage{graphicx}
 \usepackage{epsfig}
 \usepackage{amssymb}
 \usepackage{amsthm}
 \usepackage{amsmath}
 \usepackage[left=2.5cm,right=2.5cm,bottom=3.5cm,top=3.5cm]{geometry}
 \usepackage{float}
 \usepackage{lineno}
 \usepackage{textcomp}

\usepackage{epstopdf}
\usepackage{multirow}
\usepackage{color}
\usepackage{hhline}
\usepackage{booktabs}
\usepackage{lmodern}

\newtheorem{theorem}{Theorem}[section]

\newtheorem{proposition}[theorem]{Proposition}

\definecolor{nverde}{RGB}{0,61,0}
\definecolor{cr1}{RGB}{200,0,0}
\definecolor{cr2}{RGB}{0,0,200}
\definecolor{cr12}{RGB}{100,0,100}

\DeclareMathOperator{\dive}{\nabla \cdot}
\DeclareMathOperator{\curl}{\nabla \times}
\DeclareMathOperator{\gra}{\nabla}

\DeclareMathOperator{\dS}{\mathrm{dS}}
\DeclareMathOperator{\dV}{\mathrm{dV}}
\DeclareMathOperator{\dt}{\mathrm{dt}}
\newcommand{\halb}{\frac{1}{2}}

\newcommand{\Dt}{\Delta t\,}
\newcommand{\x}{\mathbf{x}}
\newcommand{\normalu}[1]{\boldsymbol{n}_{#1}}
\newcommand{\normalw}[1]{\boldsymbol{\eta}_{#1}}
\newcommand{\ci}{{\mathfrak{i}}}
\newcommand{\cell}[1]{C_{#1}}
\newcommand{\tcell}[1]{\widetilde{C}_{#1}}
\newcommand{\boundary}[1]{\Gamma_{#1}}
\newcommand{\tboundary}[1]{\widetilde{\Gamma}_{#1}}
\newcommand{\tnormal}[1]{\boldsymbol{\widetilde{n}}_{#1}}
\newcommand{\tnormalw}[1]{\boldsymbol{\widetilde{\eta}}_{#1}}
\newcommand{\tnormale}[1]{\widetilde{n}_{#1}}

\newcommand{\cj}{{\mathfrak{j}}}

\newcommand{\Pcell}[1]{T_{#1}}

\newcommand{\Bst}{\widetilde{\boldsymbol{\mathcal{B}}}}

\newcommand{\Xst}{\widetilde{\mathbf{X}}}
\newcommand{\xst}{\widetilde{\mathbf{x}}}
\newcommand{\Phist}{\widetilde{\Phi}}

\newcommand{\nablast}{\widetilde{\nabla}}
\newcommand{\nst}{\widetilde{\boldsymbol{n}}}
\newcommand{\etast}{\widetilde{\boldsymbol{\eta}}}

\newcommand{\XX}{\mathbf{X}}
\newcommand{\bV}{\mathbf{V}}
\newcommand{\dx}{\Delta \x}
\newcommand{\vel}{u}
\newcommand{\mom}{\rho\vel}
\newcommand{\bvel}{\mathbf{\vel}}
\newcommand{\bmom}{\rho\bvel}
\newcommand{\press}{p}
\newcommand{\A}{A}
\newcommand{\bA}{\mathbf{\A}}
\newcommand{\phiA}{\varphi}

\newcommand{\bpA}{\boldsymbol{\phiA}_{\bA}}
\newcommand{\bpAi}[1]{\boldsymbol{\phiA}_{\bA_{#1}}}
\newcommand{\J}{J}
\newcommand{\bJ}{\mathbf{\J}}
\newcommand{\psiJ}{\psi}

\newcommand{\bpJ}{\boldsymbol{\psiJ}_{\bJ}}
\newcommand{\tE}{\mathcal{E}}

\newcommand{\h}{h}
\newcommand{\hq}{q}

\newcommand{\devG}{\mathring{G}}
\newcommand{\devD}{\mathring{D}}
\newcommand{\g}{\mathbf{g}}
\newcommand{\q}{\mathbf{q}}
\newcommand{\Q}{\mathbf{Q}}
\newcommand{\eQ}{\overline{\Q}}

\newcommand{\p}{\mathbf{p}}
\newcommand{\Flux}{\boldsymbol{\mathcal{F}}}
\newcommand{\Fluxst}{\boldsymbol{\widetilde{\mathcal{F}}}}
\newcommand{\source}{\boldsymbol{\mathcal{S}}}
\newcommand{\BNCP}[1]{\boldsymbol{\mathcal{B}}(#1)}
\newcommand{\NCP}[1]{\BNCP{#1}\cdot \gra #1}
\newcommand{\BNCPst}[1]{\boldsymbol{\widetilde{\mathcal{B}}}(#1)}
\newcommand{\grast}{\widetilde{\gra}}
\newcommand{\NCPst}[1]{\BNCPst{#1}\cdot \grast #1}
\newcommand{\E}{E}

\newcommand{\pxk}{\frac{\partial}{\partial {x_{k}}}}
\newcommand{\pxj}{\frac{\partial}{\partial {x_{j}}}}
\newcommand{\pxi}{\frac{\partial}{\partial {x_{i}}}}

\newcommand{\pt}{\frac{\partial}{\partial t}}

\allowdisplaybreaks

\begin{document}

\thispagestyle{plain}
\begin{center}
	\textbf{ \Large{An arbitrary Lagrangian-Eulerian semi-implicit hybrid method\\[8pt] for continuum mechanics with GLM cleaning} }
	
	\vspace{0.5cm}
	{Saray Busto$^{(a,b),}$\footnote{Corresponding author. \hspace*{0.25cm} Email address: saray.busto.ulloa@usc.es}}
	
	\vspace{0.2cm}
	{\small
		\textit{$^{(a)}$ Department of Applied Mathematics, Universidade de Santiago de Compostela, 15782 Santiago de Compostela, Spain}
		
		\textit{$^{(b)}$ Galician Center for Mathematical Research and Technology, CITMAga, 15782 Santiago de Compostela, Spain}
	}
\end{center}

\begin{abstract}
This paper proposes a semi-implicit arbitrary Lagrangian-Eulerian (ALE) method for the solution of the unified Godunov-Peshkov-Romenski (GPR) model of continuum mechanics. To handle the curl free involutions arising in the solid limit of the model, the original system is augmented by adopting a thermodynamically compatible generalized Lagrangian multiplier (GLM) approach. Next, an operator splitting strategy decouples the computation of fast pressure waves from the bulk velocity of the medium yielding a transport subsystem, containing convective terms and non-conservative products, and a Poisson-type subsystem, for the pressure. A second splitting yields an ODE subsystem comprising only the potentially stiff source terms, responsible for the relaxation of the model between its fluid and solid limits.

The mesh motion can be driven by two sources: the local fluid velocity and a prescribed boundary displacement. For the spatial discretization, we employ unstructured staggered grids, with the pressure defined on the primal mesh and all remaining variables on the dual grid. The transport subsystem is advanced via an explicit finite volume method, in which integration over closed space-time control volumes ensures verification of the geometric conservation law (GCL). On the other hand, implicit continuous finite elements are used for the discretization of the pressure subsystem and an implicit DIRK scheme is employed to solve the ODE subsystem. Consequently, the proposed approach is well suited to address all Mach number flows. A comprehensive set of benchmarks is employed to assess the accuracy and robustness of the proposed methodology in tackling both fluid and solid mechanics problems.
\end{abstract}

\vspace{0.2cm}

\vspace{0.4cm}

% % % % % % % % % % % % % % % % % % % % % % % % % % % % % %
% % % % % % % % % % % % % % % % % % % % % % % % % % % % % %
%              Introduction
% % % % % % % % % % % % % % % % % % % % % % % % % % % % % %
% % % % % % % % % % % % % % % % % % % % % % % % % % % % % %
\section{Introduction} \label{sec:intro}
The development of a unified first order model for continuum mechanics has started more than sixty years ago with the pioneering work of S. Godunov and E. Romenski on the derivation of a hyperbolic model for solid mechanics \cite{God1961,GR72_GPRnnelasticity,Rom1998,GodunovRomenski72,GR03}. %that was next extended 
Contrary to classical Lagrangian descriptions of solids, this hyperbolic model was derived in the Eulerian frame of reference. This approach naturally lent itself to later extensions
including heat conduction \cite{MalyshevRom1986} and fluid dynamics \cite{PeshRom2016}. The resulting system, known as the Godunov-Peshkov-Romenski (GPR) model for continuum mechanics, relies on relaxation parameters to change across a wide range of media, from elastoplasticity to compressible fluids and porus media \cite{BDR10_GPRElastoplastic,DPR18_GPRelastovisco,JN19_GPRNNplastic,RRP21porous}. An important feature of the GPR model is its compliance with the laws of thermodynamics, placing it in the class of thermodynamically compatible models. Further, it enables a connection with Hamiltonian formulations \cite{PPRG18_GPRHamiltonian}. 
In this paper, we focus on the numerical solution of the compressible GPR model for fluids and solids given in the Godunov form. Nevertheless, it is important to notice that this modelling framework is well suited to include electromagnetic effects \cite{Rom1989Maxwell,DPRZ17_GPRelectro} and general relativity \cite{RPDF20_GPRGeneralRelativity} while recent advances regard  modelling of multiphase phenomena \cite{Rom2007TwoPhase,Rom2010TwoPhase,RRP21_GPRthreephase,RP23_GPR2Phase,RDD25_GPR2phase,FPRD25_GPRMultiphase} and reformulations of non-dispersive systems \cite{DD22_GPRNSKorterwerg}.

The hyperbolic nature of the GPR system has motivated the extension of classical methods for conservation laws to its solution. Among the different approaches proposed, we find, e.g., SPH methods \cite{KPPK23_GPR_SPH}, finite volume (FV) and discontinuous Galerkin (DG) thermodynamically compatible schemes \cite{HTCGPR,HTCA2}, and ADER methods \cite{BDL16_GPR_ALE,DPRZ16_GPRmodel,BCDGP20}. 
Most available methodologies consider a fully explicit discretisation of the model so the sound velocity appears in the system eigenvalues. This may yield to a restrictive time step condition even for slow bulk velocity phenomena.
In this case, to enhance the performance of numerical approximations, a well-established approach coming from fluid dynamics is to apply an operator splitting strategy allowing for the decoupling of the pressure waves computation  \cite{PM05,DeT11,TV12,TD17,Hybrid2,TPK20}. Then, the use of semi-implicit schemes frees the method from the sharp time step restriction while retaining an explicit scheme for the treatment of convective terms and non-conservative products. 
Moreover, this methodology eases the verification of the the asymptotic preserving (AP) property in the low Mach number limit. Further details on the development of the so-called all Mach number flow solvers can be found in \cite{BBD25_FVVEM} and references therein. Recently, this methodology has been successfully extended to the GPR model framework, using either Cartesian grids \cite{Boscheri2021SIGPR,PDBRCI2021GPRNN} or unstructured staggered meshes \cite{HybridGPR}. 

The semi-implicit method in \cite{HybridGPR} is based on the combination of explicit finite volume and implicit continuous finite element methods for the solution of the incompressible and a weakly compressible GPR models. It belongs to a family of hybrid FV/FE approaches formerly employed for the discretization of diverse hyperbolic systems such as the Navier-Stokes equations \cite{Hybrid2,LBD23_HybridImplicit}, the shallow water equations \cite{BD22} and the magnetohydrodynamics equations \cite{HybridHexa1,HybridMHD}. In this paper, we will also rely on this methodology proposing its extension to the compressible GPR model in the framework of arbitrary Lagrangian-Eulerian (ALE) methods \cite{HybridALE}. Analogously to the asymptotic preserving methods for the compressible Navier-Stokes equations, considering the low Mach limit of the proposed approach yields the incompressible formulation of the GPR model in \cite{HybridGPR}.

Let us note that, up to now, very few approaches have been developed for the GPR model in the Lagrangian \cite{BCP22_GPR_IMEXLagUnst,BDLM25_GPR_LagrangianSP} or ALE frameworks \cite{BDL16_GPR_ALE}. Nevertheless, Lagrangian schemes are the standard approach in solid mechanics \cite{Neumann1950,ShashkovCellCentered} and they excel at tracking moving interfaces and contact discontinuities \cite{Maire2007,Despres2009,ATS24_Lagrangian}. The more recent ALE methodologies can be seen as a flexible extension of Lagrangian methods where the mesh velocity can be defined independently or as a function of the local velocity of the medium reducing mesh distortion \cite{BBGKM11_ALE,Lagrange2D,Gaburro2021,BC25_ALEboundary}.
These properties make Lagrangian and ALE methods widely employed tools in the context of fluid-structure interaction \cite{LeTallecMouro,NobileVergara,Feistauer4} as well. Therefore, the development of efficient ALE schemes for the GPR model would represent an step toward its use in a wide range of practical applications.

When discretizing the GPR model, we should take into account the curl free involution constraints that arise, for the distortion and thermal impulse fields, in the solid limit of the model. To tackle these propert at the discrete level, one possibility is the development of exactly involution preserving schemes, e.g. using an adequate mesh staggering of the discrete operators and conservative variables, as the one proposed in \cite{Boscheri2021SIGPR,PDBRCI2021GPRNN} for Cartesian staggered grids. A simpler alternative to exact involution preserving methods are divergence and curl correction techniques as the generalized Lagrangian multipliers (GLM) methods devised to address the divergence free condition of the magnetic field in magnetohydrodynamics \cite{MOSSV00,HybridHexa1}. In this paper, we adopt this latter methodology and propose an augmented GLM GPR model grounded on the seminal ideas in \cite{HyperbolicDispersion,HTCMHD} concerning thermodynamically compatible GLM formulations.

Another crucial aspect of the GPR model regards the stiffness of the source terms involved in the distortion and thermal impulse equations which may lead to stability issues in the visco-plastic limit of the equations. To overcome this problem, a robust strategy, that will be followed in this paper, consists on splitting the system decoupling the source terms computation. The corresponding subsystem,  made of ordinary differential equations (ODE), can be then efficiently solved employing unconditionally stable implicit methods leading to implicit-explicit (IMEX) approaches compatible with the asymptotic limits of the model \cite{Jack17,BCP22_GPR_IMEXLagUnst}. 

The rest of this paper is organised as follows. In Section~\ref{sec:goveq}, we recall the Godunov-Peshkov-Romenski model for continuum mechanics and introduce a thermodynamically compatible Generalized Lagrangian Multiplier extension allowing for the cleaning of curl errors in the distortion and thermal impulse fields. Section~\ref{sec:numdisc} is devoted to the numerical method. We first introduce the operator splitting approach which divides the system into a transport subsystem for the conservative variables and a Poisson type system for the pressure. The unstructured staggered grid spatial discretization is introduced and the finite volume method for the solution of the transport equations is described. Next, we detail the variational formulation for the discretization of the pressure system using continuous finite element methods and the interpolation strategies to pass data between the staggered grids. The proposed methodology is assessed in Section~\ref{sec:numericalresults} by a wide set of test cases including both fluid dynamics and solid mechanics benchmarks. Finally, the conclusions and an outlook on future research are drawn in Section~\ref{sec:conclusions}.

% % % % % % % % % % % % % % % % % % % % % % % % % % % % % %
% % % % % % % % % % % % % % % % % % % % % % % % % % % % % %
%              Mathematical model
% % % % % % % % % % % % % % % % % % % % % % % % % % % % % %
% % % % % % % % % % % % % % % % % % % % % % % % % % % % % %
\section{Governing equations} \label{sec:goveq}
To address both fluids and solids within a unified partial differential system, we consider the Godunov-Peshkov-Romenski (GPR) model for continuum mechanics \cite{PeshRom2016,DPRZ16_GPRmodel}, that, using Einstein notation, reads
\begin{subequations}\label{eqn.GPR}
	\begin{eqnarray}
		\pt \rho + \pxk \left(\rho\vel_{k}\right) = 0,\label{eqn.GPR_rho}\\
		\pt \left(\rho\vel_i\right) + \pxk \left(\rho \vel_{i} \vel_{k} \right) + \pxi \press +\pxk\sigma_{ik} +\pxk\omega_{ik} = \rho g_{i},\label{eqn.GPR_mom}\\
		\pt \A_{ik} + \pxk \left(\vel_{m}\A_{im}\right) + \vel_{j} \pxj \A_{ik} -\vel_{j}\pxk \A_{ij} = -\frac{1}{\theta_{1} \left(\tau_{1}\right)} \E_{\A_{ik}}, \label{eqn.GPR_A}\\
		\pt \J_{k} + \pxk \left(\J_{m}\vel_{m} \right) +\pxk T + u_{j}\left( \pxj \J_{k} - \pxk \J_{j} \right)  =  -\frac{1}{\theta_{2} \left(\tau_{2}\right)} \E_{\J_{k}}, \label{eqn.GPR_J}\\
		\pt \left( \rho S\right) + \pxk \left(\rho S\vel_{k}\right) + \pxk \left( \rho \E_{\J_{k}}\right)  = \frac{\rho}{T}\left(\frac{1}{\theta_{1}\left(\tau_{1}\right)}\E_{\A_{ik}}\E_{\A_{ik}} +  \frac{1}{\theta_{2}\left(\tau_{2}\right)}\E_{\J_{k}}\E_{\J_{k}}\right) \geq 0,\label{eqn.GPR_S}\\
		\pt \left( \rho \E\right)  +\pxk \left(\rho\E\vel_{k}\right) + \pxk \left(\press\vel_{k}\right) + \pxk \left(\vel_{i}\sigma_{ik} \right) + \pxk \left(\vel_{i}\omega_{ik} \right) +\pxk \hq_{k} = \rho g_{i} \vel_{i}, \label{eqn.GPR_E}
	\end{eqnarray}
\end{subequations}
being $\rho$ the density, $\bvel=\left(\vel_1,\vel_2,\vel_3\right)^{T}$ the velocity field, $\press$ the pressure, $\bA$ the distortion field, $\bJ$ the thermal impulse, $S$ the entropy and $\E$ the total energy.
The total energy density $\tE:=\rho\E$ can be divided into four terms as
\begin{equation}
	\tE\left(\rho,\bvel,\bA,\bJ,S\right) = \tE_{1}\left(\rho,S\right) + \tE_{2}\left(\bvel\right)
	+\tE_{3}\left(\bA\right)  
	+\tE_{4}\left(\bJ\right)  
	.\label{eqn.energydecomp}
\end{equation}
The first one corresponds to the internal energy, related to the kinetic energy of the molecular motion. It enables the consideration of different materials depending on the selected equation of state. In case a gas is assumed, the ideal gas equation of state can be employed, then
\begin{equation}
	\tE_{1}\left(\rho,S\right) = \frac{\rho^{\gamma}}{\left( \gamma-1\right) } e^{\frac{S}{c_{v}}},\qquad  \tE_{1}\left(\rho,\press\right)= \frac{\press}{\gamma-1}
	\label{eqn.E4}
\end{equation}
with $\gamma=\dfrac{c_{\press}}{c_{v}}$ the ratio of specific heat at constant pressure, $c_{\press}$, and at constant volume, $c_{v}$.
Further, for solids and liquids, we may use the stiffened gas EOS 
\begin{gather}
	\tE_{1}\left(\rho,S\right) = \frac{c_{0}^2}{\gamma \left(\gamma-1\right)}\left(\frac{\rho}{\rho_{0}}\right)^{\gamma-1} e^{\frac{S}{c_{v}}} + \frac{\rho_{0}c_{0}^2-\gamma \press_{0}}{\gamma \rho}, \qquad
	\tE_{1}\left(\rho,\press\right)= \frac{\press}{\gamma-1} + \frac{\rho_{0} c_{0}^2 -  \gamma \press_{0}}{\gamma\left( \gamma-1\right) },
	\label{eqn.E4_stiff}
\end{gather}
where $\rho_{0}$ is the reference material density, $\press_{0}$ denotes the reference atmospheric pressure, and $c_{0}$ refers to the adiabatic soundspeed. More complex relations could be incorporated into the model as, e.g., the Mie-Grüneisen equation of state for the thermodynamical study of solid materials \cite{BCDGP20} or the EOS in \cite{FGN14_HTChypelast,GNH16_EOS} that would provide a hyperbolic model also for large deformations.
The second term in \eqref{eqn.energydecomp} corresponds to the kinetic energy per unit volume,
\begin{equation}
	\tE_{2}\left(\bvel\right) = \halb \rho\left|\bvel\right|^{2}.\label{eqn.E1}
\end{equation}
Finally, the third and fourth terms provide the contribution of the mesoscopic, non-equilibrium, part of $\tE$ related to the material deformations and the thermal impulse,
\begin{equation}
	\tE_{3}\left(\bA\right) = \frac{1}{4} c_{s}^2 \rho	\devG_{ij}	\devG_{ij}, \qquad \tE_{4}\left(\bJ\right) = \frac{1}{2} c_{h}^2 \rho	\J_{i}	\J_{i}, \label{eqn.E23}
\end{equation}
where $c_{s}^{2}$ and $c_{h}^{2}$ are the characteristic velocities for propagation of shear and thermal perturbations while $\devG_{ik}$ denotes the trace-free part of the metric tensor $G_{ik}=\A_{ji}\A_{jk}$,
\begin{equation}
	\mathring{G}_{ik} = G_{ik}- \frac{1}{3}G_{mm}\delta_{ik}.
\end{equation}
In the momentum equations \eqref{eqn.GPR_mom} the viscous and thermal diffusivity effects are taken into account thanks to the non-isotropic part of the stress tensor, containing both shear and thermal stresses,
\begin{equation}
	\sigma_{ik} = \A_{ji}\partial_{A_{jk}} \tE = \rho c_{s}^2 G_{ij} \devG_{jk}, \qquad
	\omega_{ik} =  \J_{i}\partial_{J_k} \tE =  \rho c_{h}^{2} \J_i\J_{k}
\end{equation}
Further, 
\begin{equation}
\hq_{k} = \partial_{\rho S} \tE \partial_{J_k} \tE = \rho c_{h}^2 T \J_k
\end{equation}
denotes the heat flux, with $T:=\partial_{\rho S}\tE$ the temperature.
Finally, the shear and thermal stress relaxation functions read
\begin{equation}
	\theta_{1}\left(\tau_{1}\right) = \frac{1}{3} \rho_{0}\tau_{1} c_{s}^2 \left|\bA\right|^{-\frac{5}{3}},\qquad
%\end{equation}
%\begin{equation}
	\theta_{2}\left(\tau_{2}\right) = \frac{\rho_{0} T_{0}}{T} \tau_{2} c_{h}^{2}
\end{equation}
with $\tau_{1}$ and $\tau_{2}$ the corresponding relaxation times and $\rho_0$, $T_0$ reference values for the density and the temperature. Let us recall that in the stiff limit, i.e., for $\tau_{1},\, \tau_{2} \rightarrow 0$, we retrieve the 
Navier-Stokes-Fourier equations for fluids with shear viscosity $\mu=\frac{1}{6} \rho_{0} c_{s}^{2} \tau_{1}$ and heat conductivity $\kappa = \rho_{0} T_{0} c_{h}^{2} \tau_{2}$, see \cite{DPRZ16_GPRmodel} for a formal asymptotic analysis of this property.

It is very important to remark that system \eqref{eqn.GPR} is an overdetermined system of PDEs where the energy equation \eqref{eqn.GPR_E} can be obtained as the dot product of all other equations by the corresponding main field or thermodynamical dual variables,
$\p = \left( r, \vel_{i}, \alpha_{ik}, \beta_{k}, T \right)^{T},$
$r=\partial_{\rho} \tE$, $\vel_{i}=\partial_{\rho \vel_{i}} \tE$, $\alpha_{ik}=\partial_{\A_{ik}} \tE$, $\beta_{k}=\partial_{\J_{k}} \tE$, $T=\partial_{\rho S} \tE$.
Hence, this system belongs to the class of thermodynamically compatible schemes in the sense of Godunov, \cite{God1961,GR03}.

\subsection{Augmented GLM curl cleaning GPR model}

A key aspect of the GPR model is that the distortion field, $\bA$, and the thermal impulse field, $\bJ$, are endowed with natural involution constraints of the curl-type. That is, given initially curl free fields $\bA(\x,0)$ and $\bJ(\x,0)$ then  $\curl \bA_i(\x,t)=0$ and $\curl \bJ(\x,t)=0$ for all $t\in\mathbb{R}^{+}$. To deal with this structural property several approaches have been proposed in the last decades including the Godunov-Powell approach, \cite{God72,Powell94,PRLGD99,CPGD20}, and exact involution preserving methods \cite{HS97,JT06,DCP20_curlinv,HTCGPR,BKBD23}. An alternative successful strategy to address these constraints in an approximated manner consists on the use of the so-called Generalized Lagrangian Multipliers (GLM). Augmented GLM models follow the seminal ideas put forward in \cite{MOSSV00,DKKMSW02_GLM} for the treatment of the divergence free condition of magnetohydrodynamics (MHD) equations. Then, they have been extended to further systems of conservation laws, as in \cite{DFGR20_GLMCCZ4} and \cite{HyperbolicDispersion} where GLM curl cleaning approaches have been presented for the FO-CCZ4 formulation of the Einstein field equations and the turbulent shallow water (TSW) model, respectively. Let us remark that these original GLM cleaning approaches often neglected the important thermodynamically compatible property of the underlying systems. Recently, two families of compatible GLM cleaning models have been proposed in \cite{DWGWB18_GLMentropyconsistent,HTCMHD} for the MHD equations. In this work, we follow the approach in \cite{HyperbolicDispersion,HTCA2} and introduce a compatible GLM curl cleaning for the GPR model. 

It is important to recall that the distortion field, even if written as a matrix, does not correspond to a classical tensor field but to a non-holonomic basis triad, \cite{PRD19_torsion}. Consequently, the curl free property is verified for each matrix column $\bA_i, \, i\in\left\lbrace 1,2,3 \right\rbrace,$ and we can design a GLM cleaning approach for each of these vector fields instead of considering a complete tensor. Accordingly, we introduce the cleaning fields $\bpA=\left( \bpAi{1} \mid \bpAi{2} \mid \bpAi{3}\right)$ and $\bpJ$ with associated cleaning speeds $c_{\bA}$, for the distortion field components, and $c_{\bJ}$, for the heat flux. The equations for the distortion field, the thermal impulse and the corresponding cleaning variables read
\begin{eqnarray}
	\pt \A_{ik} + \pxk \left(\vel_{m}\A_{im}\right) + \vel_{j} \pxj \A_{ik} -\vel_{j}\pxk \A_{ij} + c_{\bA} \epsilon_{ijl} \pxj \phiA_{lk} = -\frac{1}{\theta_{1} \left(\tau_{1}\right)} \E_{\A_{ik}},\label{eqn.GPRGLM_A} \\
	\pt \phiA_{ik} + \vel_{j} \pxj \phiA_{ik} - c_{\bA}\epsilon_{ijl} \pxj \A_{lk}= 0,\label{eqn.GPRGLM_PhiA}\\
	\pt \J_{k} + \pxk \left(\J_{m}\vel_{m} \right) +\pxk T + u_{j}\left( \pxj \J_{k} - \pxk \J_{j} \right)   + c_{\bJ} \epsilon_{kjl} \pxj \psiJ_{l} =  -\frac{1}{\theta_{2} \left(\tau_{2}\right)} \E_{\J_{k}}, \label{eqn.GPRGLM_J}\\
	\pt \psiJ_{k} + \vel_{j} \pxj \psiJ_{k} - c_{\bJ}\epsilon_{kjl} \pxj \J_{l}= 0 \label{eqn.GPRGLM_PsiJ}
\end{eqnarray}
with $\epsilon_{ijl}$ the Levi-Civita symbol. We observe that as $c_{\bA} \rightarrow\infty$ and $c_{\bJ} \rightarrow\infty$, we recover the sought curl free relations, $\epsilon_{ijl} \pxj \A_{lk}= 0$, $\epsilon_{kjl} \pxj \J_{l}= 0$. To ensure thermodynamical compatibility, we need to include the dependence of the total energy on the cleaning primal variables as
\begin{equation}
	\tE = \rho E = 
	\tE_{1} + \tE_{2} + \tE_{3} + \tE_{4} + \tE_{5} + \tE_{6},
	\qquad 
	\tE_{5} = \frac{1}{4} \rho c_{s}^{2} \devD_{ik}\devD_{ik},\qquad
	\tE_{6} = \frac{1}{2} \rho c_{h}^{2} \psiJ_{k}\psiJ_{k},
	\label{eqn.energyGLM}
\end{equation}
with $\devD$ the deviatoric of tensor $D_{ik}=\phiA_{ji}\phiA_{jk}$.
Then, the new vector of thermodynamical dual variables, corresponding to the primal variables
\begin{equation}
	\q_{GLM} = \left( \rho, \rho\vel_{i}, \A_{ik}, \phiA_{ik}, \J_{k}, \psiJ_{k}, S \right)^{T},
\end{equation} 
results
\begin{equation*}
	\p_{GLM} = \left( r, \vel_{i}, \alpha_{ik}, \gamma_{ik}, \beta_{k}, \xi_{k}, T \right)^{T}
\end{equation*}
with
\begin{equation}
	\gamma_{ik} = \partial_{\phiA_{ik}}\tE = \rho c_{s}^{2} \phiA_{ij} \devD_{jk},
	\qquad \xi_{i} = \partial_{\psiJ_{i}}\tE = \rho c_{h}^{2} \psiJ_{i}. \label{tederivativescleaning}
\end{equation}

\begin{proposition}
	The augmented GLM curl cleaning GPR model, \eqref{eqn.GPR_rho}-\eqref{eqn.GPR_mom}, \eqref{eqn.GPRGLM_A}-\eqref{eqn.GPRGLM_PsiJ}, \eqref{eqn.GPR_S}, is thermodynamically compatible with the extra energy conservation law
	\begin{gather}
		\pt \tE  +\pxk \left(\tE\vel_{k}\right) + \pxk \left(\press\vel_{k}\right) + \pxk \left(\vel_{i}\sigma_{ik} \right) + \pxk \left(\vel_{i}\omega_{ik} \right) +\pxk \hq_{k} 
		\notag\\
		+ \rho c_{\bA} c_{s}^2 \epsilon_{ijl} \left( \A_{im} \devG_{mk} \pxj \phiA_{lk} + \phiA_{lm}\devD_{mk}\pxj \A_{ik} \right) + c_{\bJ} c_{h}^{2}  \epsilon_{kjl} \left( \rho\psiJ_{l} \pxj \J_{k} + \rho\J_{k} \pxj \psiJ_{l} \right) 
		= \rho g_{i} \vel_{i}. \label{eqn.GPRGLM_E}
	\end{gather}	
\end{proposition}
\begin{proof}
	Let us first focus in the time derivative and the convective terms but for $\pxk \left(\vel_{m}\A_{im}\right) -\vel_{j}\pxk \A_{ij}$ and $\pxk \left(\J_{m}\vel_{m} \right) - u_{j} \pxk \J_{j}$. Then, taking into account the definition of the new energy \eqref{eqn.energyGLM}, we have
	\begin{eqnarray}
		&&r\! \left[ \pt \rho \!+ \!\pxk \left(\rho\vel_{k}\right) \right]  \!
		+\! \vel_i\! \left[ \pt \left(\rho\vel_i\right)\! +\! \pxk \left(\rho \vel_{i} \vel_{k} \right) \right] \!
		+\! \alpha_{ik\!} \left[ \pt \A_{ik}\! +\! \vel_j \pxj \A_{ik} \right] \!
		+\! \gamma_{ik}\! \left[ \pt \phiA_{ik}\! +\! \vel_{j} \pxj \phiA_{ik} \right] \notag\\ &&
		+ \beta_{k} \left[ \pt \J_{k} +\vel_{j} \pxj \J_{k}  \right]
		+ \xi_{i} \left[\pt \psiJ_{i} + \vel_{k} \pxk \psiJ_{i} \right]  
		+ T \left[ \pt \left( \rho S\right) + \pxk \left(\rho S\vel_{k}\right)  \right]  \notag\\
		&=&  \pt \tE 		 
		+ \partial_{\rho} \tE_{1} \pxk \left(\rho\vel_{k}\right) 
		+ \partial_{\rho S}\tE_{1} \pxk \left( \rho S\right) 
		+ \partial_{\rho} \tE_{2} \pxk \left(\rho\vel_{k}\right) 
		+ \partial_{\rho\vel_{i}}\tE_{2}  \pxk \left(\rho \vel_{i} \vel_{k} \right) \notag\\
		& & 		
		+ \partial_{\rho} \tE_{3} \pxk \left(\rho\vel_{k}\right)
		+ \partial_{\A_{ij}}\tE_{3} \vel_k \pxk \A_{ij} 
		+ \partial_{\rho} \tE_{4} \pxk \left(\rho\vel_{k}\right)  
		+ \partial_{\J_{i}}\tE_{4} \vel_{k} \pxk \J_{i}   
		+ \partial_{\rho} \tE_{5} \pxk \left(\rho\vel_{k}\right)  \notag\\
		& & 		
		+ \partial_{\phiA_{ij}}\tE_{5} \vel_{k} \pxk \phiA_{ij} 
		+ \partial_{\rho} \tE_{6} \pxk \left(\rho\vel_{k}\right) 
		+ \partial_{\psiJ_{i}}\tE_{6} \vel_{k} \pxk \psiJ_{i}  \notag\\
		%
%		&=&  \pt \tE 		 
%		+ \pxk \left( \tE_{1} \vel_{k} \right) 
%		+ \pxk \left( \tE_{2} \vel_{k}\right) 
%		+ \pxk \left( \tE_{3} \vel_{k}\right)  
%		+ \pxk \left( \tE_{4} \vel_{k}\right)     \notag\\
%		& & 		
%		+ \partial_{\rho} \tE_{5} \pxk \left(\rho\vel_{k}\right) 
%		+ \partial_{\phiA_{ij}}\tE_{5} \vel_{k} \pxk \phiA_{ij} 
%		+ \partial_{\rho} \tE_{6} \pxk \left(\rho\vel_{k}\right) 
%		+ \partial_{\psiJ_{i}}\tE_{6} \vel_{k} \pxk \psiJ_{i}  \notag\\
		%
		&=&  \pt \tE 		 
		+ \pxk \left[ \left( \tE_{1} + \tE_{2} + \tE_{3} + \tE_{4}\right)  \vel_{k}\right]    
		+ \partial_{\rho} \tE_{5} \pxk \left(\rho\vel_{k}\right) 
		+ \partial_{\phiA_{ij}}\tE_{5} \vel_{k} \pxk \phiA_{ij}   		
		+ \partial_{\rho} \tE_{6} \pxk \left(\rho\vel_{k}\right)  \notag\\
		& &
		+ \partial_{\psiJ_{i}}\tE_{6} \vel_{k} \pxk \psiJ_{i}  \notag\\
		%
%		&=&  \pt \tE 		 
%		+ \pxk \left[ \left( \tE_{1} + \tE_{2} + \tE_{3} + \tE_{4}\right)  \vel_{k}\right]   
%		+ \partial_{\rho} \left( \frac{1}{4} \rho c_{s}^{2} \devD_{il}\devD_{il}\right)  \pxk \left(\rho\vel_{k}\right)  	
%		+ \partial_{\phiA_{ij}} \left( \frac{1}{4} \rho c_{s}^{2} \devD_{il}\devD_{il}\right) \vel_{k} \pxk \phiA_{ij}      \notag\\
%		& & 	
%		+ \partial_{\rho} \left( \frac{1}{2} \rho c_{h}^{2} \psiJ_{i}\psiJ_{i}\right)  \pxk \left(\rho\vel_{k}\right) 
%		+ \partial_{\psiJ_{i}}\left( \frac{1}{2} \rho c_{h}^{2} \psiJ_{l}\psiJ_{l}\right)  \vel_{k} \pxk \psiJ_{i}  \notag\\
		%
		&=&  \pt \tE 		 
		+ \pxk \left[ \left( \tE_{1} + \tE_{2} + \tE_{3} + \tE_{4}\right)  \vel_{k}\right]     
		+ \frac{1}{4} c_{s}^{2} \devD_{il}\devD_{il}  \pxk \left(\rho\vel_{k}\right)
		+    \rho c_{s}^{2} \phiA_{il} \devD_{lj} \vel_{k} \pxk \phiA_{ij}		 \notag\\
		& & 		 
		+ \frac{1}{2}  c_{h}^{2} \psiJ_{i}\psiJ_{i} \pxk \left(\rho\vel_{k}\right)   
		+  \rho c_{h}^{2} \psiJ_{i}  \vel_{k} \pxk \psiJ_{i}  \notag\\
		&=&  \pt \tE 		 
		+ \pxk \left[ \left( \tE_{1} + \tE_{2} + \tE_{3} + \tE_{4}\right)  \vel_{k}\right]   
		+ \pxk \left( \frac{1}{4} \rho c_{s}^{2} \devD_{il}\devD_{il} \vel_{k}\right)     
		+ \pxk \left( \frac{1}{2} \rho c_{h}^{2} \psiJ_{i}\psiJ_{i} \vel_{k}\right)  \notag\\
		&=&  \pt \tE 		 
		+ \pxk \left[ \left( \tE_{1} + \tE_{2} + \tE_{3} + \tE_{4}\right)  \vel_{k}\right]    
		+ \pxk \left( \tE_{5} \vel_{k}\right)     
		+ \pxk \left( \tE_{6}  \vel_{k}\right)  \notag\\
		&=& \pt \tE + \pxk \left(  \tE\vel_{k} \right).
		\label{HTCdtflux}
	\end{eqnarray}	
	Next, the product of $(\alpha_{ik},\gamma_{ik},\beta_{k},\xi_{k})$ by the curl terms in \eqref{eqn.GPRGLM_A}-\eqref{eqn.GPRGLM_PsiJ}, substitution of $\alpha_{ik}= \rho c_{s}^{2} \A_{ij}\devG_{jk} $, $\beta_{k} = \rho c_{h}^{2}\J_{k}$  and \eqref{tederivativescleaning}, and reordering terms gives
	\begin{eqnarray}
	&&\alpha_{ik}  c_{\bA} \epsilon_{ijl} \pxj \phiA_{lk}
	- \gamma_{ik} c_{\bA}\epsilon_{ijl} \pxj \A_{lk}
	+ \beta_{k} c_{\bJ} \epsilon_{kjl} \pxj \psiJ_{l}
	- \xi_{k} c_{\bJ}\epsilon_{kjl} \pxj \J_{l} \notag \\
%	%
	&=& \rho c_{s}^2 \A_{im}\devG_{mk} c_{\bA} \epsilon_{ijl} \pxj \phiA_{lk}
	- \rho c_{s}^2 \phiA_{im}\devD_{mk} c_{\bA}\epsilon_{ijl} \pxj \A_{lk}
	+ \rho c_{h}^{2} \J_{k} c_{\bJ} \epsilon_{kjl} \pxj \psiJ_{l}
	- \rho c_{h}^{2} \psiJ_{k} c_{\bJ}\epsilon_{kjl} \pxj \J_{l} \notag \\
	&=& \rho c_{s}^2 c_{\bA} \left( \A_{im}\devG_{mk} \epsilon_{ijl} \pxj \phiA_{lk}
	- \phiA_{lm}\devD_{mk} \epsilon_{lji} \pxj \A_{ik}\right) 
	+ \rho c_{h}^{2} c_{\bJ}\left( \J_{k}  \epsilon_{kjl} \pxj \psiJ_{l}
	-  \psiJ_{l}\epsilon_{ljk} \pxj \J_{k}\right)  \notag \\
	&=& \rho c_{s}^2 c_{\bA} \left( \A_{im}\devG_{mk} \epsilon_{ijl} \pxj \phiA_{lk}
	+ \phiA_{lm}\devD_{mk} \epsilon_{ijl} \pxj \A_{ik}\right) 
	+ \rho c_{h}^{2} c_{\bJ}\left( \J_{k}  \epsilon_{kjl} \pxj \psiJ_{l}
	+  \psiJ_{l}\epsilon_{kjl} \pxj \J_{k}\right)  \notag \\
	&=& \rho c_{\bA} c_{s}^2 \epsilon_{ijl} \left( \A_{im} \devG_{mk} \pxj \phiA_{lk} + \phiA_{lm}\devD_{mk}\pxj \A_{ik} \right) + c_{\bJ} c_{h}^{2}  \rho \epsilon_{kjl} \pxj \left(\psiJ_{l} \J_{k} \right) .
	\end{eqnarray}
	
	The remaining terms exactly correspond to those in the original GPR system and are therefore thermodynamically compatible. We briefly recall their relations. Considering the remaining convective terms of \eqref{eqn.GPRGLM_A} and \eqref{eqn.GPRGLM_J} and the terms on $\sigma_{ik}$, $\omega_{ik}$, it results
	\begin{equation}
		  \vel_i\! \left[ \pxk\sigma_{ik}\! +\!\pxk\omega_{ik}\right] 
		+ \alpha_{ik}\! \left[ \pxk \left(\vel_{m}\A_{im}\right)\! -\!\vel_{j}\pxk \A_{ij} \right] 
		+ \beta_{k}\! \left[ \pxk \left(\J_{m}\vel_{m} \right) \!-\! u_{j} \pxk \J_{j}  \right] \!%\notag\\
		= \pxk\! \left( \vel_{i}\sigma_{ik}\! +\! \vel_{i}\omega_{ik} \right) .
	\end{equation}
	Now, from the temperature term in \eqref{eqn.GPRGLM_J} and the term on $\tE_{J_{k}}$ in \eqref{eqn.GPR_S}, we get
	\begin{eqnarray}
		 \beta_{k} \pxk T + T \pxk\tE_{\J_{k}}  = \pxk \hq_{k}.
	\end{eqnarray}
	Finally, the dot product of the source terms multiplied by the corresponding dual variables yields
	\begin{eqnarray}
		\vel_i \rho g_i
		-\alpha_{ik} \frac{1}{\theta_{1} \left(\tau_{1}\right)} \E_{\A_{ik}}
		-\beta_{k} \frac{1}{\theta_{2} \left(\tau_{2}\right)} \E_{\J_{k}}
		T\left[ \frac{\rho}{T}\left(\frac{1}{\theta_{1}\left(\tau_{1}\right)}\E_{\A_{ik}}\E_{\A_{ik}} +  \frac{1}{\theta_{2}\left(\tau_{2}\right)}\E_{\J_{k}}\E_{\J_{k}}\right)\right] 
		= \vel_i \rho g_i. \label{HTCsources}
	\end{eqnarray}
	
	Hence, from \eqref{HTCdtflux}-\eqref{HTCsources}, we conclude that the dot product of equations \eqref{eqn.GPR_rho}-\eqref{eqn.GPR_mom}, \eqref{eqn.GPRGLM_A}-\eqref{eqn.GPRGLM_PsiJ}, \eqref{eqn.GPR_S} by the thermodynamically dual variables $\p_{GLM}$ yields \eqref{eqn.GPRGLM_E}.
\end{proof}

Let us note that the resulting augmented GLM GPR system \eqref{eqn.GPR_rho}-\eqref{eqn.GPR_mom}, \eqref{eqn.GPRGLM_A}-\eqref{eqn.GPRGLM_PsiJ}, \eqref{eqn.GPR_S}, \eqref{eqn.GPRGLM_E} is overdetermined. In what follows, we will address the submodel based on the total energy conservation law, 
\begin{subequations}\label{eqn.GPRGLM}
	\begin{eqnarray}
		\pt \rho + \pxk \left(\rho\vel_{k}\right) = 0,\label{eqn.GRPGLMs_rho}\\
		\pt \left(\rho\vel_i\right) + \pxk \left(\rho \vel_{i} \vel_{k} \right) + \pxi \press +\pxk\sigma_{ik} +\pxk\omega_{ik} = \rho g_{i},\label{eqn.GPRGLMs_mom}\\
		\pt \A_{ik} + \pxk \left(\vel_{m}\A_{im}\right) + \vel_{j} \pxj \A_{ik} -\vel_{j}\pxk \A_{ij} + c_{\bA} \epsilon_{ijl} \pxj \phiA_{lk} = -\frac{1}{\theta_{1} \left(\tau_{1}\right)} \E_{\A_{ik}},\label{eqn.GPRGLMs_A} \\
		\pt \phiA_{ik} + \vel_{j} \pxj \phiA_{ik} - c_{\bA}\epsilon_{ijl} \pxj \A_{lk}= 0,\label{eqn.GPRGLMs_PhiA}\\
		\pt \J_{k} + \pxk \left(\J_{m}\vel_{m} \right) +\pxk T + u_{j}\left( \pxj \J_{k} - \pxk \J_{j} \right)   + c_{\bJ} \epsilon_{kjl} \pxj \psiJ_{l} =  -\frac{1}{\theta_{2} \left(\tau_{2}\right)} \E_{\J_{k}}, \label{eqn.GPRGLMs_J}\\
		\pt \psiJ_{k} + \vel_{j} \pxj \psiJ_{k} - c_{\bJ}\epsilon_{kjl} \pxj \J_{l}= 0,\label{eqn.GPRGLMs_PsiJ}\\
		\pt \tE  +\pxk \left(\tE\vel_{k}\right) + \pxk \left(\press\vel_{k}\right) + \pxk \left[\vel_{i}\left( \sigma_{ik} + \omega_{ik} \right)\right] +\pxk \hq_{k} \notag\\
		+ \rho c_{\bA} c_{s}^2 \epsilon_{ijl} \left( \A_{im} \devG_{mk} \pxj \phiA_{lk} + \phiA_{lm}\devD_{mk}\pxj \A_{ik} \right) + c_{\bJ} c_{h}^{2}\epsilon_{kjl} \left( \rho \psiJ_{l} \pxj \J_{k} + \rho \J_{k} \pxj \psiJ_{l} \right) 
		= \rho g_{i} \vel_{i}. \label{eqn.GPRGLMs_E}
	\end{eqnarray}
\end{subequations}
i.e. we neglect the entropy inequality. For thermodynamically compatible methods discretizing the original GPR system containing the entropy relation, we refer to \cite{HTCGPR,HTCAbgrall,HTCA2}.

% % % % % % % % % % % % % % % % % % % % % % % % % % % % % %
% % % % % % % % % % % % % % % % % % % % % % % % % % % % % %
%              Numerical discretization
% % % % % % % % % % % % % % % % % % % % % % % % % % % % % %
% % % % % % % % % % % % % % % % % % % % % % % % % % % % % %
\section{Numerical discretization} \label{sec:numdisc}
The discretization of the GPR system will be performed in the framework of the hybrid finite volume/finite element approach put forward in 
\cite{BFSV14,BFTVC17} for incompressible flows and then extended to solve the compressible Navier-Stokes equations  \cite{Hybrid1,Hybrid2}, the shallow water equations  \cite{HybridMPI} and  the incompressible and weakly compressible GPR model \cite{HybridGPR}. Further, following the seminal ideas in \cite{HybridALE}, we propose an extension of the methodology to the arbitrary-Lagrangian-Eulerian framework. % allowing for small mesh deformations. 

This family of hybrid methods relies on an operator splitting approach which splits the system into a transport subsystem focused on the conservative variables and a Poisson type subsystem for the pressure. Consequently, we decouple the fast moving pressure waves from the bulk velocity of the medium. Then, an explicit discretization of the transport equations provides an accurate solution of discontinuities while the pressure subsystem can be efficiently solved using an implicit approach. Hence, the CFL time step condition does not depend on pressure waves that would greatly restrict the time step size. Further, the corresponding semi-discretization in time verifies the asymptotic preserving property for decreasing Mach numbers, which would yield to the so called all Mach number methods \cite{PM05,TV12}. 

In what follows, we first present the operator splitting of the augmented GLM GPR model and the overall methodology both in the purely Eulerian framework and in the ALE context. Then, we introduce the spatial discretization detailing the mesh notation and describe the algorithm's stages.

\subsection{Operator splitting and semi-discretization in time}
To ease the presentation of the operator splitting, we first focus on a purely Eulerian approach and introduce a semi-discretization in time of system \eqref{eqn.GPRGLM}, yielding
\begin{subequations}\label{eqn.GPR_semidiscrete}
	\begin{eqnarray}
		\frac{1}{\Dt} \left( \rho^{n+1}\!-\rho^{n} \right) + \pxk \left(\mom_{k}^{n}\right) = 0, \label{eqn.GPR_semidiscrete_rho}\\
		\frac{1}{\Dt} \left(\mom_i^{n+1}\!-\mom_i^{n}\right) + \pxk \left(\mom_{i}^{n} \vel_{k}^{n} \right) + \pxi \press^{n+1} +\pxk\sigma_{ik}^{n} +\pxk\omega_{ik}^{n}= \rho^{n} g_{i}, \label{eqn.GPR_semidiscrete_mom}\\
		\frac{1}{\Dt} \left( \A_{ik}^{n+1}\!-\A_{ik}^{n} \right) + \pxk \left(\vel_{m}^{n}\A_{im}^{n}\right) + \vel_{j}^{n} \pxj \A_{ik}^{n} -\vel_{j}^{n}\pxk \A_{ij}^{n} + c_{\bA} \epsilon_{ijl} \pxj \phiA_{lk}^{n} = -\frac{1}{\theta_{1}^{n} \left(\tau_{1}\right)} \E_{\A_{ik}}^{n}, \label{eqn.GPR_semidiscrete_A}\\
		\frac{1}{\Dt} \left(\phiA_{ik}^{n+1}\! - \phiA_{ik}^{n}\right)  + \vel_{j}^{n} \pxj \phiA_{ik}^{n} - c_{\bA}\epsilon_{ijl} \pxj \A_{lk}^{n}= 0 , \label{eqn.GPR_semidiscrete_pA}\\
		\frac{1}{\Dt} \left(\J_{k}^{n+1}\!-\J^{n}_{k}\right) + \pxk \left(\rho\J_{m}^{n}\vel_{m}^{n} \right) +\pxk T^{n} + \vel_{j}^{n} \left(\! \pxj \J_{k}^{n}\! -\! \pxk \J_{j}^{n}\! \right)  + c_{\bJ} \epsilon_{kjl} \pxj \psiJ_{l}^{n} = -\frac{1}{\theta_{2}^{n} \left(\tau_{2}\right)} \E_{\J_{k}}^{n}, \label{eqn.GPR_semidiscrete_J}\\
		\frac{1}{\Dt} \left(\psiJ_{k}^{n+1}\! -\psiJ_{k}^{n}\right) + \vel_{j}^{n} \pxj \psiJ_{k}^{n} - c_{\bJ}\epsilon_{kjl} \pxj \J_{l}^{n}= 0, \label{eqn.GPR_semidiscrete_pJ}\\
		\frac{1}{\Dt} \left( \tE^{n}\!-\tE^{n+1}\right)  
		+\pxk \left(\tE_{1}^{n+1}\vel_{k}^{n+1}\right)  
		+\pxk \left(\tE_{2}^{n}\vel_{k}^{n}\right) 
		+\pxk \left(\tE_{3}^{n}\vel_{k}^{n}\right) 
		+\pxk \left(\tE_{4}^{n}\vel_{k}^{n}\right)			
		+\pxk \left(\tE_{5}^{n}\vel_{k}^{n}\right) \notag\\	 
		+\pxk \left(\tE_{6}^{n}\vel_{k}^{n}\right)
		+ \pxk \left(\press^{n+1}\vel_{k}^{n+1}\right) 
		+ \pxk \left(\vel_{i}^{n}\sigma_{ik}^{n} \right) 
		+ \pxk \left(\vel_{i}^{n}\omega_{ik}^{n} \right) 
		+\pxk \hq_{k}^{n} 
		=  g_{i} \mom_{i}^{n}. \label{eqn.GPR_semidiscrete_E}
	\end{eqnarray}
\end{subequations}

Next, extending the methodology proposed in \cite{Hybrid2}, for solving the compressible Navier-Stokes equations, to the solution of the GPR model, we consider a TV-splitting of the momentum \eqref{eqn.GPR_semidiscrete_mom} and energy equations \eqref{eqn.GPR_semidiscrete_E},  which  leads to a transport and a pressure subsystem, \cite{TV12}. 
Introducing the intermediate auxiliary variable $\bmom^{\ast}$, we get
\begin{gather}
	\frac{1}{\Dt} \left(\rho\vel_i^{\ast}-\rho^{n}\vel_i^{n}\right) + \pxk \left(\rho^{n} \vel_{i}^{n} \vel_{k}^{n} \right) + \pxi \press^{n} +\pxk\sigma_{ik}^{n} +\pxk\omega_{ik}^{n} = \rho^{n} g_{i}, \label{eqn.sdmom_td}\\
	\frac{1}{\Dt} \left(\rho^{n+1}\vel_i^{n+1}-\rho\vel_i^{\ast}\right) + \pxi \left( \press^{n+1} -  \press^{n}\right)  = 0.
\end{gather}
Hence, the momentum at the new time step can be computed from
\begin{equation}
	\mom_k^{n+1}=\mom_k^{\ast} - \Dt\pxk \delta \press^{n+1}, \qquad \delta \press^{n+1}= \press^{n+1}-\press^{n},\label{eqn.sd_mom_update}
\end{equation}
once the transport-diffusion equation \eqref{eqn.sdmom_td} is solved.

Similarly, focusing on the energy equation, \eqref{eqn.GPR_semidiscrete_E}, decomposing the energy density within the flux term into its six components and introducing the intermediate auxiliary energy density $\tE^{\ast}$, we obtain
\begin{eqnarray}
	\frac{1}{\Dt} \left( \tE^{\ast}-\tE^{n}\right)  
	+ \pxk \left(\tE^{n}_{2}\vel_{k}^{n}\right)
	+ \pxk \left(\tE^{n}_{3}\vel_{k}^{n}\right)  
	+ \pxk \left(\tE^{n}_{4}\vel_{k}^{n}\right)
	+ \pxk \left(\tE^{n}_{5}\vel_{k}^{n}\right) 
	+ \pxk \left(\tE^{n}_{6}\vel_{k}^{n}\right)\notag\\
	+ \pxk \left(\vel_{i}^{n}\sigma_{ik}^{n} \right) 
	+ \pxk \left(\vel_{i}^{n}\omega_{ik}^{n} \right) 
	+\pxk \hq_{k}^{n} 
	=g_{i} \mom_{i}^{n}, \label{eqn.sdE_td}\\
	\frac{1}{\Dt} \left( \tE^{n+1}-\tE^{\ast}\right)  
	+ \pxk \left(\tE^{n+1}_{1}\vel_{k}^{n+1}\right)
	+ \pxk \left(\press^{n+1}\vel_{k}^{n+1}\right)
	= 0. \label{eqn.sdE_ntd}
\end{eqnarray}
Therefore, from \eqref{eqn.sdE_ntd}, 
\begin{equation}
	 \tE^{n+1}= \tE^{\ast}- \Dt \pxk \left(\tE^{n+1}_{1}\vel_{k}^{n+1}\right)
	- \Dt\pxk \left(\press^{n+1}\vel_{k}^{n+1}\right) \label{eqn.sd_E_update}
\end{equation}
with $\tE^{\ast}$ the solution of \eqref{eqn.sdE_td}. Substituting the expression for $\tE_{1}$ in terms of the pressure and gathering terms yields
\begin{equation}
	\tE^{n+1}= \tE^{\ast}- \Dt \pxk \left(\h^{n+1} \mom_{k}^{n+1}\right), \label{eqn.sd_En1}
\end{equation}
where the enthalpy reads 
\begin{equation}\h=\dfrac{\gamma}{\rho\left( \gamma-1\right) }\, \press \label{eqn.enthalpy_ideal}
\end{equation}
if the ideal gas EOS \eqref{eqn.E4} is employed. Meanwhile, for the stiffened gas EOS \eqref{eqn.E4_stiff}, we have
\begin{equation}
	\h=\dfrac{\gamma \press + \rho_{0}c_{0}^{2}-\gamma\press_{0}}{\rho\left( \gamma-1\right) } .
	\label{eqn.enthalpy_stiff}
\end{equation}
Taking into account \eqref{eqn.energydecomp}, the energy equation is then rewritten in terms of the pressure unknown as
\begin{gather}
    \frac{\press^{n+1}}{\left( \gamma-1\right) } - \frac{\press^{n}}{\left( \gamma-1\right) }= -\tE_{1}^{n}	-\tE_{2}^{n+1} -\tE_{3}^{n+1}-\tE_{4}^{n+1} -\tE_{5}^{n+1}-\tE_{6}^{n+1}+\tE^{\ast} -\Dt \pxk \left(  \h^{n+1} \mom_{k}^{n+1}\right),
\end{gather}
where we have subtracted the term $\tE_{1}$ at the two sides of \eqref{eqn.sd_En1}. Finally, substitution of \eqref{eqn.sd_mom_update} yields the pressure equation
\begin{eqnarray}
	\frac{1}{\left( \gamma-1\right) }\delta \press^{n+1}
	-\Dt^2 \pxk \left(  \h^{n+1} \pxk \delta \press^{n+1}\right) 
	=\tE^{\ast}  
	-\tE_{1}^{n}
	-\tE_{2}^{n+1} 
	-\tE_{3}^{n+1}
	-\tE_{4}^{n+1}
	-\tE_{5}^{n+1} \notag \\
	-\tE_{6}^{n+1}
	-\Dt \pxk \left(  \h^{n+1} \mom_k^{\ast} \right) \label{eqn.sd_press}
\end{eqnarray}
whose solution is employed in \eqref{eqn.sd_mom_update} and \eqref{eqn.sd_E_update} to obtain the momentum and total energy at the new time step.

The highly nonlinear system \eqref{eqn.sd_press} involves the term $\tE_{1}^{n+1}=\halb\rho^{n+1}\left|\bvel^{n+1}\right|^2$ so it cannot be solved independently of \eqref{eqn.sd_mom_update}, i.e., the system to be solved involves \eqref{eqn.GPR_semidiscrete_rho}, \eqref{eqn.sdmom_td}, \eqref{eqn.GPR_semidiscrete_A}-\eqref{eqn.GPR_semidiscrete_pJ}, \eqref{eqn.sd_press}, \eqref{eqn.sdE_td}, \eqref{eqn.sd_mom_update}, \eqref{eqn.sd_E_update}. To deal with these crossed $\press^{n+1}$ and $\bmom^{n+1}$ terms and with the non linearity introduced by the presence of the enthalpy in the stiffness matrix, we apply a Picard procedure getting the final system
\begin{subequations}\label{eqn.GPRGLM_disc}
	\begin{align}
		&\rho^{n+1}  = \rho^{n} -\Dt \pxk \left(\mom_{k}^{n}\right), \label{eqn.GPRGLM_disc_rho}\\
		&\mom_i^{\ast} = \mom_i^{n} -\Dt\left(  \pxk \left(\rho^{n} \vel_{i}^{n} \vel_{k}^{n} \right) + \pxi \press^{n} +\pxk\sigma_{ik}^{n} +\pxk\omega_{ik}^{n} - \rho^{n} g_{i}\right) ,\label{eqn.GPRGLM_disc_mom}\\
		&\A_{ik}^{n+1} = \A_{ik}^{n}- \Dt \left( \pxk \left(\vel_{m}^{n}\A_{im}^{n}\right)+ \vel_{j}^{n} \pxj \A_{ik}^{n} -\vel_{j}^{n}\pxk \A_{ij}^{n} + c_{\bA} \epsilon_{ijl} \pxj \phiA_{lk}^{n} +\frac{1}{\theta_{1}^{n} \left(\tau_{1}\right)} \E_{\A_{ik}}^{n}\right) , \label{eqn.GPRGLM_disc_A}\\
		&\phiA_{ik}^{n+1} = \phiA_{ik}^{n}  -\Dt\left(  \vel_{j}^{n} \pxj \phiA_{ik}^{n} - c_{\bA}\epsilon_{ijl} \pxj \A_{lk}^{n}\right) ,\label{eqn.GPRGLM_disc_phiA}\\
		&\J_{k}^{n+1} = \J_{k}^{n} -\Dt\left(  \pxk \left(\J_{m}^{n}\vel_{m}^{n} \right) 
		+\vel_{j} \left( \pxj \J_{k}^{n} - \pxk \J_{j}^{n} \right)
		+\pxk T^{n} + c_{\bJ} \epsilon_{kjl} \pxj \psiJ_{l}^{n}  +\frac{1}{\theta_{2}^{n} \left(\tau_{2}\right)} \E_{\J_{k}}^{n}\right) ,\label{eqn.GPRGLM_disc_J}\\
		& \psiJ_{k}^{n+1} = \psiJ_{k}^{n}-\Dt\left(  \vel_{j}^{n} \pxj \psiJ_{k}^{n} - c_{\bJ}\epsilon_{kjl} \pxj \J_{l}^{n}\right) , \label{eqn.GPRGLM_disc_psiJ}\\
		&\tE^{\ast} = \tE^{n} 
		-\Dt\!\left( \pxk \left(\tE^{n}_{2}\vel_{k}^{n}\right)
		+ \pxk \left(\tE^{n}_{3}\vel_{k}^{n}\right)
		+ \pxk \left(\tE^{n}_{4}\vel_{k}^{n}\right)
		+ \pxk \left(\tE^{n}_{5}\vel_{k}^{n}\right)
		+ \pxk \left(\tE^{n}_{6}\vel_{k}^{n}\right) \right.
		\notag\\&\hspace*{7cm}  \left.
		+ \pxk \left[\vel_{i}^{n}\left( \sigma_{ik}^{n}+\omega_{ik}^{n} \right)\right] 
		+\pxk \hq_{k}^{n} 
		-g_{i} \mom_{i}^{n}\right) , \label{eqn.GPRGLM_disc_E}\\
		&\frac{\delta \press^{n+1,\ell+1}}{ \gamma-1 }
		-\Dt^2 \pxk \left(  \h^{n+1,\ell} \pxk \delta \press^{n+1,\ell+1}\right) 
		=\!\tE^{\ast}\!\! 
		-\tE^{n}_{1}
		-\!\tE_{2}^{n+1,\ell} \!
		-\!\tE_{3}^{n+1}\!
		\notag\\
		&\hspace*{8cm} 
		-\!\tE_{4}^{n+1}\!
		-\!\tE_{5}^{n+1}\!
		-\!\tE_{6}^{n+1}\!
		-\Dt \pxk \left(  \h^{n+1,\ell} \mom_k^{\ast} \right), \label{eqn.GPRGLM_disc_press}\\
		&\mom_k^{n+1,\ell+1}=\mom_k^{\ast} - \Dt\pxk \delta \press^{n+1,\ell+1}, \label{eqn.GPRGLM_disc_monf}\\
		&\tE^{n+1}= \tE^{\ast}- \Dt \pxk \left(\h^{n+1,\ell+1} \mom_{k}^{n+1,\ell+1}\right) \label{eqn.GPRGLM_disc_Ef}
	\end{align}
\end{subequations}
with $\ell$ the Picard iteration index, $\ell=1,\dots,N_{\mathrm{Pic}}$ and
\begin{equation}
	\press^{n+1,\ell+1} = \press^{n} + \delta\press^{n+1,\ell+1},
\end{equation}
needed to compute $\h^{n+1,\ell+1}$. Note that the Picard iteration does not affect the computation of the terms related to the distortion field and the thermal impulse in \eqref{eqn.GPRGLM_disc_press} since their value at the new time step is directly obtained from \eqref{eqn.GPRGLM_disc_A}-\eqref{eqn.GPRGLM_disc_psiJ}. 

Besides, the source terms related to the distortion tensor and thermal impulse may become very stiff, so employing an explicit approach for their discretization may require very small time steps. To circumvent this issue, we split \eqref{eqn.GPRGLM_A} and \eqref{eqn.GPRGLM_J} into two subsystems \cite{BCP22_GPR_IMEXLagUnst,HybridGPR}. The first one accounts for the contributions of the flux term and the non-conservative products,
\begin{align}
	% Distortion
	&\A_{ik}^{\ast}=\A_{ik}^{n} - \Dt \left( 
		 \pxk \left(\vel_{m}^{n}\A_{im}^{n}\right) 
	+ \vel_{j}^{n} \pxj \A_{ik}^{n} 
	-\vel_{j}^{n}\pxk \A_{ij}^{n}
	+ c_{\bA} \epsilon_{ijl} \pxj \phiA_{lk}^{n}\right) , \\
	% thermal impulse
	& \J_{k}^{\ast} = \J_k^{n} - \Dt \left( 
	 \pxk \left(\J_{m}^{n}\vel_{m}^{n} \right) 
	 +\pxk T^{n} 
	 + \vel_{j}^{n}\left(\pxj J_{k}^{n} + \pxk J_{j}^{n}\right)
	 + c_{\bJ} \epsilon_{kjl} \pxj \psiJ_{l}^{n}\right).
\end{align}
Meanwhile, the second subsystem contains only the algebraic source terms and can be implicitly discretized:
\begin{subequations}\label{eqn.sourcesystem_sd}
	\begin{align}
		% Distortion
		&\frac{1}{\Delta t}\left( \A_{ik}^{n+1}\!-\!\A_{ik}^{\ast}\right) 
		= -\frac{1}{\theta_{1}^{n+1} \left(\tau_{1}\right)} \E_{\A_{ik}}^{n+1}, \\
		% thermal impulse
		&\frac{1}{\Delta t}\!\left( \J_{k}^{n+1}\!-\!\J_k^{\ast}\right) 
		=  -\frac{1}{\theta_{2}^{\star} \left(\tau_{2}\right)} \E_{\J_{k}}^{n+1}.
	\end{align}
\end{subequations}

\subsubsection{Eulerian hybrid FV/FE method}
Attending to the different nature of the equations on \eqref{eqn.GPRGLM_disc}, we split the Eulerian hybrid finite volume/finite element algorithm into four main stages: 

\begin{enumerate}
	\item Transport stage. A finite volume method is employed to solve system \eqref{eqn.GPRGLM_disc_rho}-\eqref{eqn.GPRGLM_disc_E} yielding the solution at the new time step for the density, distortion and thermal impulse fields, $\rho^{n+1}$, $\bA^{n+1}$, $\bJ^{n+1}$, and intermediate values for the momentum and total energy unknowns, $\bmom^{\ast}$, $\tE^{\ast}$. In case the GLM GPR model is solved, also the cleaning variables are computed at this stage, $\bpA^{n+1}$, $\bpJ^{n+1}$.
	\item Intermediate stage. The intermediate values obtained for the momentum and total energy are interpolated between the staggered grids allowing for the computation of the source terms in the next stage.
	\item Pressure stage. A Picard iteration procedure is employed to approximate the pressure at the new time step, $\press^{n+1}$. Within each Picard iteration a $\mathbb{P}^{1}$ implicit finite element method is applied to solve \eqref{eqn.GPRGLM_disc_press}.  
	Then, the enthalpy and momentum unknowns are updated following \eqref{eqn.enthalpy_stiff} and \eqref{eqn.GPRGLM_disc_monf}, getting $h^{n+1,\ell+1}$ and $\bmom^{n+1,\ell+1}$, that allows the computation of $\tE_{2}^{n+1,\ell+1}$.
	\item Correction stage. Once the final momentum, $\bmom^{n+1}$, is calculated substituting the pressure correction $\delta \press^{n+1}$ in \eqref{eqn.GPRGLM_disc_monf}, the obtained velocity is employed in \eqref{eqn.GPRGLM_disc_Ef} to get the total energy at time $t^{n+1}$, $\tE^{n+1}$.
\end{enumerate}
A detailed description of each of these stages will be provided in Sections~\ref{sec:transport}-\ref{sec:intermediate}.

\subsubsection{Extension to the ALE framework}
The former methodology assumes a purely Eulerian framework. In order to deal also with small geometry deformations, we extend the hybrid approach to the context of Arbitrary-Lagrangian-Eulerian methods. Following \cite{HybridALE}, we introduce a new stage at the beginning of each time step: 
\begin{itemize}
	\item[0.] Mesh motion stage. In this stage, we move the mesh attending to the local fluid velocity at the previous time step, $\bvel^{n}$.  Further, the mesh may also be reshaped attending to prescribed displacements of the boundaries of the computational domain. 
\end{itemize}
Next, in the transport stage, we consider the space-time divergence form of the transport equations which are integrated over space-time control volumes. As a result, the intermediate approximations reside on the deformed mesh. Then, the remaining stages can be performed as for the purely Eulerian method only considering the new grid configuration. 

\subsection{Spatial discretization: staggered unstructured grids}
\begin{figure}[h]
	\centering
	\includegraphics[width=0.45\linewidth]{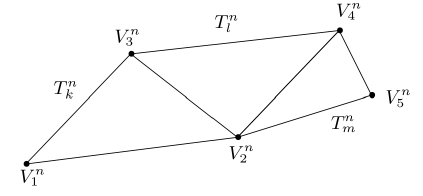}\hspace{0.25cm}
	\includegraphics[width=0.45\linewidth]{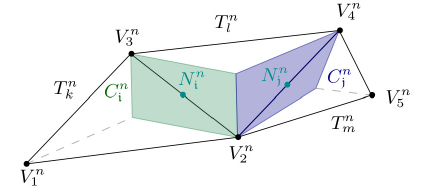}
	\caption{Construction of two interior dual elements $\cell{\ci}$ and $\cell{\cj}$ from primal elements $\Pcell{k}$, $\Pcell{l}$, $\Pcell{m}$.}
	\label{fig.mesh}
\end{figure}

For the discretization of the spatial domain $\Omega$,  we employ unstructured staggered grids of the face type, \cite{BDDV98,BFSV14,TD16}, also known as diamond grids. We consider a primal mesh made of non-overlapping triangles in 2D and tetrahedra in 3D denoted by $\Pcell{k}$, $k=1,\dots,N_{\mathrm{el}}$, $\Omega=\bigcup_{k} \Pcell{k}$. Focusing on one interior edge/face of the domain, $\Gamma_{i}$, we can compute the barycentres, $B_{i_L}$,  $B_{i_R}$, of the two neighbouring primal elements related to the edge/face  and construct a dual cell $\cell{\ci}$ by merging the two subtriangles/subtetrahedra with vertex each of the barycentres and the vertex of $\Gamma_{i}$. Similarly, for a boundary edge/face,  we construct a boundary dual cell as the subelement with vertex those of $\Gamma_{i}$ and the barycentre of the related primal element. A sketch on the dual mesh construction in 2D is provided in Figure~\ref{fig.mesh}. Further details on the staggered grid generation can be found in \cite{BFSV14,Hybrid2}. 
To describe the proposed hybrid methodology, we furthermore need to introduce the following notation related to the mesh:
\begin{itemize}
\item $V_{j}$, $j\in\left\lbrace 1,\dots,N_{\mathrm{ver}}\right\rbrace$, are the vertex of the primal grid, $N_{\mathrm{ver}}$ the total number of vertex.
\item $\left|\Pcell{k}\right|$ is the area/volume of the primal element $\Pcell{k}$.
\item $\mathcal{K}_{k}$ is the set of indexes of primal elements adjacents to $\Pcell{k}$.
\item $\Gamma_{kl}$ is the boundary shared by the primal elements $\Pcell{k}$ and  $\Pcell{l}$.
\item $\normalu{kl}$ is the unit normal vector of $\Gamma_{kl}$ exterior to $\Pcell{k}$.
\item $\normalw{kl}$ is the normal vector of $\Gamma_{kl}$ exterior to $\Pcell{k}$ and weighted the length/area of $\Gamma_{kl}$, $\|\normalw{kl}\|$.
\item $\mathcal{\hat{K}}_{k}$ is the set of indexes of dual cells generated using $\Pcell{k}$.
\item $\left|\Pcell{k\ci}\right|$ is the area/volume of the intersection of $\Pcell{k}$ with $\cell{\ci}$.
\item $\mathcal{\hat{K}}_{\ci}$ is the set of indexes of primal elements generating the dual cell $\cell{\ci}$.
\item $N_{\ci}$, $\ci\in\left\lbrace 1,\dots,N_{\mathrm{cell}}\right\rbrace$, are the vertex of the primal grid, $N_{\mathrm{cell}}$ the total number of dual cells.
\item $\left|\cell{\ci}\right|$ is the area/volume of $\cell{\ci}$.
\item $\mathcal{K}_{\ci}$ is the set of indexes of dual elements adjacents to cell $\cell{\ci}$.
\item $\boundary{\ci}=\partial \cell{\ci}=\bigcup_{\cj\in \mathcal{K}_{\ci}} \boundary{\ci\cj}$ is the boundary of $\cell{\ci}$ and $\boundary{\ci\cj}$ represents the edge/face shared with cell $\cell{\cj}$. 
\item $\normalu{\ci\cj}$ is the outward unit normal vector of $\boundary{\ci\cj}$.
\item $\normalw{\ci\cj}$ is the outward normal vector of $\boundary{\ci\cj}$ exterior to $\cell{\ci}$ weighted with the  length/area of $\boundary{\ci\cj}$, $\|\normalw{\ci\cj}\|$.
\item $ \cell{\ci}^{\circ} = \partial \cell{\ci}\setminus \boundary{\ci}$ is the interior of $\cell{\ci}$.
\end{itemize}

On the other hand, the ALE approach, at each time step, requires the update of the coordinates of the primal and dual grids as well as related data; e.g. areas, volumes and boundary normals. It is important to remark that topology changes are not allowed so the connectivities of the mesh structures do not change. Thus, the work load lowers with respect to the generation of a completely new mesh in the modified computational domain. For alternative ALE approaches allowing large deformations and topology changes without the need of remeshing the whole domain, we refer to \cite{GBCKSD19,Gab25_ALEDGtopo}.

\begin{figure}[H]
	\centering
	\includegraphics[width=0.4\linewidth]{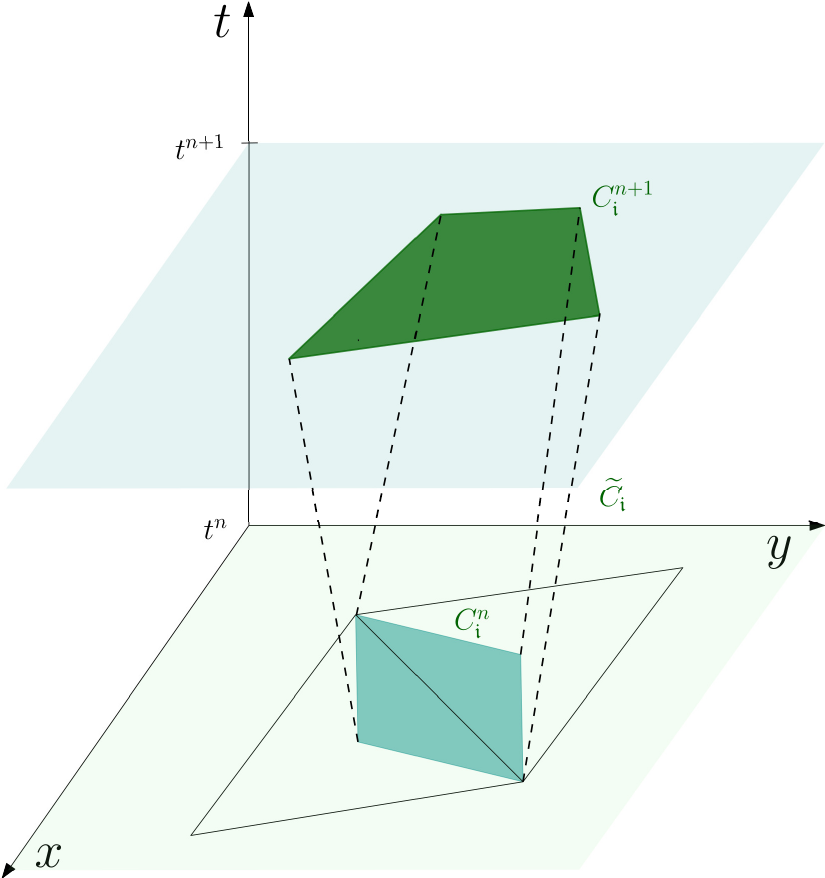}
	\caption{Construction of a space-time control volume $\tcell{\ci}$ generated by the deformation of the dual cell  $\cell{\ci}$ between times $t^{n}$ and $t^{n+1}$. $\cell{\ci}^{n}$ and $\cell{\ci}^{n+1}$ are shadowed in green while the space time boundary is indicated with dashed lines. The black triangles at time $t^{n}$ plane correspond to the primal elements generating the dual cell $\cell{\ci}^{n}$.}
	\label{fig.spacetimecontrolvolume}
\end{figure}

\begin{figure}[h]
	\centering
	\includegraphics[width=0.9\linewidth]{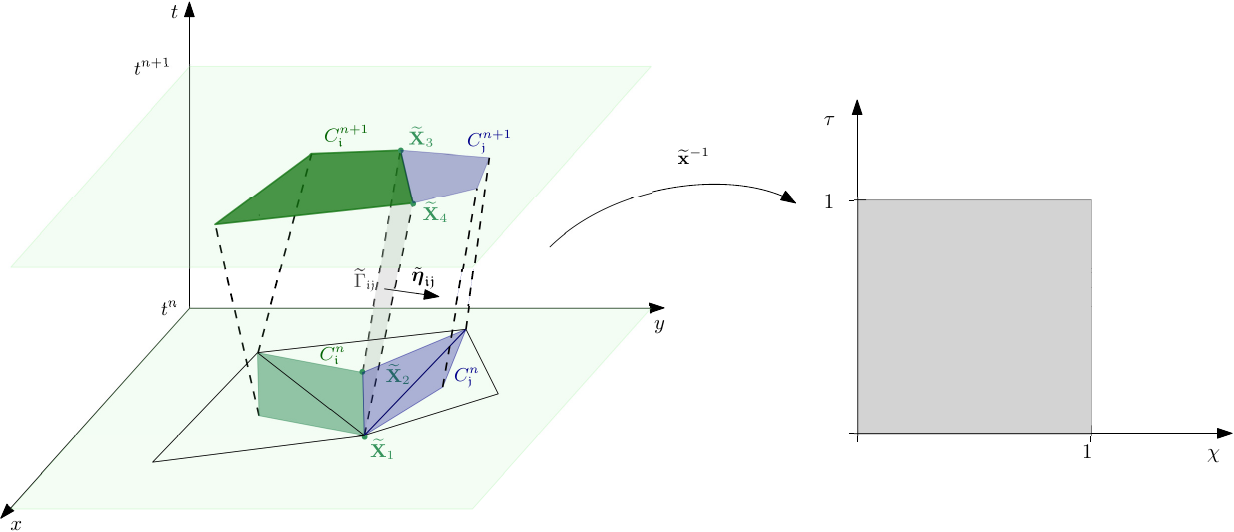}
	\caption{Map of a space-time face $\tboundary{\ci\cj}$ to the reference element in the local coordinate system $\chi-\tau$. The face $\tboundary{\ci\cj}$ (shadowed in gray) shared by the space-time control volumes $\tcell{\ci}$  and $\tcell{\cj}$ is determined by vertex $\Xst_1$, $\Xst_2$, $\Xst_3$ and $\Xst_4$.}
	\label{fig.referencemap}
\end{figure}

Given a dual cell at two subsequent time steps, $\cell{\ci}^{n}$ and $\cell{\ci}^{n+1}$, we design a space-time control volume $\tcell{\ci}$ joining the vertex of each dual element at both times so $\tcell{\ci\mid_{t^{n}}}=\cell{\ci}^{n}$ and $\tcell{\ci\mid_{t^{n+1}}}=\cell{\ci}^{n+1}$. This control volume, depicted in Figure~\ref{fig.spacetimecontrolvolume}, will be then employed within the transport stage to solve the transport subsystem. Analogously with the notation of steady grids, we label $\widetilde{N}_{\ci}$, $\tboundary{\ci} = \partial \tcell{\ci}$ and $ \tcell{\ci}^{\circ} = \partial \tcell{\ci}\setminus \tboundary{\ci}$ the node, the boundary and the interior of $\tcell{\ci}$, respectively. Further, the space-time faces of the control volumes located between two neighbouring elements $\tcell{\ci}$ and $\tcell{\cj}$ are denoted by $\tboundary{\ci\cj}$ and its corresponding outward unit normal with respect to $\tcell{\ci}$ is $\tnormal{\ci\cj}=(\tnormale{t},\tnormal{\x})$. Consequently, we have 
\begin{equation}
	\tboundary{i} = \cell{\ci}^{n} \cup \left( \bigcup_{\cj\in \mathcal{K}_{\ci}} \tboundary{\ci\cj} \right) \cup  \cell{\ci}^{n+1}  \label{tboundary}
\end{equation}
and the outward pointing unit normal vectors to $\cell{\ci}^{n}$ and $\cell{\ci}^{n+1}$ are $\tnormal{}=(-1,\boldsymbol{0})$ and $\tnormal{}=(1,\boldsymbol{0})$, respectively. 
Besides, the space-time faces $\tboundary{\ci\cj}$ can be mapped into a reference element in a local coordinate system. For instance, in the bidimensional case, we consider the coordinate system $\chi-\tau$ with basis functions $\Phist_{\ell}=\Phist_{\ell}(\chi,\tau)$ given by
\begin{equation}
	\Phist_1 = (1-\chi)(1-\tau), \quad
	\Phist_2 = \chi(1-\tau),     \quad
	\Phist_3 = (1-\chi)\tau,     \quad
	\Phist_4 = \chi\tau,         \qquad
	0\le\chi\le 1,\quad
	0\le\tau\le 1.
\end{equation}
Then, denoting $\Xst_{\ell}$ the space-time coordinate vectors of the four vertices which form the face, 
\begin{align*}
	\Xst_1 &= (t^{n}, \XX_1^{n})=(t^{n}, x_1, y_1),& \qquad
	&\Xst_2 = (t^{n}, \XX_2^{n})=(t^{n}, x_2, y_2),&     \\
	\Xst_3 &= (t^{n+1}, \XX_2^{n+1})=(t^{n+1}, x_3, y_3), &    \qquad
	&\Xst_4 = (t^{n+1}, \XX_1^{n+1})=(t^{n+1}, x_4, y_4),&        
\end{align*}
we define the map from the reference configuration to $\tboundary{\ci\cj}$ as
\begin{equation*}
	\xst_{\ci\cj}(\chi,\tau) = \Phist_1\Xst_1+\Phist_2\Xst_2+\Phist_3\Xst_3+\Phist_4\Xst_4,
\end{equation*}
see Figure~\ref{fig.referencemap}. Furthermore, the space-time unit normal vector on the space-time face $\tboundary{\ci\cj}$ reads $$\nst_{\ci\cj}=(\tilde{n}_{\ci\cj}^t,\tilde{n}_{\ci\cj}^x,\tilde{n}_{\ci\cj}^y) = \frac{\left(\frac{\partial \xst}{\partial \chi}\times \frac{\partial \xst}{\partial \tau}\right)}{\left\|\frac{\partial \xst}{\partial \chi}\times \frac{\partial \xst}{\partial \tau}\right\|}.$$ 
Then, denoting its integral 
\begin{equation}
	\etast_{\ci\cj} := \int_{\tboundary{\ci\cj}} \nst_{\ci\cj} \, \dt \dS,  
\end{equation}
we observe that
\begin{equation}
	\int_{\tboundary{\ci}}\nst \dS = \sum_{\cj \in \mathcal{K}_{\ci}}  \etast_{\ci\cj} = 0,
	\label{eqn.GCL2}
\end{equation}
i.e. the geometric conservation law (GCL) \cite{Farhat2001} will be verified when using closed space-time control volumes for the integration of the equations at each time step \cite{Lagrange2D,BD14}.

\subsection{Mesh motion}
Given a vertex on the primal grid at time $t^{n}$, $\XX^{n}_{\ell}$, we compute its position at the new time step $\XX_{\ell}^{n+1}$ as
\begin{equation}
	\XX_{\ell}^{n+1} = \XX_{\ell}^{n} + \Dt \bV_{\ell}^{n}
\end{equation}
with $\bV_{\ell}^{n}$ the mesh velocity. This velocity can be set in two different ways:
\begin{itemize}
	\item According to a prescribed boundary velocity, $\bV^{n}_{BC}$. 
	Assuming that we have the description of the boundary motion, the mesh displacement at the interior of the domain is computed by solving a Laplace equation on the velocity with Dirichlet boundary conditions,
	\begin{equation}
		\left\lbrace
		\begin{array}{ll}
			\nabla^{2} \bV^{n} = 0,  &\forall\, \XX\in\Omega^{n},\\
			\bV^{n}(\XX) = \bV^{n}_{BC}, & \forall\, \XX\in\partial\Omega^{n}.
		\end{array}
		\right.
	\end{equation}

	\item As a smoothed local fluid velocity. Since we are working in the ALE framework, we can arbitrarily choose the mesh velocity. For instance, we can set it to match the local fluid velocity, as in classical Lagrangian schemes, or apply a smoothing operator in order to reduce mesh distortion at the interior of the computational domain \cite{BD14,GBCKSD19}. Accordingly, we define a regularization parameter $\zeta\in[0,\infty)$ and compute the mesh velocity solving
	\begin{equation}
		\left\lbrace
		\begin{array}{ll}
			\bV^{n}(\XX) - \Dt \zeta \nabla^2\bV^{n}(\XX) = \bvel^{n} ,  &\forall\, \XX\in\Omega^{n},\\
			\bV^{n}(\XX) = \bvel^{n}_{BC}, & \forall\, \XX\in\partial\Omega^{n},
		\end{array}
		\right.
	\end{equation}
	at the aid of $\mathbb{P}^{1}$ continuous finite element methods.
	Following this approach, for $\zeta = 0$ we retrieve the pure Lagrangian scheme while as $\zeta\rightarrow \infty$ we have a Laplacian smoothing of the mesh velocity.
	
\end{itemize} 
Let us note that, as it will be seen in the numerical results section, both types of motion may be combined by setting some of the boundaries to have a specific boundary velocity while others freely move with the local fluid velocity. This methodology allows, e.g., for the simulation of free surface flows in confined moving domains. Further details on the mesh motion within the hybrid FV/FE methodology can be found in \cite{HybridALE}.

\subsection{Transport stage: a finite volume scheme on the space-time divergence form of the transport subsystem}\label{sec:transport}
The transport subsystem \eqref{eqn.GPRGLM_disc_rho}-\eqref{eqn.GPRGLM_disc_E} is solved at the aid of an explicit finite volume method.
Introducing the vector of conservative variables $\Q = \left(\rho, \bvel, \bA, \bpA,  \bJ, \bpJ,\tE \right)^{T}$, the continuous counterpart of the transport equations can be recast into a general system of the form
\begin{equation}
	\partial_{t} \Q + \dive \Flux\left(\Q\right) + \NCP{\Q} = \source\left(\Q\right),\label{eqn.generaltransportncp}
\end{equation}
with $\Flux\left(\Q\right)$ the convective flux terms, $\NCP{\Q}$ the non-conservative products, and $\source\left(\Q\right)$ the source terms. 
Hence, its space-time divergence formulation reads
\begin{equation}
	\nablast \cdot \Fluxst (\Q) + \NCPst{\Q} = \mathbf{S}\left( \Q \right) 	\label{eqn.stdiv}  
\end{equation} 
with 
\begin{equation*}
	\nablast = \left( \frac{\partial}{\partial t}, \nabla \right), \qquad 
	\Fluxst = \left( \Q, \Flux (\Q) \right), \qquad 
	\BNCPst{\Q} = \left( 0, 
	\boldsymbol{\mathcal{B}}_1(\Q), 
	\boldsymbol{\mathcal{B}}_2(\Q), 
	\boldsymbol{\mathcal{B}}_3(\Q)  
	\right).
\end{equation*}
Integration of \eqref{eqn.stdiv} on a space-time control volume $\tcell{\ci}$, leads to
\begin{equation}
	 \int\limits_{\tcell{\ci}} \nablast \cdot \Fluxst\left(\Q\right) \dt \dV 
	+  \int\limits_{\tcell{\ci}} \NCPst{\Q}\dt\dV 
	= \int\limits_{\tcell{\ci}}  \source\left(\Q\right)\dt\dV.
	\label{eqn.generaltransportncp_disc}
\end{equation}
Now, Gauss theorem is applied to transform the space-time integral of the flux into the integral over the space-time boundary $\tboundary{\ci}$.  Moreover, we employ a path conservative approach to deal with the non-conservative products splitting their contribution into a smooth part within the control volume and the jumps at its boundaries \cite{Par06}. Hence, we have
\begin{equation}
	\int\limits_{\tboundary{\ci}} \Fluxst\left(\Q\right) \cdot \tnormal{\ci}\dt \dS
	+  \int\limits_{\tboundary{\ci}} \boldsymbol{\widetilde{\mathcal{D}}} \left(\Q\right) \cdot \tnormal{\ci} \dt\dS 
	+  \int\limits_{\tcell{\ci}^{\circ}} \NCPst{\Q} \dt\dV 
	= \int\limits_{\tcell{\ci}}  \source\left(\Q\right)\dt\dV.
	\label{eqn.generaltransportncp_disc2}
\end{equation}
From \eqref{tboundary}, we can split the integrals over the faces into the integral over the space-time boundary surfaces $\tboundary{\ci\cj}$ and the two boundaries of $\tcell{\ci}$ orthogonal to the time axis, i.e. $\cell{\ci}^{n}$ and $\cell{\ci}^{n+1}$, so
\begin{eqnarray}
	\left|\tcell{\ci}\right| \Q_{\ci}^{\ast}  =
	\left|\tcell{\ci}\right| \Q_{\ci}^{n} 
	- \sum_{\cj\in\mathcal{K}_{\ci}} \int\limits_{\tboundary{\ci\cj}} \left(  \Fluxst\left(\Q\right) +  \boldsymbol{\widetilde{\mathcal{D}}} \left(\Q\right)\right)  \cdot \tnormal{\ci\cj} \dt\dS 
	-  \int\limits_{\tcell{\ci}^{\circ}} \NCPst{\Q} \dt\dV 
	+ \int\limits_{\tcell{\ci}^{\circ}}  \source\left(\Q\right)\dt\dV,\label{eqn.spacetimedisc}
\end{eqnarray}
where, for the GLM GPR model, we define $\Q^{\ast}=\left(\rho^{n+1}, \bmom^{\ast}, \bA^{n+1}, \bpA^{n+1}, \bJ^{n+1}, \bpJ^{n+1}, \tE^{\ast} \right)^{T}$.

\subsubsection{Convective terms}
Denoting $\Fluxst\left(\Q_{\ci}^{n},\Q_{\cj}^{n},\tnormalw{\ci\cj}\right)$ a numerical flux function, we approximate the flux terms in \eqref{eqn.spacetimedisc} as
\begin{equation}
	\sum_{\cj\in\mathcal{K}_{\ci}} \int\limits_{\tboundary{\ci\cj}}  \Fluxst\left(\Q\right) 
	 \cdot \tnormal{\ci\cj} \dt\dS 
	= \sum_{\cj\in\mathcal{K}_{\ci}}  \Fluxst\left(\Q_{\ci}^{n},\Q_{\cj}^{n},\tnormalw{\ci\cj}\right). 
\end{equation}
For instance, for the Rusanov numerical flux function \cite{Rus62}, we have 
\begin{equation}
	\Fluxst^{R}\left(\eQ^{n}_{\ci,R},\eQ^{n}_{\cj,L},\tnormalw{\ci\cj}\right) = \halb\left(\Fluxst\left(\eQ_{\ci,R}^{n}\right)+\Fluxst\left(\eQ_{\cj,L}^{n}\right)\right)\cdot \etast_{\ci\cj}
	-\halb \alpha_{\ci\cj}^{n} \left(\eQ_{\cj,L}^{n}-\eQ_{\ci,R}^{n}\right)
	\label{eqn.fluxRS}
\end{equation}
with the maximum signal speed on the edge
\begin{gather}
	\alpha_{\ci\cj}^{n} = 
	\max\left\lbrace
	\left|\overline{\bvel}_{\ci,R}^{n}\cdot \normalw{\ci\cj} \pm c_{\ci,R}^{n} \right|, 
	\left| \overline{\bvel}_{\cj,L}^{n}\cdot \normalw{\ci\cj}\pm c_{\cj,L}^{n}\right|,
	\left| \overline{\bvel}_{\ci,R}^{n}\cdot \normalw{\ci\cj}\pm c_{\bA, \ci}\right|,
	\left| \overline{\bvel}_{\cj,L}^{n}\cdot \normalw{\ci\cj}\pm c_{\bA, \cj}\right|, \right. 
	\notag\\\left.
	\left| \overline{\bvel}_{\ci,R}^{n}\cdot \normalw{\ci\cj}\pm c_{\bJ, \ci}\right|,
	\left| \overline{\bvel}_{\cj,L}^{n}\cdot \normalw{\ci\cj}\pm c_{\bJ, \cj}\right|
	\right\rbrace,\qquad
%	\notag\\
	c_{\ci,R}^{n} = \sqrt{\frac{4}{3} c_s^2 + \frac{2 c_{h}^{2}  \overline{T}_{\ci,R}^{n}}{\left( \overline{\rho}^{n}_{\ci,R}\right)^2 c_v}  } ,
	\label{eq:eig_gpr}
\end{gather}
approximated from the eigenvalues obtained for the unidimensional GPR model. To get a second order scheme, the half in time evolved and boundary extrapolated states, $\eQ^{n}_{\ci,R}$, $\eQ^{n}_{\cj,L}$, are computed following a local ADER approach \cite{TMN01,BTVC16,BFTVC17}:
\begin{enumerate}	
	\item Calculate the gradients in space, $\nabla \Q^{n}$, and time, $\nabla_{t} \Q^{n}_{\ci}$. To address the spacial gradient, the Crouzeix-Raviart basis functions on the primal elements are used, assuming the averaged values of $\Q^{n}$ approximate the values at the barycentres of the primal element faces. Second order is achieved by choosing the gradient to be the one of the primal element containing the dual face. For the time derivative, a Cauchy-Kovalevskaya procedure is applied.
	
	\item Employ the reference element in order to calculate the coordinates of the barycentre of the space-time face and compute the displacement with respect to the barycentres of cells  $\cell{\ci}^{n}$ and  $\cell{\cj}^{n}$:
	$$
	\dx_{\ci} =  \xst_{\ci\cj}(0.5,0.5) - \x_{\ci},\qquad
	\dx_{\cj} =  \xst_{\ci\cj}(0.5,0.5) - \x_{\cj}.
	$$
	
	\item Perform the half in time evolution and extrapolate the data at the boundary:  
	 $$
	 \eQ_{\ci,R} = \Q^{n}_{\ci} + \frac{\Delta t}{2} \nabla_{t} \Q^{n}_{\ci} + \nabla \Q^{n} \cdot \dx_{\ci}, \qquad
	 \eQ_{\cj,L} = \Q^{n}_{\cj} + \frac{\Delta t}{2} \nabla_{t} \Q^{n}_{\cj} + \nabla \Q^{n} \cdot \dx_{\cj}.
	 $$
\end{enumerate}
In the numerical results section, to circumvent Godunov's theorem and guarantee a stable scheme, the former ADER approach is modified including the ENO strategy \cite{Toro} or the min-mod limiter of Roe \cite{Roe81}. Further, we observe that, in most simulations, employing ENO limiting based on the physical variables instead of the conservative ones would lead to smoother solutions. Besides, adding a small artificial viscosity depending on a coefficient $c_{\alpha}$ can improve stability of the overall algorithm when small velocities are observed in comparison with the magnitude of the pressure field \cite{Hybrid2}. 

\subsubsection{Non-conservative terms: a path-conservative approach}
As aforementioned, the non-conservative products are approximated employing a path conservative scheme. Accordingly, we approximate
\begin{equation}
	\int\limits_{\tboundary{\ci}} \boldsymbol{\widetilde{\mathcal{D}}} \left(\Q^{n}\right) \cdot \tnormal{\ci} \dt\dS 
	= \sum_{\cj\in\mathcal{K}_{\ci}}  \boldsymbol{\widetilde{\mathcal{D}}}\left(\Q_{\ci}^{n},\Q_{\cj}^{n},\tnormalw{\ci\cj}\right).
\end{equation}
Considering the straight line segment path
\begin{equation*}
	\psi = \psi\left(\Q^{n}_{\ci},\Q^{n}_{\cj},s\right) = \Q^{n}_{\ci} + s \left( \Q^{n}_{\cj} - \Q^{n}_{\ci} \right), \qquad s\in [0,1],
\end{equation*}
we have
\begin{equation*}
	\boldsymbol{\widetilde{\mathcal{D}}}(\Q^{n}_{\ci},\Q^{n}_{\cj},\tnormalw{\ci\cj})  =\frac{1}{2}\Bst_{\ci\cj} \left(\Q^{n}_{\cj}-\Q^{n}_{\ci}\right), \qquad 
	\Bst_{\ci\cj} = \int_{0}^{1} \Bst\left(\psi\left(\Q^{n}_{\ci},\Q^{n}_{\cj},s\right)\right) \cdot\tnormalw{\ci\cj} \dS. 
\end{equation*}
On the other hand, we compute the half in time evolved conservative variables at each space-time control volume $\tcell{\ci}$ as 
\begin{equation}
	\overline{\Q}_{\ci}^{n} = \Q_{\ci}^{n} + \frac{\Dt}{4} \left(\nabla_{t} \Q^{n}_{\ci_{1}}  + \nabla_{t} \Q^{n}_{\ci_{2}}  \right),
\end{equation}
where the subindex $\ci_{1}$, $\ci_{2}$ refer to the two primal elements from which the dual cell is built, $\widetilde{T}_{\ci_{1}}$, $\widetilde{T}_{\ci_{2}}$. Then, using a weighted average for the computation of the gradient at the control volume, we have
\begin{equation}
	\int\limits_{\tcell{\ci}^{\circ}} \NCPst{\overline{\Q}} \dt\dV = 
		\left| \tcell{\ci} \right|
	\BNCP{\eQ_{\ci}^{n}} \left(\frac{\left|\cell{\ci_1}^{n}\right|}{\left|\cell{\ci}^{n}\right|} \gra \eQ_{\ci_{1}}^{n} +\frac{\left|\cell{\ci_2}^{n}\right|}{\left|\cell{\ci}^{n}\right|} \gra \eQ_{\ci_{2}}^{n}\right),
	\qquad \left| \tcell{\ci} \right|=\frac{1}{2} \left( \left|\cell{\ci}^{n}\right|+\left|\cell{\ci}^{n+1}\right|\right)  \Dt.
\end{equation}
Above, $\left|\cell{\ci_{1}}^{n}\right|$, $\left|\cell{\ci_{2}}^{n}\right|$ denote the areas of the two halves of the dual element $\cell{\ci}^{n}$ corresponding to  $\Pcell{\ci_{1}}^{n}$ and $\Pcell{\ci_{2}}^{n}$, respectively.

\subsubsection{Source terms}
The source terms of the momentum, \eqref{eqn.GPRGLM_disc_mom}, and total energy, \eqref {eqn.GPRGLM_disc_E}, equations are simply integrated on the control volume employing the half in time evolved density and momentum as
\begin{equation}
		\int\limits_{\tcell{\ci}^{\circ}} \rho_{\ci} \g \dt\dV = 	\left| \tcell{\ci} \right| \overline{\rho}_{\ci}^{n} \g,
		\qquad
		\int\limits_{\tcell{\ci}^{\circ}} \rho_{\ci} \g \cdot \bvel\dt\dV = 	\left| \tcell{\ci} \right|  \g \cdot \overline{\bmom}_{\ci}^{n}.
\end{equation}
As mentioned previously, in the visco-plastic limit of the model, the source terms related to the distortion tensor and thermal impulse can become very stiff. Therefore, \eqref{eqn.GPRGLM_A} and \eqref{eqn.GPRGLM_J} are split into two subsystems. The first subsystem accounts for the contributions of the flux term and the non-conservative products and is discretized as described in the previous sections. This yields intermediate values of the distortion and thermal impulse fields, 
and $\bJ^{\ast}_{\ci}$. The second subsystem, corresponding to \eqref{eqn.sourcesystem_sd}, is given by
\begin{subequations}\label{eqn.sourcesystem}
	\begin{align}
		% Distortion
		&\partial_{t} \A_{ik} = -\frac{1}{\theta_{1} \left(\tau_{1}\right)} \E_{\A_{ik}}, \\
		% thermal impulse
		& \partial_{t}  \J_{k}=  -\frac{1}{\theta_{2} \left(\tau_{2}\right)} \E_{\J_{k}}.
	\end{align}
\end{subequations}
Following \cite{HybridGPR}, this system of ordinary differential equations is locally solved at each dual cell using the Runge-Kutta DIRK scheme \cite{PR05} combined with an inexact Newton algorithm. As a consequence, we get the updated fields $\bA^{n+1}_{\ci}$, $\bJ^{n+1}_{\ci}$.

\subsection{Pressure stage: solution of the pressure subsystem} \label{sec:pressure}
To solve the pressure system, we first derive its weak formulation. Multiplication of \eqref{eqn.GPRGLM_disc_press} by a test function $z\in H^{1}_{0}= \left\lbrace z\in H^{1}(\Omega) \mid \int_{\Omega} z \dV = 0\right\rbrace,$ and integration over the computational domain gives
\begin{eqnarray}
	\frac{1}{ \gamma-1 }\int\limits_{\Omega}\delta \press^{n+1,\ell+1} z \dV
	-\Dt^2 \int\limits_{\Omega}\pxj \left(  \h^{n+1,\ell} \pxj \delta \press^{n+1,\ell+1}\right)z \dV
	= 
	-\Dt \int\limits_{\Omega}\pxj \left(  \h^{n+1,\ell} \mom_j^{\ast} \right) z\dV
	\notag\\ 
	+ \int\limits_{\Omega}\left( \tE^{\ast} 
	- \tE^{n}_{1}
	- \tE_{2}^{n+1,\ell}
	- \tE_{3}^{n+1}
	- \tE_{4}^{n+1}
	- \tE_{5}^{n+1}
	- \tE_{6}^{n+1}\right) z\dV.
\end{eqnarray}
Then, applying Green's theorem and taking into account the relation \eqref{eqn.GPRGLM_disc_monf}, we obtain the weak problem
\begin{eqnarray}
	\frac{1}{ \gamma-1 }\int\limits_{\Omega} \delta \press^{n+1,\ell+1} z \dV
	+\Dt^2 \int\limits_{\Omega}  \h^{n+1,\ell} \pxj \delta \press^{n+1,\ell+1}  \pxj z \dV
	=
	\Dt \int\limits_{\Omega} \h^{n+1,\ell} \mom_j^{\ast}  \pxj z \dV \notag\\
	+\int\limits_{\Omega} \left( \tE^{\ast} 
	-\tE^{n}_{1}
	-\tE_{2}^{n+1,\ell} 
	-\tE_{3}^{n+1}
	-\tE_{4}^{n+1}
	-\tE_{5}^{n+1}
	-\tE_{6}^{n+1}\right) z \dV
	- \Dt \int\limits_{\Gamma} \h^{n+1,\ell} \mom_j^{n+1}  n_{j} z \dS. \label{eqn.weakproblem}
\end{eqnarray}
This problem can be discretized using  $\mathbb{P}^{1}$ Lagrange finite element methods in the updated primal grid.
Let us remark that the use of Picard iterations increase one order in time per iteration. Hence, two Picard iterations are enough to preserve the accuracy of the overall scheme. 
Further, the resulting pressure system is symmetric and positive definite so a classical conjugate gradient method can be employed for its solution. 

We observe that the right hand side of \eqref{eqn.weakproblem} involves intermediate and updated values of several variables that have been computed in the dual grid. Therefore, previous to the pressure stage, we perform an intermediate stage to interpolate the needed data between the staggered grids. Moreover, $\tE_{2}^{n+1,\ell}$ is updated at each Picard iteration approximating the momentum according to \eqref{eqn.GPRGLM_disc_monf} as
\begin{equation}
	\bmom_{\ci}^{n+1,\ell} = \bmom_{\ci}^{\ast} - \Dt \frac{1}{\left|\Pcell{k}\right|}\sum\limits_{k \in\mathcal{\hat{K}}_{\ci}} 
	\left|\Pcell{k\ci}\right| \left( \gra \delta\press^{n+1,\ell}\right)_{k}
	\label{eqn.moupdate}
\end{equation}
with $ \left( \gra \press\right)_{k}$ the pressure gradient computed using the $\mathbb{P}^{1}$ basis functions in the primal element $\Pcell{k}$,
and computing
\begin{equation}
	\tE_{2}^{n+1,\ell} = \frac{1}{2\rho^{n+1}} \left|\bmom^{n+1,\ell}\right|^{2}.
\end{equation}
Finally, the enthalpy $ \h^{n+1,\ell}$ is calculated at each primal vertex by substituting \mbox{$\press^{n+1,\ell} = \press^{n} + \delta \press^{n+1,\ell} $} in \eqref{eqn.enthalpy_stiff}. The enthalpy at each primal face, which is also assumed to be the updated average value in the corresponding dual cell, is then computed as the average of its values on the vertex of the face.

\subsection{Correction stage: update of the intermediate velocity and total energy} \label{sec:correction}
Once the pressure correction $\delta \press^{n+1}$ is computed, we update the velocity and total energy according to \eqref{eqn.moupdate} and \eqref{eqn.GPRGLM_disc_Ef}. In particular, the energy correction is first performed by primal element
\begin{equation}
	\tE_{k}^{n+1} = \tE_{k} - \frac{\Delta t}{\left|\Pcell{k}\right|} \sum_{l\in\mathcal{K}_{k}} \int_{\Gamma_{kl}}  \h^{n+1}\bmom^{n+1} \cdot \normalu{kl} dS.\label{eq:WEupdatek}
\end{equation}
Next, the obtained values are interpolated to the dual cells as
\begin{equation}
	\tE^{n+1}_{\ci} = \frac{1}{\left|\cell{\ci}\right|}\sum_{k\in\mathcal{\hat{K}}_{i}} \left|\Pcell{ki}\right|\tE_{k}^{n+1}.
\end{equation}

\subsection{Intermediate stage: interpolation between grids} \label{sec:intermediate}
It is important to recall that once the transport stage has been performed, all variables live in the new mesh configuration.
Therefore, mesh interpolation between the dual and primal grids for ALE is carried out in the same way as for the purely Eulerian approach \cite{HybridGPR}. More precisely, given a scalar field at the dual cells, $q_{\ci}$, we approximate the solution at each primal element, $q_{k}$, as a weighted average of the form
\begin{equation}
	q_{k} = \sum_{\ci \in \mathcal{\hat{K}}_{k}} q_{\ci} \frac{\left|T_{k\ci}\right|}{\left| T_{k}\right|}. \label{eq:interpolation_dual2primal}
\end{equation}
This interpolation needs to be done for the density, $\rho^{n+1}$, and the intermediate momentum, $\bmom^{\ast}$, to compute the initial kinetic energy $\tE_{2}^{n+1,1}$. 
Further, we also apply \eqref{eq:interpolation_dual2primal} to get the intermediate total energy density, $\tE^{\ast}$, the total energy components $\tE_{3}^{n+1}$, $\tE_{4}^{n+1}$, $\tE_{5}^{n+1}$ and $\tE_{6}^{n+1}$, and the initial enthalpy which is first approximated at each dual cell as
\begin{equation}
\h^{n+1,1}_{\ci} = 
\dfrac{\gamma \press^{n}_{\ci} + \rho_{0}c_{0}^{2}-\gamma\press_{0}}{\rho^{n+1}_{\ci}\left( \gamma-1\right) }.
\end{equation}
Then, the related volume integrals in the right hand side of \eqref{eqn.weakproblem} are calculated assuming a constant value of the involved variables by primal element. Meanwhile, for the term depending on $\tE_{1}^{n}$, we employ the values of the pressure at each primal vertex computed in the previous iteration. Hence, its contribution is calculated using the mass matrix. The boundary integral in the right hand side of \eqref{eqn.weakproblem} is performed setting the enthalpy at each primal boundary to be the averaged value at the dual cell containing it.

% % % % % % % % % % % % % % % % % % % % % % % % % % % % % %
% % % % % % % % % % % % % % % % % % % % % % % % % % % % % %
%              Numerical Results
% % % % % % % % % % % % % % % % % % % % % % % % % % % % % %
% % % % % % % % % % % % % % % % % % % % % % % % % % % % % %
\section{Numerical results} \label{sec:numericalresults}

In this section, we analyse several test cases both in the fluid and solid limits of the GPR model to assess the proposed methodology. In what follows, the international system of units is employed and all test cases are initialized considering $\bA=\mathbf{I}$ and $\bJ=\boldsymbol{0}$. Further, unless stated the contrary, the time step is dynamically computed at each time iteration attending to a CFL stability condition based on the transport-diffusion subsystem, i.e.
\begin{equation*}
	\Delta t = \min_{\cell{\ci}} \left\lbrace \Delta t_{\ci}\right\rbrace, \qquad \Delta t_{\ci} =\frac{\mathrm{CFL}\, r_{\ci}}{\left|\lambda_{\ci}\right|_{\max}},
\end{equation*}
where $\left|\lambda_{\ci}\right|_{\max}$ is de maximum of the absolute approximated eigenvalues of the subsystem at cell $\cell{\ci}$ and $r_{\ci}$ is the incircle diameter of the cell. The eigenvalues are approximated from those calculated in the 1D case as
\begin{equation}
	\lambda \in \left\lbrace \left| \bvel \right|-c, \left| \bvel \right|, \left| \bvel \right|+c \right\rbrace,
	\quad
	c = \sqrt{\frac{4}{3} c_{s}^2 + 2 \frac{c_{h}^2\, T}{\rho^2\, c_{v}} },	
\end{equation}
for the original GPR model, \eqref{eqn.GPR}, and
\begin{equation}
	\lambda \in\left\lbrace \left| \bvel \right|-c_{\bJ}, \left| \bvel \right|-c_{\bA}, \left| \bvel \right|-c,\left| \bvel \right|, \left| \bvel \right|+c, \left| \bvel \right|+c_{\bA},\left| \bvel \right|+c_{\bJ} \right\rbrace,
\end{equation}
if we consider the augmented GLM curl cleaning GPR model, \eqref{eqn.GPR_rho}-\eqref{eqn.GPR_mom}, \eqref{eqn.GPRGLM_A}-\eqref{eqn.GPRGLM_PsiJ}, \eqref{eqn.GPR_S}.
Generally, for the second order scheme, we take $\mathrm{CFL}=0.5$. 

%For a discussion of the parallelization versatility of hybrid methods, we refer to \cite{HybridMPI,HybridGPR}.

\subsection{Convergence test and low Mach number asymptotic behaviour}
% TGV
As first test case, we consider the 2D Taylor Green vortex, a classical benchmark for assessing the order of convergence of numerical methods in fluid dynamics. The initial condition %in the computational $\Omega=[0,2\pi]^{2}$ 
reads
\begin{gather*}
	\bvel \left(\mathbf{x},0\right) = \left( \begin{array}{r} 
		\sin(x)\cos(y) \\
		-\cos(x)\sin(y) \end{array} \right), \quad 
	\press \left(\mathbf{x},0\right) = \frac{\press_{0}}{\gamma} + \frac{1}{4} \left(\cos(2x)+\cos(2y) \right).
\end{gather*}
Meanwhile, the model parameters, leading to an low Mach inviscid flow with $M\approx 3.2 \times 10^{-3}$, are $c_{s} = c_{h} = 0$, $\mu =\kappa = 0$, $c_{v}=2.5$, $c_{p}=3.5$, $p_{0}=10^{5}$. As computational domain we set $\Omega=[0,2\pi]^{2}$ with periodic boundary conditions everywhere. The errors and orders of convergence for a set of successively refined grids described in Table~\ref{tab:TGV_mesh} are reported in Table~\ref{tab:TGV_errors}. %We observe that the expected order of accuracy is reached.
\begin{table}[H]
	\renewcommand{\arraystretch}{1.2}
	\begin{center}
		\begin{tabular}{cccc}
			\hline 
			Mesh & Elements & Vertices & Dual elements \\\hline
			$M_1$ & $128 $ & $81 $ & $208 $ \\ %8
			$M_2$ & $512 $ & $289 $ & $800 $ \\ %16
			$M_3$ & $2048 $ & $1089 $ & $3136 $ \\ %32
			$M_4$ & $8192 $ & $4225 $ & $12416 $ \\ %64
			$M_5$ & $32768 $ & $16641 $ & $49408 $ \\ %128
			$M_6$ & $131072 $ & $66049 $ & $197120 $ \\ %256
			$M_7$ & $524288 $ & $263169 $ & $787456 $ \\ %512
			\hline 
		\end{tabular}
		\caption{2D Taylor-Green vortex. Main features of the primal triangular grids used to run the convergence table.} \label{tab:TGV_mesh}
	\end{center}
\end{table}

\begin{table}[h]
	\renewcommand{\arraystretch}{1.2}
	\begin{center}
		\begin{tabular}{c}
			\hline 
			Mesh  
			\\ \hline
			M1 \\
			M2 \\
			M3 \\
			M4  \\
			M5  \\
			M6  \\
			M7 \\
			\hline 
		\end{tabular}
		\begin{tabular}{cccccc}
			\hline 
			$L^{2}_{\Omega}\left(\rho\right)$ & $\mathcal{O}\left(\rho\right)$                  
			&$L^{2}_{\Omega}\left(\bmom\right)$ & $\mathcal{O}\left(\bmom\right)$ & $L^{2}_{\Omega}\left(\press\right)$ & $\mathcal{O}\left(\press\right)$ \\ \hline
			$2.43\cdot 10^{-2} $ & $     $ &$1.05\cdot 10^{-1} $ & $     $ & $5.03\cdot 10^{-1} $ & $ $ \\
			$3.51\cdot 10^{-3} $ & $ 2.79$ &$2.95\cdot 10^{-2} $ & $1.83 $ & $1.26\cdot 10^{-1} $ & $2.00 $ \\
			$5.15\cdot 10^{-4} $ & $ 2.77$ &$7.62\cdot 10^{-3} $ & $1.95 $ & $3.04\cdot 10^{-2} $ & $2.05 $ \\
			$9.16\cdot 10^{-5} $ & $ 2.49$ &$1.90\cdot 10^{-3} $ & $2.00 $ & $7.56\cdot 10^{-3} $ & $2.01 $ \\
			$1.96\cdot 10^{-5} $ & $ 2.23$ &$4.75\cdot 10^{-4} $ & $2.00 $ & $1.89\cdot 10^{-3} $ & $2.00 $ \\
			$4.50\cdot 10^{-6} $ & $ 2.12$ &$1.19\cdot 10^{-4} $ & $2.00 $ & $1.02\cdot 10^{-3} $ & $0.90 $ \\
			$1.68\cdot 10^{-6} $ & $ 1.42$ &$2.97\cdot 10^{-5} $ & $2.00 $ & $1.21\cdot 10^{-4} $ & $3.07 $ \\
			\hline 
		\end{tabular}
		\caption{Taylor-Green vortex. Spatial $L_{2}$ error norms and convergence rates at time $t=0.1$.} \label{tab:TGV_errors}
	\end{center}
\end{table}

This test case is also run for a set of decreasing Mach numbers to analyse the asymptotic behaviour of the proposed methodology in the low Mach number limit. Table~\ref{tab:TGV_Mach} confirms that the order of convergence is preserved for very low Mach number regimes.  Let us note that the limit scheme for the incompressible Navier-Stokes equations in which the proposed hybrid method is based \cite{BFTVC17}, is not exactly divergence free. Hence, when $M\rightarrow 0$ the errors are not expected to decrease quadratically but to be preserved. Due to the magnitude of the pressure variable, for very low Mach numbers quadruple precision is required.

\begin{table}[H]
	\renewcommand{\arraystretch}{1.2}
	\setlength{\tabcolsep}{3.8pt}
	\begin{center}
		\begin{tabular}{ccccccccc}
			\hline 
			Mesh & $L^{2}_{\Omega}\left(\rho\right)$ & $\mathcal{O}\left(\rho\right)$     & $L^{2}_{\Omega}\left(\bmom\right)$ & $\mathcal{O}\left(\bmom\right)$                 
			&$L^{2}_{\Omega}\left(\tE\right)$ & $\mathcal{O}\left(\tE\right)$ & $L^{2}_{\Omega}\left(\press\right)$ & $\mathcal{O}\left(\press\right)$ \\ \hline
			\multicolumn{9}{c}{$M=10^{-2}$, $p_{0}=10^{4}$ {\scriptsize (double precision)}}\\ \hline
			M1 & $2.43\cdot 10^{-2} $ & $  $ & $1.05\cdot 10^{-1}  $ & $ $ &$8.57\cdot 10^{1} $ & $ $ & $5.07\cdot 10^{-1}  $ & $ $ \\
			M2 & $3.50\cdot 10^{-3} $ & $2.79  $ & $2.95\cdot 10^{-2}  $ & $1.83 $ &$5.59\cdot 10^{0} $ & $3.94 $ & $1.26\cdot 10^{-1}  $ & $2.01 $ \\
			M3 & $5.09\cdot 10^{-4} $ & $2.78  $ & $7.62\cdot 10^{-3}  $ & $1.95 $ &$3.99\cdot 10^{-1} $ & $3.81 $ & $3.04\cdot 10^{-2}  $ & $2.05 $ \\
			M4 & $8.97\cdot 10^{-5} $ & $2.50  $ & $1.91\cdot 10^{-3}  $ & $2.00 $ &$3.04\cdot 10^{-2} $ & $3.71 $ & $7.55\cdot 10^{-3}  $ & $2.01 $ \\
			\hline
			\multicolumn{9}{c}{$M=10^{-3}$, $p_{0}=10^{6}$ {\scriptsize (double precision)}}\\ \hline
			M1 & $2.43\cdot 10^{-2} $ & $  $ & $1.05\cdot 10^{-1}  $ & $ $ &$8.53\cdot 10^{3} $ & $ $ & $5.03\cdot 10^{-1}  $ & $ $ \\
			M2 & $3.51\cdot 10^{-3} $ & $2.79  $ & $2.95\cdot 10^{-2}  $ & $1.83 $ &$5.51\cdot 10^{2} $ & $3.95 $ & $1.26\cdot 10^{-1}  $ & $2.00 $ \\
			M3 & $5.15\cdot 10^{-4} $ & $2.77  $ & $7.62\cdot 10^{-3}  $ & $1.95 $ &$3.83\cdot 10^{1} $ & $3.85 $ & $3.04\cdot 10^{-2}  $ & $2.05 $ \\
			M4 & $9.16\cdot 10^{-5} $ & $2.49  $ & $1.90\cdot 10^{-3}  $ & $2.00 $ &$2.62\cdot 10^{0} $ & $3.87 $ & $7.56\cdot 10^{-3}  $ & $2.01 $ \\
			\hline
			\multicolumn{9}{c}{$M=10^{-4}$, $p_{0}=10^{8}$ {\scriptsize (quadruple precision)}}\\ \hline
			M1 & $2.43\cdot 10^{-2} $ & $  $ & $1.05\cdot 10^{-1}  $ & $ $ &$8.53\cdot 10^{5} $ & $ $ & $5.03\cdot 10^{-1}  $ & $ $ \\
			M2 & $3.51\cdot 10^{-3} $ & $2.79  $ & $2.95\cdot 10^{-2}  $ & $1.83 $ &$5.51\cdot 10^{4} $ & $3.95 $ & $1.26\cdot 10^{-1}  $ & $2.00 $ \\
			M3 & $5.15\cdot 10^{-4} $ & $2.77  $ & $7.62\cdot 10^{-3}  $ & $1.95 $ &$3.83\cdot 10^{3} $ & $3.85 $ & $3.04\cdot 10^{-2}  $ & $2.05 $ \\
			M4 & $9.16\cdot 10^{-5} $ & $2.49  $ & $1.90\cdot 10^{-3}  $ & $2.00 $ &$2.62\cdot 10^{2} $ & $3.87 $ & $7.56\cdot 10^{-3}  $ & $2.01 $ \\
			\hline
			\multicolumn{9}{c}{$M=10^{-5}$, $p_{0}=10^{10}$ {\scriptsize (quadruple precision)}}\\ \hline
			M1 & $2.43\cdot 10^{-2} $ & $  $ & $1.05\cdot 10^{-1} $ & $    $ & $8.53\cdot 10^{7}  $ & $    $ & $5.03\cdot 10^{-1}  $ & $ $ \\
			M2 & $3.51\cdot 10^{-3} $ & $2.79$ & $2.95\cdot 10^{-2} $ & $1.83$ & $5.51\cdot 10^{6}  $ & $3.95$ & $1.26\cdot 10^{-1}  $ & $2.00$ \\
			M3 & $5.15\cdot 10^{-4} $ & $2.77$ & $7.62\cdot 10^{-3} $ & $1.95$ & $3.83\cdot 10^{5}  $ & $3.85$ & $3.04\cdot 10^{-2}  $ & $2.05$ \\
			M4 & $9.58\cdot 10^{-5} $ & $2.43$ & $1.91\cdot 10^{-3} $ & $2.00$ & $5.72\cdot 10^{4}  $ & $2.74$ & $7.56\cdot 10^{-3}  $ & $2.01$ \\
			\hline 
		\end{tabular}
		\caption{2D Taylor-Green vortex. Spatial $L_{2}$ error norms and convergence rates of the density, momentum, total energy and pressure at time $t=0.1$ for $M\in\left\lbrace 10^{-2}, 10^{-3}, 10^{-4}, 10^{-5} \right\rbrace$. Results computed using the second order approach.} \label{tab:TGV_Mach}
	\end{center}
\end{table}

% Isentropic vortex 

To analyse also the behaviour of the ALE hybrid approach, we consider the isentropic vortex and let the mesh freely move according to the fluid velocity \cite{HybridALE}. The solution is given by a rotating vortex defined as
\begin{equation}
	\rho\left(\x,0\right) = 1,\quad 
	\bvel \left(\mathbf{x},0\right) = \left( \begin{array}{r} 
		-r e^{-\frac{1}{2}\left(r^2-1 \right)}\sin \varphi  \\
		r e^{-\frac{1}{2}\left(r^2-1 \right)}\cos \varphi 
		\end{array} \right), \quad 
	\press \left(\x,0\right) = \frac{\press_{0}}{\gamma} - \frac{1}{2} e^{-r^2+1},
	\quad \varphi = \arctan \left( y-5,x-5\right),
\end{equation}
with $ r=\sqrt{(x-5)^2+(y-5)^2}$ the radius to the centre of the computational domain $\Omega=[0,10]^{2}$. Since this benchmark studies an inviscid incompressible fluid, the parameters for the GPR model are taken as $c_{h}=c_{s}=\mu=\kappa=0$, $c_{v}=2.5$, $c_{\press}=3.5$. A set of reference pressures is considered to analyse the behaviour of the method for different Mach numbers, $M\in\left\lbrace 10^{-2}, 10^{-3}, 10^{-4}, 10^{-5} \right\rbrace$. Dirichlet boundary conditions are imposed everywhere. The errors and convergence rates obtained at time $t_{\mathrm{e}}=0.1$ are reported in Table~\ref{tab.IV_Mach}. The expected second order is achieved for all studied Mach numbers. 

\begin{table}[H]
	\renewcommand{\arraystretch}{1.2}
	\setlength{\tabcolsep}{3.8pt}
	\begin{center}
		\begin{tabular}{ccccccccc}
			\hline 
			Mesh & $L^{2}_{\Omega}\left(\rho\right)$ & $\mathcal{O}\left(\rho\right)$     &                
			$L^{2}_{\Omega}\left(\rho\bvel\right)$ & $\mathcal{O}\left(\rho\bvel\right)$ & 
			$L^{2}_{\Omega}\left(\tE\right)$ & $\mathcal{O}\left(\tE\right)$ &
			$L^{2}_{\Omega}\left(\press\right)$ & $\mathcal{O}\left(\press\right)$ \\
			\hline
			\multicolumn{7}{c}{$M=10^{-2}$, $\press_{0}=10^{4}$ {\scriptsize (double precision)}}\\ \hline
			M2 & $1.25\cdot 10^{-2} $ & $     $ & $6.58 \cdot 10^{-2}  $ & $     $ & $5.58\cdot 10^{1} $ & $ $ & $2.54\cdot 10^{-1} $ & $ $ \\
			M3 & $2.05\cdot 10^{-3} $ & $2.61 $ & $1.54 \cdot 10^{-2}  $ & $2.10 $ & $2.72\cdot 10^{0} $ & $4.36 $ & $6.03\cdot 10^{-2} $ & $2.07 $ \\
			M4 & $3.40\cdot 10^{-4} $ & $2.59 $ & $3.73 \cdot 10^{-3}  $ & $2.05 $ & $1.80\cdot 10^{-1} $ & $3.92 $ & $1.50\cdot 10^{-2} $ & $2.01 $ \\
			M5 & $6.03\cdot 10^{-5} $ & $2.50 $ & $9.27 \cdot 10^{-4}  $ & $2.01 $ & $1.55 \cdot 10^{-2} $ & $3.54 $ & $3.63\cdot 10^{-3} $ & $2.05 $ \\
			\hline
			\multicolumn{7}{c}{$M=10^{-3}$, $\press_{0}=10^{6}$ {\scriptsize (double precision)}}\\ \hline
			M2 & $1.26\cdot 10^{-2} $ & $     $ & $6.59\cdot 10^{-2} $ & $     $ &$5.58\cdot 10^{3} $ & $ $ & $2.63\cdot 10^{-1} $ & $ $ \\
			M3 & $2.07\cdot 10^{-3} $ & $2.60 $ & $1.54\cdot 10^{-2} $ & $2.09 $ &$2.70\cdot 10^{2} $ & $4.37 $ & $5.91\cdot 10^{-2} $ & $2.15 $ \\
			M4 & $3.47\cdot 10^{-4} $ & $2.58 $ & $3.74\cdot 10^{-3} $ & $2.04 $ &$1.74\cdot 10^{1} $ & $3.95 $ & $1.49\cdot 10^{-2} $ & $1.99 $ \\
			M5 & $6.17\cdot 10^{-5} $ & $2.49 $ & $9.29\cdot 10^{-4} $ & $2.01 $ &$1.32 \cdot 10^{0} $ & $3.73 $ & $3.78\cdot 10^{-3} $ & $1.98 $ \\
			\hline
			\multicolumn{7}{c}{$M=10^{-4}$, $\press_{0}=10^{8}$ {\scriptsize (quadruple precision)}}\\ \hline
			M2 & $1.26\cdot 10^{-2} $ & $      $ & $6.59\cdot 10^{-2}  $ & $     $ & $5.58\cdot 10^{5} $ & $ $ & $5.63\cdot 10^{0} $ & $ $ \\
			M3 & $2.07\cdot 10^{-3} $ & $2.60  $ & $1.54\cdot 10^{-2}  $ & $2.09 $ & $2.70\cdot 10^{4} $ & $4.37 $ & $6.17\cdot 10^{-2} $ & $6.51 $ \\
			M4 & $3.47\cdot 10^{-4} $ & $2.58  $ & $3.74\cdot 10^{-3}  $ & $2.04 $ & $1.74\cdot 10^{3} $ & $3.96 $ & $1.50\cdot 10^{-2} $ & $2.04 $ \\
			M5 & $6.17\cdot 10^{-5} $ & $2.49  $ & $9.30\cdot 10^{-4}  $ & $2.01 $ & $1.31\cdot 10^{2} $ & $3.73 $ & $4.18\cdot 10^{-3} $ & $1.85 $ \\
			\hline
			\multicolumn{7}{c}{$M=10^{-5}$, $\press_{0}=10^{10}$ {\scriptsize (quadruple precision)}}\\ \hline
			M2 & $1.26\cdot 10^{-2} $ & $     $ & $6.59\cdot 10^{-2}  $ & $     $ & $5.58\cdot 10^{7} $ & $ $ & $5.63\cdot 10^{2}  $ & $ $ \\
			M3 & $2.07\cdot 10^{-3} $ & $2.60 $ & $1.54\cdot 10^{-2}  $ & $2.09 $ & $2.70\cdot 10^{6} $ & $4.37 $ & $1.80\cdot 10^{0}  $ & $8.29 $ \\
			M4 & $3.47\cdot 10^{-4} $ & $2.58 $ & $3.74\cdot 10^{-3}  $ & $2.04 $ & $1.74\cdot 10^{5} $ & $3.96 $ & $2.04\cdot 10^{-1} $ & $3.14 $ \\
			M5 & $6.17\cdot 10^{-5} $ & $2.49 $ & $9.30\cdot 10^{-4}  $ & $2.01 $ & $1.31\cdot 10^{4} $ & $3.73 $ & $2.35\cdot 10^{-1} $ & $-0.02 $ \\
			\hline 
		\end{tabular}
		\caption{2D isentropic vortex. Spatial $L_{2}$ error norms and convergence rates of the density, momentum, total energy and pressure at time $t=0.1$ for $M\in\left\lbrace 10^{-2}, 10^{-3}, 10^{-4}, 10^{-5} \right\rbrace$. Results computed using the second order ALE hybrid method.} \label{tab.IV_Mach}
	\end{center}
\end{table}

\subsection{Riemann problems}
We now analyse a set of Riemann problems both in the fluid and solid limits of the model. In particular, we consider the computational domain $\Omega=[-0.5,0.5]\times[-0.05,0.05]$ and define the initial condition for the primitive variables $\mathbf{V}=\left(\rho,\bvel,\bA,\bJ,\press\right)$ as
\begin{equation*}
	\mathbf{V}\left(\mathbf{x},0\right) = \left\lbrace \begin{array}{lc}
		\mathbf{V}^{L} & \mathrm{ if } \; x \le  x_{d},\\
		\mathbf{V}^{R} & \mathrm{ if } \; x>  x_{d},
	\end{array}\right. 
\end{equation*}
with the remaining left and right states, $\rho^{L}$, $\bvel^{L}$, $\press^{L}$, $\rho^{R}$, $\bvel^{R}$, $\press^{R}$, the initial position of the discontinuity, $x_{d}$, the final simulation time, $t_{\mathrm{e}}$, and the number of mesh divisions along the $x$-axis given in Table~\ref{tab.RPIC}. The material parameters for each test case are provided in Table~\ref{tab.RPparameters}.
\begin{table}[H]
	\renewcommand{\arraystretch}{1.2}
	\begin{center}
		\begin{tabular}{ccccccccccc}
			Test &  $\rho^{L}$ &  $\rho^{R}$  &  $u_{1}^{L}$ &  $u_{1}^{R}$ &  $u_{2}^{L}$ &  $u_{2}^{R}$ &  $\press^{L}$ &  $\press^{R}$ & $x_{c}$ & $t_{\mathrm{e}}$ \\ \hline
			RP1 & $ 1 $ & $ 0.125 $ & $ 0 $ & $ 0 $ & $ 0 $ & $ 0 $ & $ 1 $ & $ 0.1 $& $ 0 $ &  $ 0.2 $  \\  
			RP2 & $ 1 $ & $ 1 $ & $ -1 $ & $ 1 $& $ 0 $ & $ 0 $  & $ 0.4 $ & $ 0.4 $& $ 0 $ &  $ 0.15$ \\
			RP3 & $ 1 $ & $ 0.125 $ & $ 0.5 $ & $ 0 $ & $ 0 $ & $ 0 $ & $ 1 $ & $ 1 $& $ 0 $ &  $ 0.1 $  \\
			RP4 & $ 5.99924 $ & $ 5.99242 $ & $ 19.5975 $ & $ -6.19633 $ & $ 0 $ & $ 0 $ & $ 460.894 $ & $ 46.095 $& $ -0.2 $ &  $ 0.035 $  \\
			RP5 &  $ 1 $ & $ 1 $ & $ -19.59745 $ & $ -19.59745 $ & $ 0 $ & $ 0 $ & $ 1000 $ & $ 0.01 $& $ 0.3 $ &  $ 0.01 $  \\				
			RP6 &  $ 1 $ & $ 1 $ & $ 2 $ & $ -2 $ & $ 0 $ & $ 0 $ & $ 0.1 $ & $ 0.1 $& $ 0 $ &  $ 0.8 $ \\
			RP7 &  $ 1 $ & $ 0.125 $ & $ 0.5 $ & $ 0 $ & $ 0 $ & $ 0 $ & $ 1 $ & $ 1 $& $ 0 $ &  $ 0.1 $  \\
			RP$8$--$9$ & $ 1 $ & $ 0.5 $ & $ 0 $ & $ 0 $ & $ -0.2 $ & $ 0.2 $ & $ 1 $ & $ 0.5 $& $ 0 $ &  $ 0.2 $  \\
%			RP8 &  $ 1 $ & $ 0.5 $ & $ 0 $ & $ 0 $ & $ -0.2 $ & $ 0.2 $ & $ 1 $ & $ 0.5 $& $ 0 $ &  $ 0.2 $  \\
%			RP9 &  $ 1 $ & $ 0.5 $ & $ 0 $ & $ 0 $ & $ -0.2 $ & $ 0.2 $ & $ 1 $ & $ 0.5 $& $ 0 $ &  $ 0.2 $  \\
		\end{tabular} 
	\end{center}
	\caption{Riemann problems. Initial condition, location of the initial discontinuity, $x_{d}$, and final time, $t_{\mathrm{end}}$.
		%, and number of mesh divisions on $x$-direction, $N_{x}$.
	}
	\label{tab.RPIC}
\end{table}
\begin{table}[H]
	\renewcommand{\arraystretch}{1.2}
	\begin{center}
		\begin{tabular}{ccccccc}
			Test &  $c_{s}$ &  $c_{h}$ &  $c_{v}$ &  $c_{p}$  &  $\mu$ &  $\kappa$ \\ \hline
			RP1-7 & $ 0 $ & $ 0 $ & $ 2.5 $ & $ 3.5 $ & $ 0 $ & $ 0 $\\  
%			RP8 & $ 1 $ & $ 1 $ & $ 2.5 $ & $ 3.5 $ & $ 10^{-5} $ & $ 10^{-5} $\\
			RP8 & $ 1 $ & $ 0 $ & $ 1 $ & $ 1.4 $ & $ 10^{20} $ & $ 10^{20} $\\
		    RP9 & $ 1 $ & $ 1 $ & $ 1 $ & $ 1.4 $ & $ 10^{20} $ & $ 10^{20} $\\
		\end{tabular} 
	\end{center}
	\caption{Riemann problems. Material parameters for the GPR model.}
	\label{tab.RPparameters}
\end{table}

\begin{figure}
	\centering
	\includegraphics[width=0.45\linewidth]{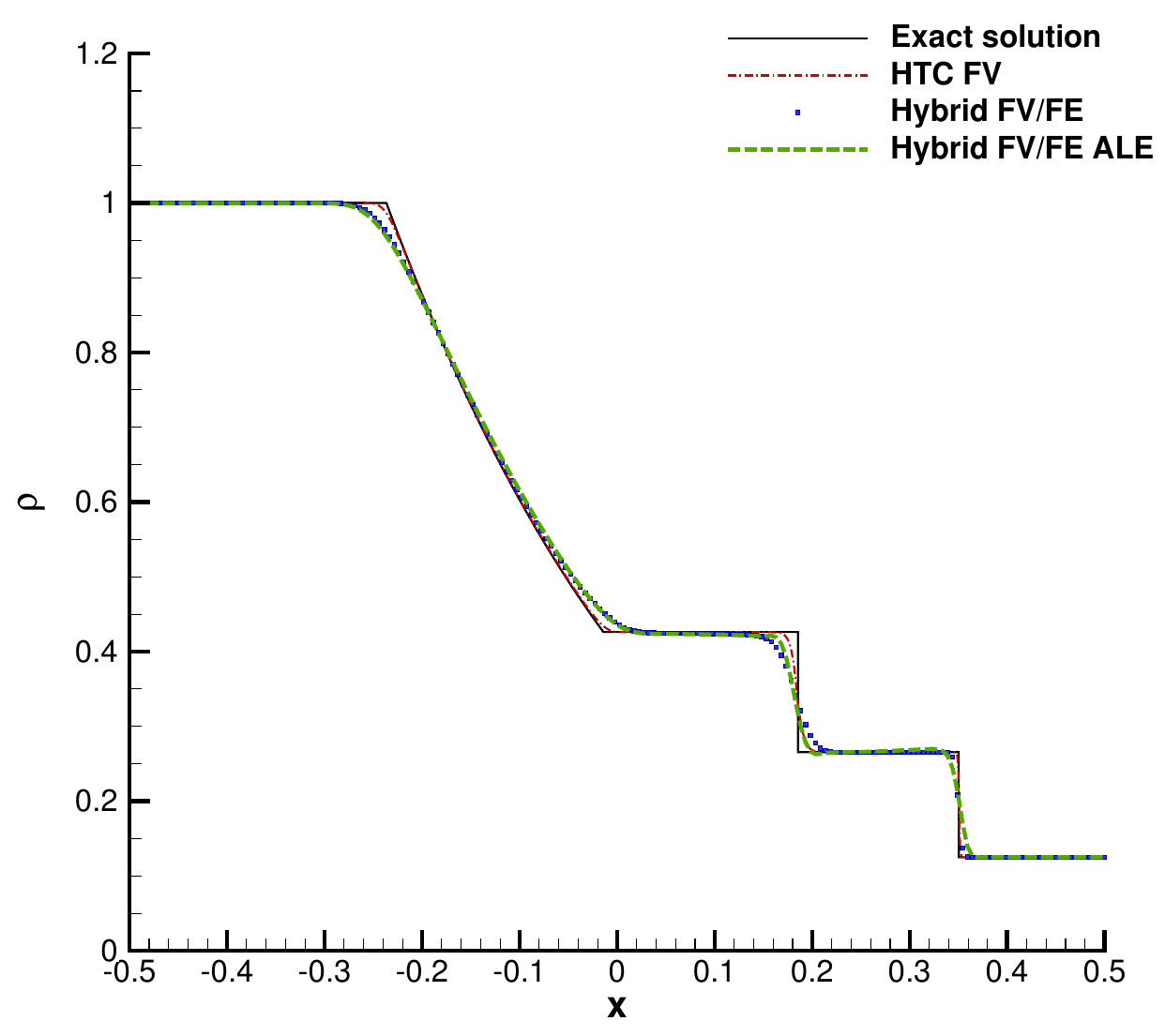}
	\includegraphics[width=0.45\linewidth]{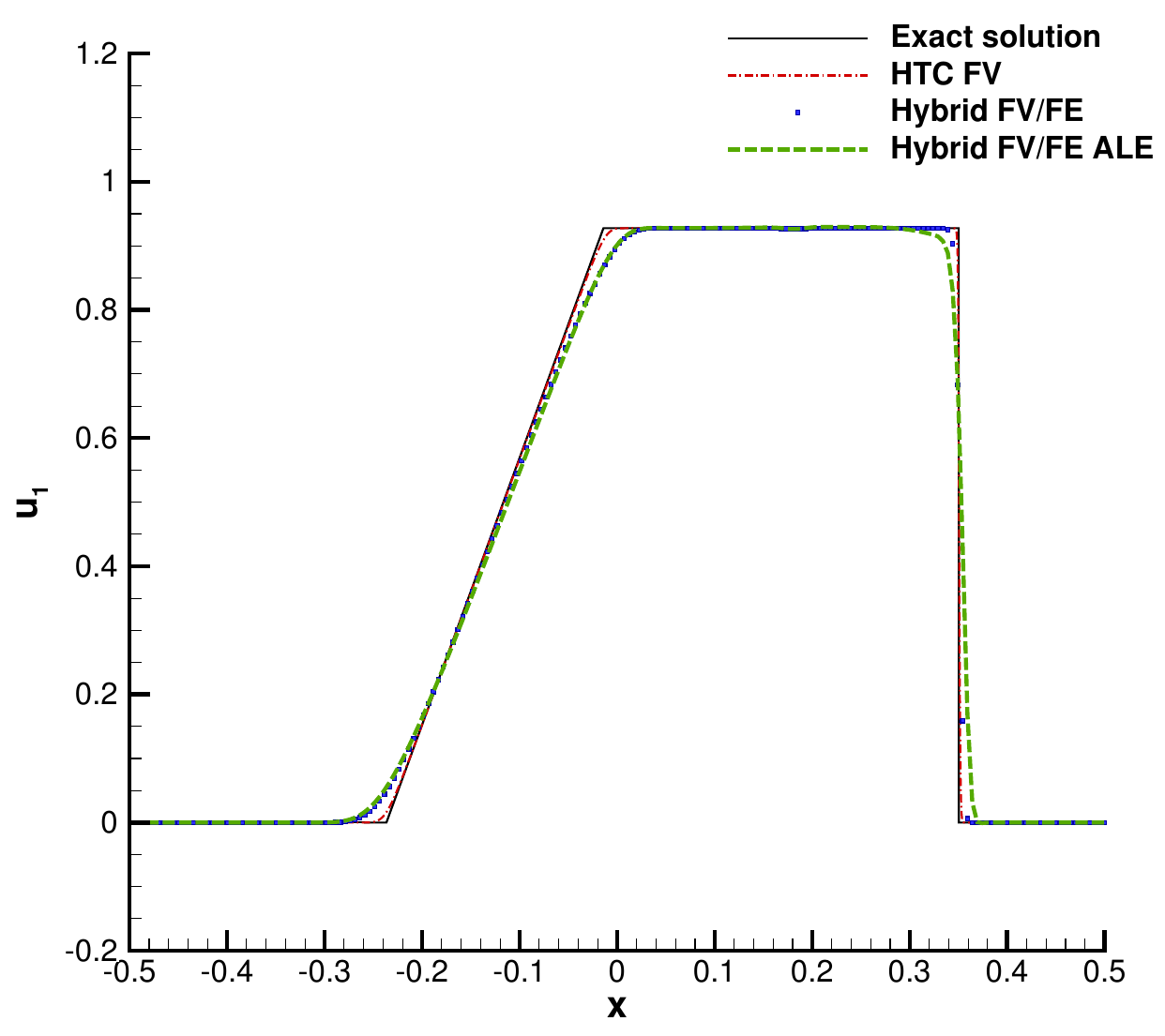}\\
	\includegraphics[width=0.45\linewidth]{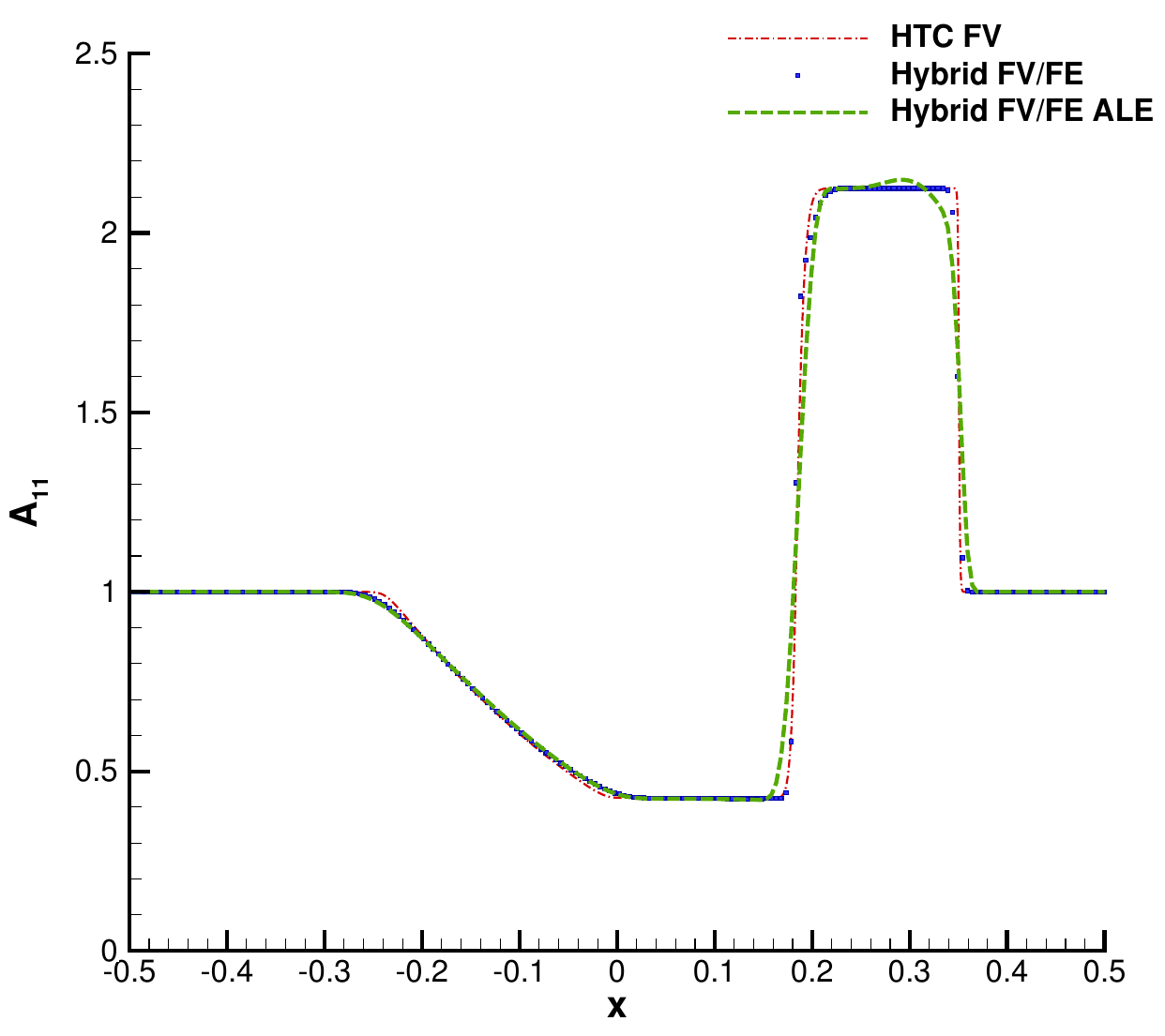}
	\includegraphics[width=0.45\linewidth]{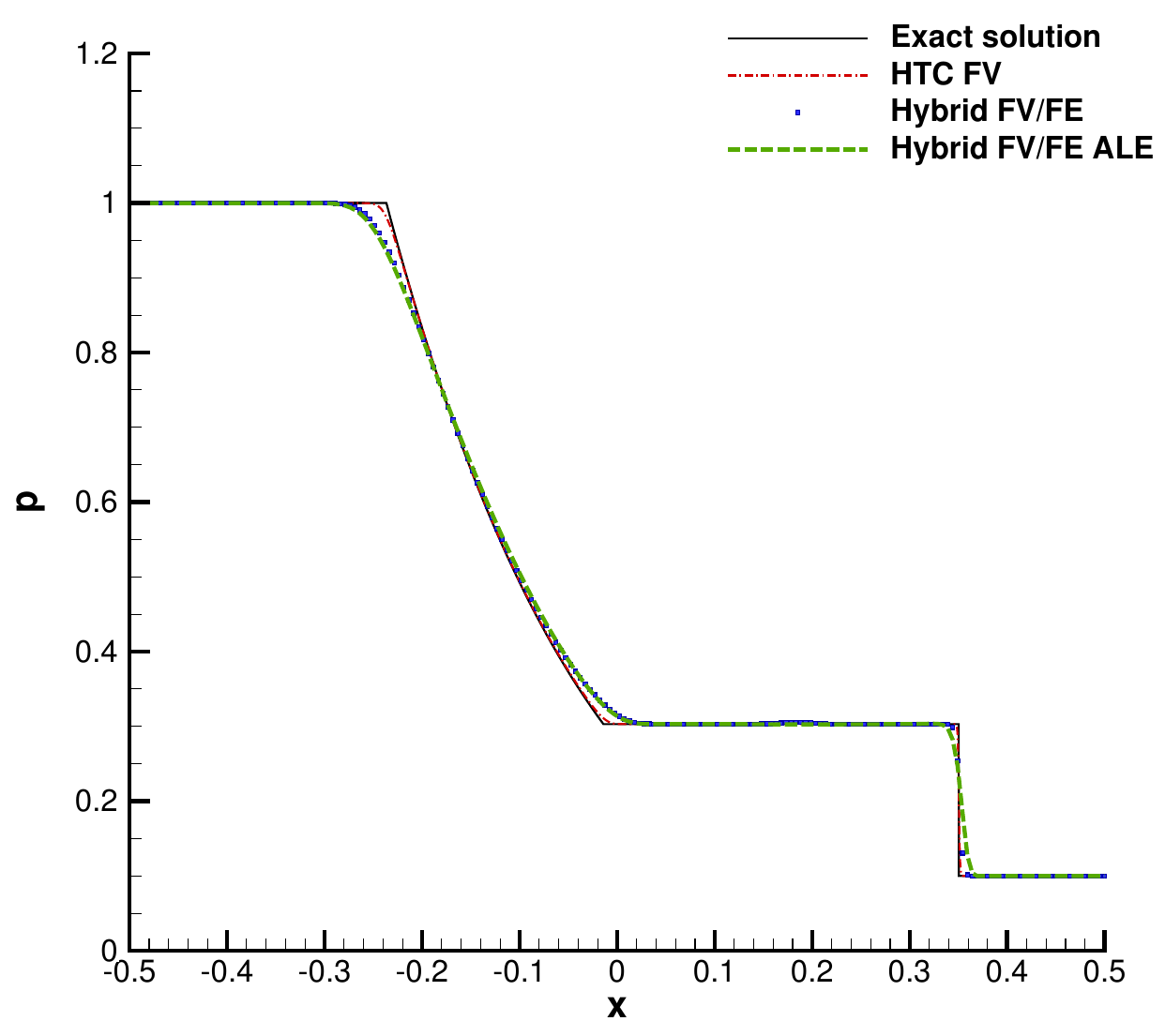}
	\caption{RP1 Sod. 1D cuts of the solution obtained using the hybrid FV/FE method for the compressible GPR model with the fully Eulerian code (blue squares) and the ALE scheme (green dashed line) compared against the exact solution for the Navier-Stokes equations (black line) and the solution of the GPR model obtained using the HTC FV scheme in \cite{HTCA2} (red dashed line). From left-top to right-bottom: density, horizontal velocity component, distortion field component $A_{11}$, and pressure.}
	\label{fig:rp1SodALE}
\end{figure}

\begin{figure}
	\centering
	\includegraphics[width=0.32\linewidth]{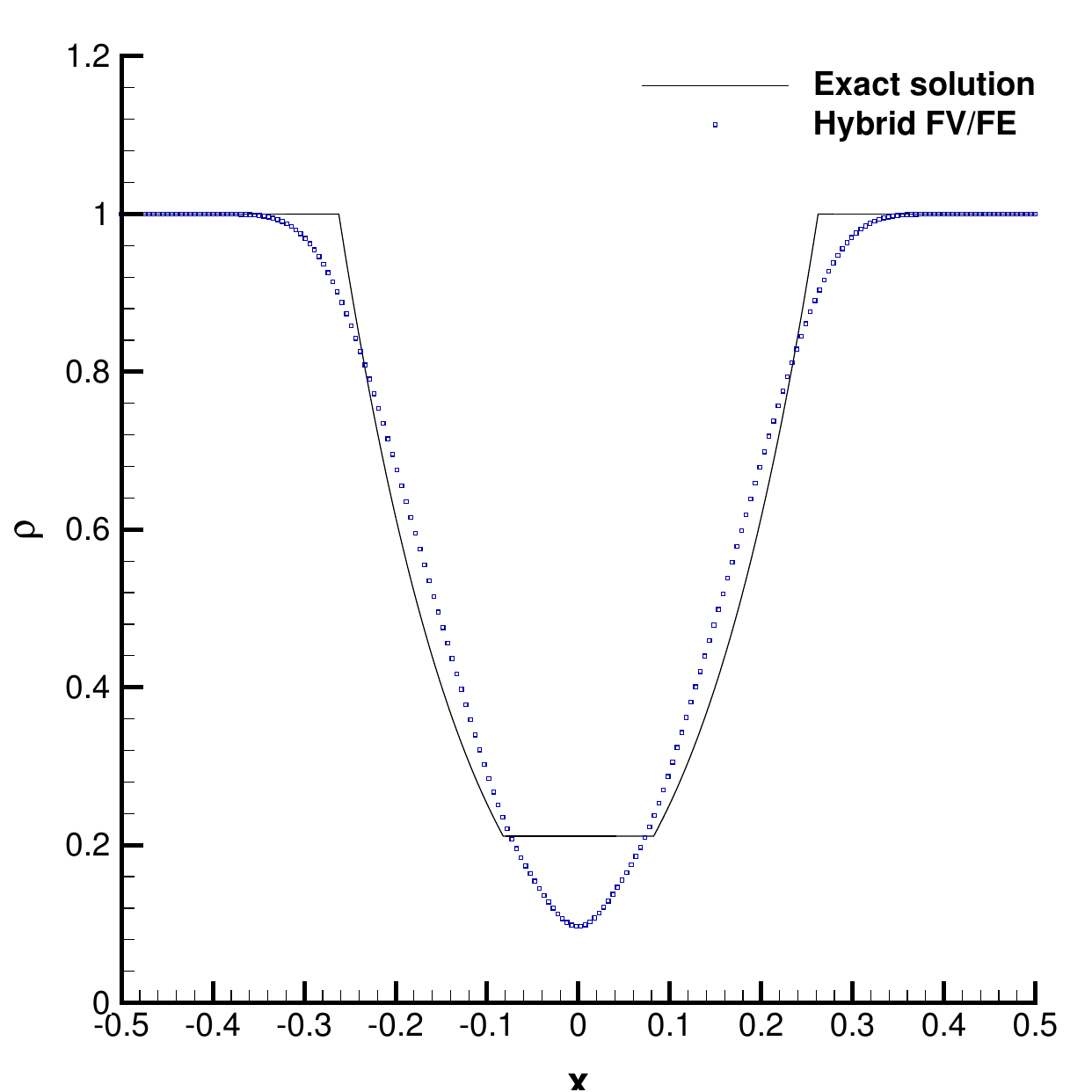}
	\includegraphics[width=0.32\linewidth]{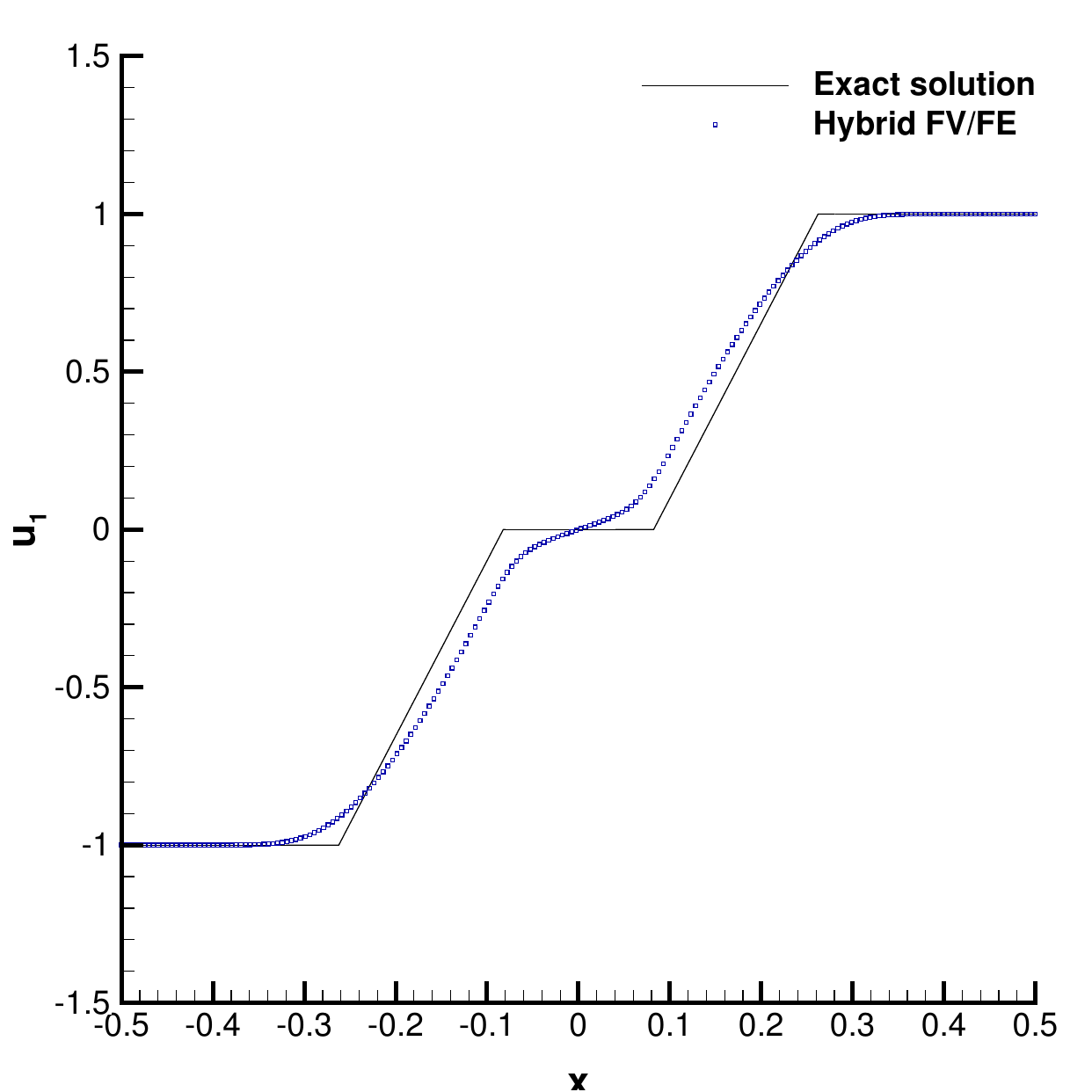}
	\includegraphics[width=0.32\linewidth]{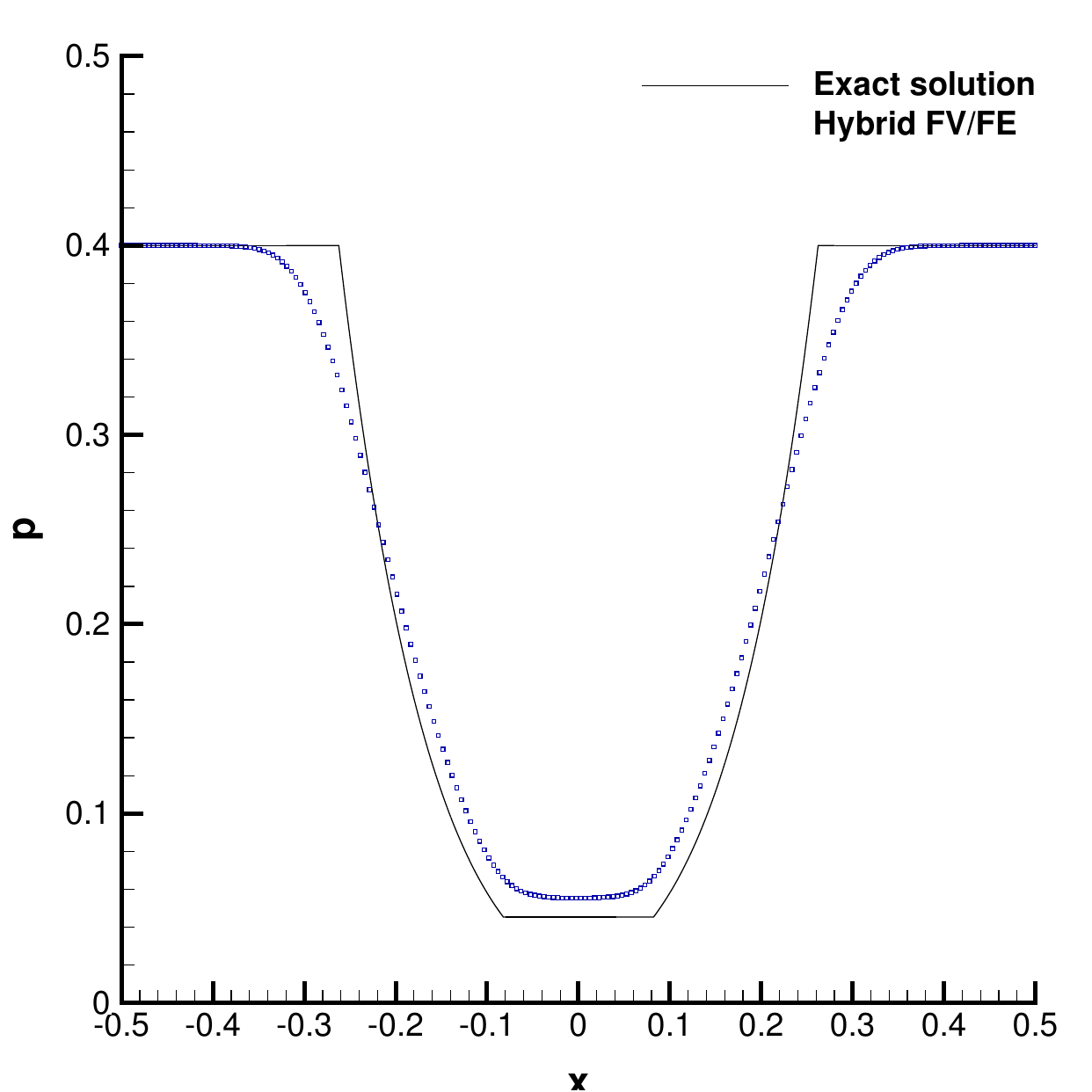}
	\caption{RP2 smooth double rarefaction. 1D cuts of the solution obtained using the hybrid FV/FE method for the compressible GPR model (blue squares) compared against the exact solution (black line). From left to right: density, horizontal velocity component and pressure.}
	\label{fig:rp2DR}
\end{figure}

\begin{figure}
	\centering
	\includegraphics[width=0.32\linewidth]{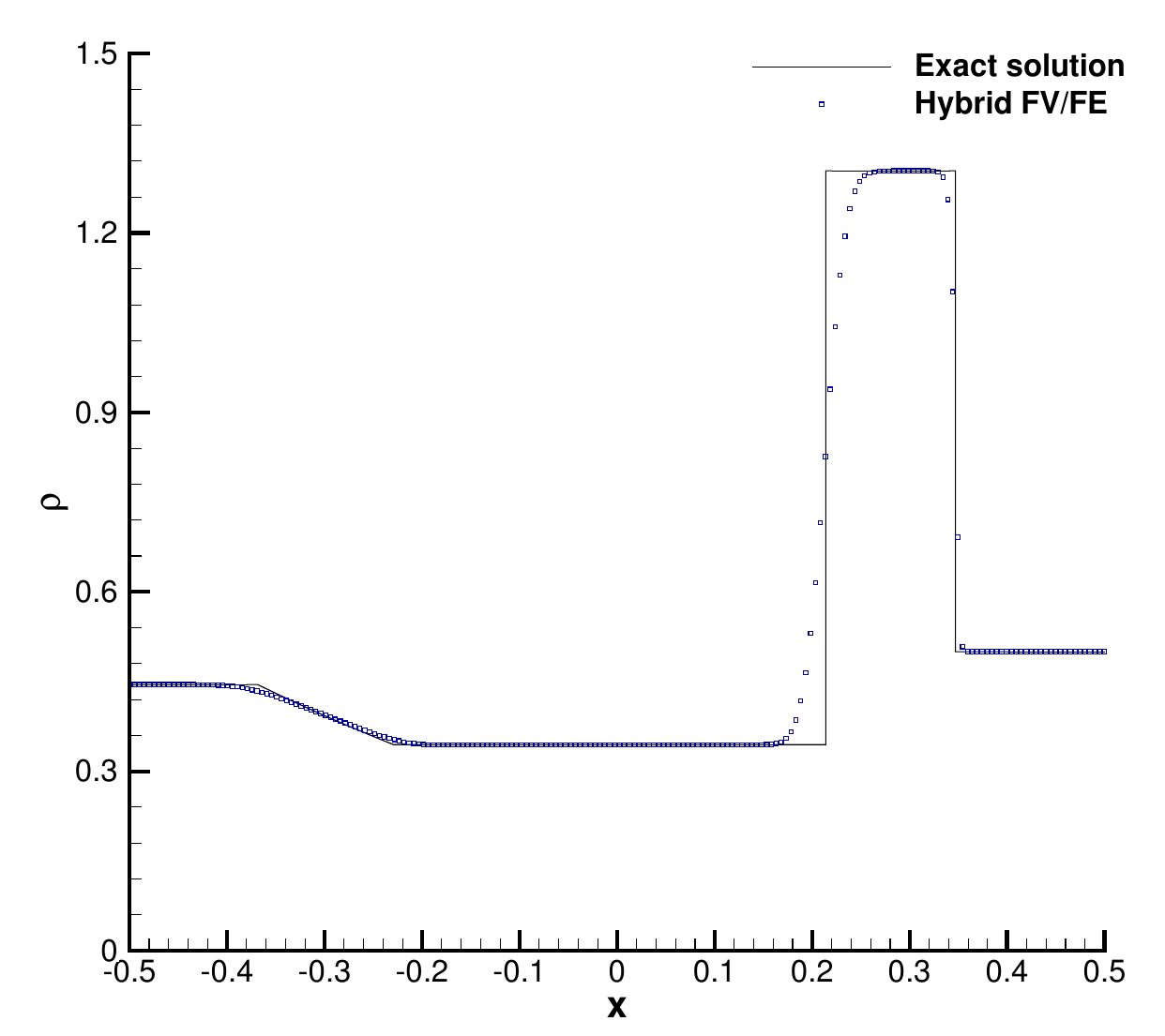}
	\includegraphics[width=0.32\linewidth]{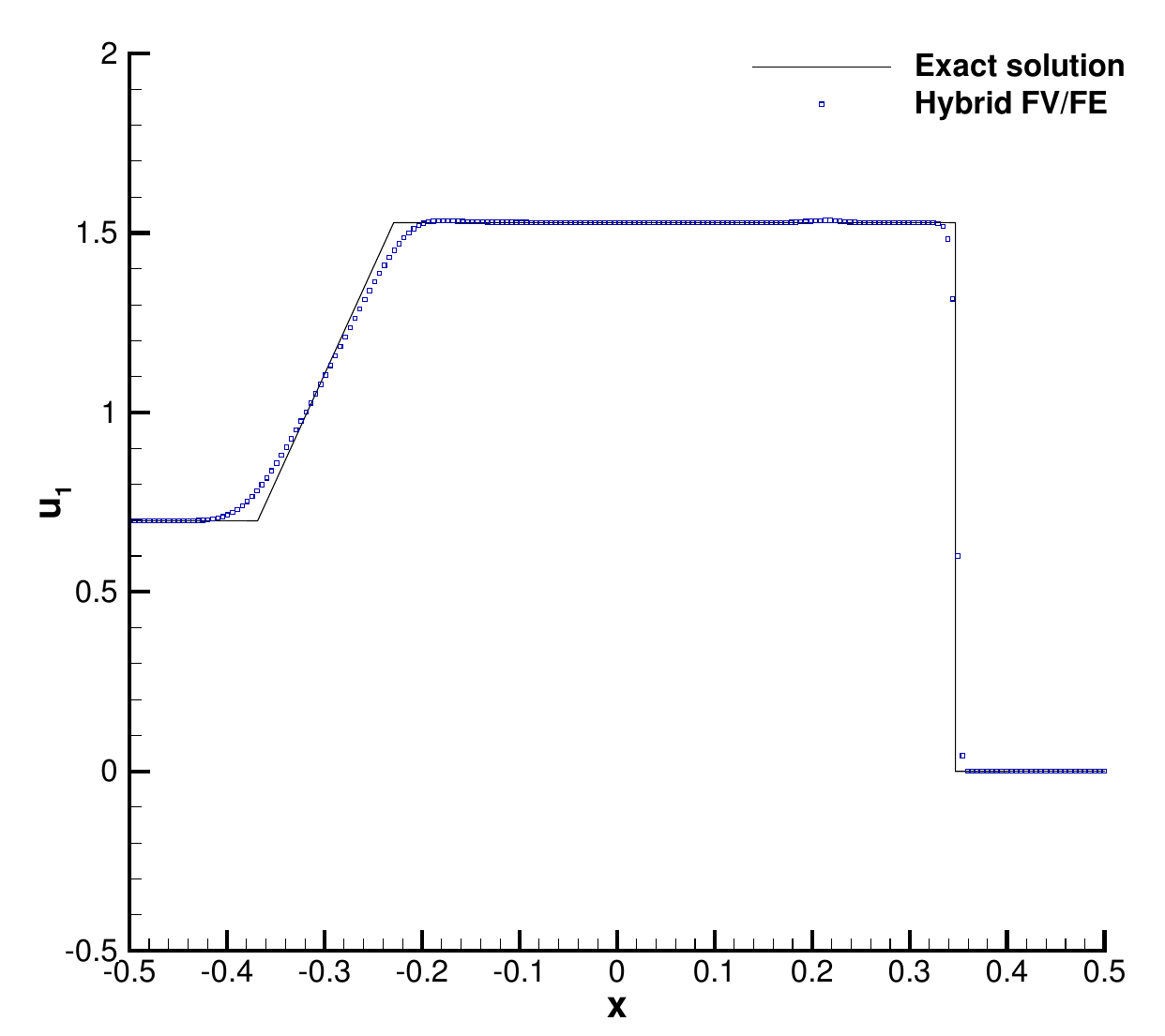}
	\includegraphics[width=0.32\linewidth]{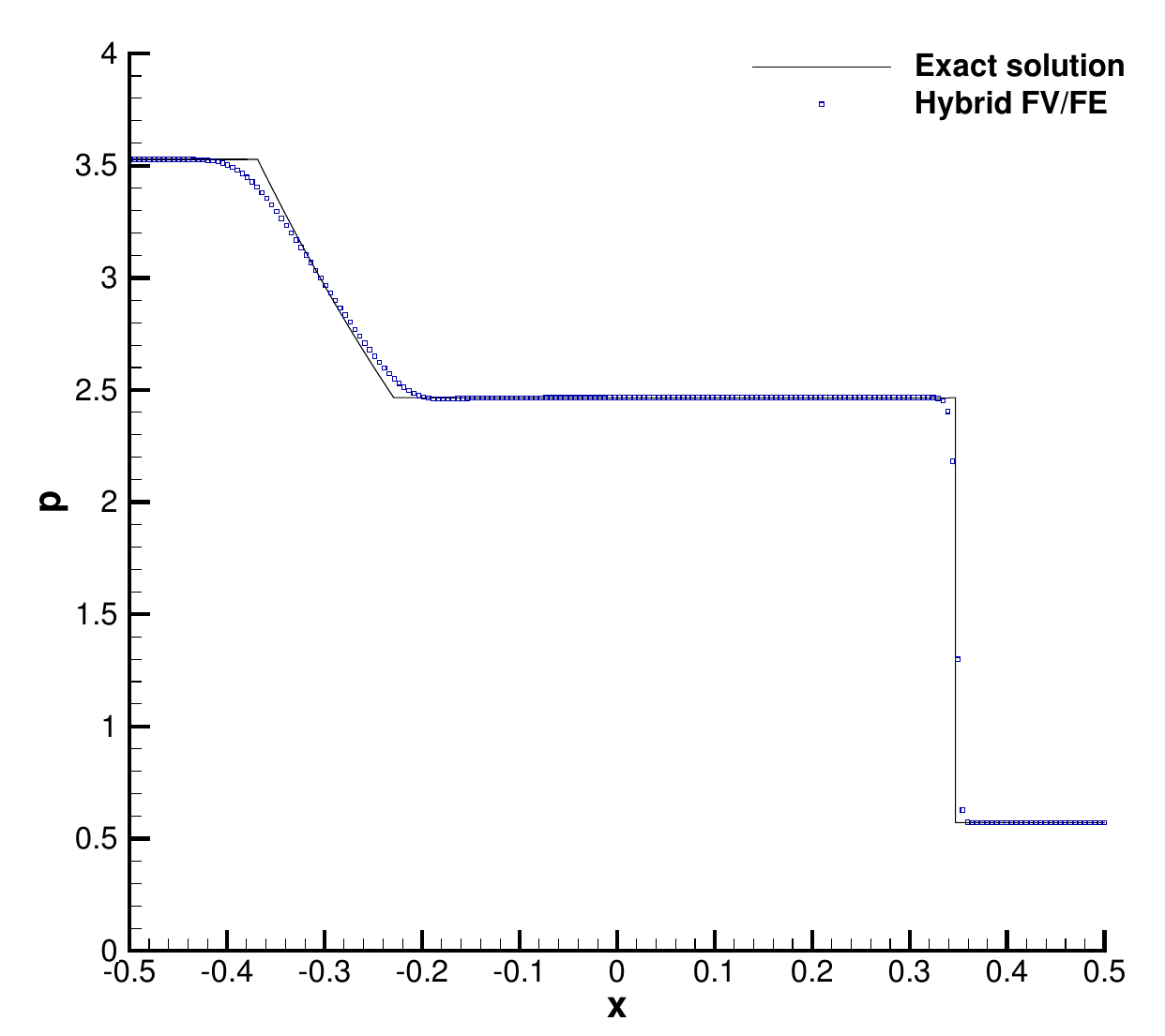}
	\caption{RP3 Lax. 1D cuts of the solution obtained using the hybrid FV/FE method for the compressible GPR model (blue squares) compared against the exact solution (black line). From left to right: density, horizontal velocity component and pressure.}
	\label{fig:rp3Lax}
\end{figure}

\begin{figure}
	\centering
	\includegraphics[width=0.32\linewidth]{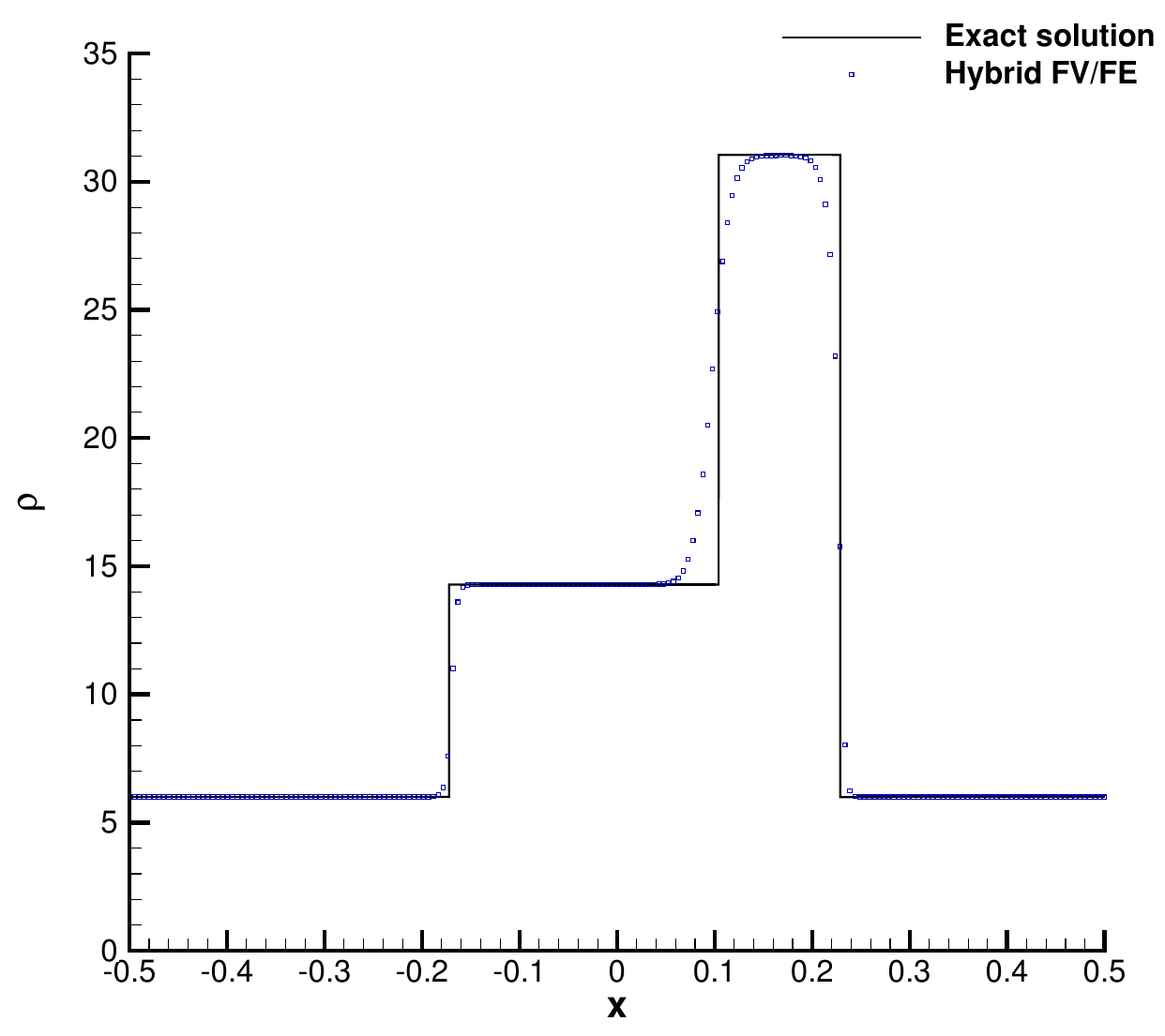}
	\includegraphics[width=0.32\linewidth]{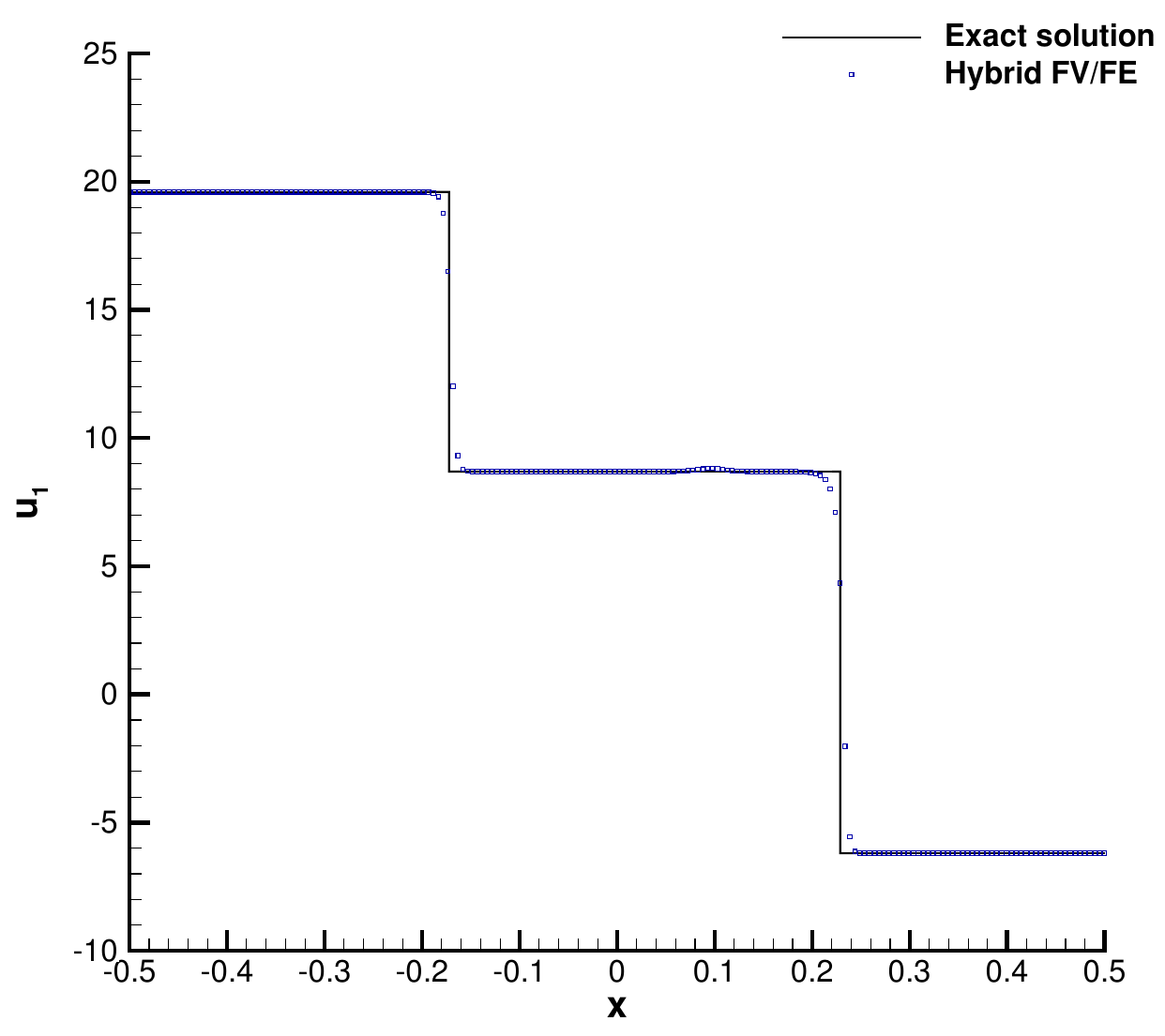}
	\includegraphics[width=0.32\linewidth]{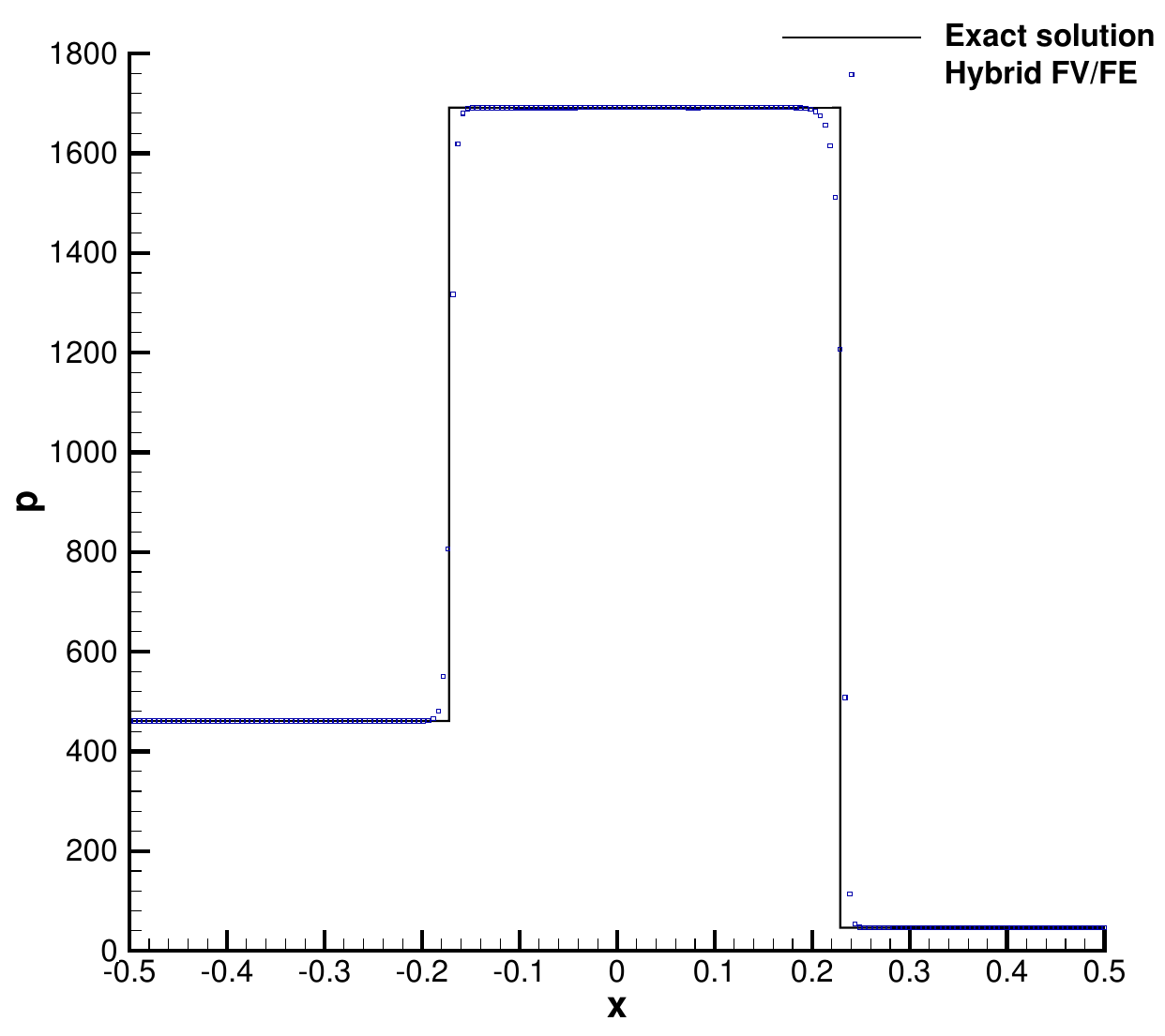}
	\caption{RP4. 1D cuts of the solution obtained using the hybrid FV/FE method for the compressible GPR model (blue squares) compared against the exact solution (black line). From left to right: density, horizontal velocity component and pressure.}
	\label{fig:rp4}
\end{figure}

\begin{figure}
	\centering
	\includegraphics[width=0.32\linewidth]{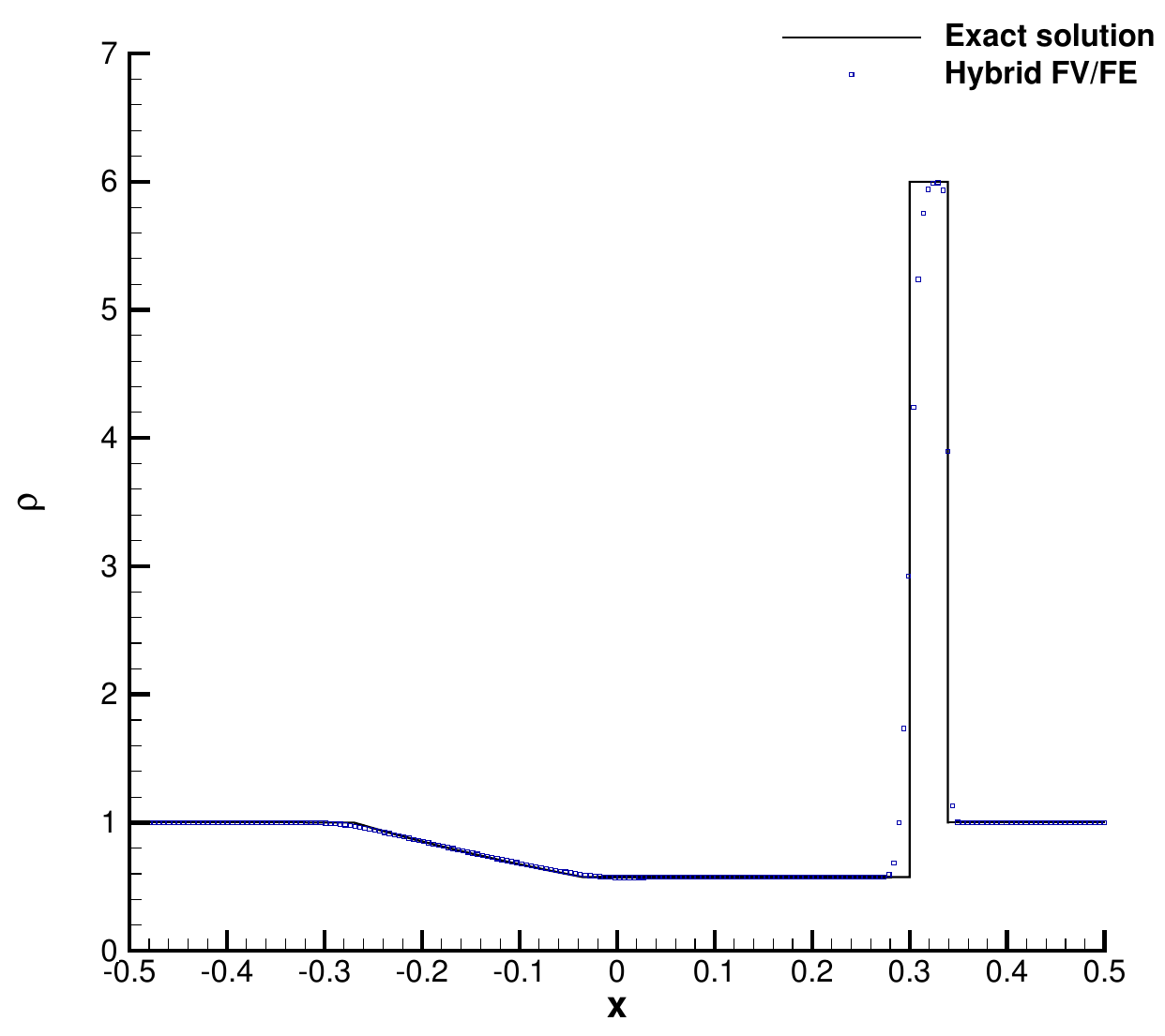}
	\includegraphics[width=0.32\linewidth]{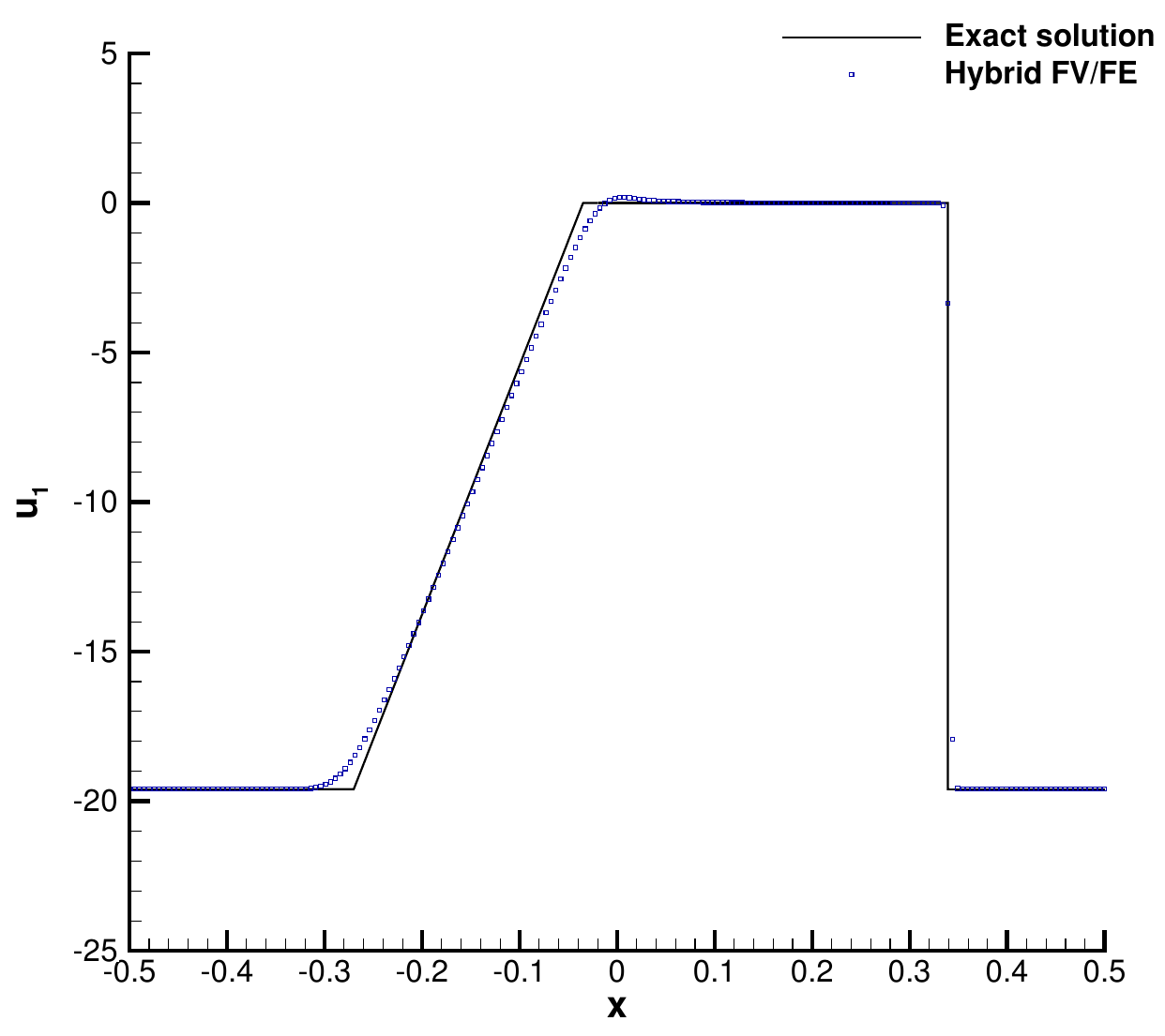}
	\includegraphics[width=0.32\linewidth]{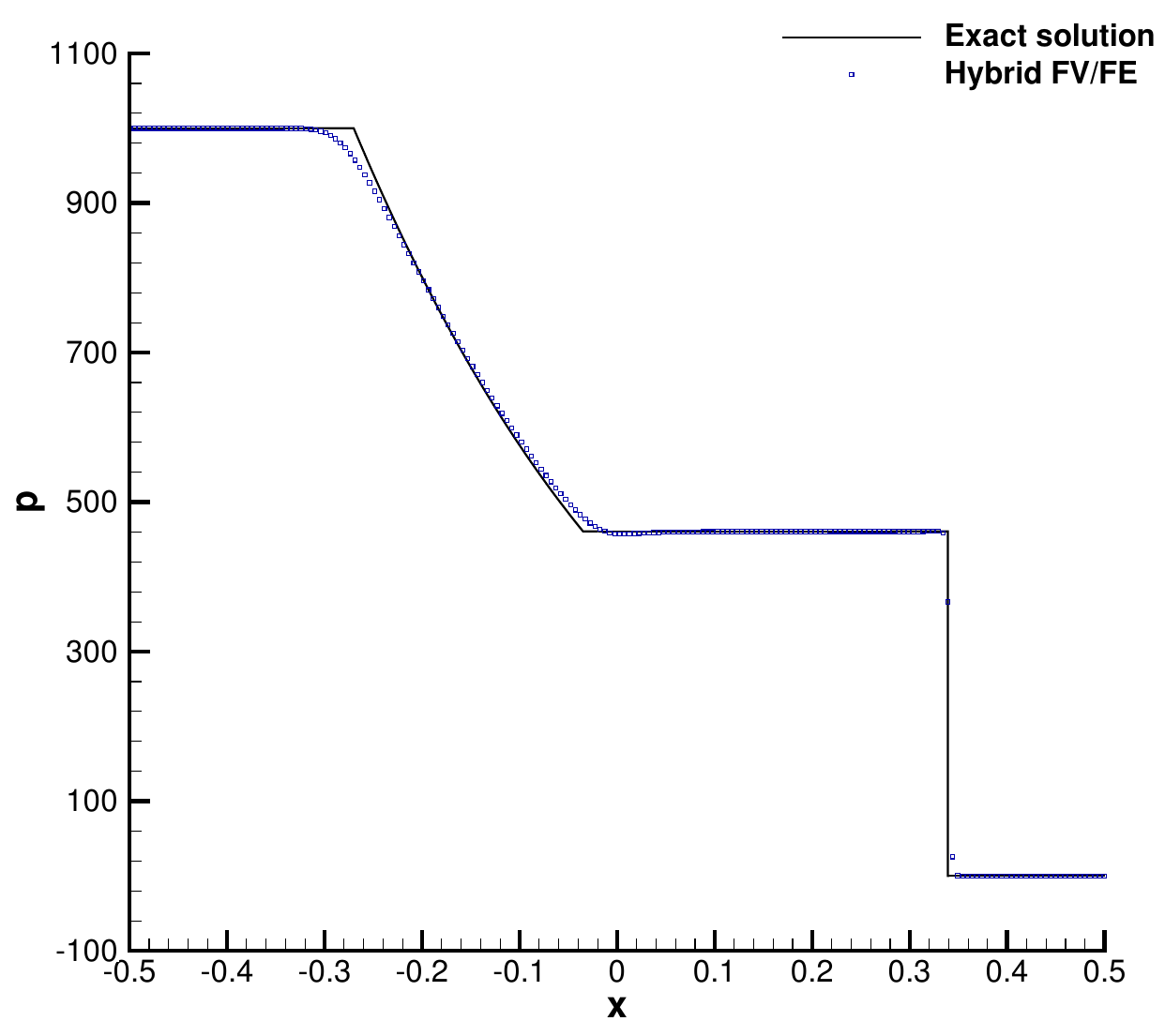}
	\caption{RP5 left blast problem. 1D cuts of the solution obtained using the hybrid FV/FE method for the compressible GPR model (blue squares) compared against the exact solution (black line). From left to right: density, horizontal velocity component  and pressure.}
	\label{fig:rp5blast}
\end{figure}

\begin{figure}
	\centering
	\includegraphics[width=0.32\linewidth]{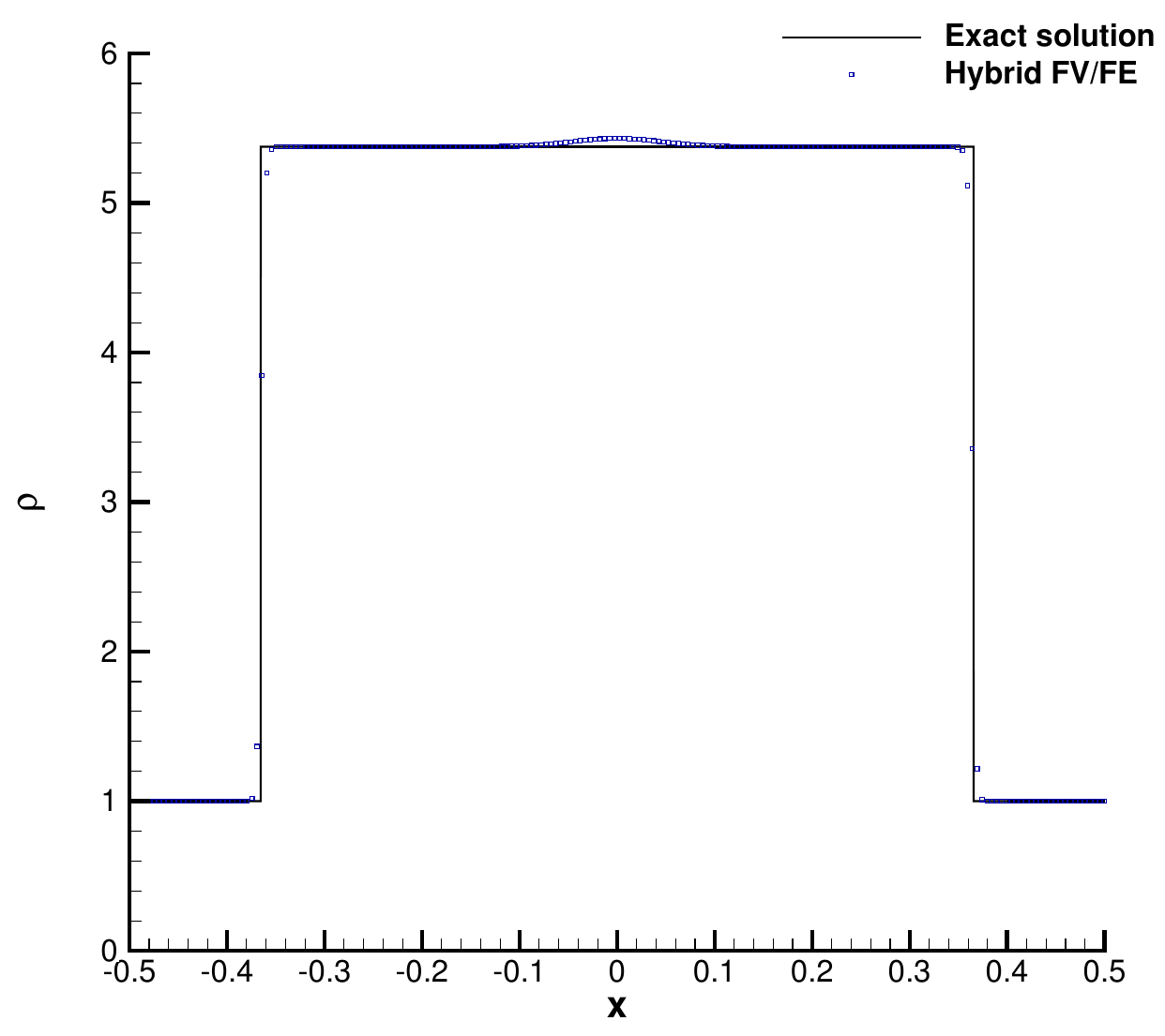}
	\includegraphics[width=0.32\linewidth]{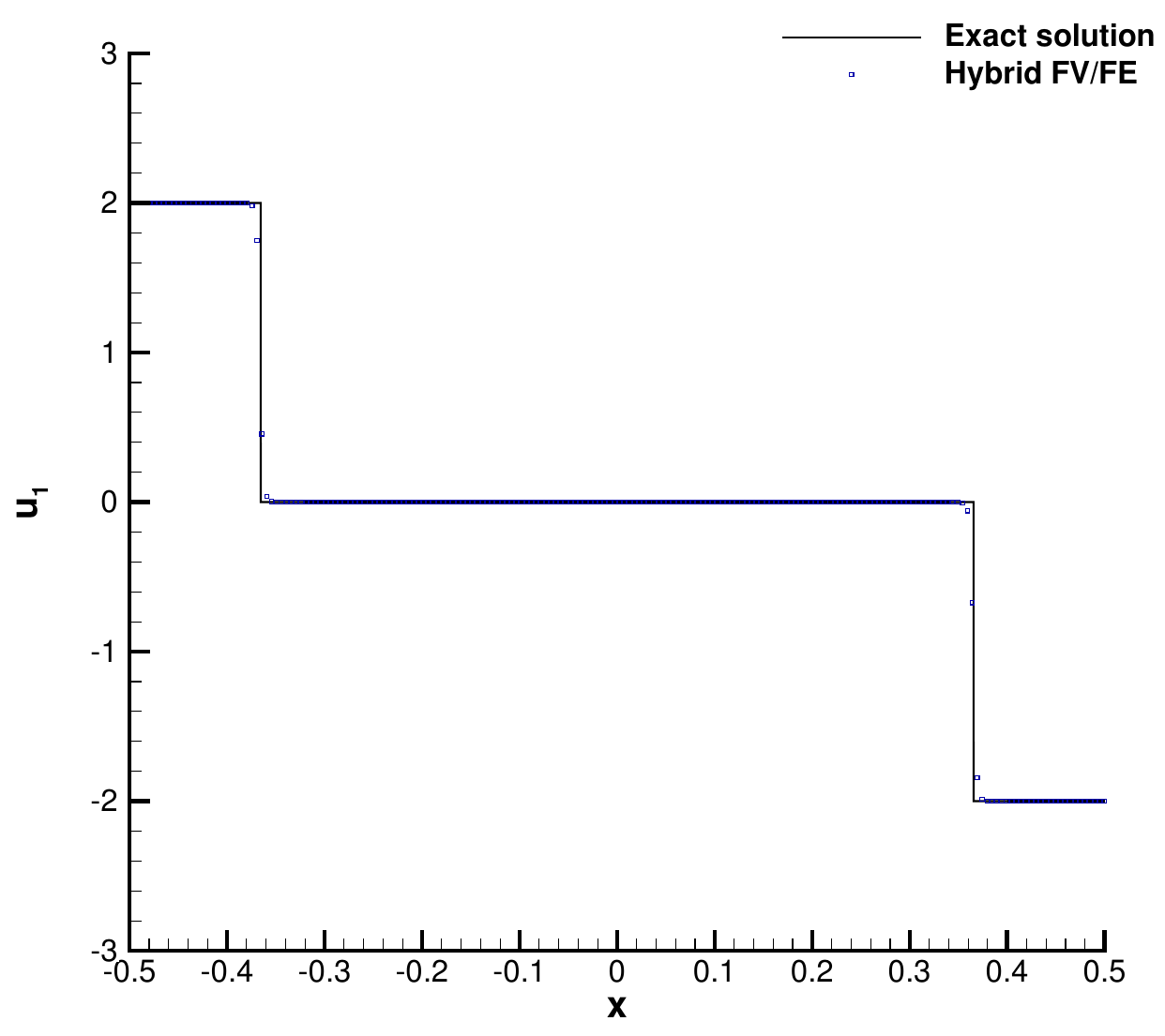}
	\includegraphics[width=0.32\linewidth]{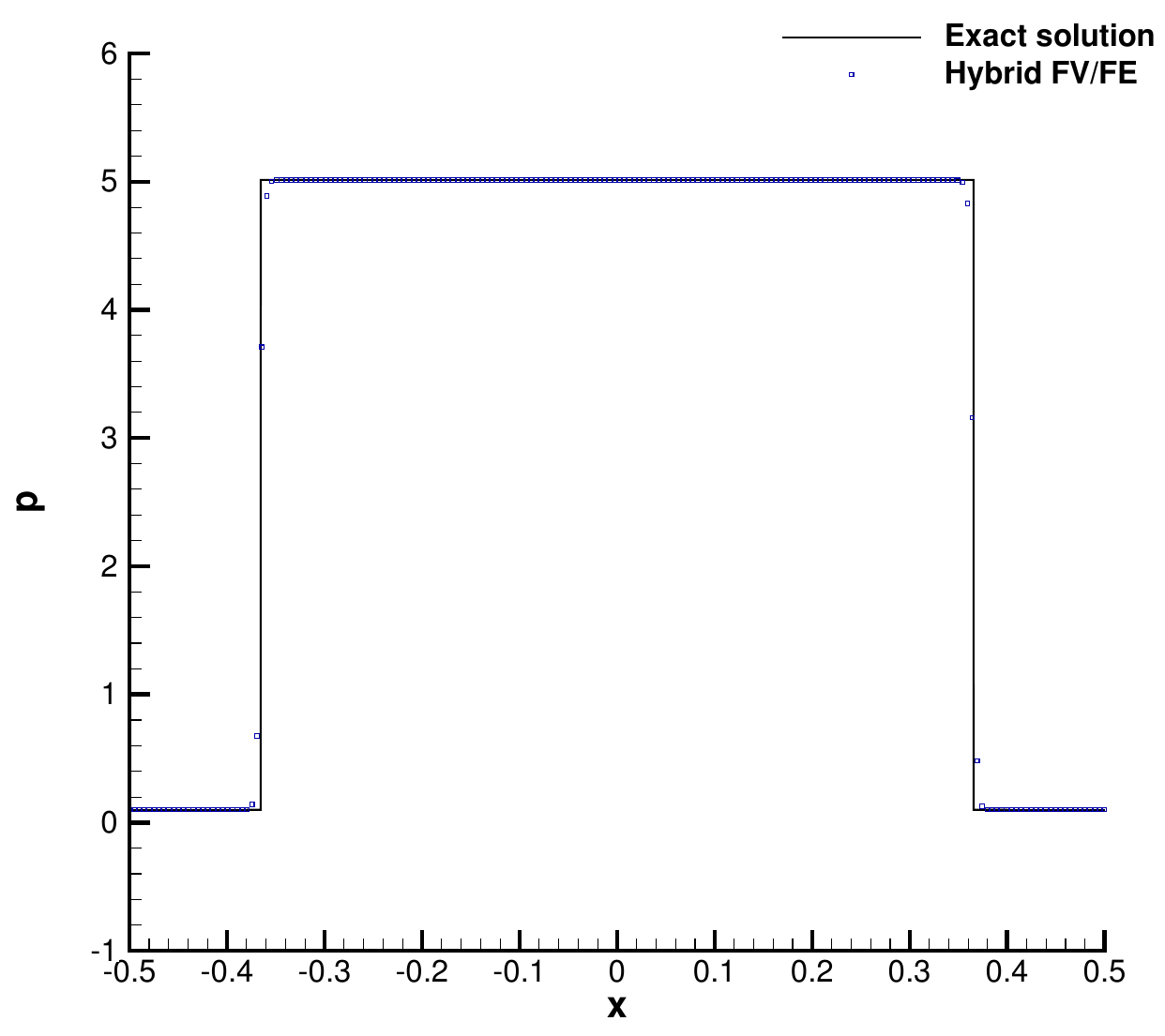}
	\caption{RP6 double shock. 1D cuts of the solution obtained using the hybrid FV/FE method for the compressible GPR model (blue squares) compared against the exact solution (black line). From left to right: density, horizontal velocity component and pressure.}
	\label{fig:rp6DS}
\end{figure}

\begin{figure}
	\centering
	\includegraphics[width=0.32\linewidth]{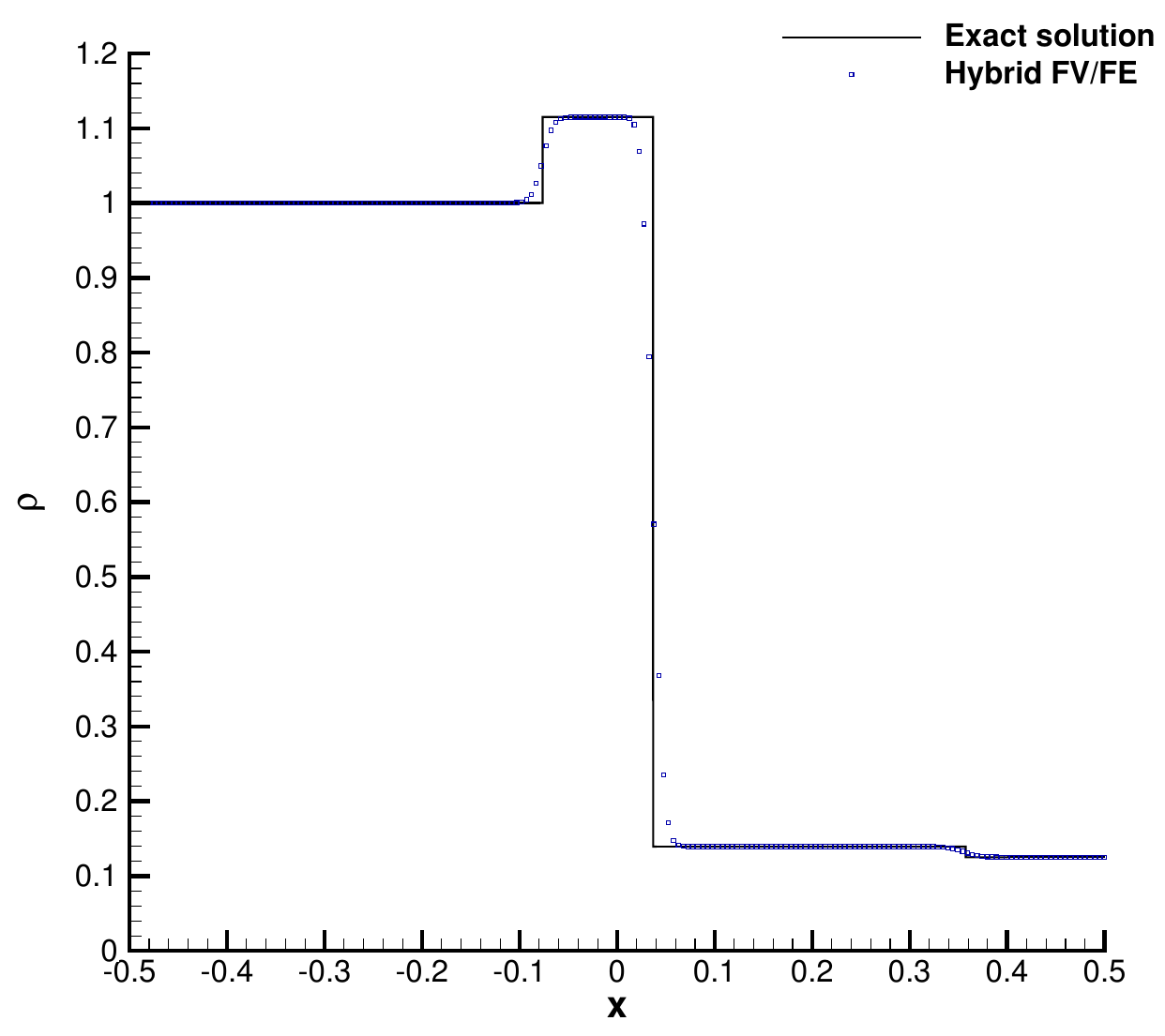}
	\includegraphics[width=0.32\linewidth]{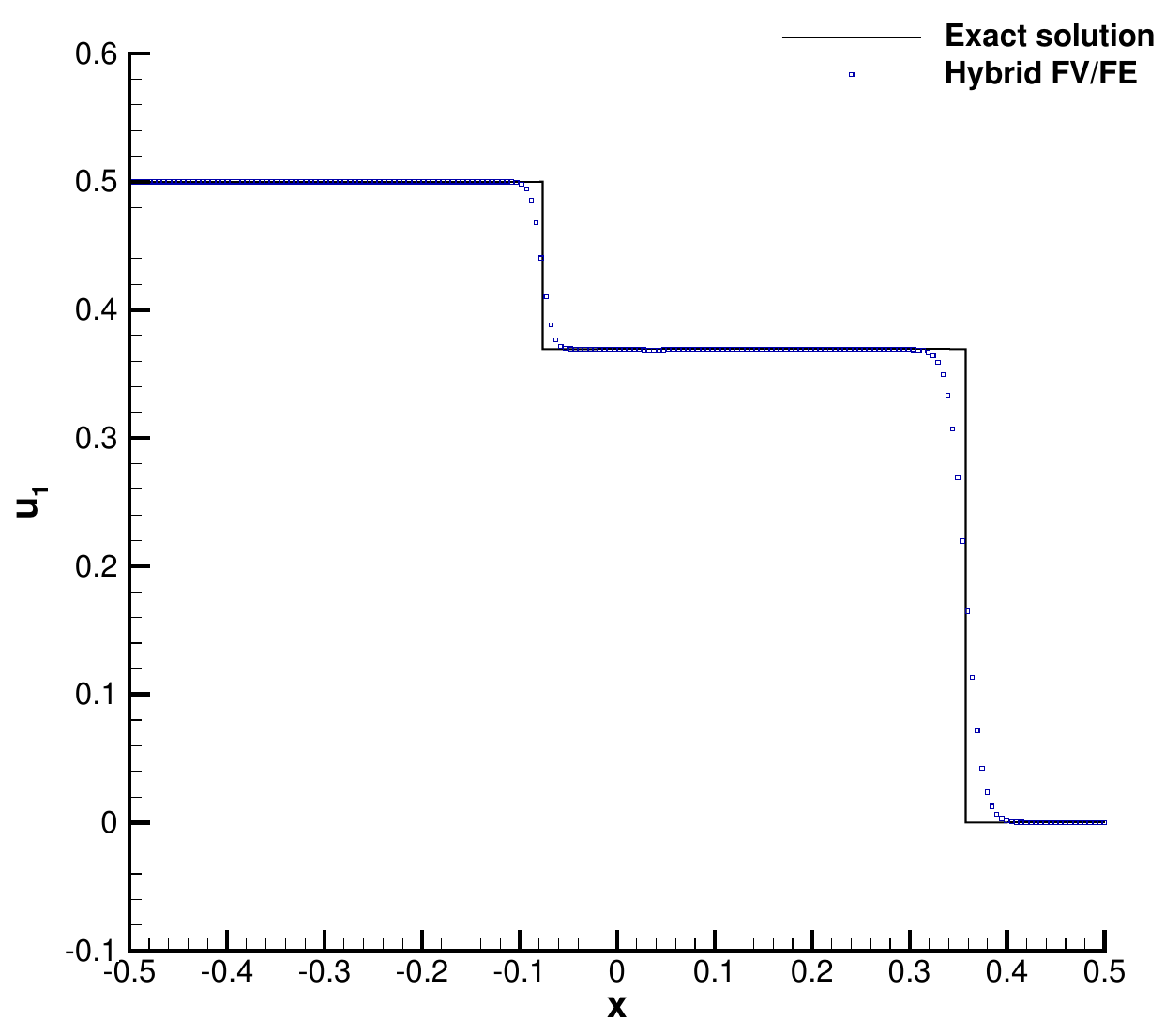}
	\includegraphics[width=0.32\linewidth]{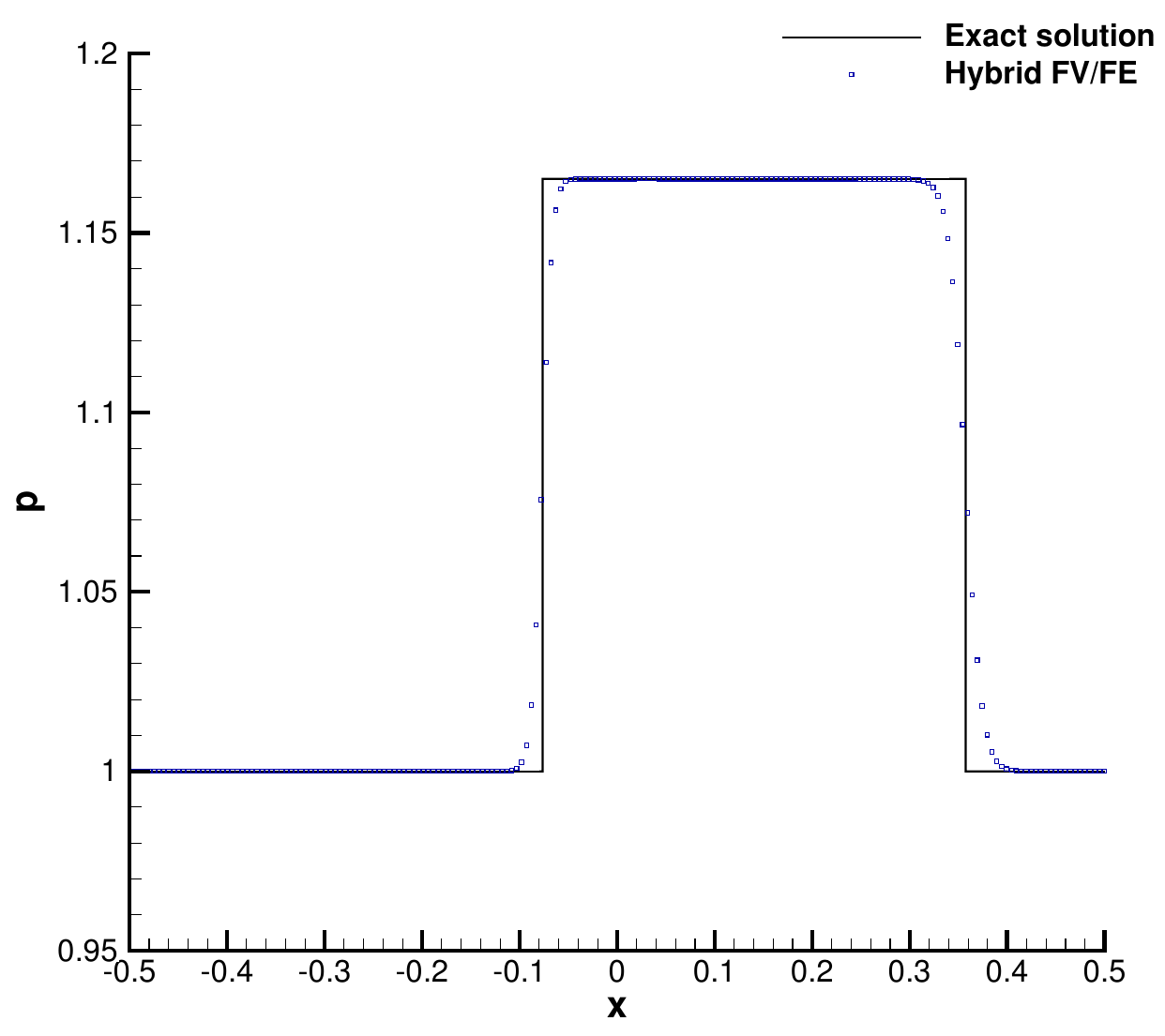}
	\caption{RP7. 1D cuts of the solution obtained using the hybrid FV/FE method for the compressible GPR model (blue squares) compared against the exact solution (black line). From left to right: density, horizontal velocity component and pressure.}
	\label{fig:rp7}
\end{figure}

Riemann problems RP$1$--RP$7$ correspond to classical tests for inviscid flows and, therefore, an exact Riemann solver for the Euler equations has been employed to get the exact solution of the density, pressure and velocity fields \cite{Toro}. In Figure~\ref{fig:rp1SodALE}, we report the results obtained for the Sod shock tube problem which presents a shock, a contact discontinuity and a rarefaction wave. For this Riemann problem both the purely Eulerian and the ALE scheme with $\varsigma=5$ have been run. The numerical solution computed using the HTC FV scheme in \cite{HTCA2} has also been included for comparison. 
Figures~\ref{fig:rp2DR} and \ref{fig:rp6DS} show the solutions obtained for a double rarefaction, RP2, and a double shock, RP6, tests. 
Combined waves are analysed through the Lax benchmark (RP3), RP4, where two shock colliding waves lead to three strong discontinuities, and RP7, where a shock presenting a small wave in the density field is generated. The obtained results are reported in Figures~\ref{fig:rp3Lax}, \ref{fig:rp4} and \ref{fig:rp7}. Finally, Figure~\ref{fig:rp5blast} depicts the solution obtained for the severe left blast problem introduced in \cite{WC84}.

\begin{figure}
	\centering
	\includegraphics[width=0.45\linewidth]{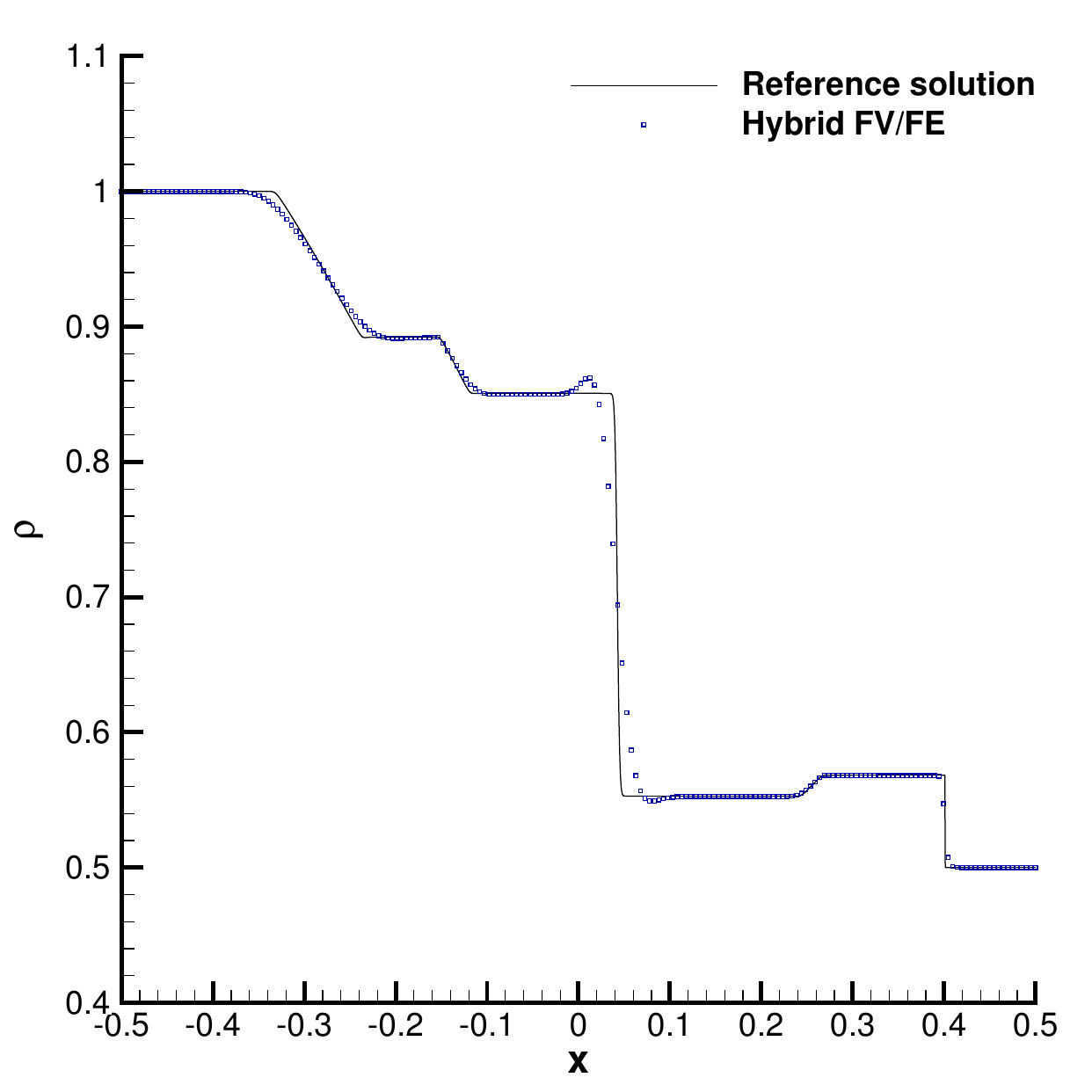}
	\includegraphics[width=0.45\linewidth]{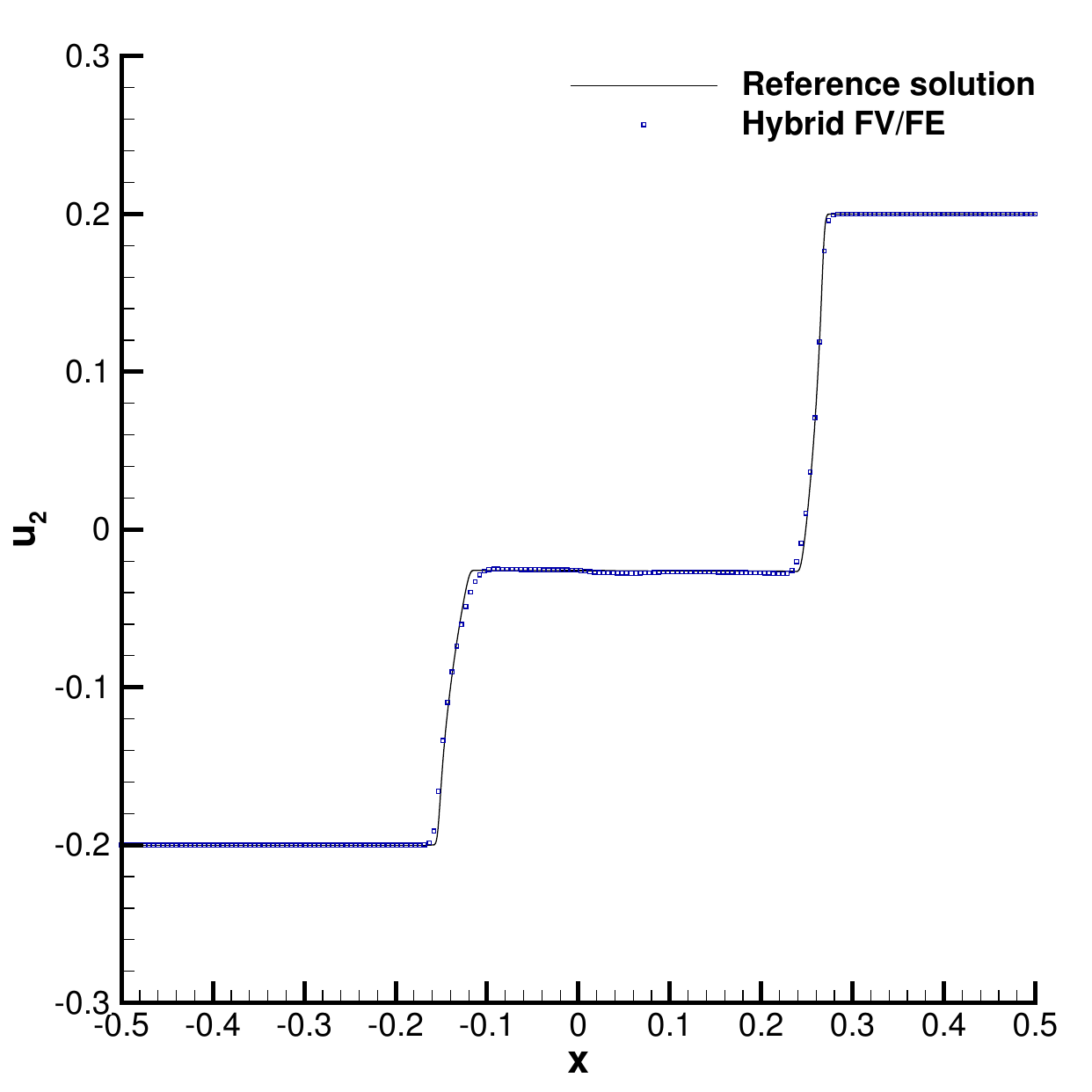}
	\includegraphics[width=0.45\linewidth]{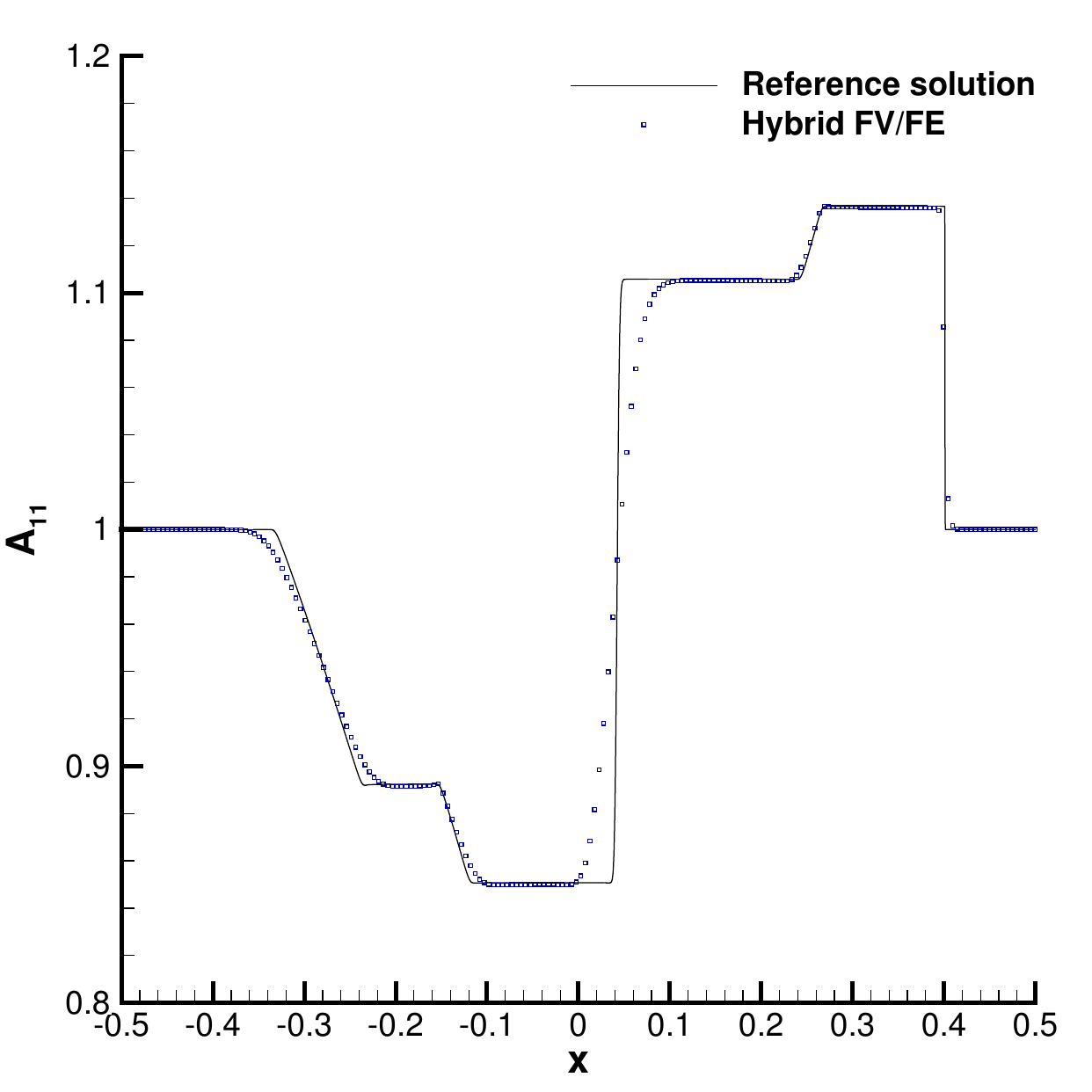}
	\includegraphics[width=0.45\linewidth]{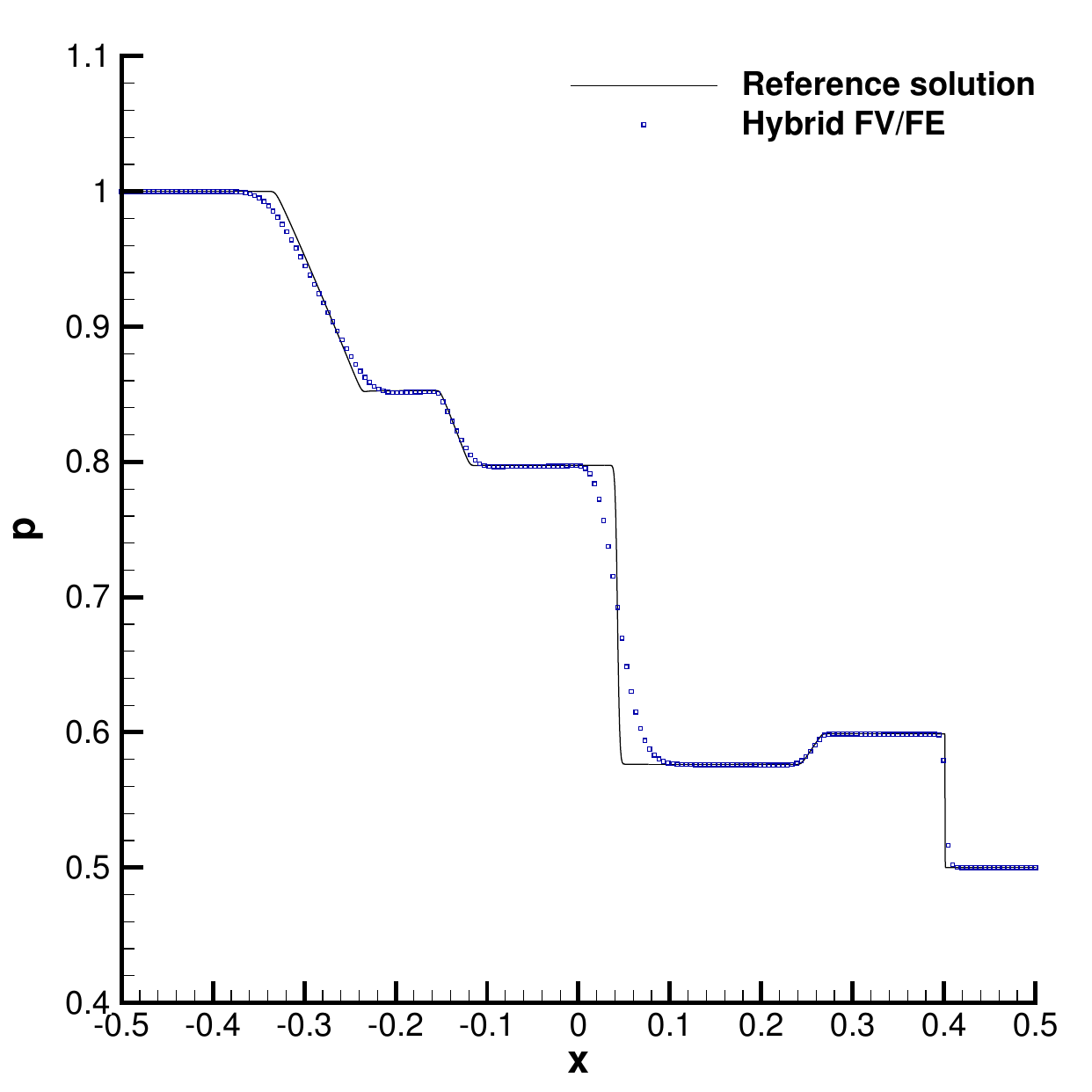}
	\caption{RP8. 1D cuts of the solution obtained using the hybrid FV/FE method for the compressible GPR model (blue squares) compared against the reference solution computed using the HTC method proposed in \cite{HTCA2}  (black line). From left-top to right-bottom: density, velocity component $\vel_{2}$, distortion field component $A_{11}$ and pressure.}
	\label{fig:rp9}
\end{figure}

\begin{figure}
	\centering
	\includegraphics[width=0.45\linewidth]{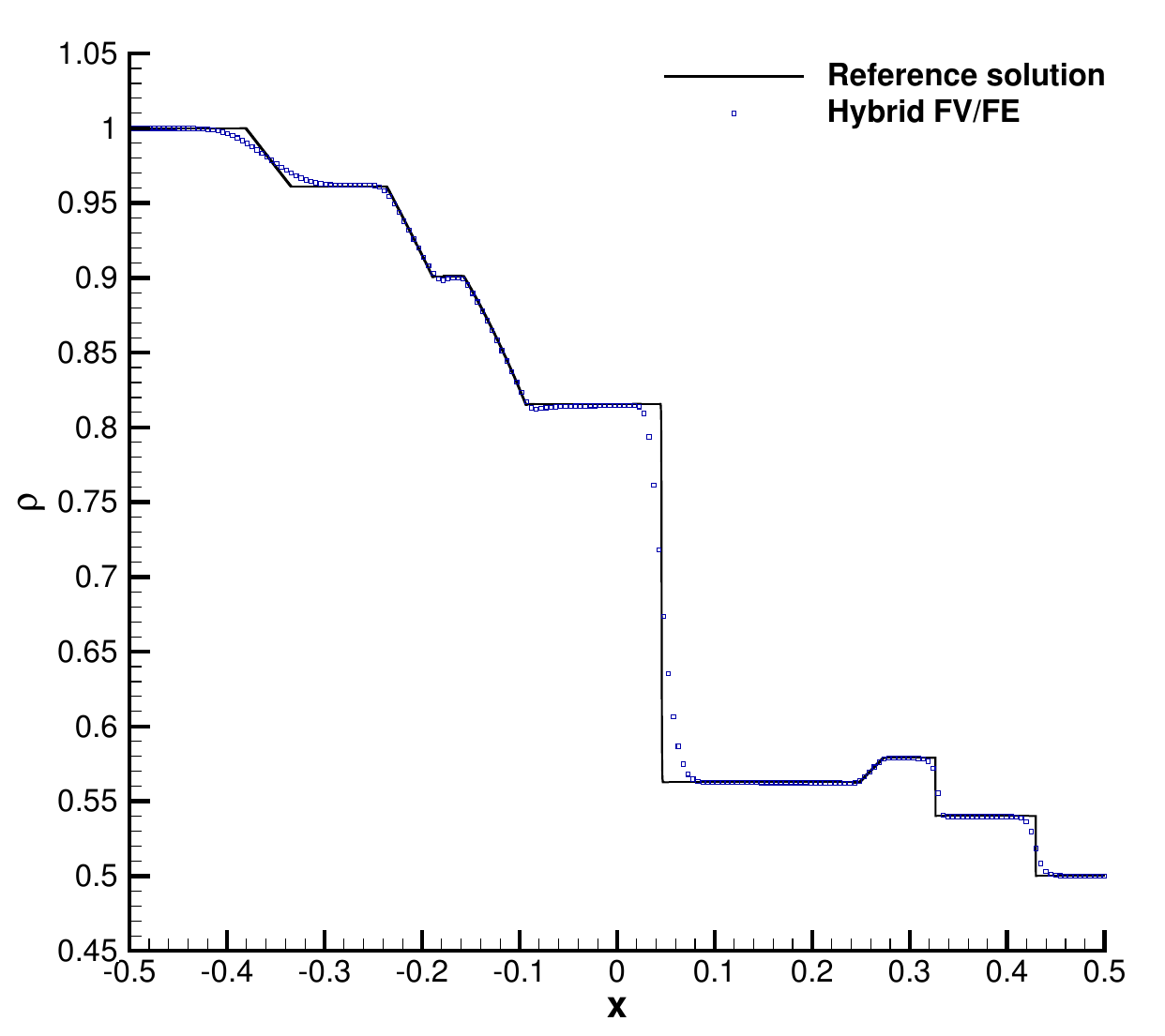}
	\includegraphics[width=0.45\linewidth]{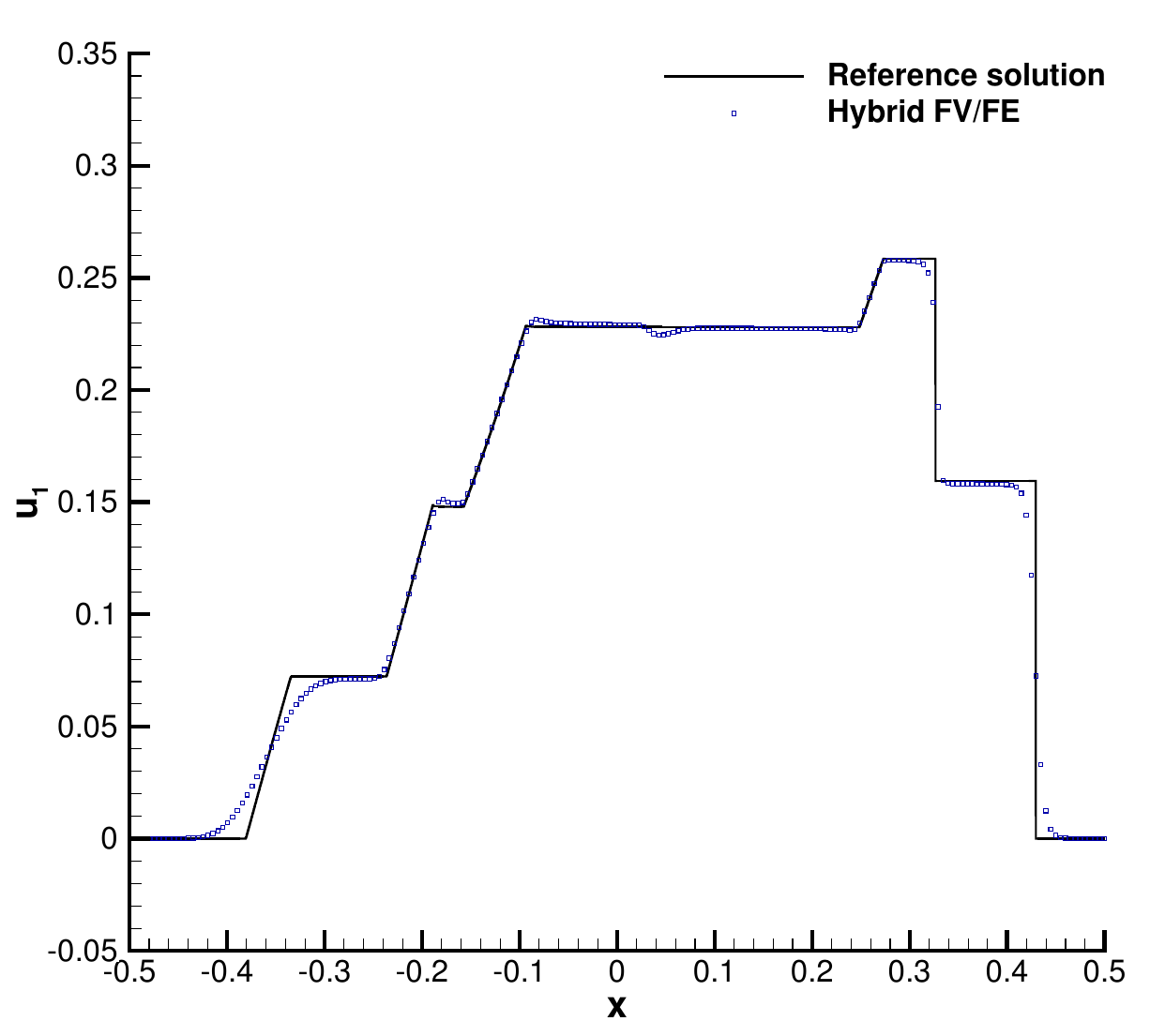}
	\includegraphics[width=0.45\linewidth]{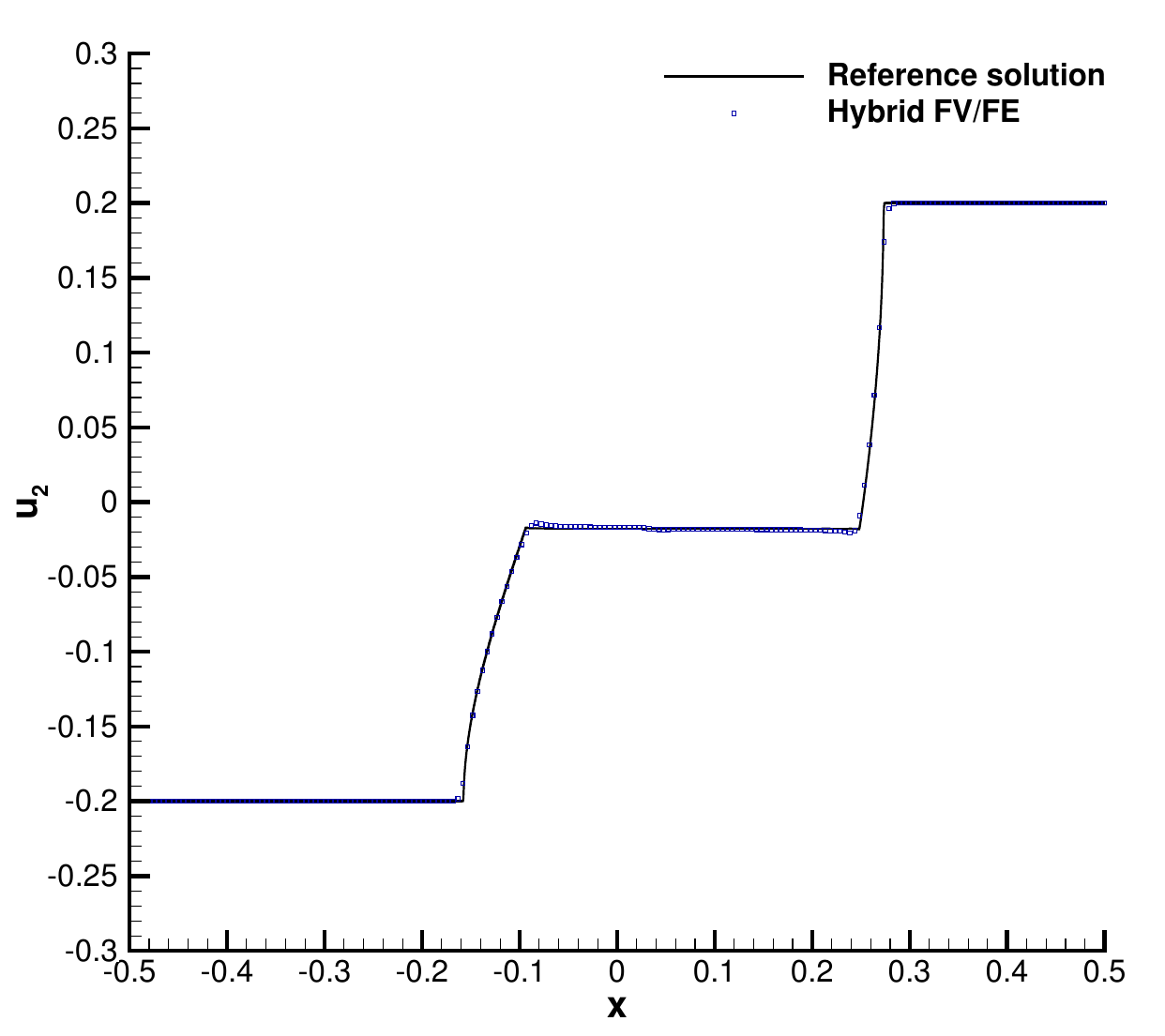}
	\includegraphics[width=0.45\linewidth]{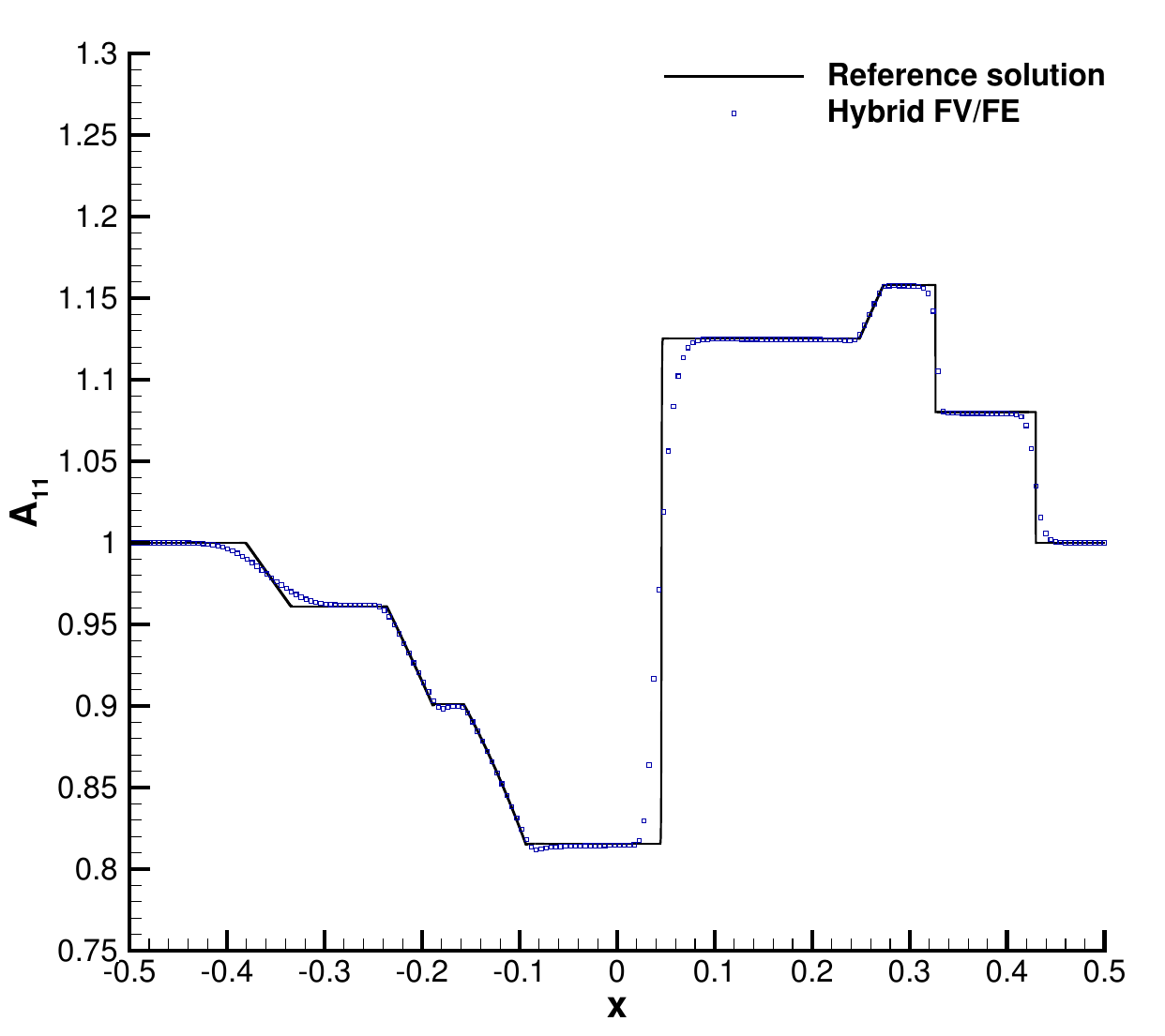}
	\includegraphics[width=0.45\linewidth]{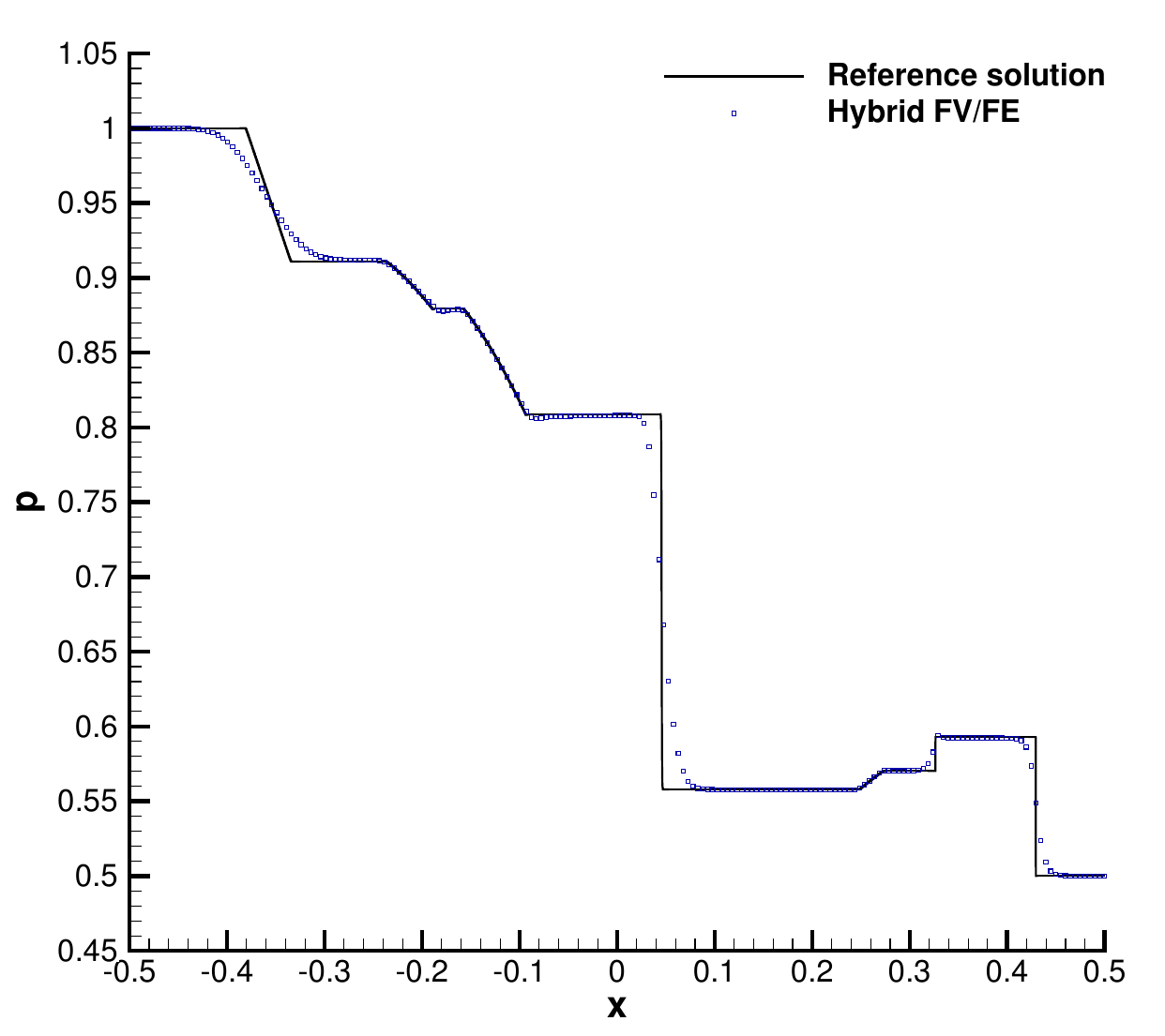}
	\includegraphics[width=0.45\linewidth]{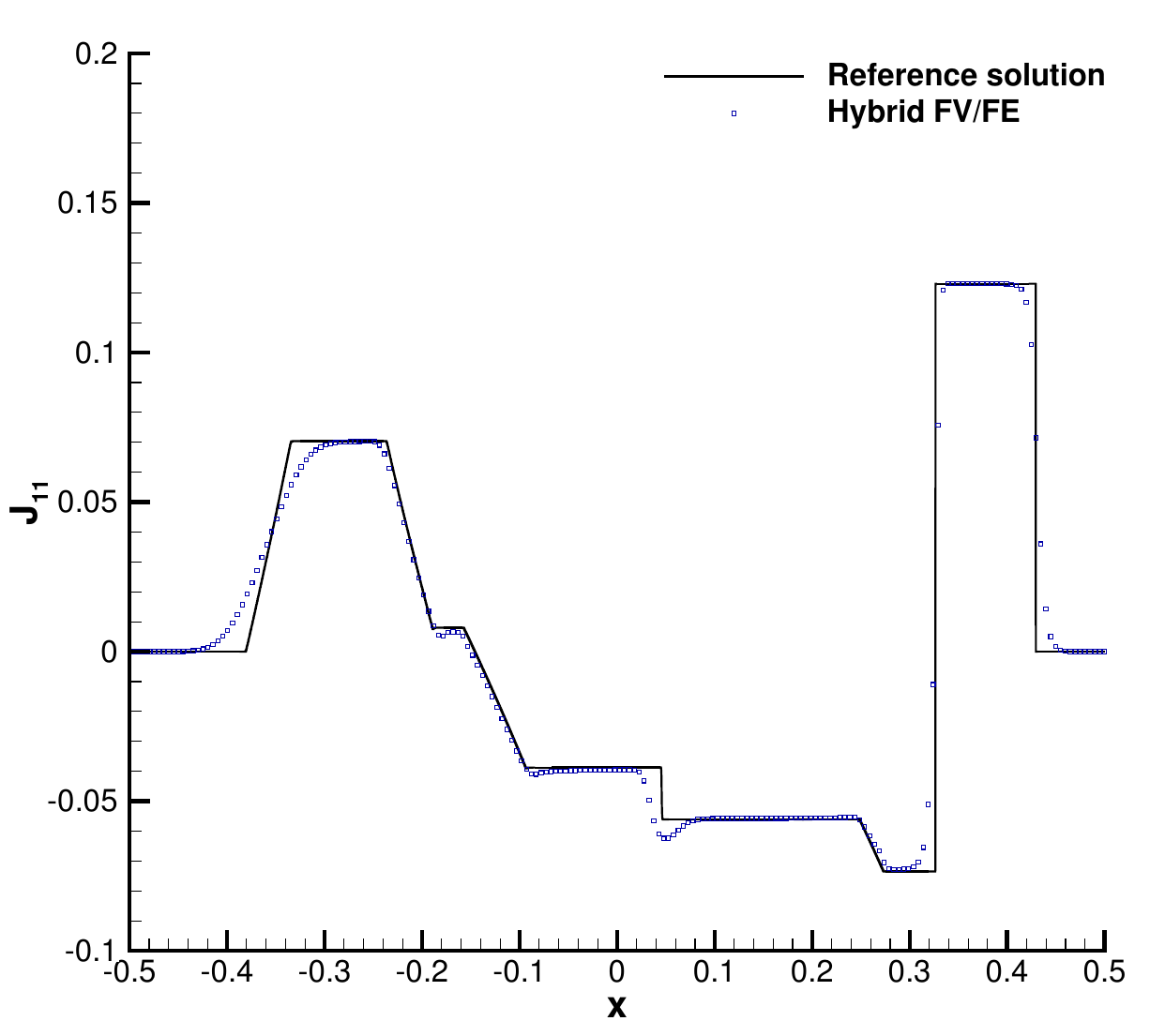}
	\caption{RP9. 1D cuts of the solution obtained using the hybrid FV/FE method for the compressible GPR model (blue squares) compared against the reference solution computed using the HTC method proposed in \cite{HTCA2} (black line). From left-top to right-bottom: density, velocity components $\vel_{1}$, $\vel_{2}$, distortion field component $A_{11}$, pressure and thermal impulse $J_{1}$.}
	\label{fig:rp10}
\end{figure}

On the other hand, test cases RP$8$, RP$9$ are classically employed to specifically assess numerical methods in the solid limit of the GPR model \cite{Boscheri2021SIGPR,HTCGPR}. The initial condition coincides for both tests but different media are considered: RP8 corresponds to an ideal elastic solid without heat conduction and RP9 incorporates heat conduction. As expected, a contact discontinuity, two acoustic and two shear waves plus two additional thermo-acoustic waves, when accounting for heat conduction, are generated, see Figures~\ref{fig:rp9}--\ref{fig:rp10}. The provided reference solutions have been computed using the thermodynamically compatible scheme (HTC) proposed in \cite{HTCA2}, which solves the entropy equation instead of the total energy conservation law. 

All solutions of the semi-implicit hybrid methodology have been computed using an ENO limiting strategy on the primitive variables. 
The artificial viscosity has been activated with values 
$c_{\alpha}=0.2$ for RP1 and RP7,
$c_{\alpha}=2$ for RP2 and RP6,
$c_{\alpha}=0.5$ for RP3,
$c_{\alpha}=5$ for RP4,
$c_{\alpha}=10$ for RP5 and
$c_{\alpha}=1$ for RP8-9. 
A structured mesh of $N_{x}=400$ divisions along the $x$-axis has been employed for all test cases except for RP5 and RP10, which have been run on grids of $N_{x}=800$ and  $N_{x}=2000$ divisions, respectively. 
Overall, the results obtained show a good agreement with the exact and reference solutions in presence of strong waves.

\subsection{Fluid and solid circular explosions}
\begin{figure}[H]
	\centering
	\includegraphics[width=0.38\linewidth]{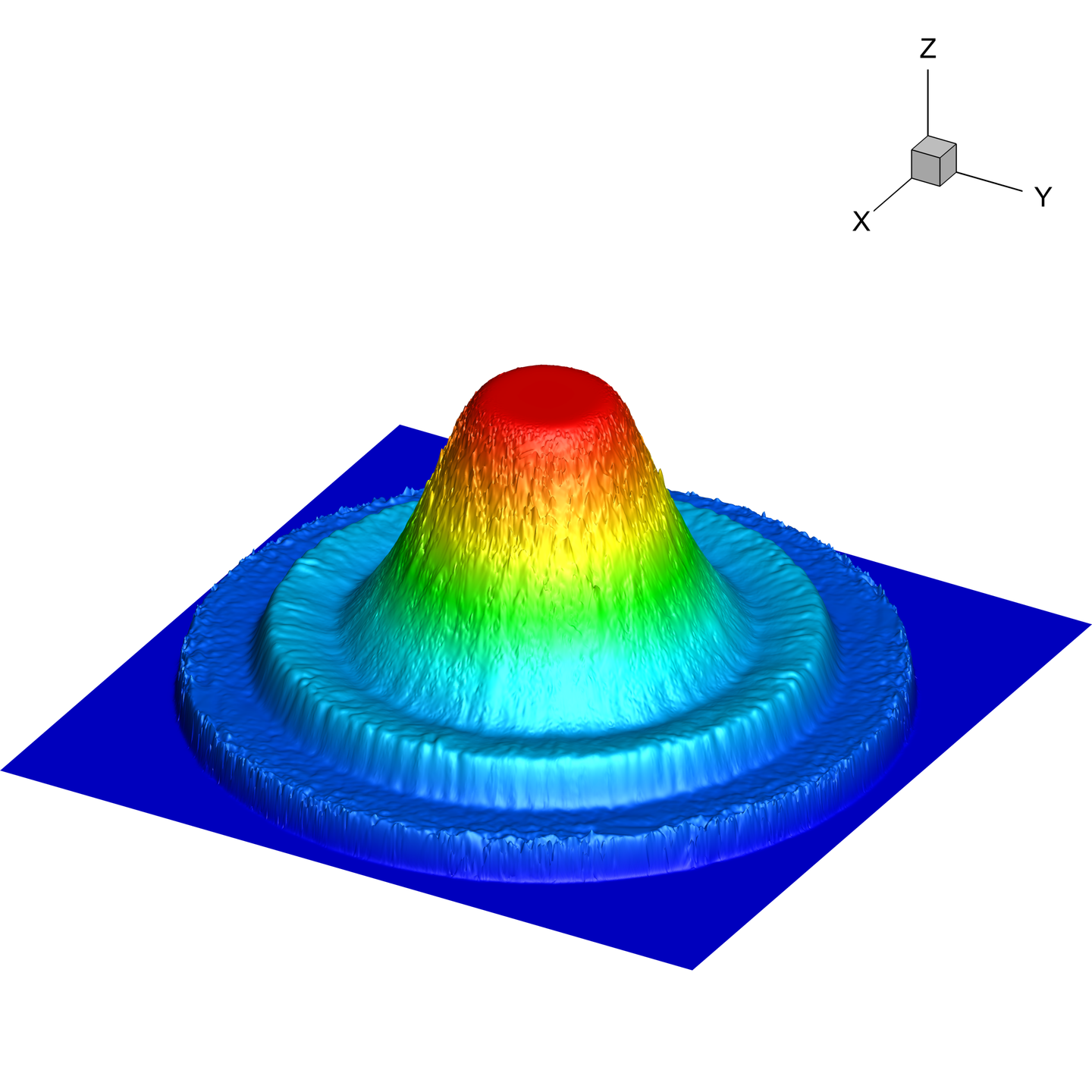}
	\includegraphics[width=0.38\linewidth]{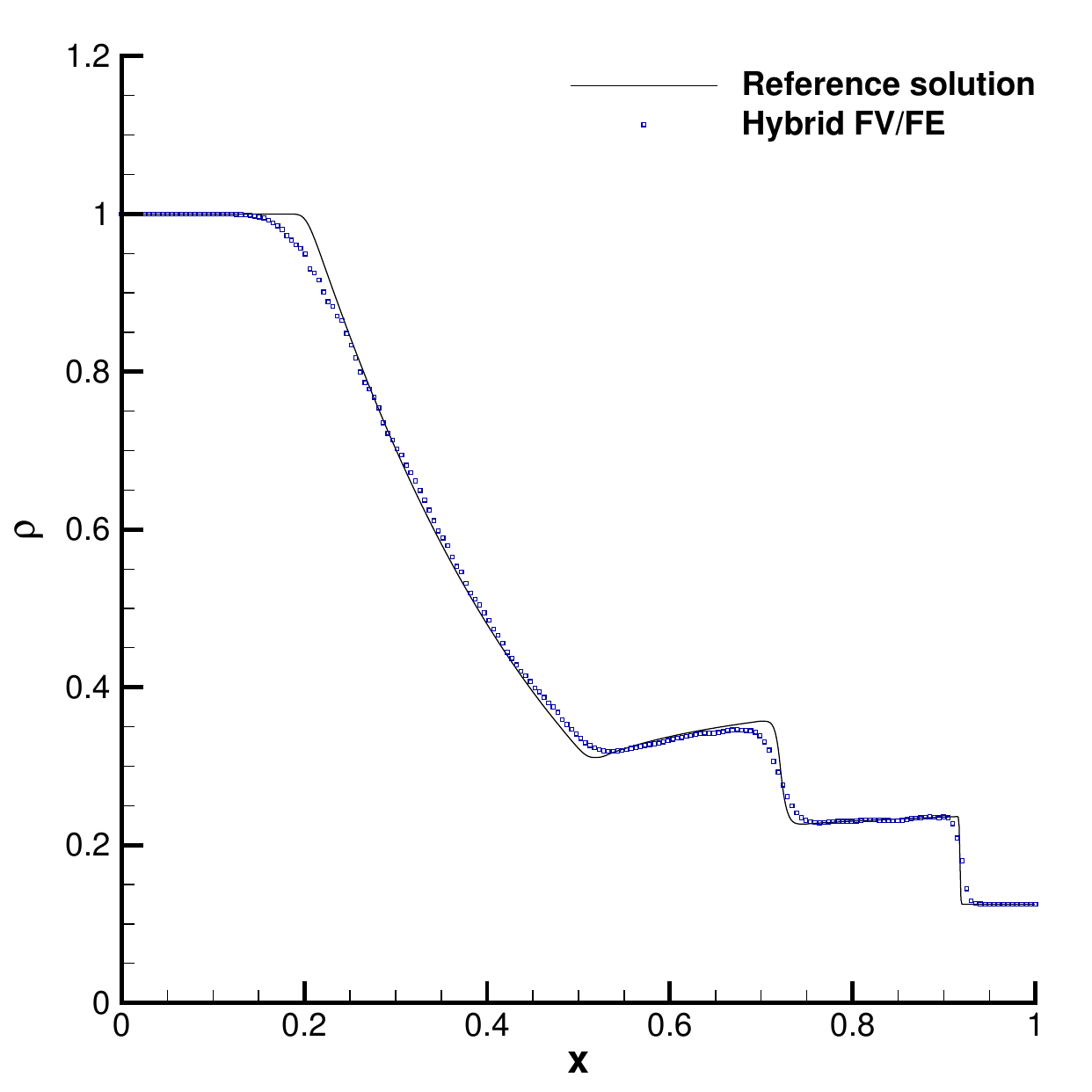}
	\includegraphics[width=0.38\linewidth]{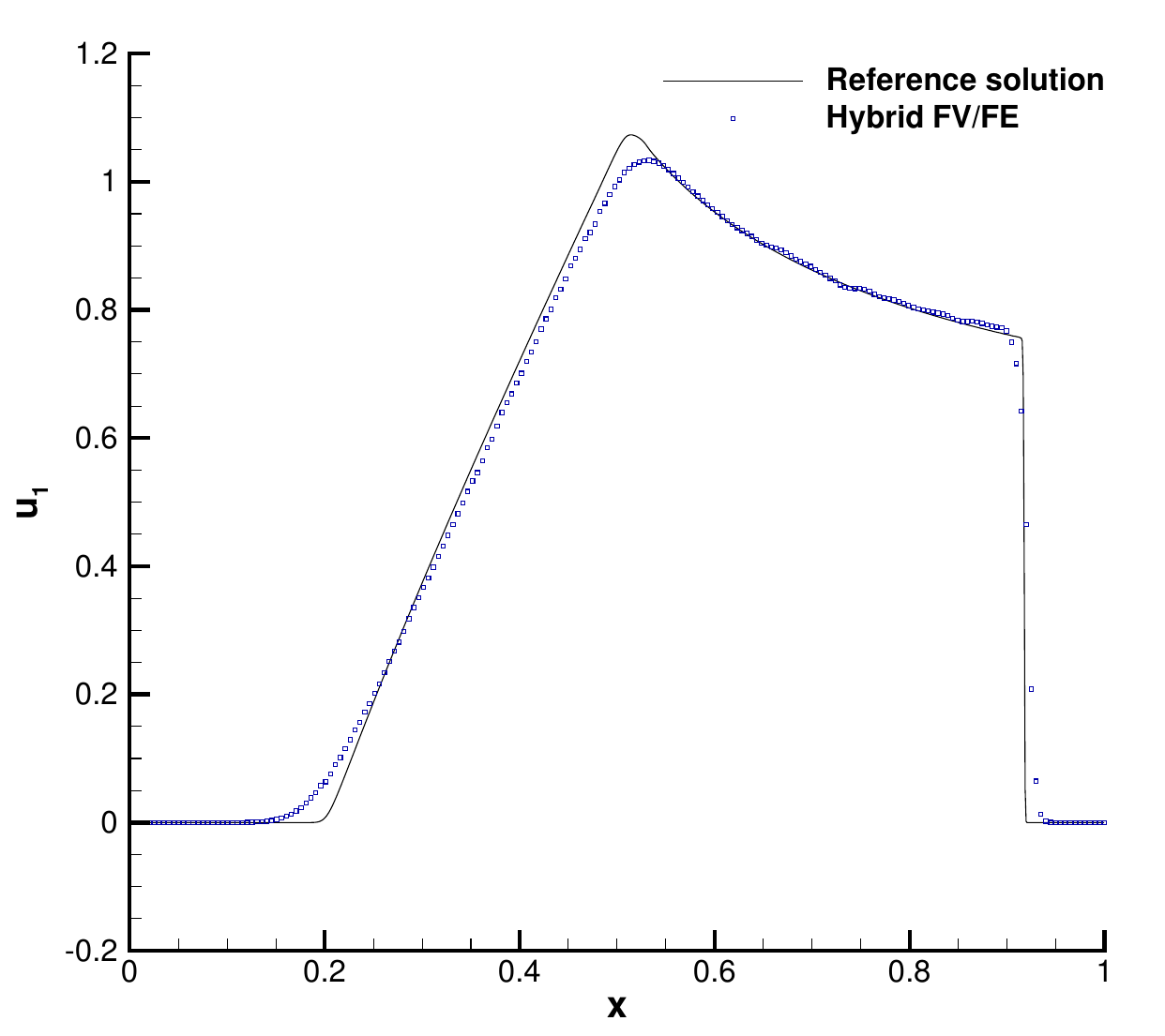}
	\includegraphics[width=0.38\linewidth]{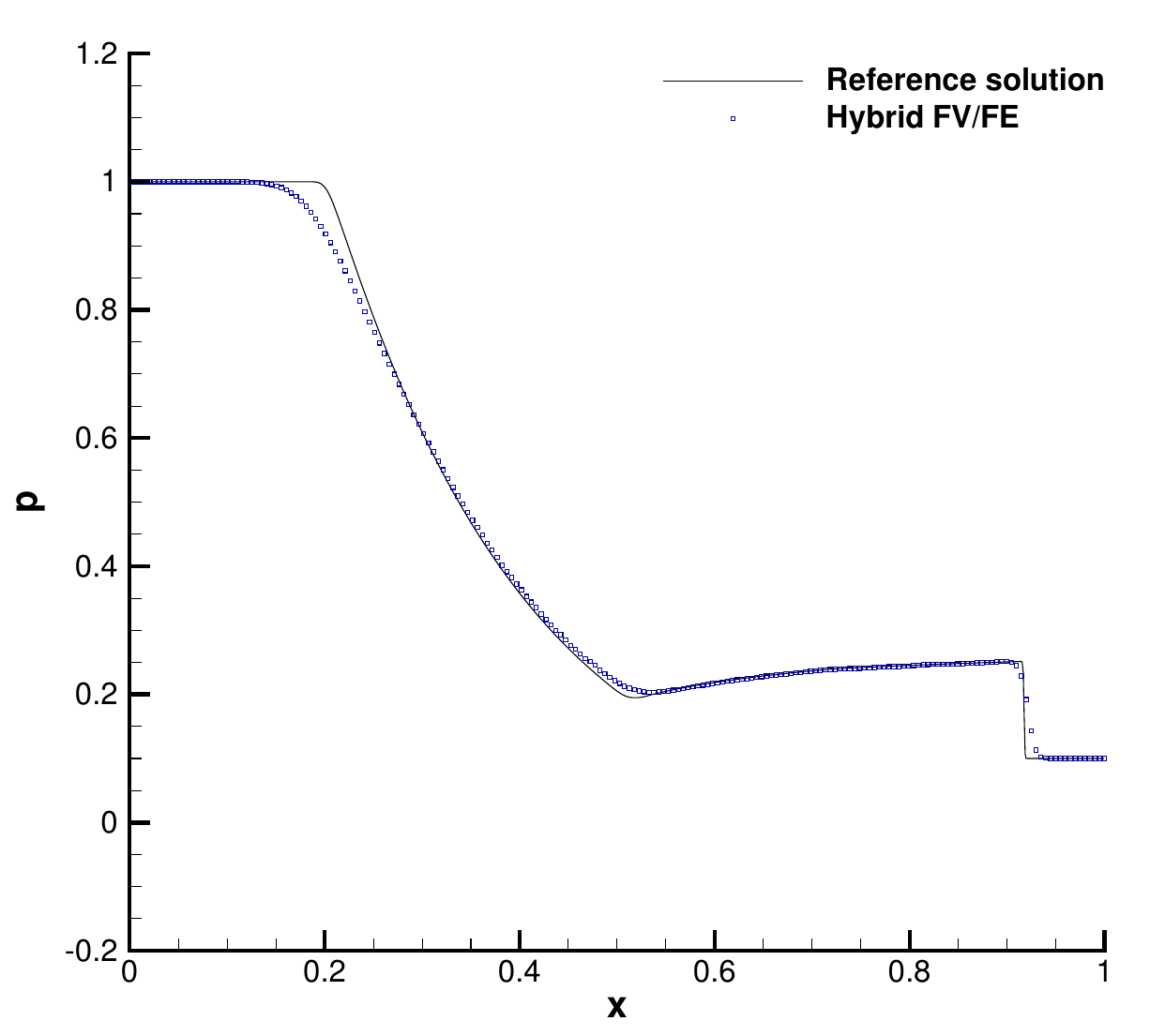}
	\caption{Fluid circular explosion 2D. 3D extrusion of the density field and 1D cuts of the density, $\vel_{1}$ and pressure for $\left\lbrace (x,y)\in\mathbb{R}^{2}\,  |\,x \in[0,1],\, y=0 \right\rbrace$ using the hybrid method (blue squares) compared against the reference solution (black line).}
	\label{fig:CEFluid}
\end{figure}

We now consider two circular explosions. The first one corresponds to the radial extension of the Sod shock tube problem \cite{Lagrange2D}, i.e. an inviscid flow with $c_{s}=c_{h}=0$ and $\mu =\kappa = 0$, and initial condition given by
%one in the inviscid flow limit of the model and another one for a solid medium. 
%The initial conditions are given by
\begin{gather*}
	\rho\left(\x,0\right) =  \left\lbrace \begin{array}{lr}
		1 & \mathrm{ if } \; r \le 0.5,\\
		0.125 & \mathrm{ if } \; r > 0.5,
	\end{array}\right. \qquad
	\mathbf{u}\left(\x,0\right) =  \boldsymbol{0}, \qquad
	\press\left(\x,0\right) = \left\lbrace \begin{array}{lr}
		1 & \mathrm{ if } \; r \le 0.5,\\
		0.1 & \mathrm{ if } \; r > 0.5,
	\end{array}\right.
	\qquad r = \sqrt{x^2+y^2}.
\end{gather*}
The results obtained at $t_{\mathrm{e}}=0.25$ using the hybrid FV/FE method with ENO limiting are depicted in Figure~\ref{fig:CEFluid}. For comparison, we include a reference solution computed solving, with a TVD-FV scheme, the 1D radial PDE with source terms equivalent to the compressible Euler system \cite{DT11}.

\begin{figure}[H]
	\centering
	\includegraphics[clip,trim= 0 5 0 15,width=0.42\linewidth]{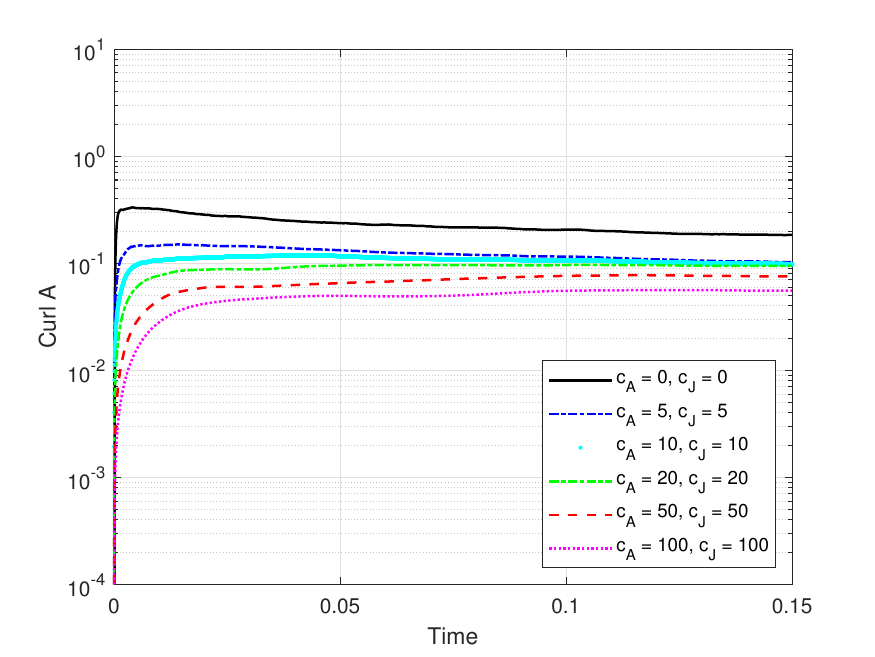}
	\includegraphics[clip,trim= 0 5 0 15,width=0.42\linewidth]{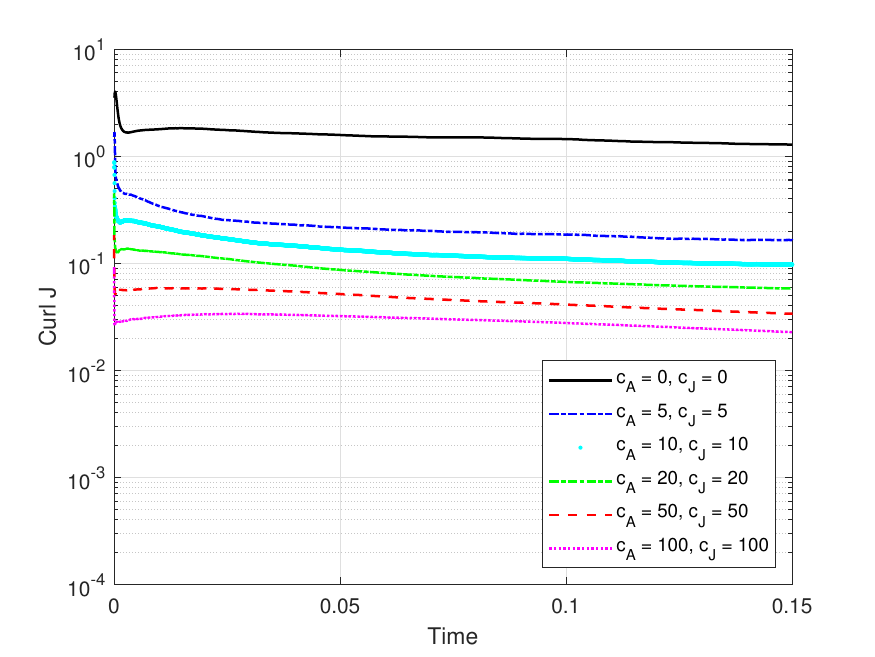}
	\caption{Solid circular explosion. Time evolution of the curl errors for the distortion field, $\bA$, (left) and the thermal impulse, $\bJ$, (right) for cleaning speeds $c_{\bA},\, c_{\bJ} \in \left\lbrace0, 5, 10, 20, 50, 100 \right\rbrace$. }
	\label{fig_ce85r0curla}
\end{figure}

\begin{figure}[H]
	\centering
	\includegraphics[clip,trim= 0 250 0 50,width=0.4\linewidth]{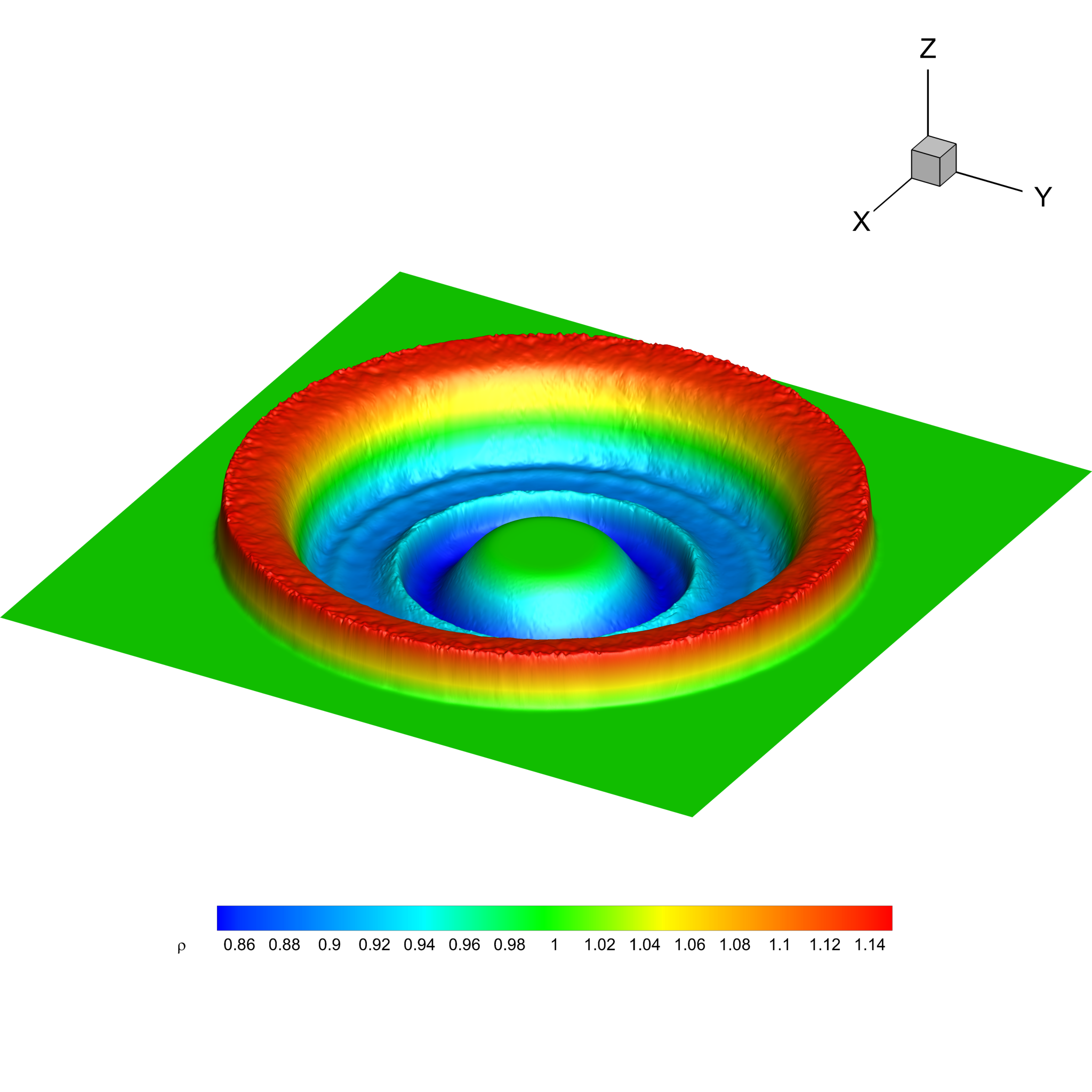}
	\includegraphics[clip,trim= 0 250 0 50,width=0.4\linewidth]{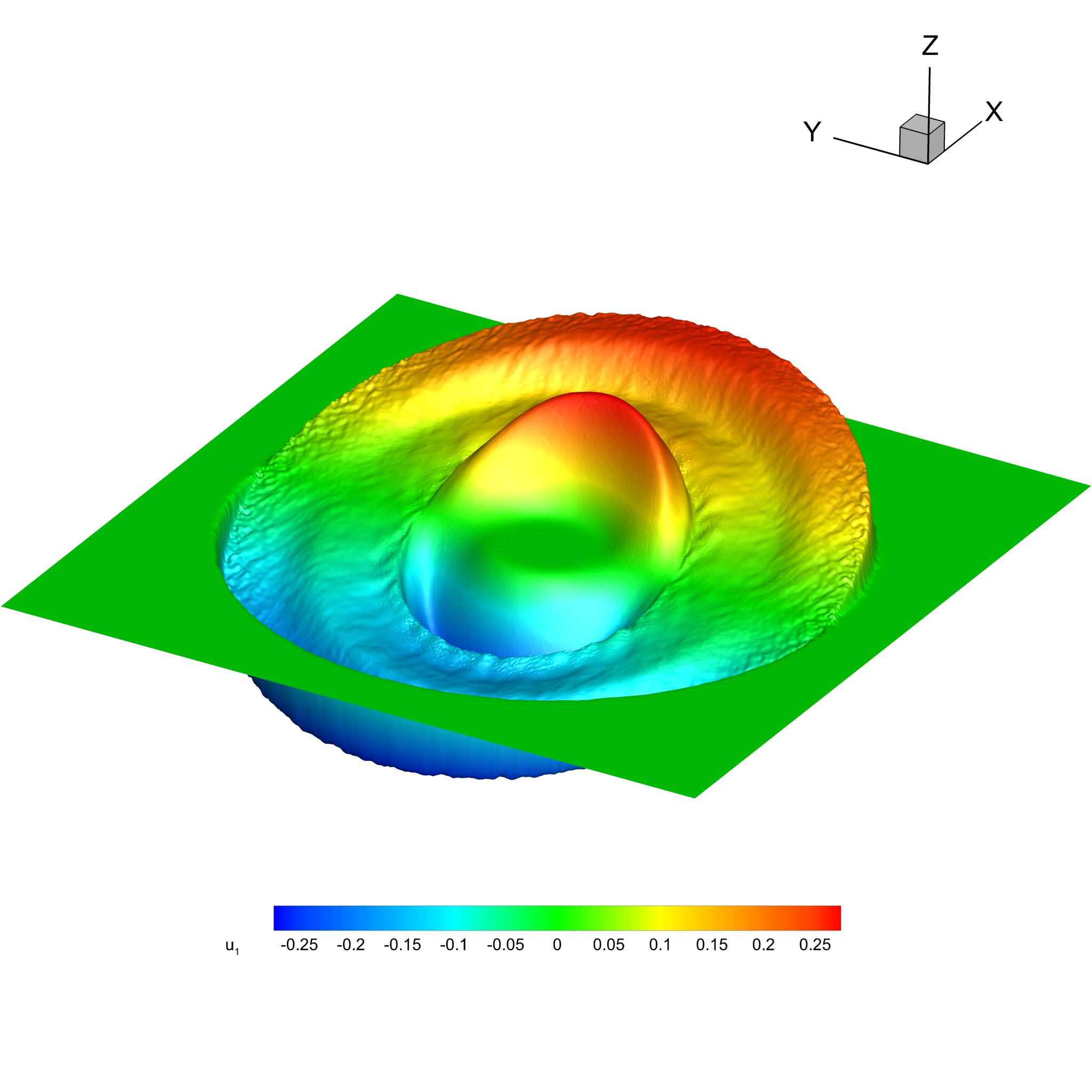}
	\includegraphics[clip,trim= 0 250 0 50,width=0.4\linewidth]{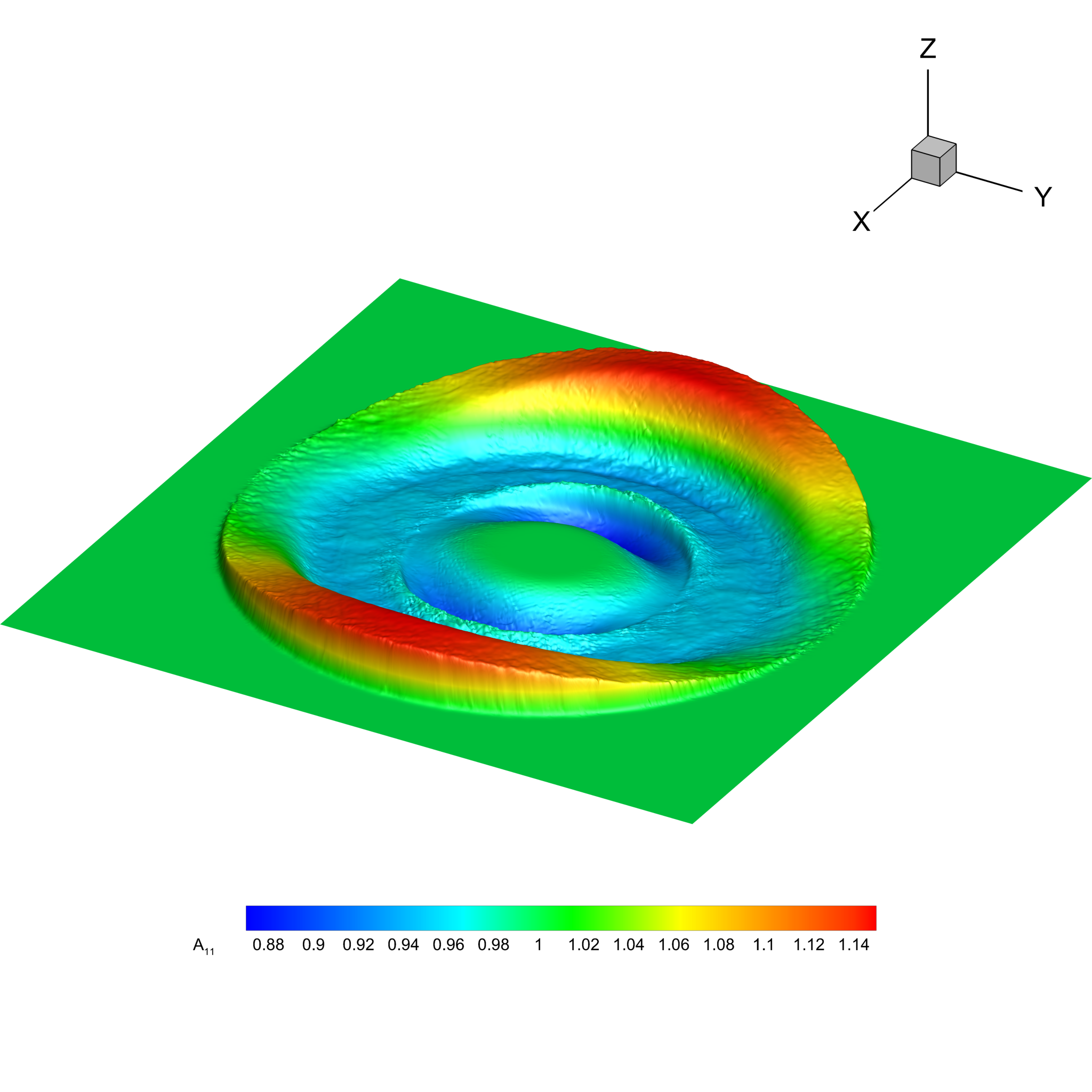}
	\includegraphics[clip,trim= 0 250 0 50,width=0.4\linewidth]{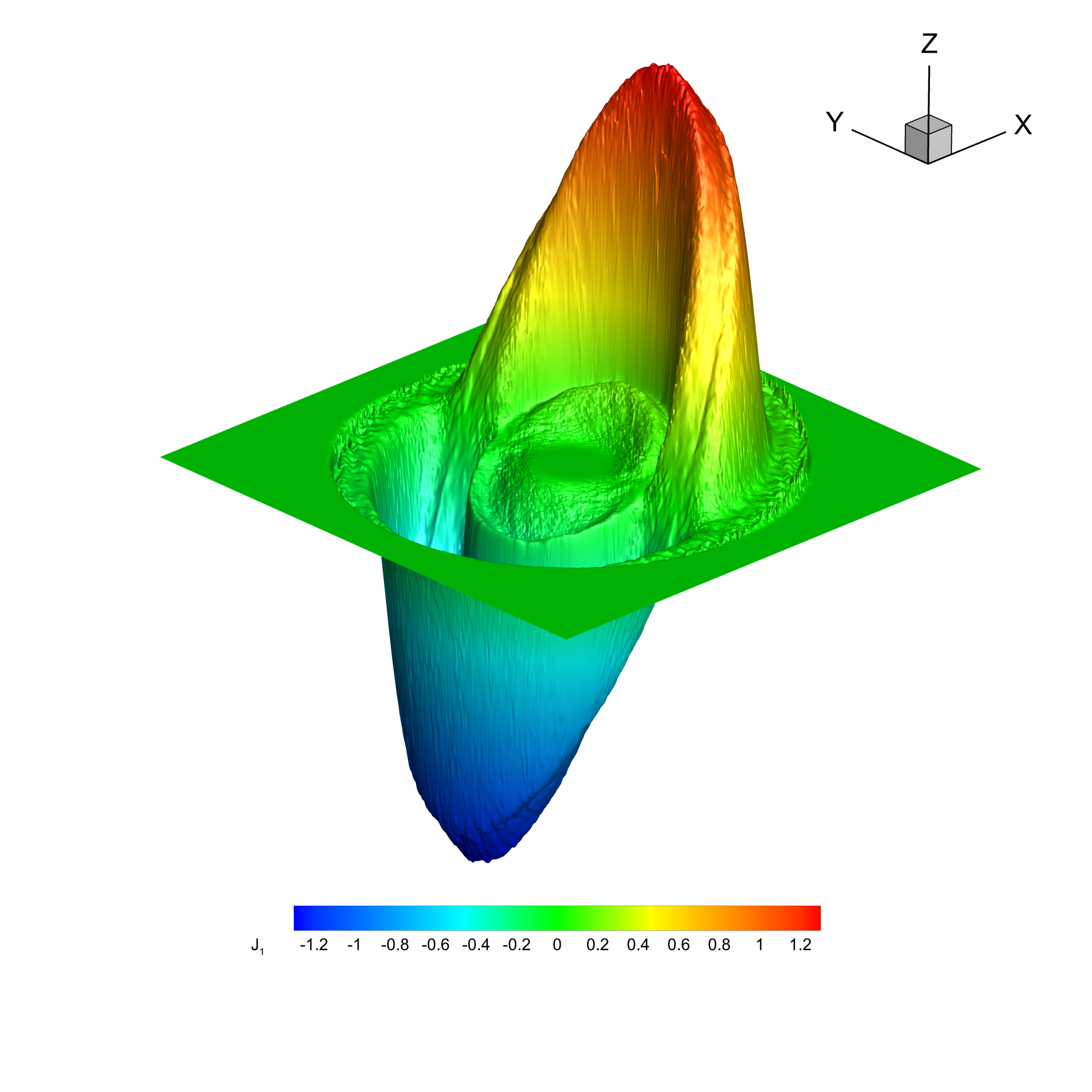}
	\caption{Solid circular explosion 2D. Elevated contour plots of the density, $\vel_{1}$, $A_{11}$ and $J_1$ for $\left\lbrace (x,y)\in\mathbb{R}^{2}\,  |\,x \in[0,1],\, y=0 \right\rbrace$ using the hybrid FV/FE method for the compressible GPR model.}
	\label{fig:CESolid3D}
\end{figure}

\begin{figure}[H]
	\centering
	\includegraphics[width=0.42\linewidth]{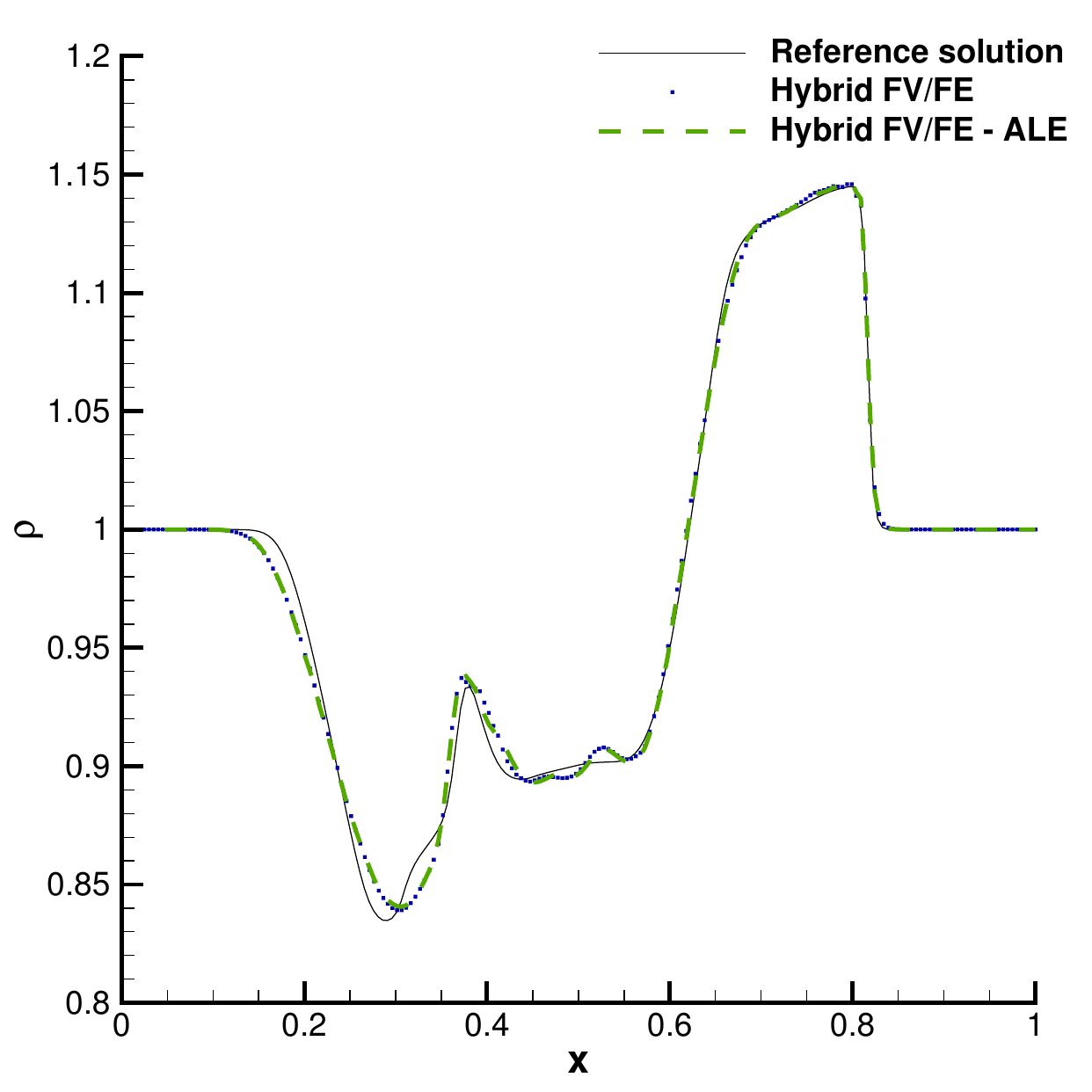}
	\includegraphics[width=0.42\linewidth]{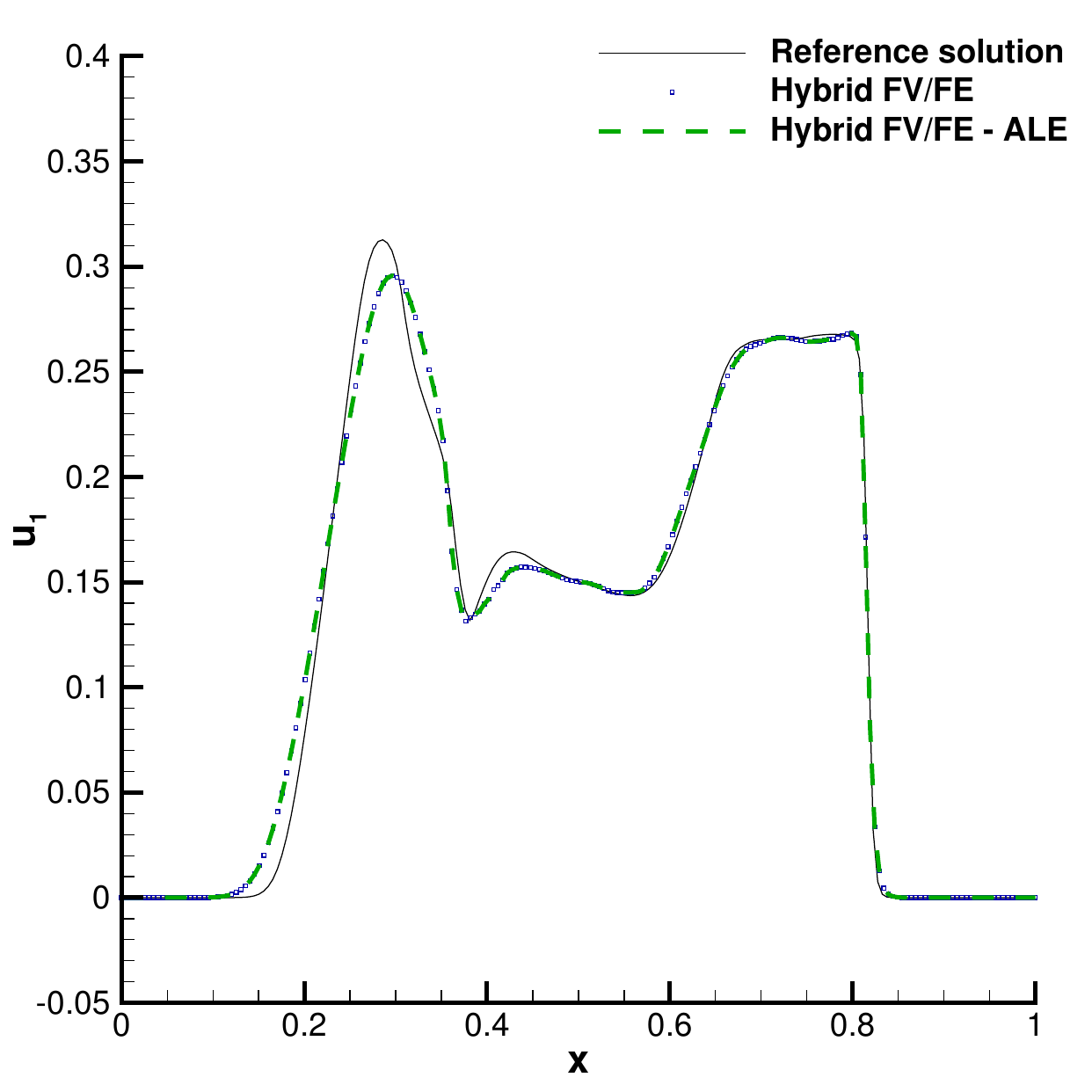}
	\includegraphics[width=0.42\linewidth]{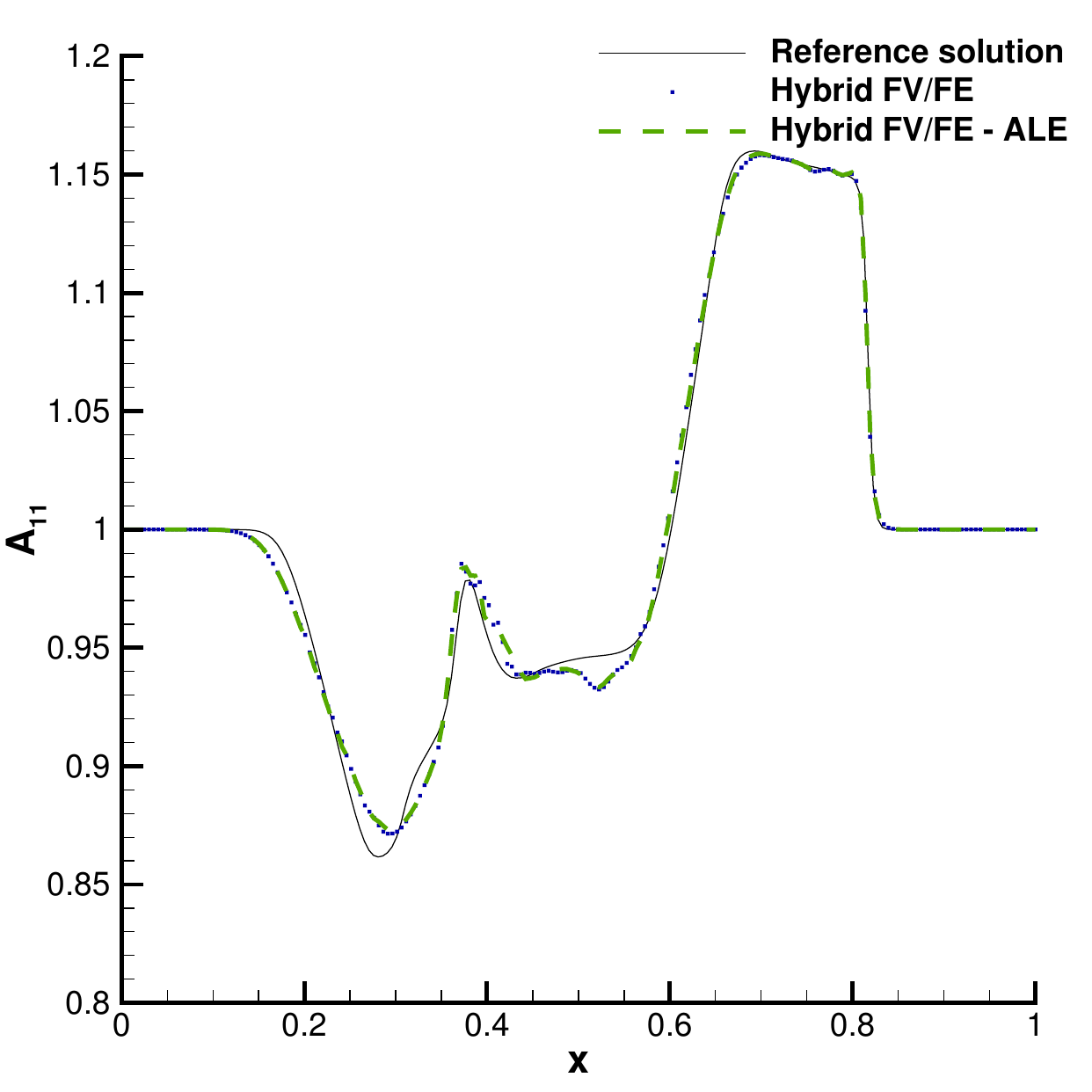}
	\includegraphics[width=0.42\linewidth]{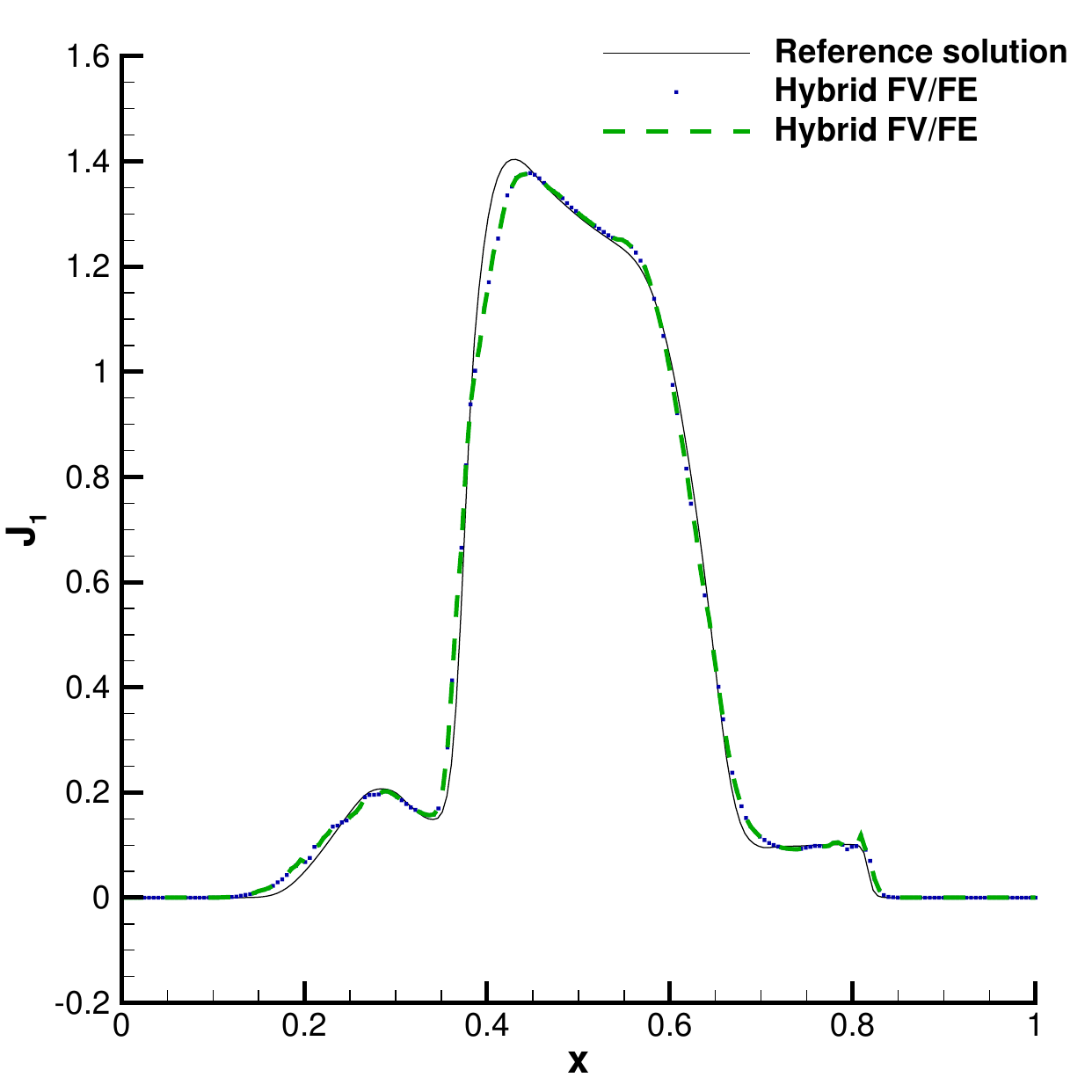}
	\caption{Solid circular explosion 2D. 1D cuts of the density, $\vel_{1}$, $A_{11}$ and $J_1$ for $\left\lbrace (x,y)\in\mathbb{R}^{2}\,  |\,x \in[0,1],\, y=0 \right\rbrace$ using the hybrid FV/FE method for the compressible GPR model with the Eulerian scheme (blue squares) and with the ALE method (green dashed line) compared against the reference solution (black line).}
	\label{fig:CESolid1D_ALE}
\end{figure}

The second circular Riemann problem regards a solid medium \cite{Boscheri2021SIGPR}. We set $\tau_1=\tau_2=10^{20}$, $\rho_{0}=1$, $c_v=1.0$, $c_s=1$, $c_h=0.5$, $\gamma=1.4$, and the initial condition
\begin{equation*}
	\rho\left(\x,0\right) = 1, \;\;  \bvel\left(\x,0\right)=  \boldsymbol{0},\;\; \press \left(\x,0\right) = \left\lbrace \begin{array}{lr}
		2 & \mathrm{ if } \; r \le 0.5,\\
		1 & \mathrm{ if } \; r > 0.5.
	\end{array}\right.
\end{equation*}
As aforementioned, in the solid limit the GPR model presents curl-free involution constraints for the distortion and heat flux fields. To analyse the behaviour of the proposed GLM cleaning strategy a set of simulations for $c_{\bA},\, c_{\bJ} \in \left\lbrace0, 5, 10, 20, 50, 100 \right\rbrace$ has been run for the solid circular explosion test case. The obtained curl errors are depicted in Figure~\ref{fig_ce85r0curla}. As expected, we observe the decay of $\curl \A$ and $\curl \J$ for increasing values of the cleaning speeds.
Besides, the contour plots of density, $\A_{11}$ and $\J_{1}$ obtained at time $t_{\mathrm{e}}=0.15$ are depicted in Figure~\ref{fig:CESolid3D}. To ease comparison with available reference data, the solution for a 1D cut along the positive part of the $x$-axis are reported in Figure~\ref{fig:CESolid1D_ALE}. The solutions are presented for both the fully Eulerian scheme and the ALE method with a smoothing parameter for the mesh velocity $\varsigma=5$.

\subsection{3D explosion}
\begin{figure}[H]
	\centering
	\includegraphics[width=0.32\linewidth]{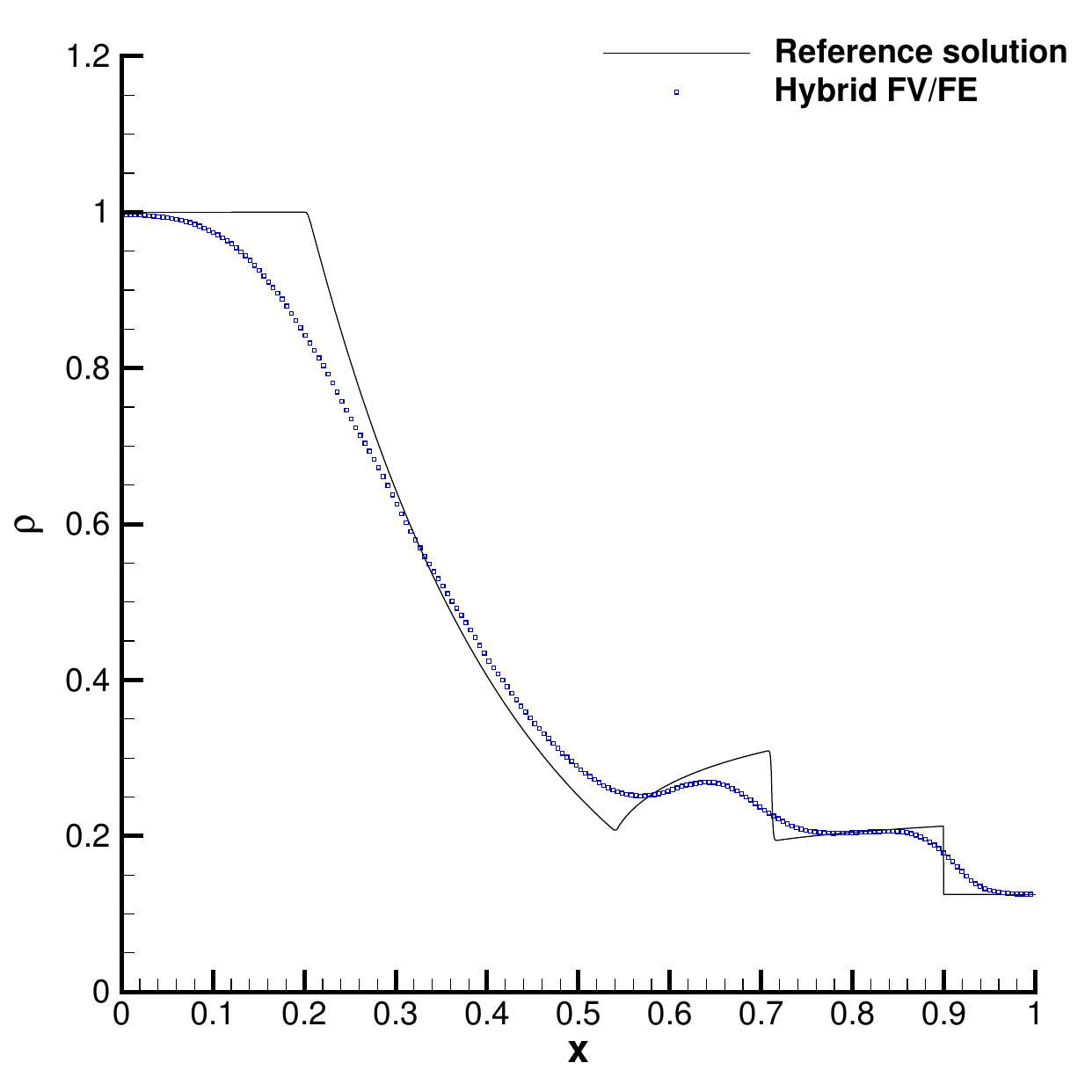}
	\includegraphics[width=0.32\linewidth]{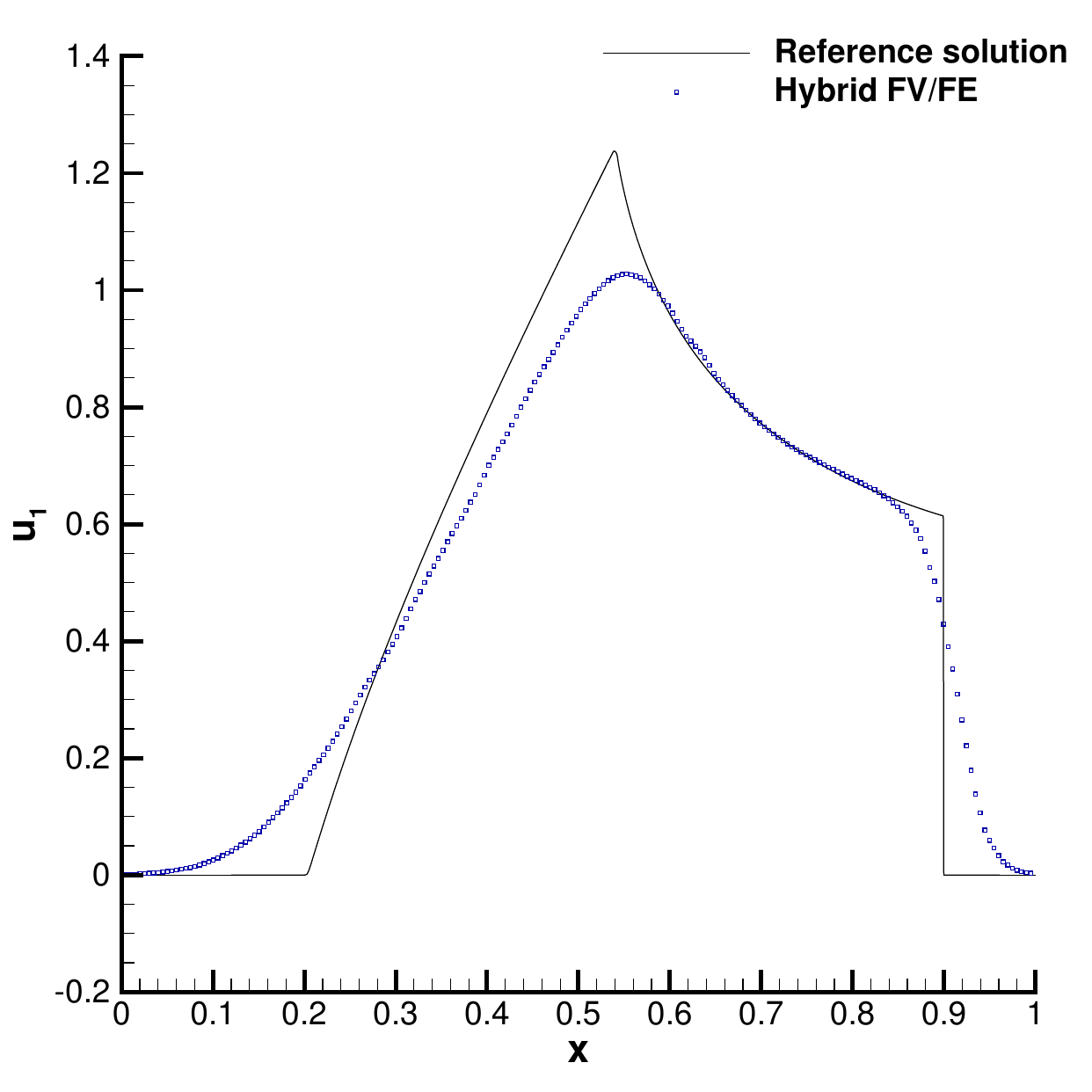}
	\includegraphics[width=0.32\linewidth]{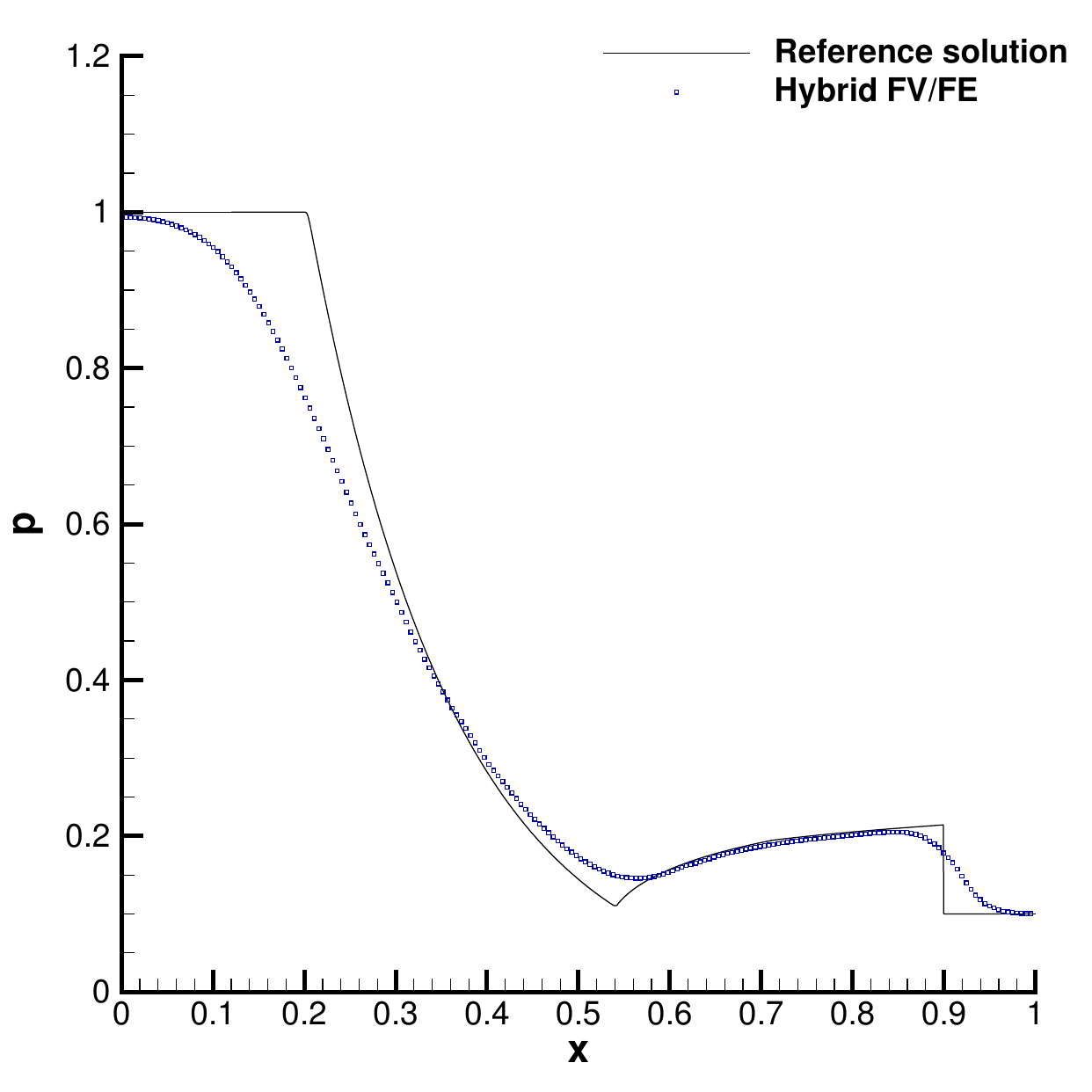}
	\caption{3D Sod explosion. 1D cuts of the density, $\vel_{1}$, and pressure for $\left\lbrace (x,y,z)\in\mathbb{R}^{3}\,  |\,x \in[0,1],\, y=0,\, z=0 \right\rbrace$ using the hybrid FV/FE method for the compressible GPR model (blue squares) compared against the reference solution (black line).}
	\label{fig:CE3Dcuts}
\end{figure}

\begin{figure}[H]
	\centering
	\includegraphics[width=0.32\linewidth]{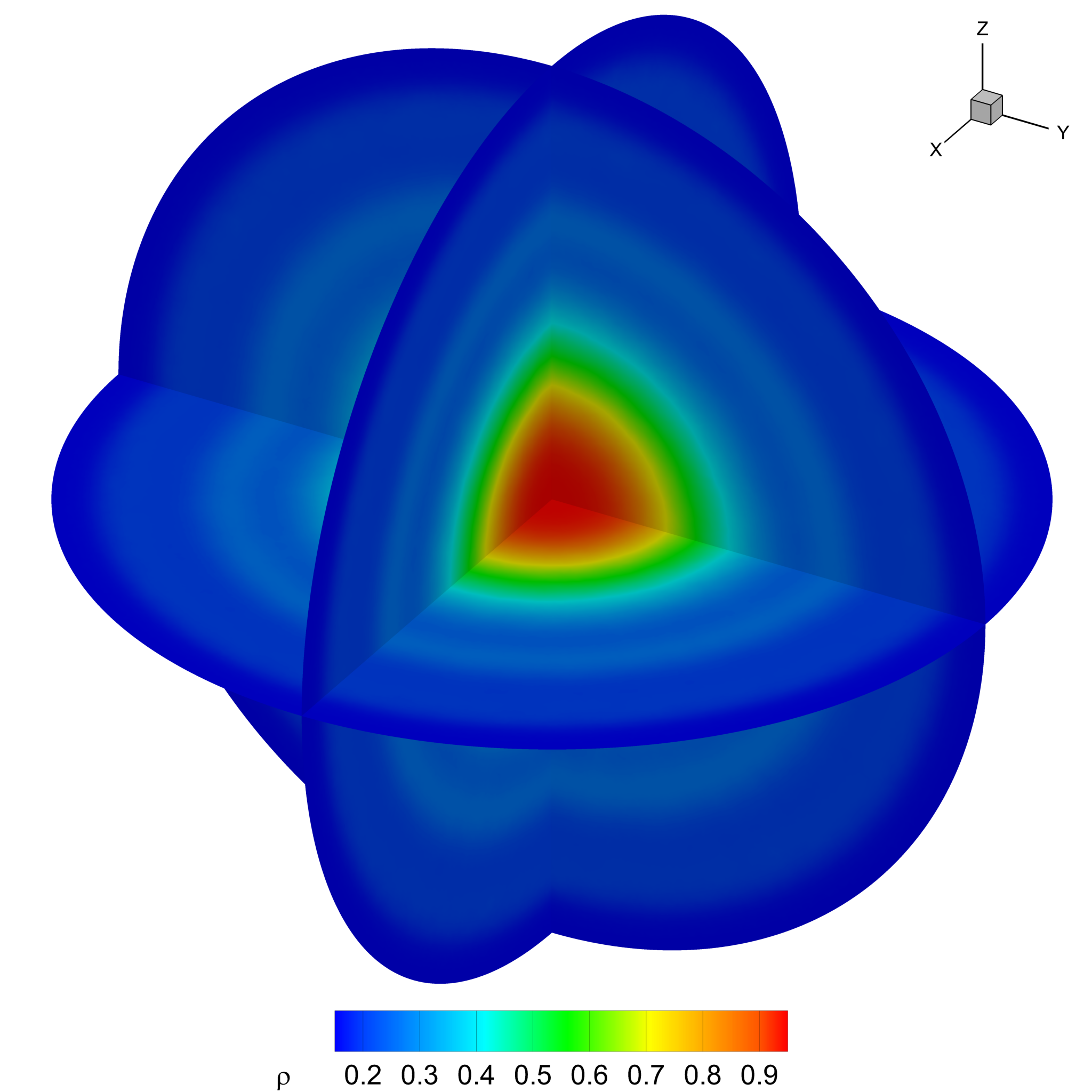}
	\includegraphics[width=0.32\linewidth]{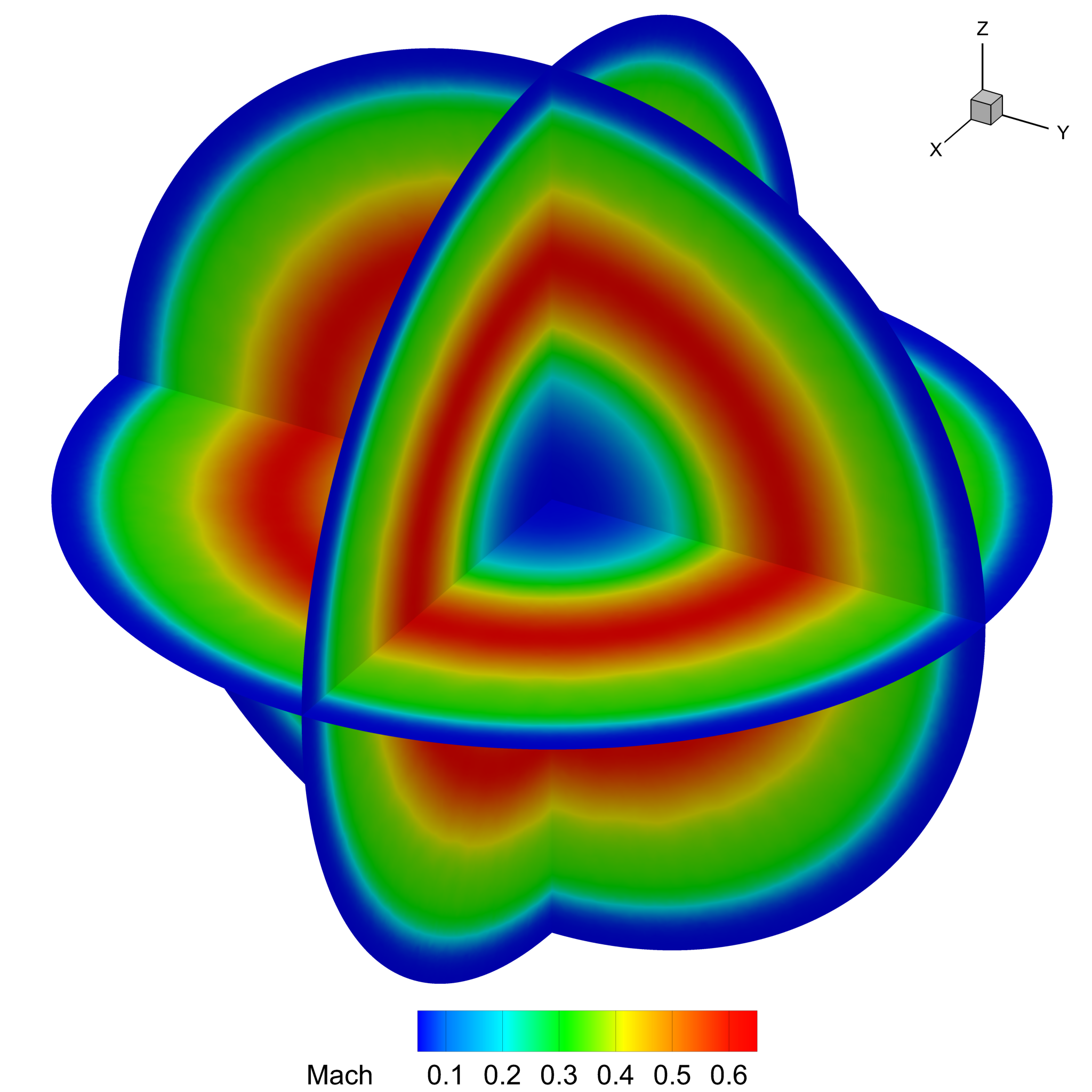}
	\includegraphics[width=0.32\linewidth]{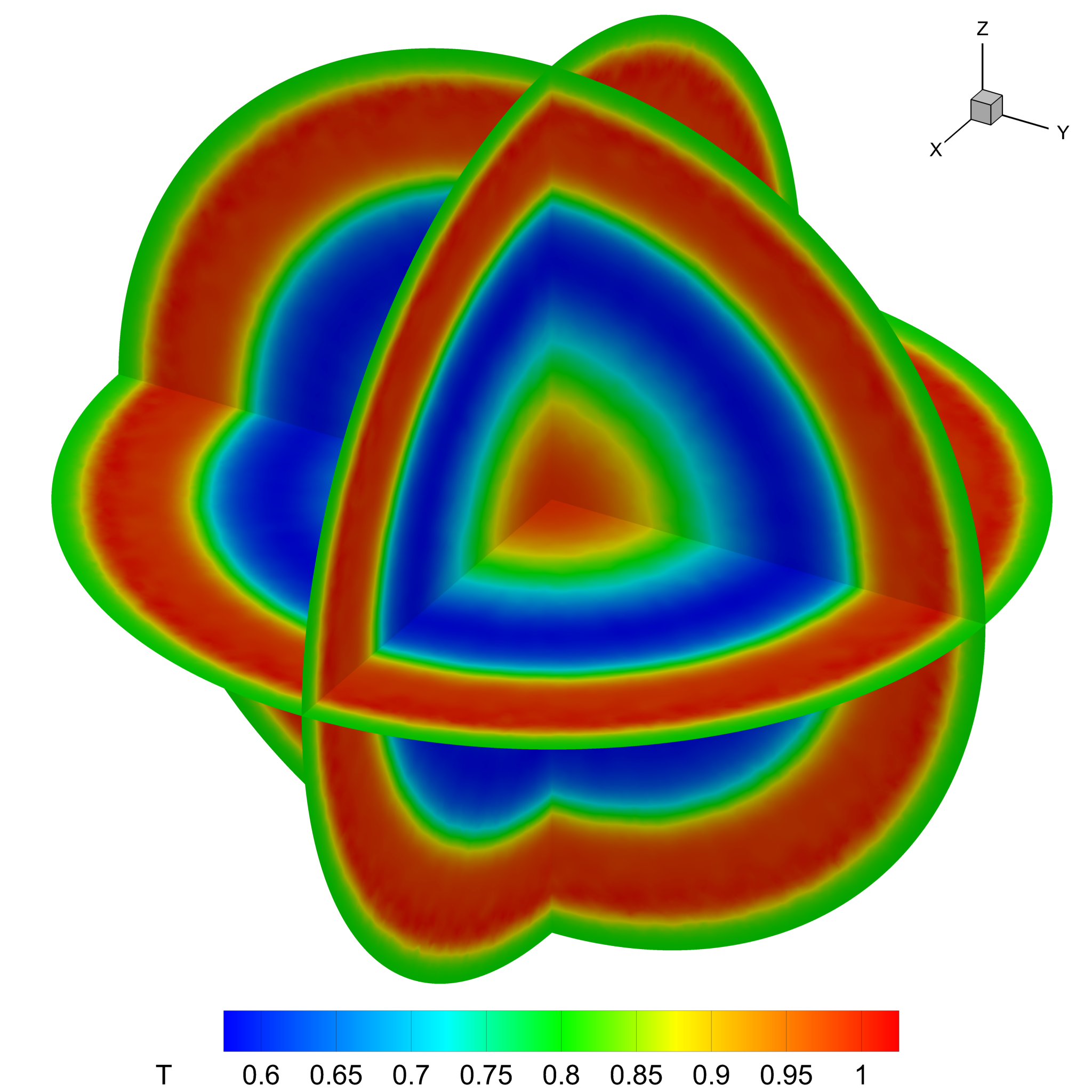}
	\includegraphics[width=0.32\linewidth]{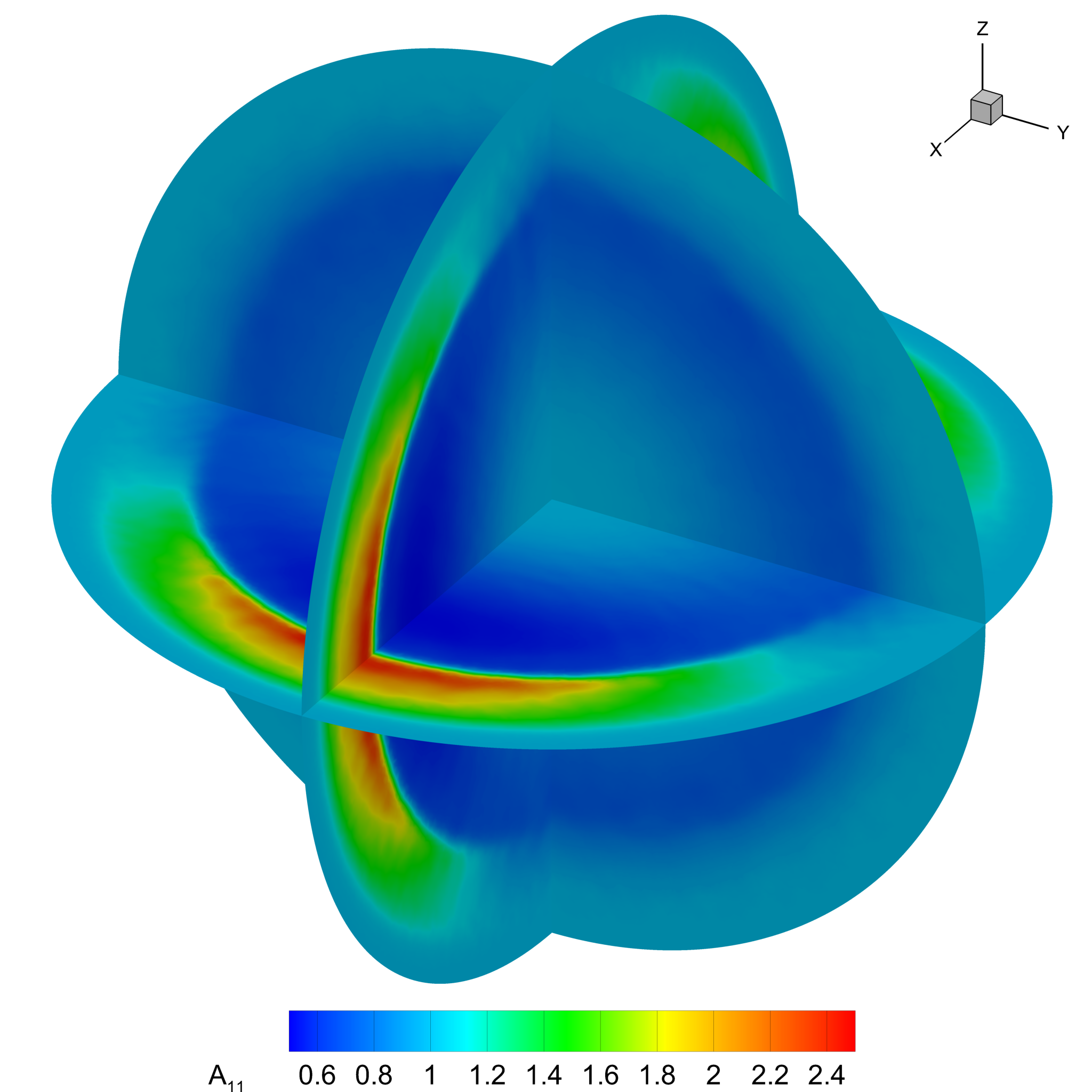}
	\includegraphics[width=0.32\linewidth]{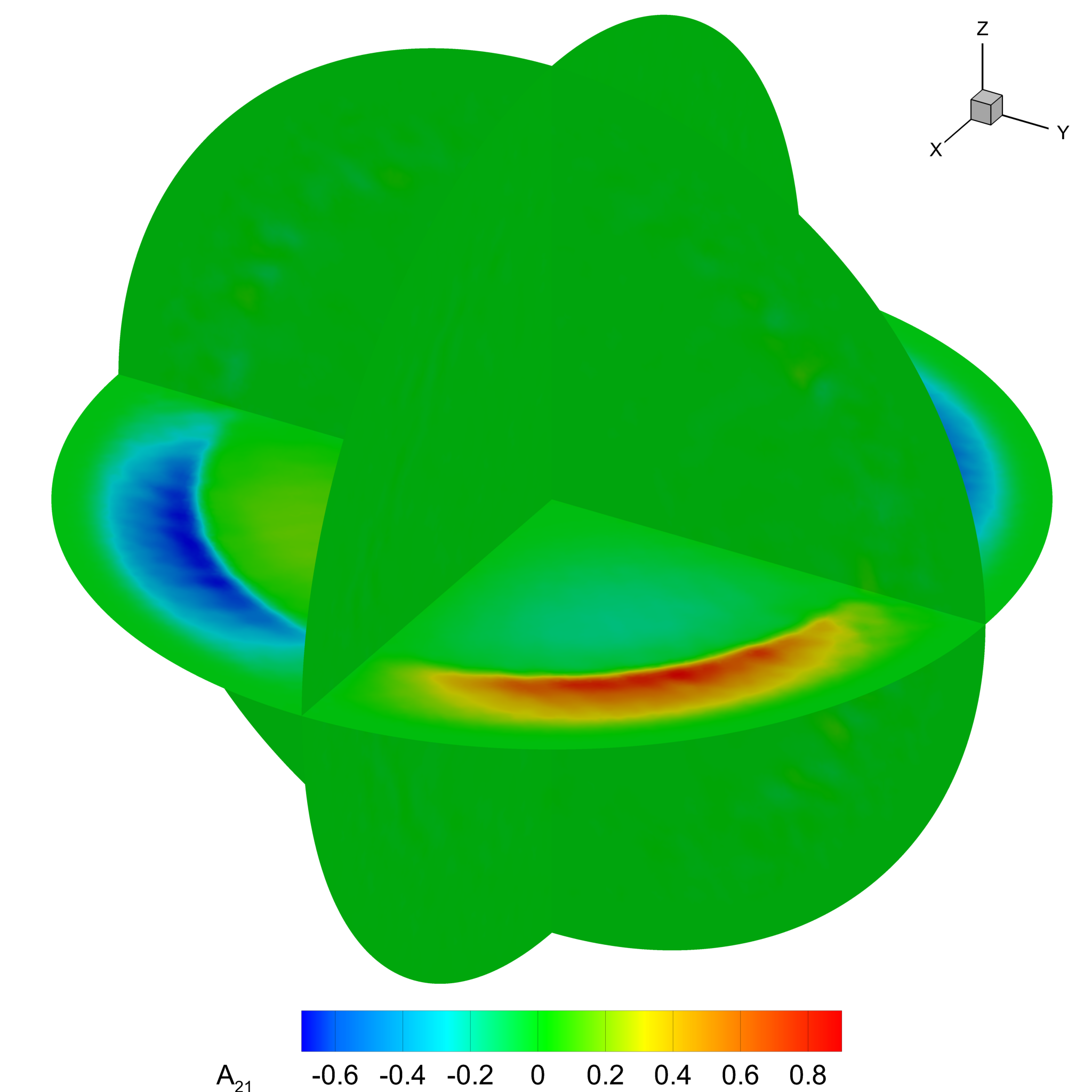}
	\includegraphics[width=0.32\linewidth]{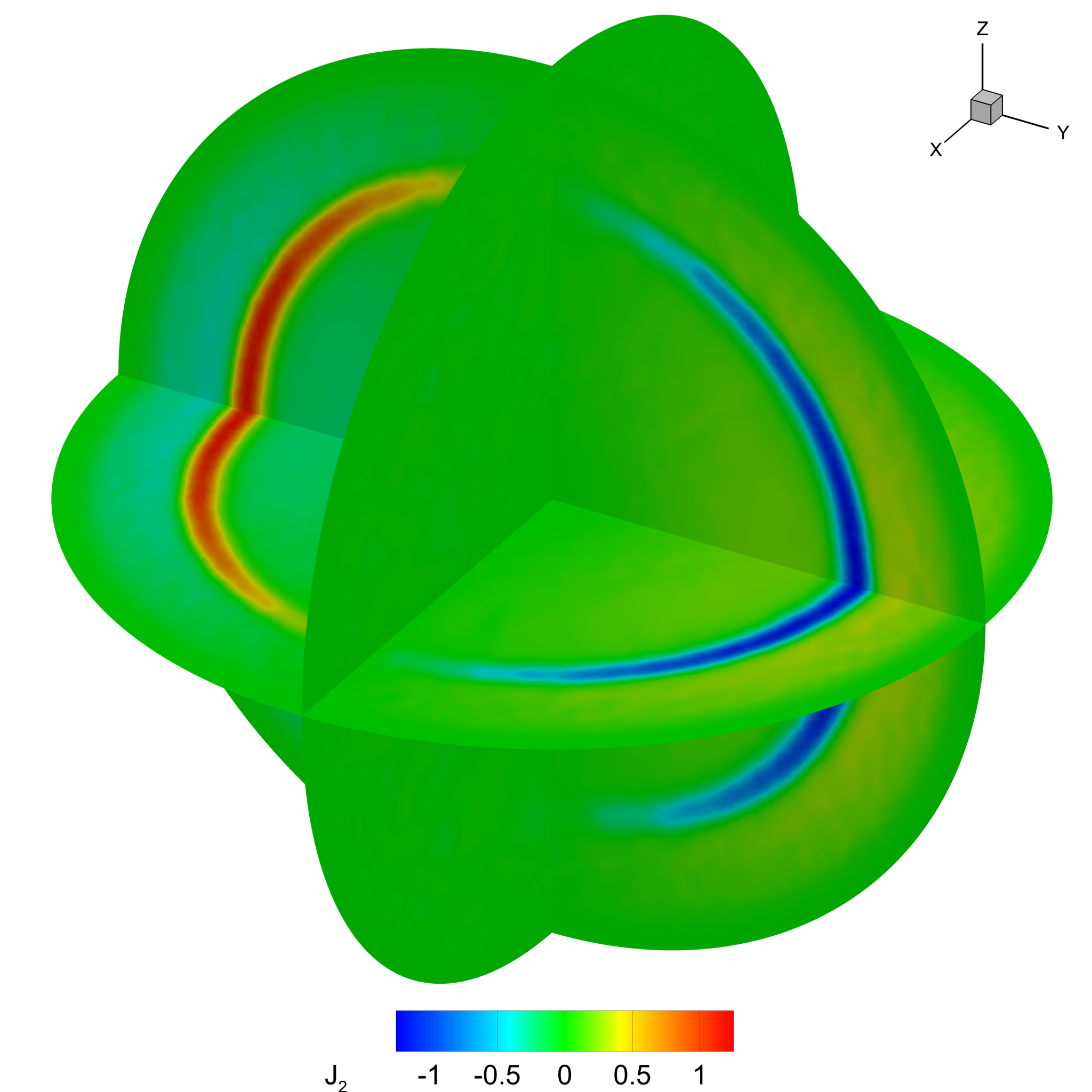}
	\caption{3D Sod explosion. Contour plots of the density, Mach number, temperature, $A_{11}$, $A_{21}$ and $J_{2}$ for the 2D slices $x=0$, $y=0$ and $z=0$.
	}
	\label{fig:CE3D}
\end{figure}

In this section, the 3D extension of the Sod problem is studied. As computational domain, we define the unit sphere centred at the origin which is discretized employing a primal mesh made of $8082535$ tetrahedral elements. The initial condition is given by 
\begin{gather*}
	\rho\left(\x,0\right) =  \left\lbrace \begin{array}{lr}
		1 & \mathrm{ if } \; r \le 0.5,\\
		0.125 & \mathrm{ if } \; r > 0.5,
	\end{array}\right. \qquad
	\press \left(\x,0\right) = \left\lbrace \begin{array}{lr}
		1 & \mathrm{ if } \; r \le 0.5,\\
		0.1 & \mathrm{ if } \; r > 0.5,
	\end{array}\right. \qquad
	\bvel \left(\x,0\right) =\boldsymbol{0}, \qquad r=\sqrt{x^2+y^2+z^2},
\end{gather*}
while the model parameters read $c_{s}=c_{h}=0$, $\mu =\kappa = 0$, $c_v=2.5$ and $\gamma=1.4$. Dirichlet boundary conditions are defined in the outer boundary. The simulation is run employing the ENO limiting strategy based on physical variables and setting $c_{\alpha}=1$. The solution obtained at $t_{\mathrm{e}}=0.25$ is depicted in \mbox{Figures~\ref{fig:CE3Dcuts}--\ref{fig:CE3D}.} 
We observe a good agreement with the reference solution computed with a TVD-FV solver of the corresponding 1D radial PDE. 

\subsection{Shear motion}
\begin{figure}[H]
	\centering
	\includegraphics[width=0.45\linewidth]{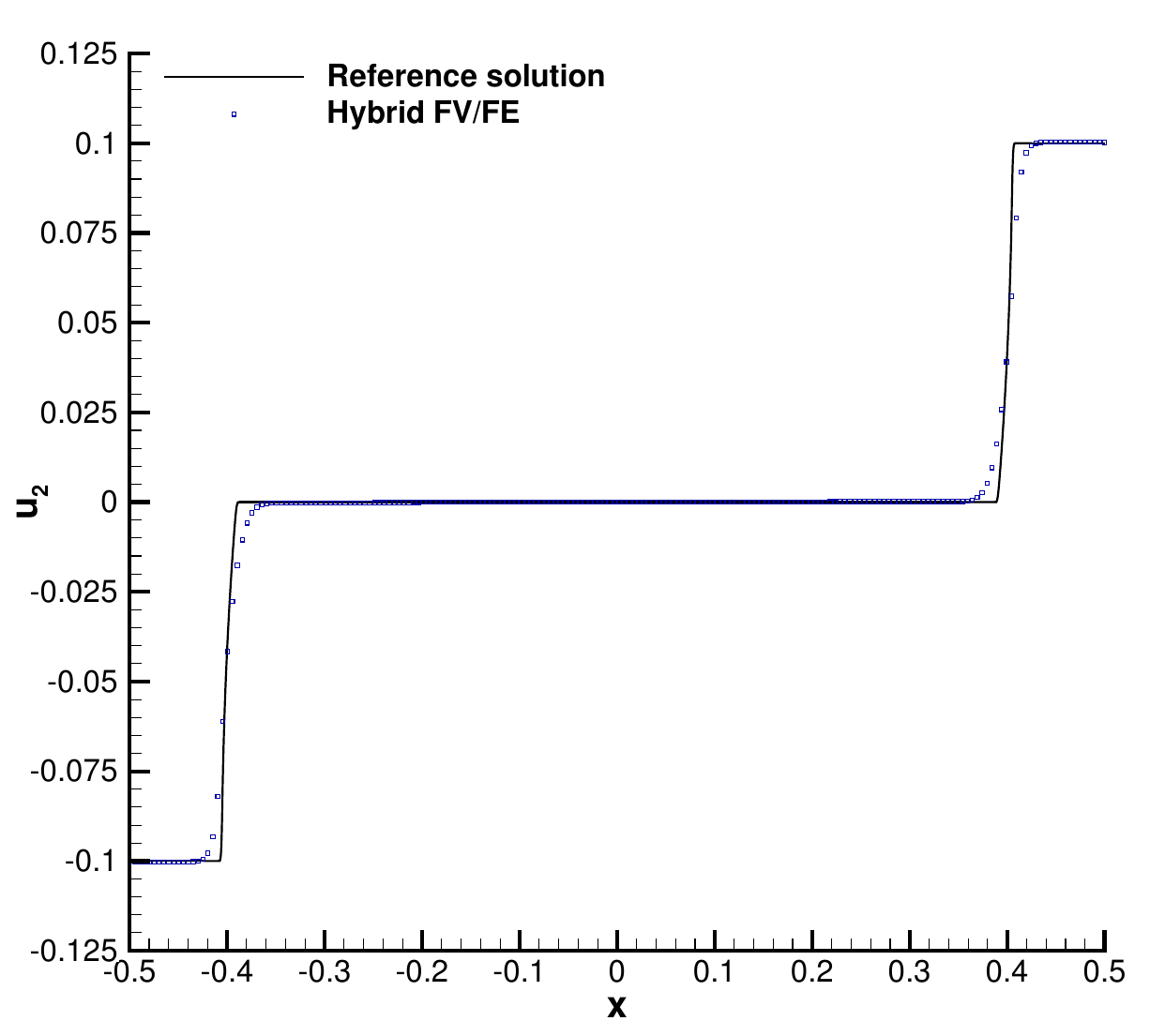}
	\includegraphics[width=0.45\linewidth]{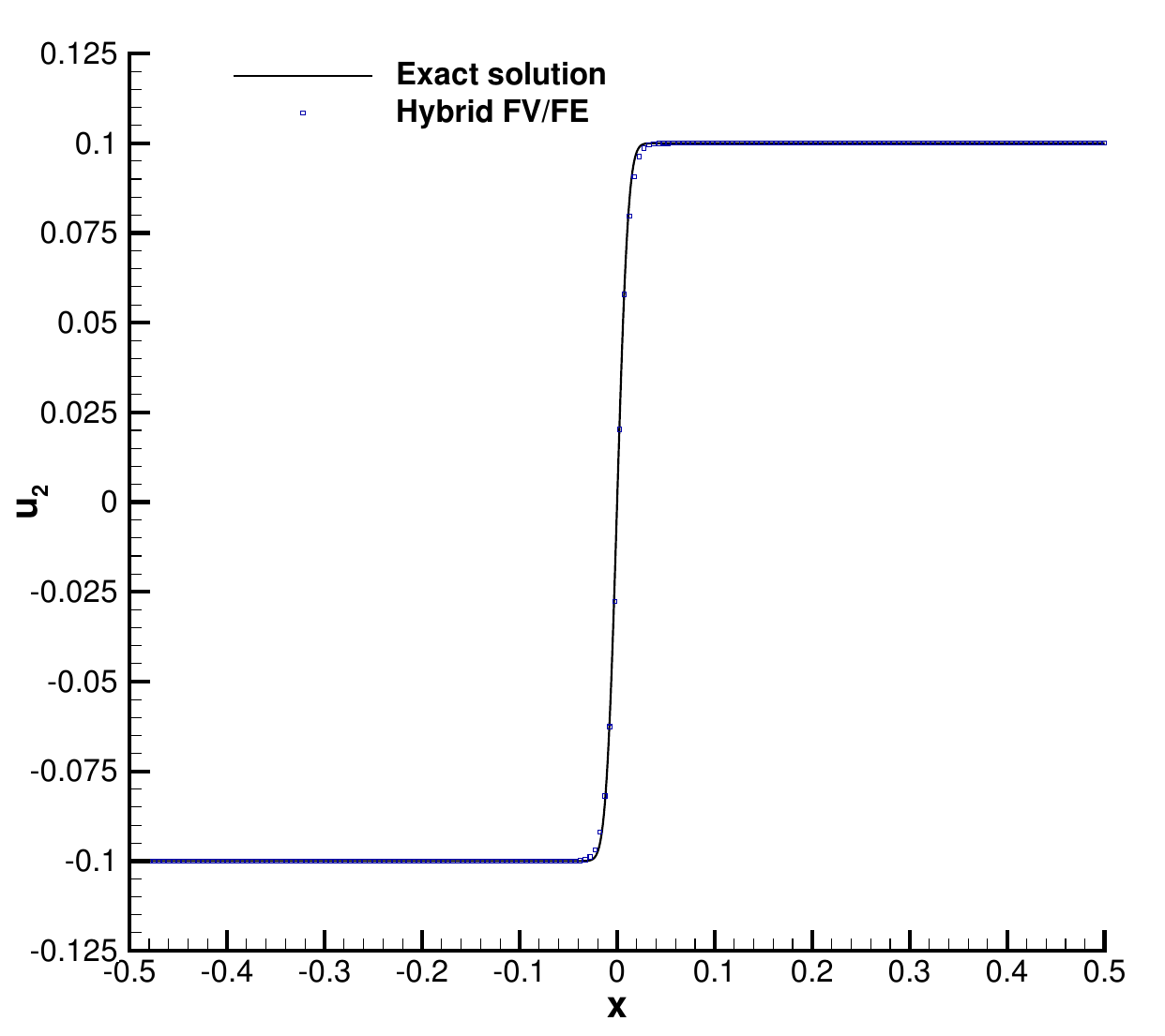}
	\includegraphics[width=0.45\linewidth]{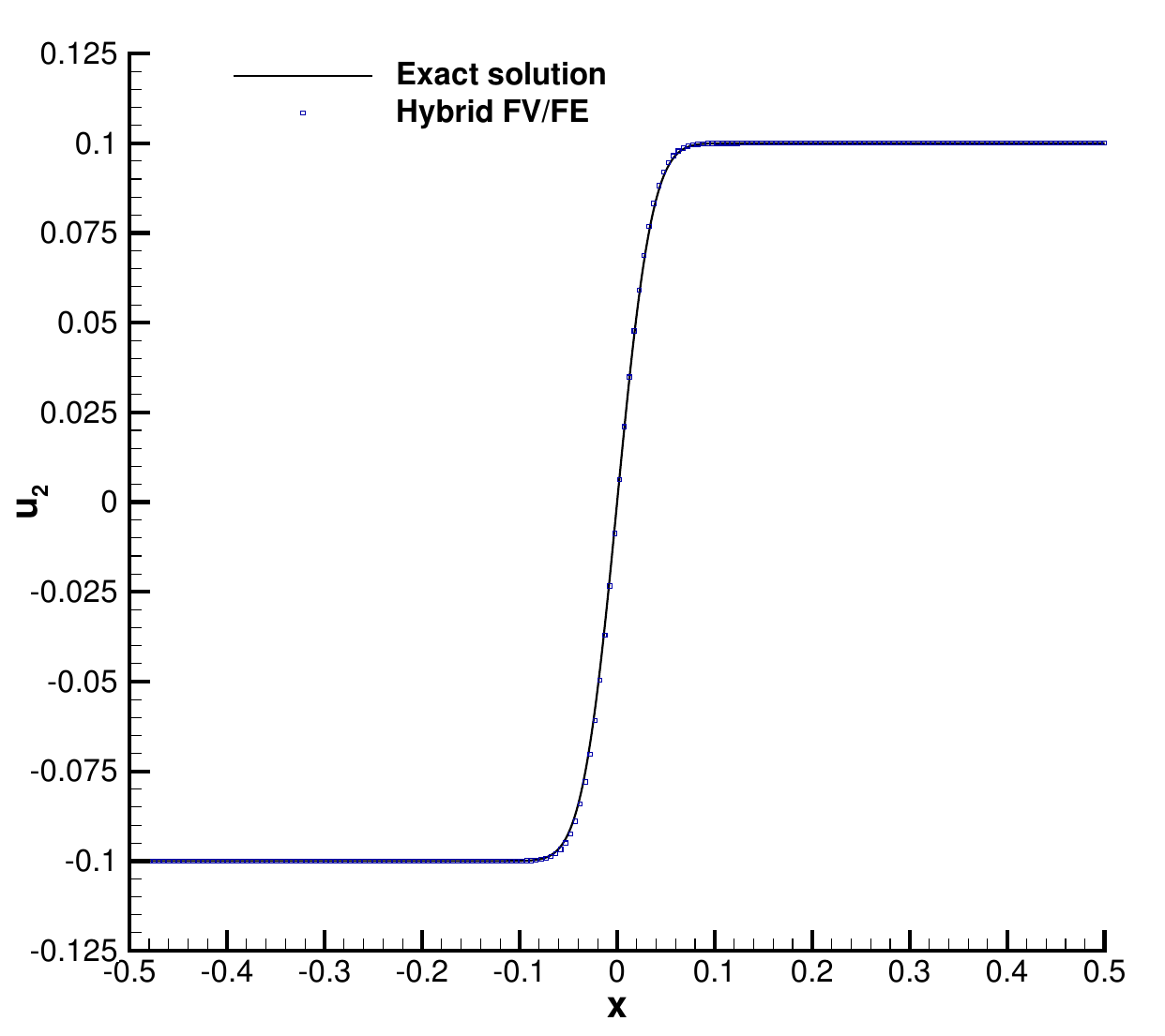}
	\includegraphics[width=0.45\linewidth]{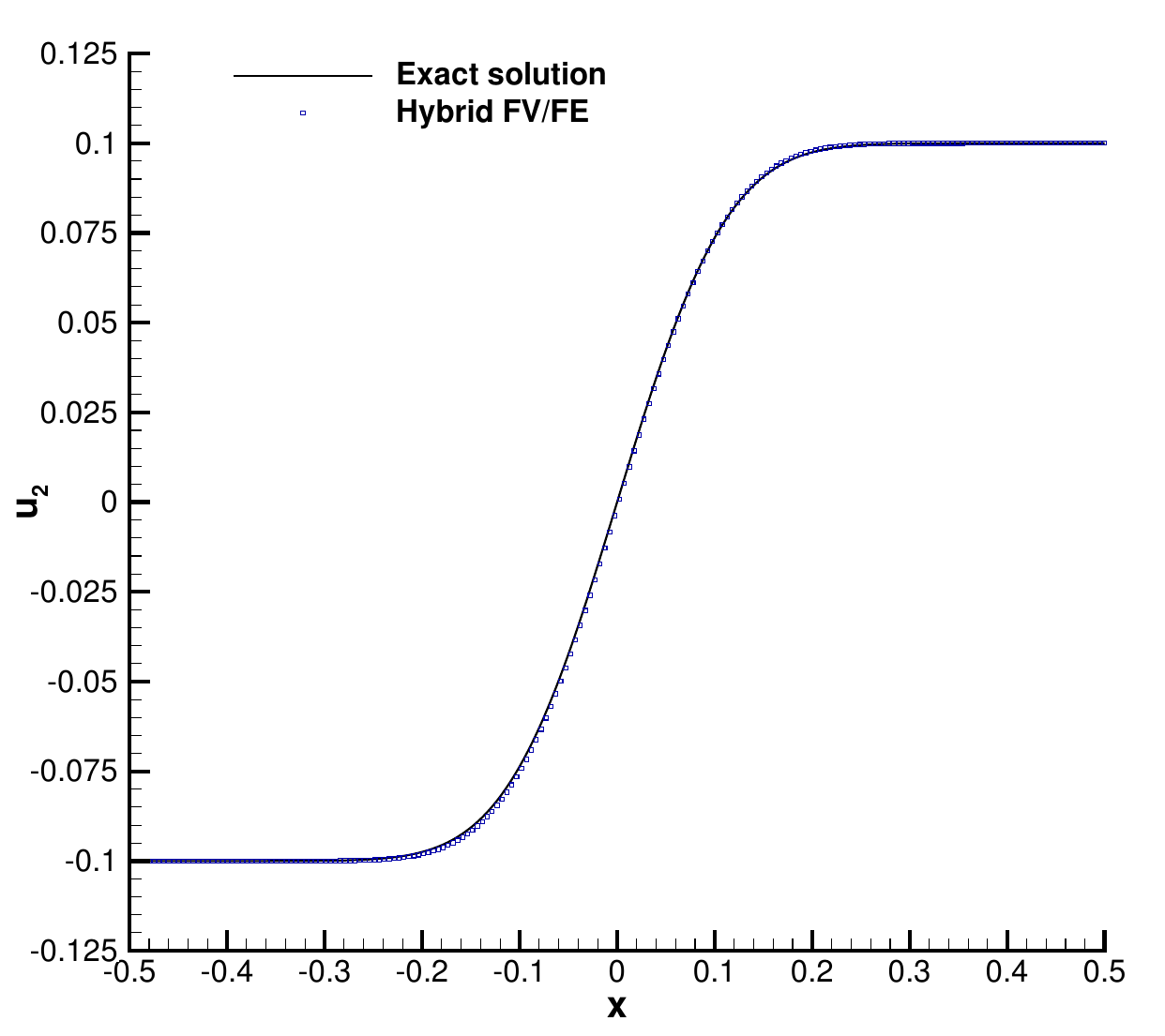}
	\caption{First Stokes. 1D cuts along the $x$-direction of the velocity component $\vel_{2}$ computed using the hybrid FV/FE method for the compressible GPR model (blue squares) compared against the exact and reference solutions (black line). From top left to right bottom: shear solid and first Stokes with $\mu\in \left\lbrace  10^{-4},  10^{-3}, 10^{-2} \right\rbrace$.}
	\label{fig:FS}
\end{figure}

To analyse the behaviour of the proposed approach for different shears, we consider a shear test in solid mechanics and the First Stokes problem for fluid dynamics with viscosities $\mu\in\left\lbrace 10^{-4}, 10^{-3}, 10^{-2} \right\rbrace$. Further, we take $c_{s}=c_{h}=1$, $c_v=2.5$ and $\gamma=1.4$ and $\mu=\kappa = 10^{20}$ for the solid test case. The initial conditions are given by a Riemann-type problem with 
\begin{gather*}
	\rho\left(\x,0\right) = 1,\quad
	\press \left(\x,0\right) = \frac{1}{\gamma}, \quad
	{u}_{1} \left(\x,0\right) = 0, \quad
	{u}_{2} \left(\x,0\right) = \left\lbrace \begin{array}{lc}
		-0.1 & \mathrm{ if } \; y \le 0,\\
		0.1 & \mathrm{ if } \; y > 0.
	\end{array}\right.
\end{gather*}
We define the computational domain $ \Omega=[-0.5,0.5]\times[-0.05,0.05]$ which is discretized employing a regular Cartesian triangular grid made of $N_{x}=400$ divisions along the $x$-axis. We define strong Dirichlet boundary conditions at the left and right boundaries of the domain while periodic conditions are imposed in $y$-direction. The 1D cut of the velocity field component $\vel_{2}$ approximated at time $t_{\mathrm{e}}=0.4$ is depicted in Figure~\ref{fig:FS}.
For the fluid test cases, the known analytical solution 
\begin{equation*}
	\vel_{2} \left(\x,t\right) = \frac{1}{10} \mathrm{erf}\left( \frac{x}{2\sqrt{\mu t}}\right),
\end{equation*}
is also reported while a reference solution, computed using a TVD-FV scheme on a very fine one-dimensional grid, is included for validation of the shear solid test case. An excellent agreement is observed in all setups.

\subsection{Double shear layer}
\begin{figure}[H]
	\begin{center}
		\includegraphics[trim=5 5 5 5,clip,width=0.48\textwidth]{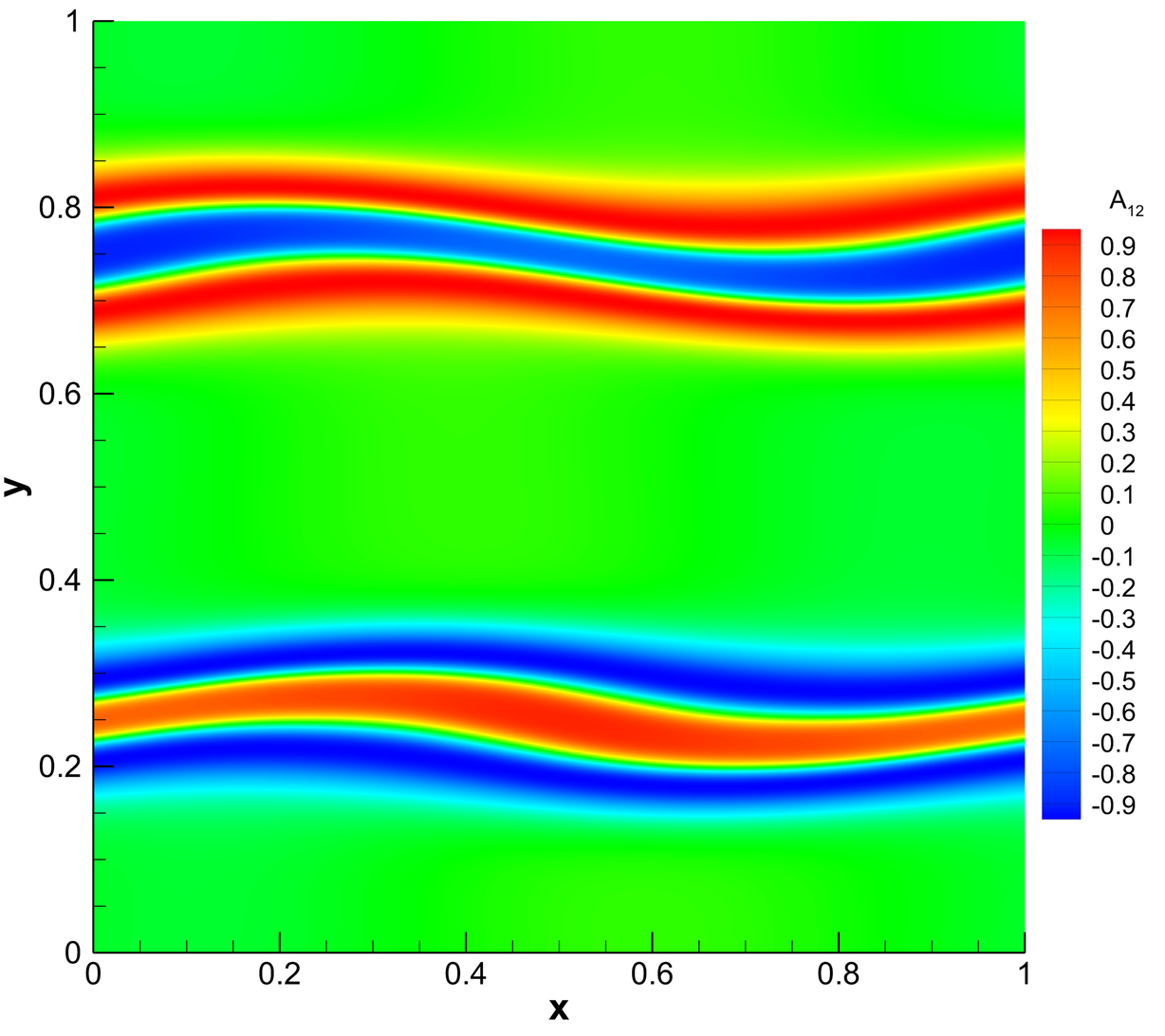} 
		\includegraphics[trim=5 5 5 5,clip,width=0.48\textwidth]{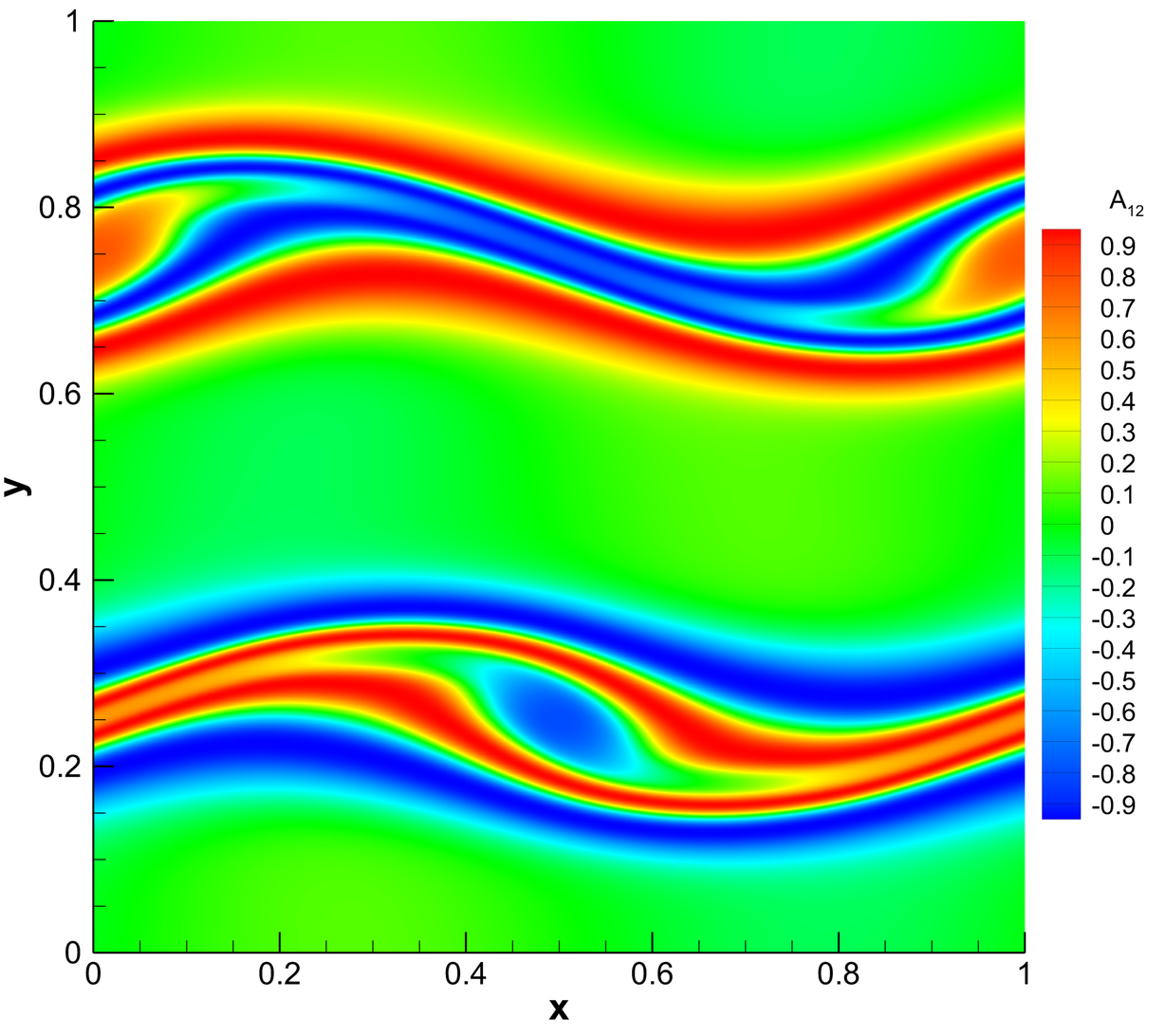}
		\includegraphics[trim=5 5 5 5,clip,width=0.48\textwidth]{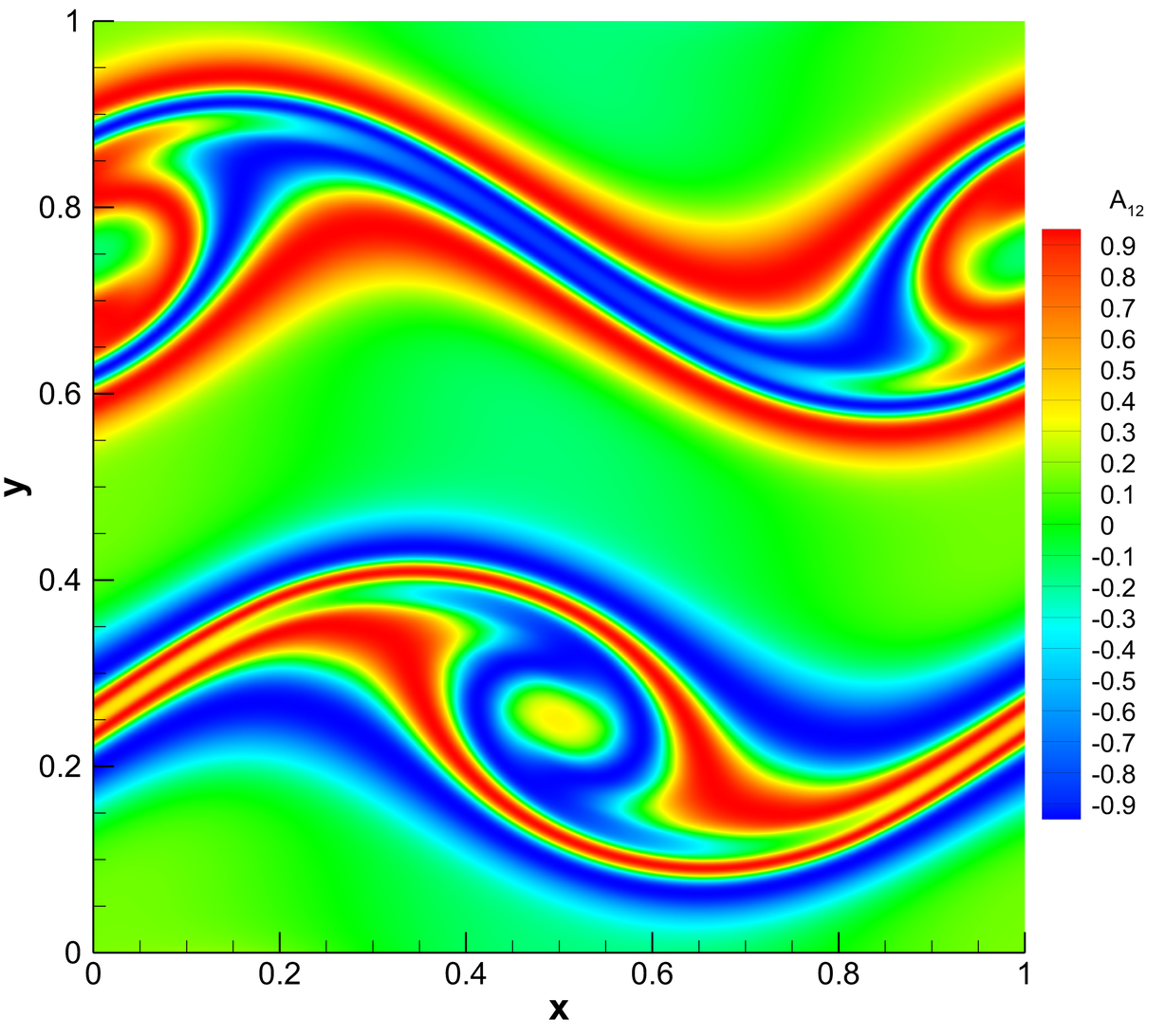}
		\includegraphics[trim=5 5 5 5,clip,width=0.48\textwidth]{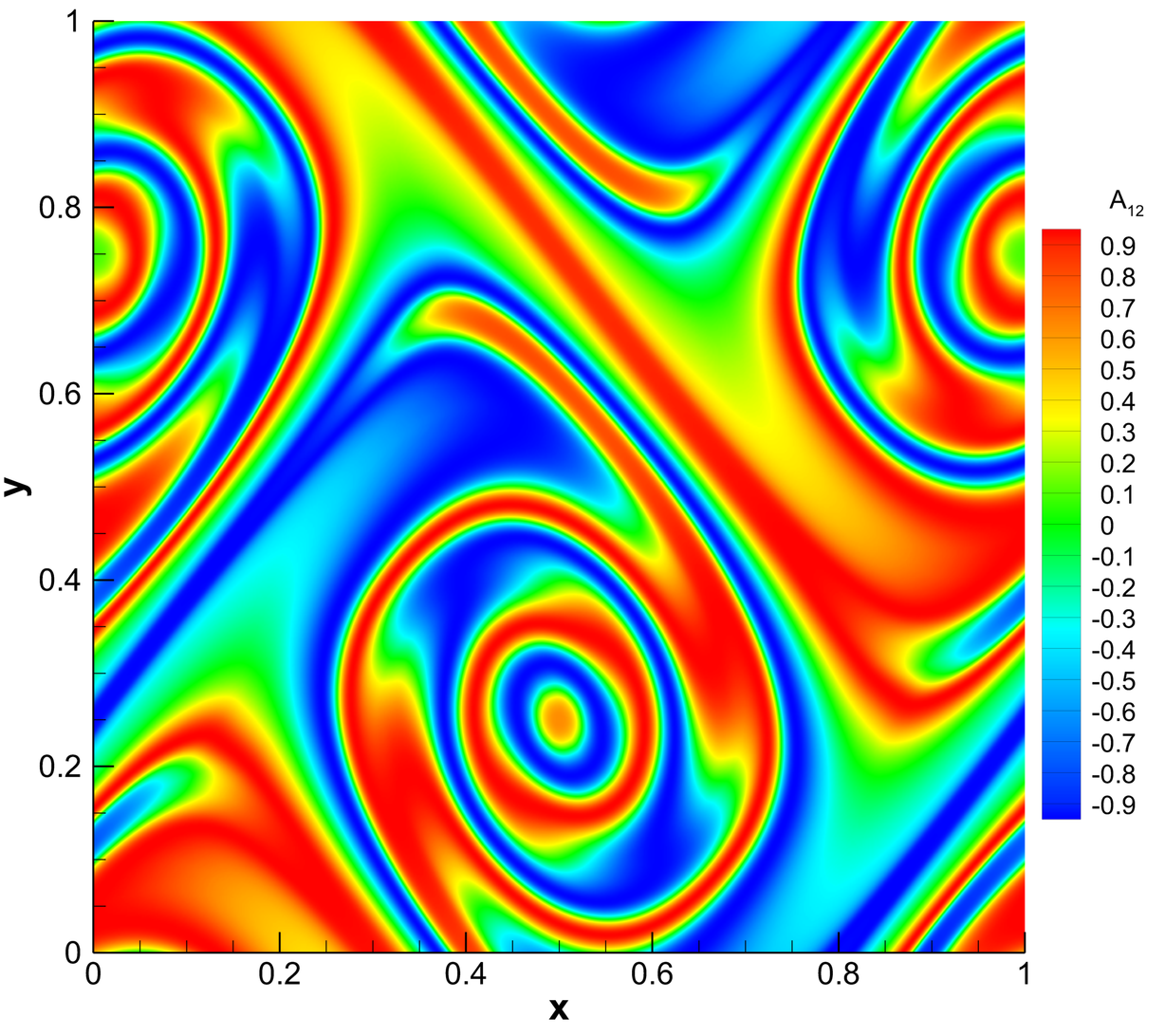}        
		\caption{Double shear layer. Contour plots of the distortion field component $A_{12}$ for $t\in\left\lbrace 0.4, 0.8, 1.2, 1.8\right\rbrace$.}  
		\label{fig:DSLdistortion}
	\end{center}
\end{figure}

The double shear layer is a classical test in fluid dynamics yielding complex flow patterns in the distortion field components when the GPR model is solved \cite{DPRZ16_GPRmodel,HTCA2}. Hence it allows verification of numerical methods in the presence of thin sharp flow structures. In particular, we consider the initial condition in the computational domain $\Omega=[0,1]^2$ given by 
\begin{gather*}
	\rho\left(\mathbf{x},0\right)  = 1,   \quad
	\vel_1\left(\mathbf{x},0\right) =\left\lbrace \begin{array}{ll}
		\tanh\left( 30 (y-0.25) \right) & \textrm{if } y \leq 0.5,\\
		\tanh\left( 30 (0.75-y) \right)  & \textrm{if } y > 0.5,
	\end{array}
	\right. \quad
	\vel_2\left(\mathbf{x},0\right)= 0.05 \sin(2\pi x), \qquad
	p\left(\mathbf{x},0\right) = 0
\end{gather*}
and we set the model parameters as $c_{p}=3.5$, $c_{s}=\sqrt{2}$, $c_{s}=8$, $c_{v}=2.5$, $\mu=\kappa=2\cdot 10^{-3}$. The second order Eulerian hybrid scheme is run up to time $t_{\mathrm{e}}=1.8$ with $\mathrm{CFL}=0.25$ using a mesh made of $524288$ primal elements. The contour plots of the distortion field component $\A_{12}$, reported in Figure~\ref{fig:DSLdistortion}, agree well with the solutions obtained for the incompressible GPR model and the weakly compressible approach in \cite{HybridGPR}.

\subsection{Smooth acoustic wave}
\begin{figure}[H]
	\centering
	\includegraphics[clip,trim = 0 200 0 100,width=0.32\linewidth]{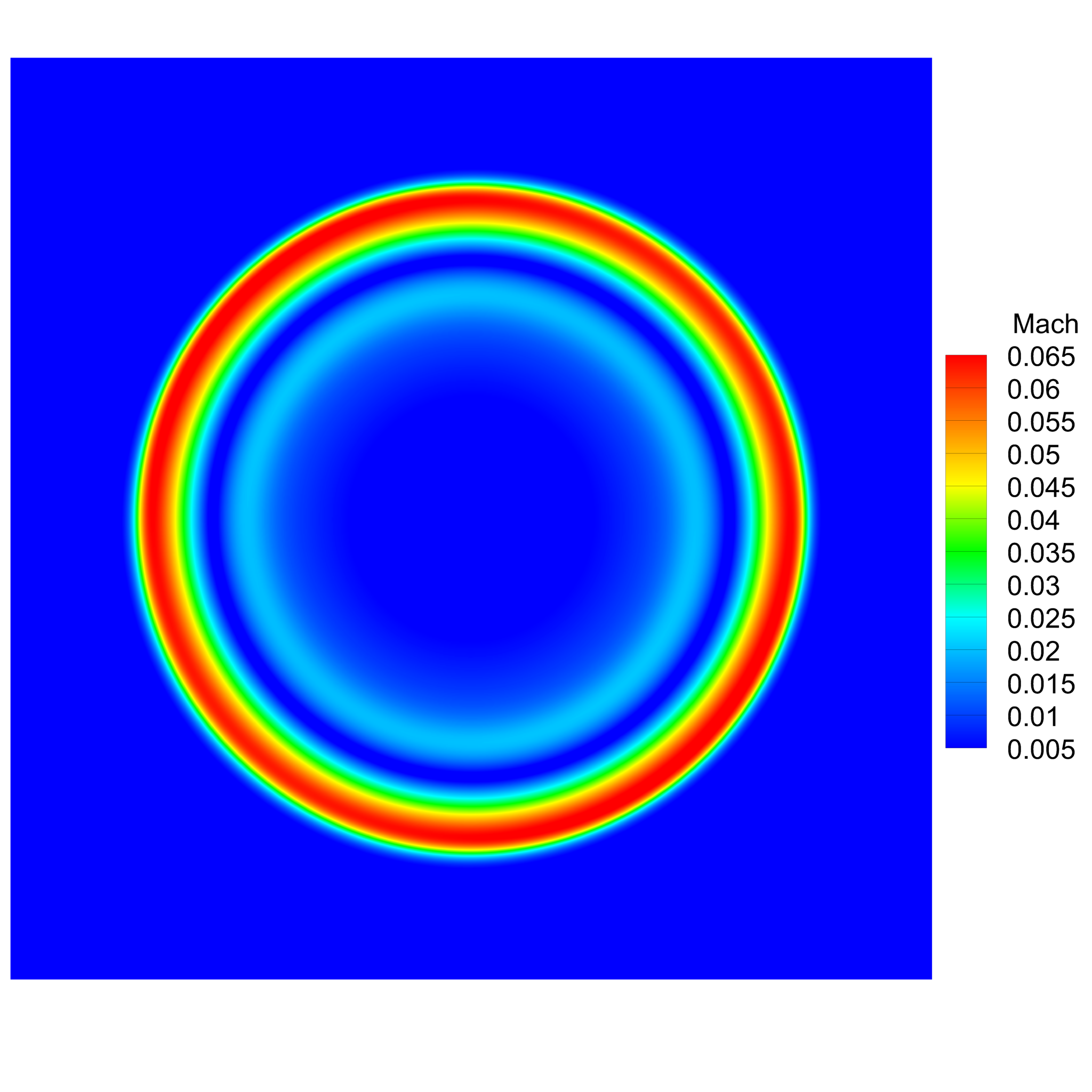}
	\includegraphics[clip,trim = 0 200 0 100,width=0.32\linewidth]{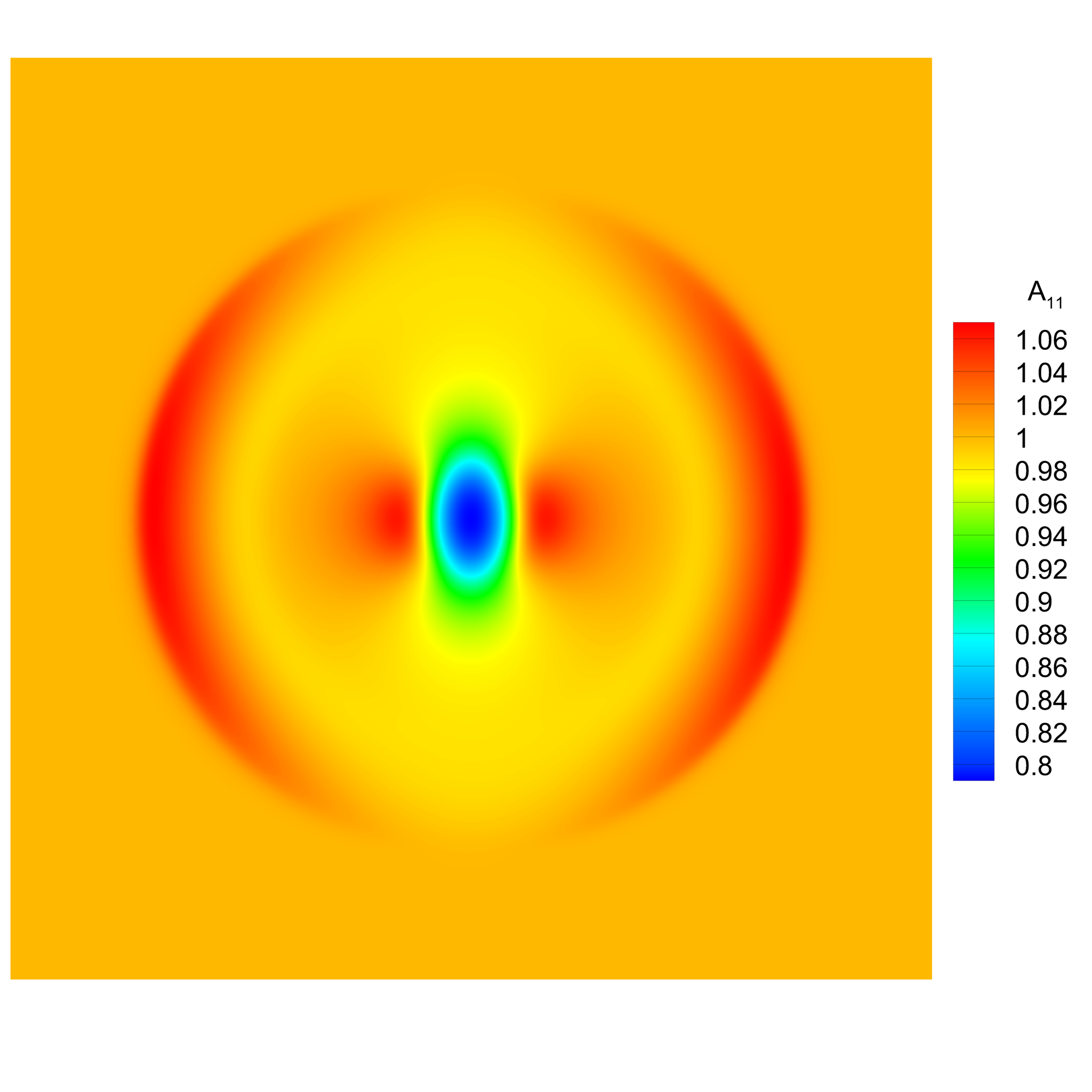}
	\includegraphics[clip,trim = 0 200 0 100,width=0.32\linewidth]{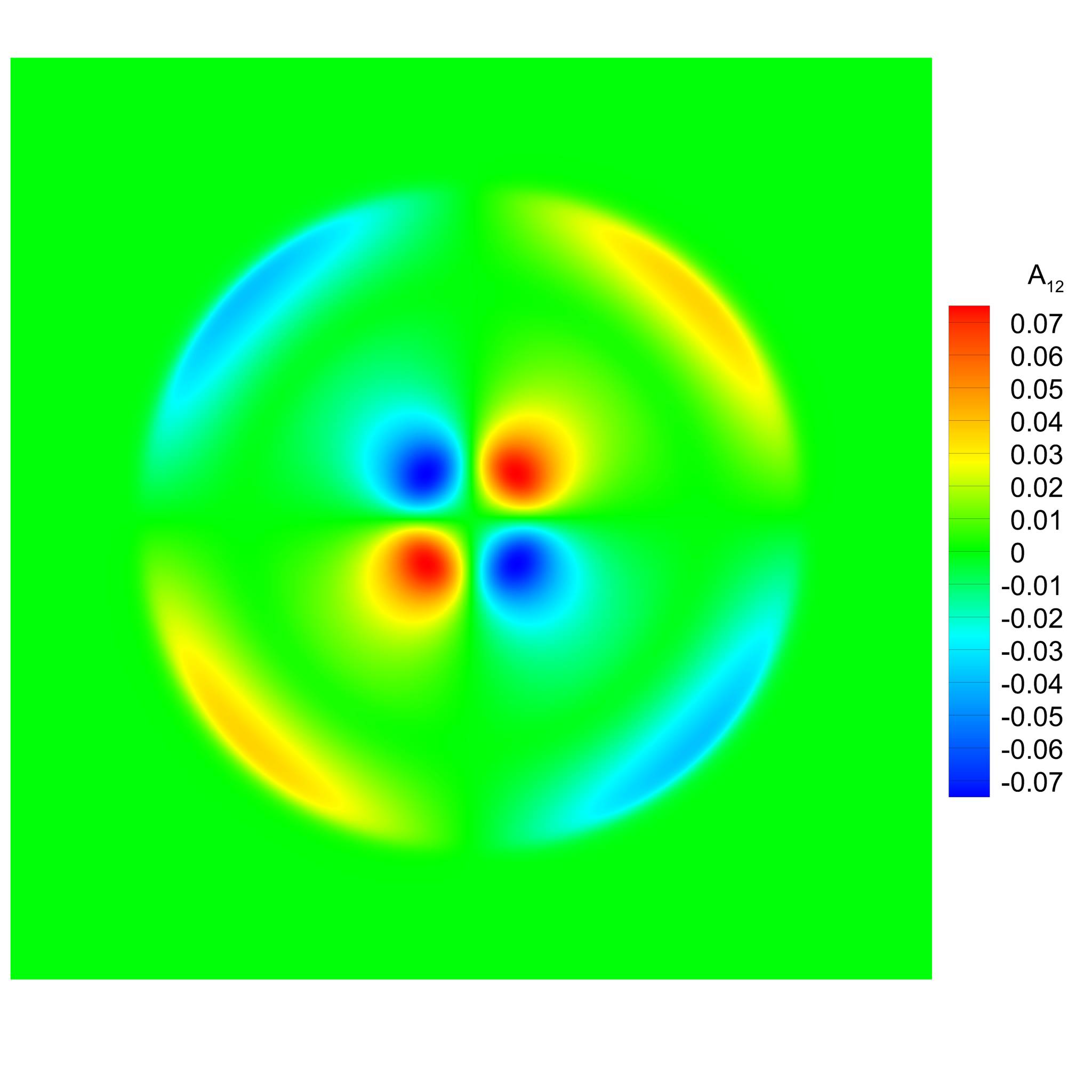}
	\caption{Smooth acoustic wave. Contour plots of the Mach number, $A_{11}$ and $A_{12}$ fields computed using the hybrid FV/FE method for the compressible GPR model with the Eulerian scheme.}
	\label{fig:SAW}
\end{figure}

\begin{figure}[H]
	\centering
	\includegraphics[width=0.32\linewidth]{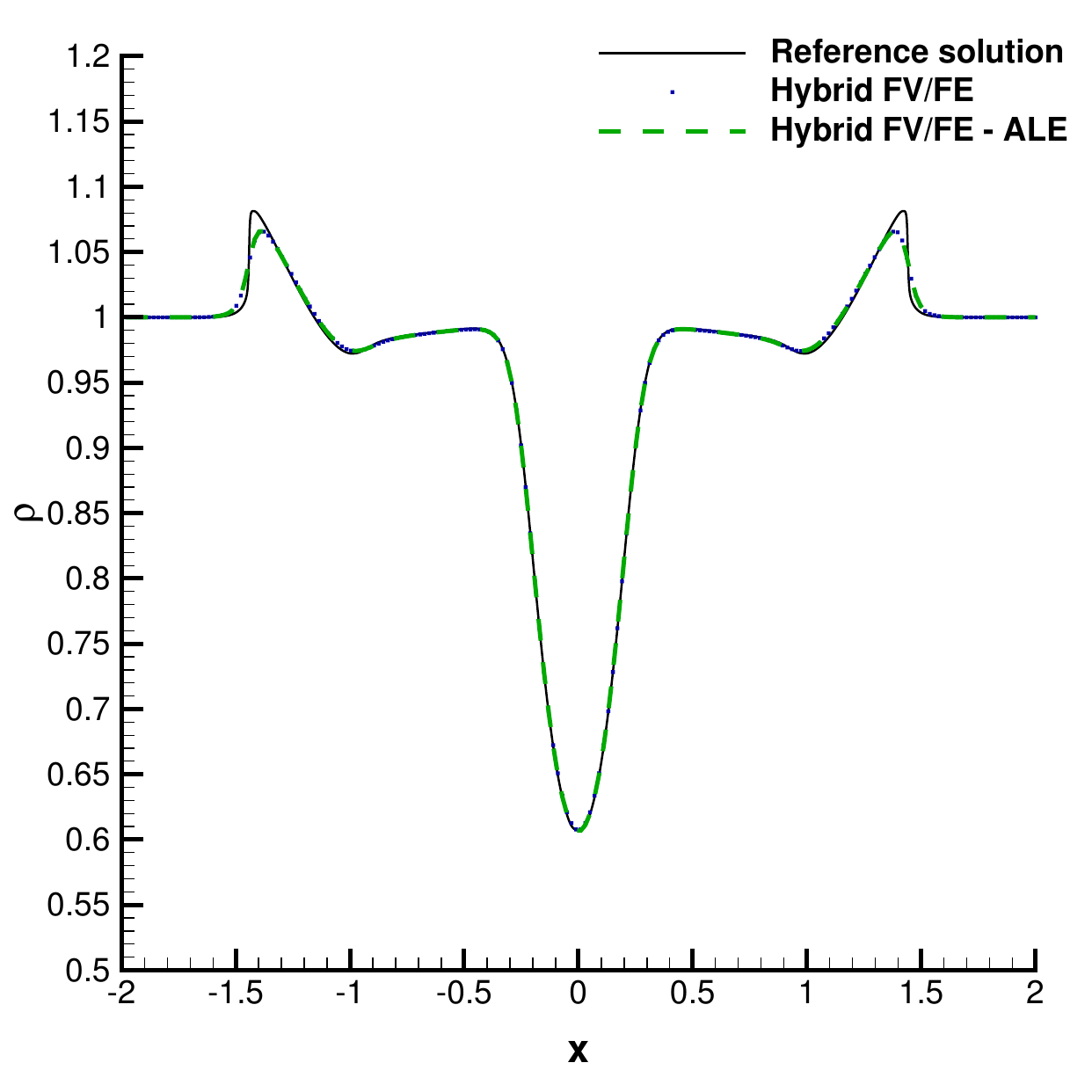}
	\includegraphics[width=0.32\linewidth]{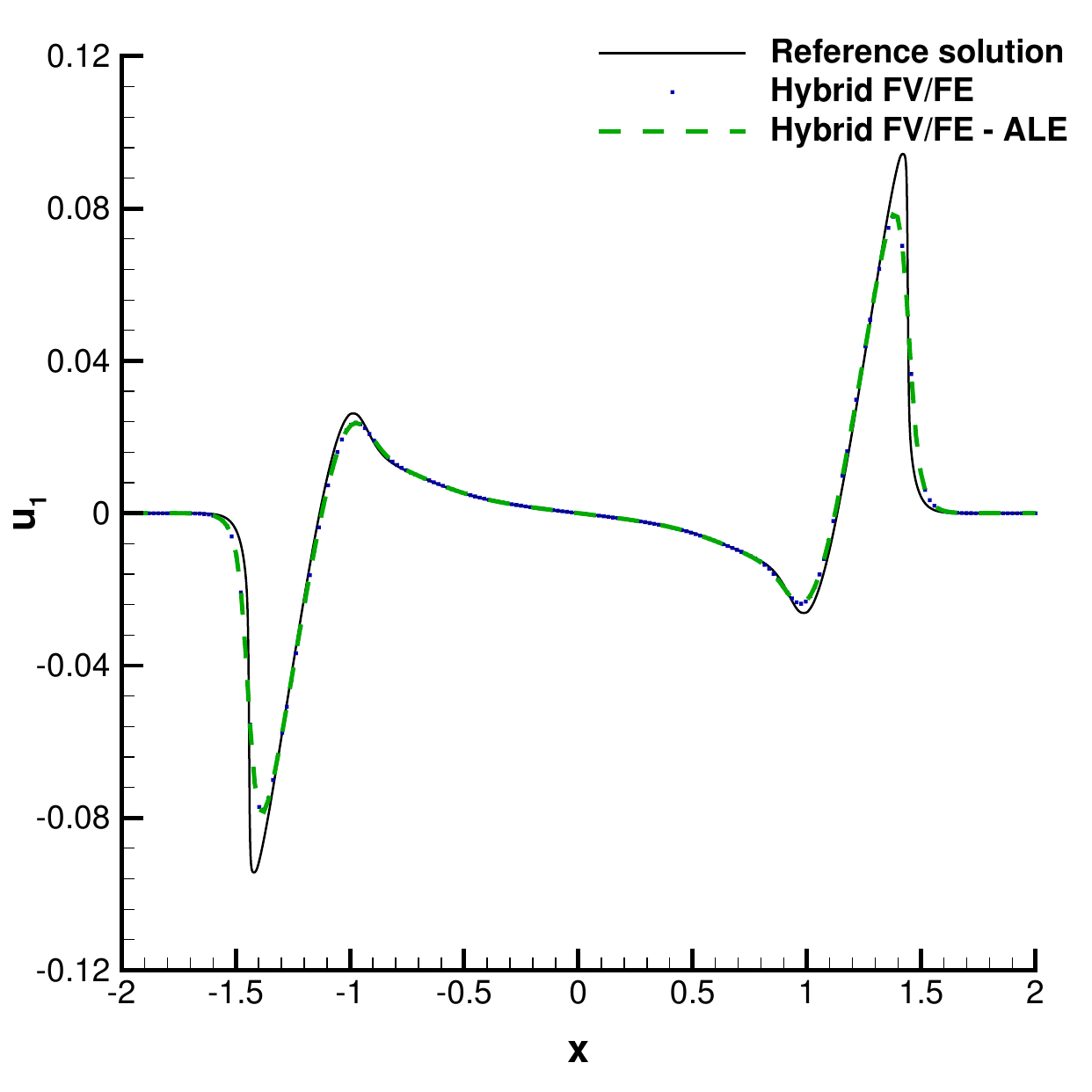}
	\includegraphics[width=0.32\linewidth]{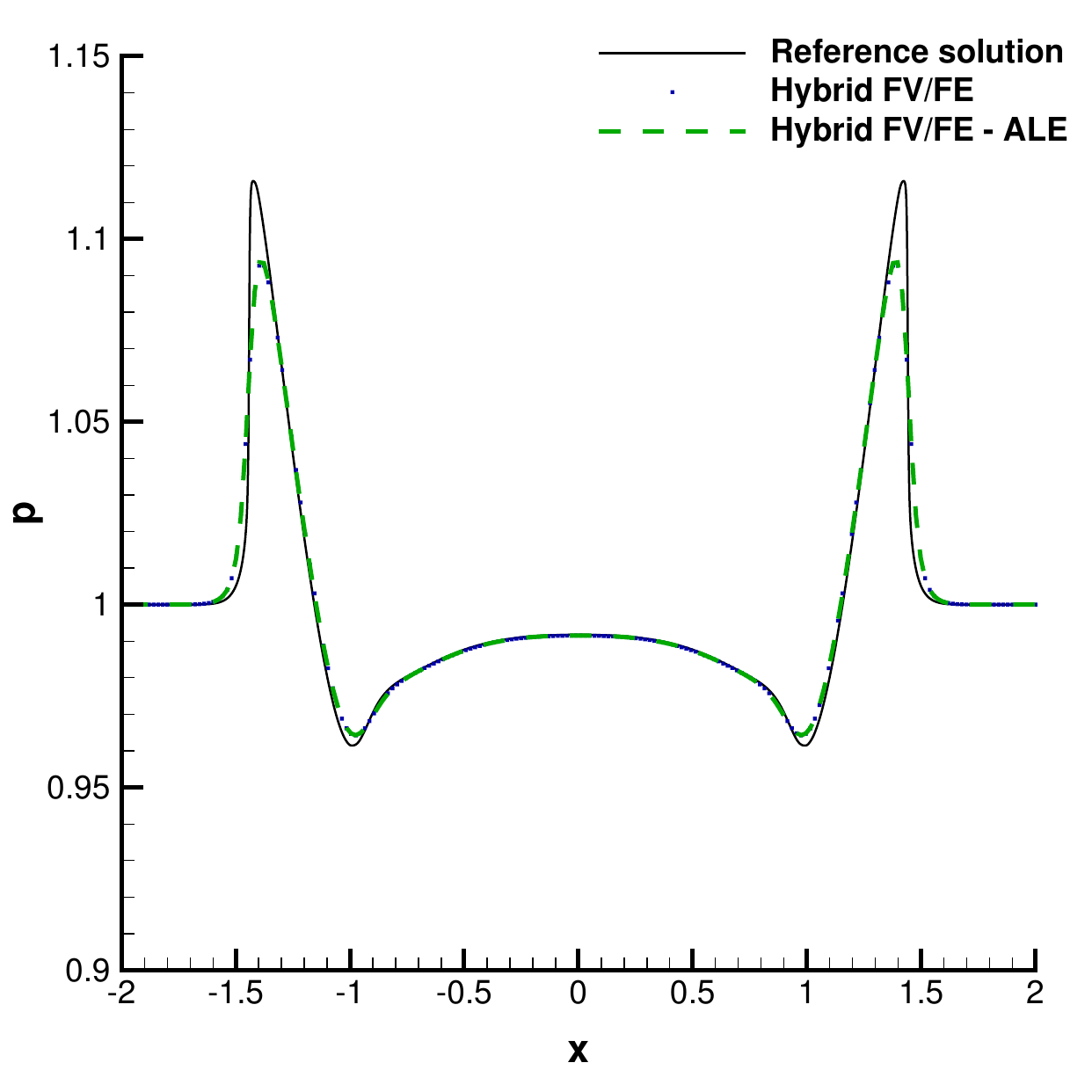}
	\caption{Smooth acoustic wave. 1D cuts of the density, $\vel_{1}$,  and pressure fields for $\left\lbrace (x,y)\in\mathbb{R}^{2}\,  |\,x \in[-1,1],\, y=0 \right\rbrace$ using the hybrid FV/FE method for the compressible GPR model with the Eulerian scheme (blue squares) and with the ALE method (green dashed line) compared against the reference solution (black line).}
	\label{fig:SAW1D}
\end{figure}

To analyse the behaviour of the proposed approach for flows dominated by acoustic waves, we consider the acoustic wave benchmark \cite{TD17}. The initial condition defined in $\Omega=[-2,2]^{2}$ is given by
\begin{equation*}
	\rho\left(\x,0\right) = 1,\quad
	\bvel \left(\x,0\right) = \boldsymbol{0}, \quad
	\press \left(\x,0\right) = 1+e^{-40 r^2}
\end{equation*}
with $r = \sqrt{x^2+y^2}$. As model parameters, we set $\mu=\kappa = c_{s}=0$, $c_{h}=10^{-10}$, $c_{v}=1$, $c_{p}=1.4$. 
A first simulation is run using the hybrid FV/FE method without limiters on a computational grid made of $2097152$ primal elements. The contour plots of the Mach number , $A_{11}$ and $A_{12}$, and 1D cuts of the density, pressure and velocity component $\vel_{1}$ at time $t_{\mathrm{e}}=1$ are depicted in Figures~\ref{fig:SAW}--\ref{fig:SAW1D}. The 1D plots also report the solution obtained using the ALE method with mesh speed smoothing $\varsigma=5$. To ease validation, a reference solution computed with a second order TVD FV scheme on  a very fine grid is included. We observe that the steep acoustic wave front is well captured with both schemes.

\subsection{Lid driven cavity}
\begin{figure}[H]
	\begin{center}
		\includegraphics[trim=15 -70 15 100,clip,width=0.45\textwidth]{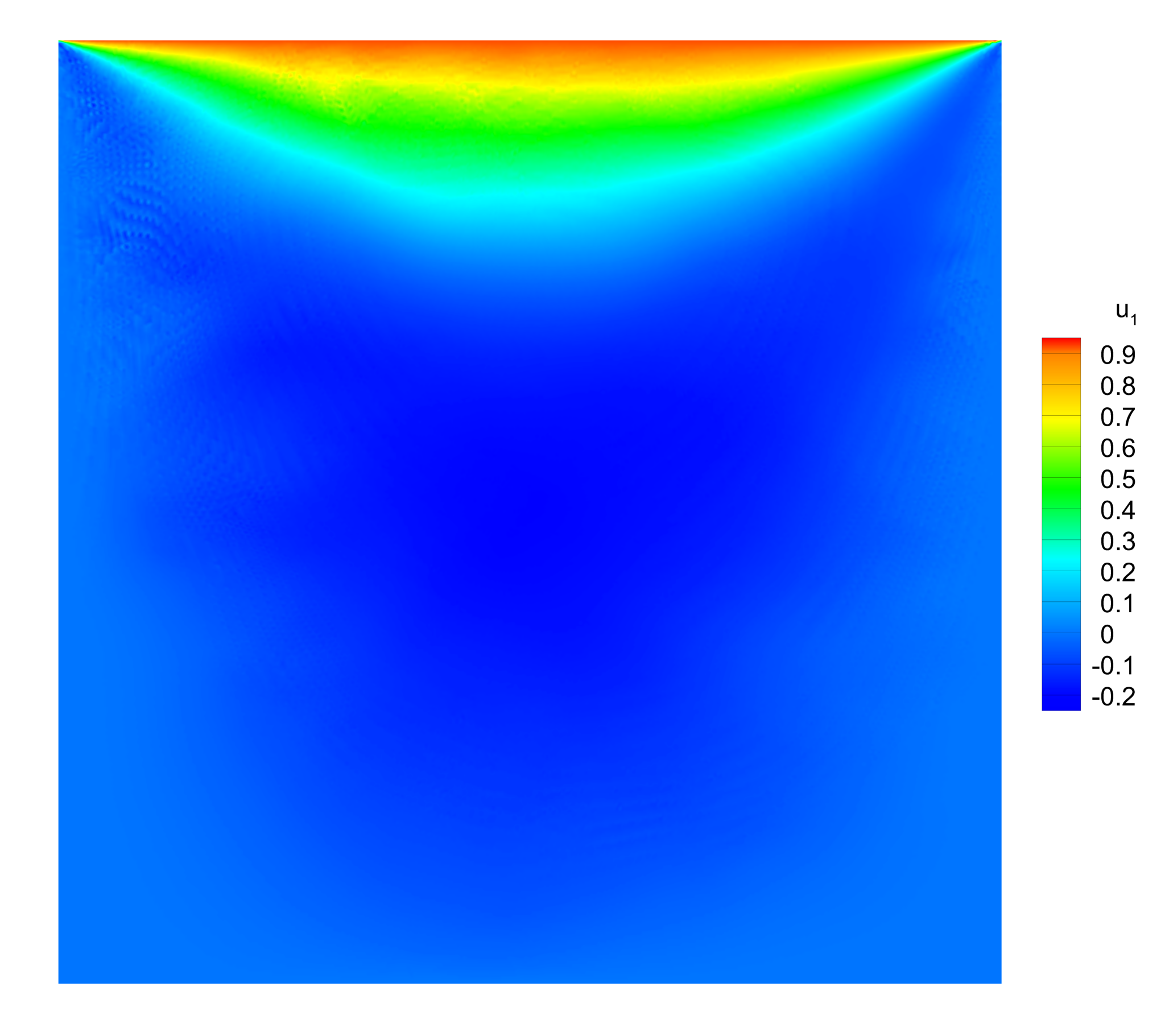}  
		\includegraphics[width=0.45\textwidth]{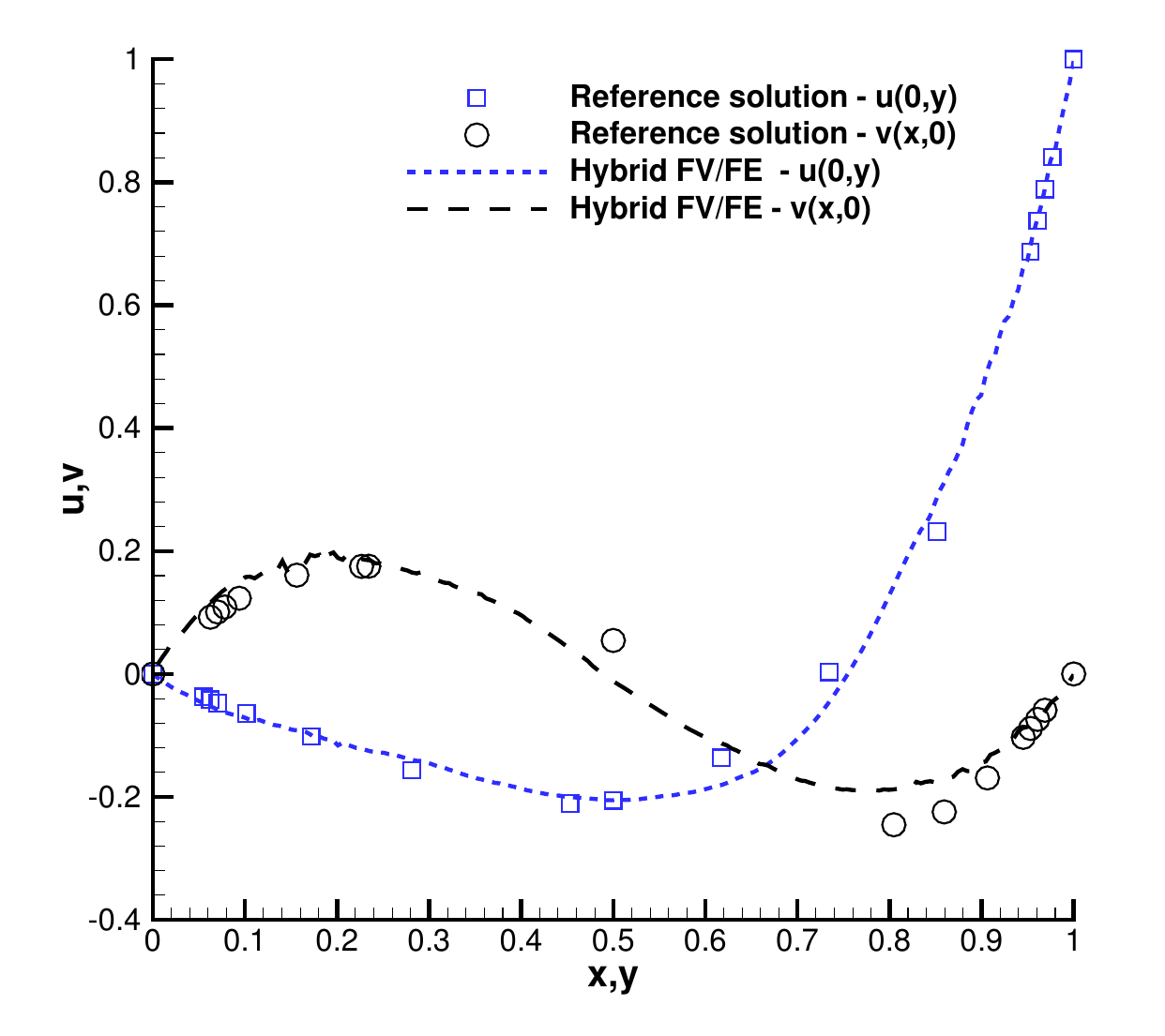} 	
%		\caption{Lid-driven cavity. Contour plot of the distortion field $A_{12}$ (left) and 1D cut in $x-$ and $y-$directions of the velocity components (right). The provided reference solution has been taken from \cite{GGS82}.}  
		\caption{Lid-driven cavity. Contour plot of the velocity field $\vel_{1}$ (left) and 1D cut in $x-$ and $y-$directions of the velocity components (right). The provided reference solution has been taken from \cite{GGS82}.}  
		\label{fig.LDC}
	\end{center}
\end{figure}

A classical test case in fluid dynamics is the lid driven cavity benchmark, for which we set $c_{s}=8$, $c_{h}=0$, $\mu = 10^{-2}$, $\kappa = 0$, $c_v=1$ and $\gamma=1.4$. In this test case, the fluid is assumed to be confined within a unit square cavity with a sliding lid at the top moving at velocity $\bvel_{\mathrm{lid}}=(1,0)^{T}$ and an initial fluid at rest. The computational domain is discretized employing $185984$ primal cells and strong wall boundary conditions are imposed in all boundaries. The results, computed at $t_{\mathrm{e}}=10$ using the asymptotic preserving method with ENO limiting and $c_{\alpha}=0.5$, are reported in Figure~\ref{fig.LDC}. A good agreement is observed for the 1D cuts in comparison with the reference data in  \cite{GGS82}.

\subsection{Solid rotor}
\begin{figure}[H]
	\centering
	\includegraphics[width=0.32\linewidth]{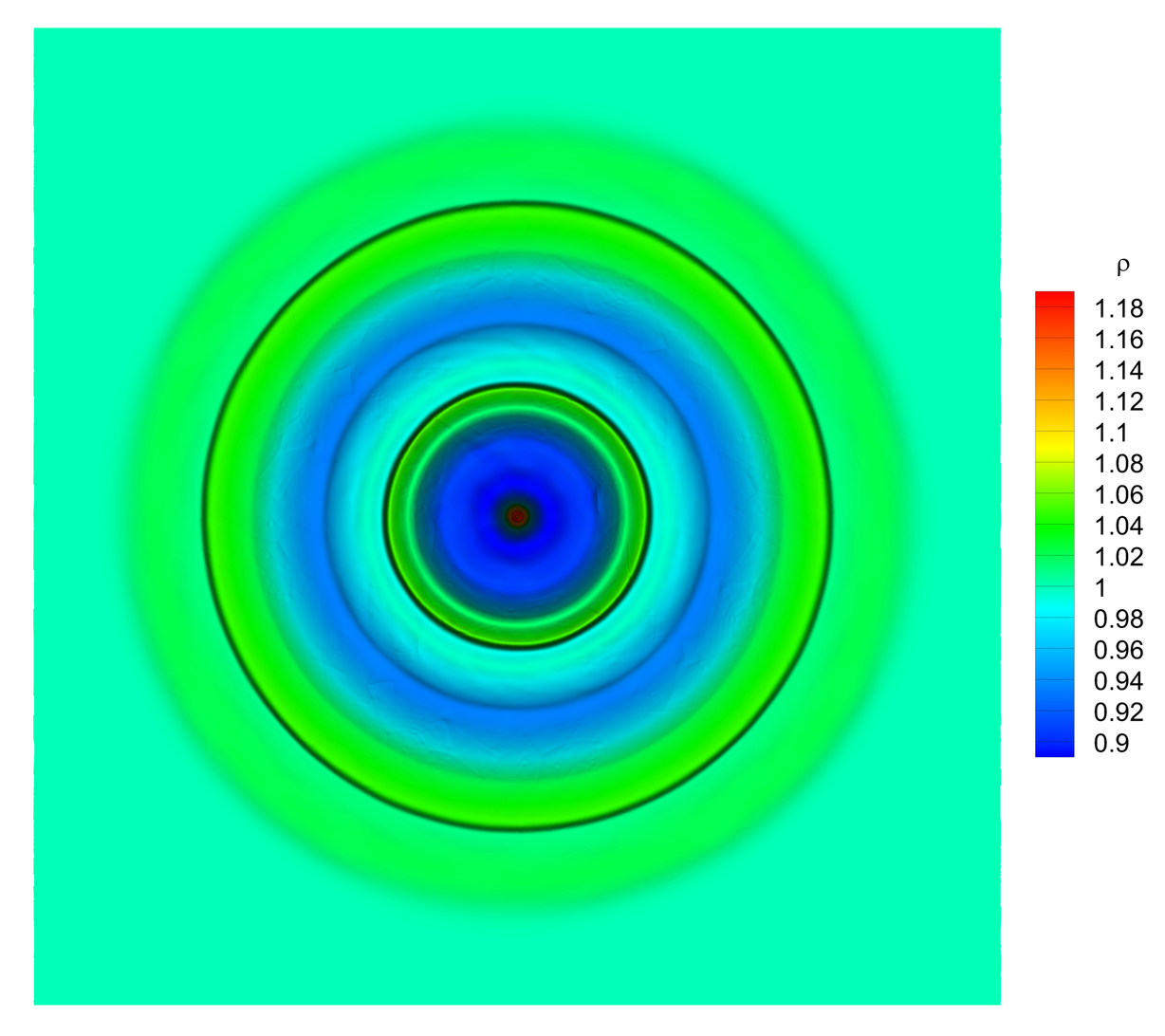}
	\includegraphics[width=0.32\linewidth]{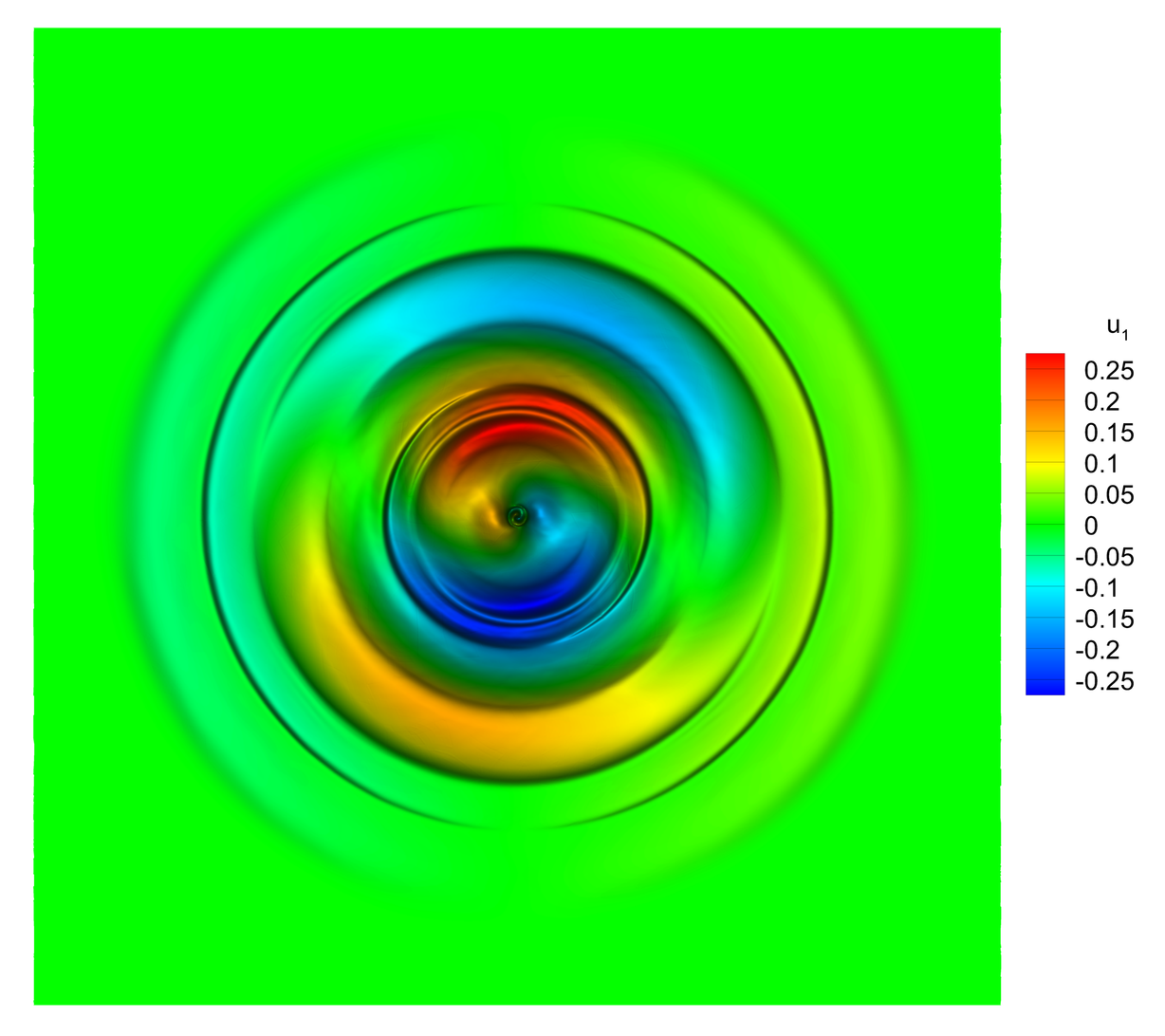}
	\includegraphics[width=0.32\linewidth]{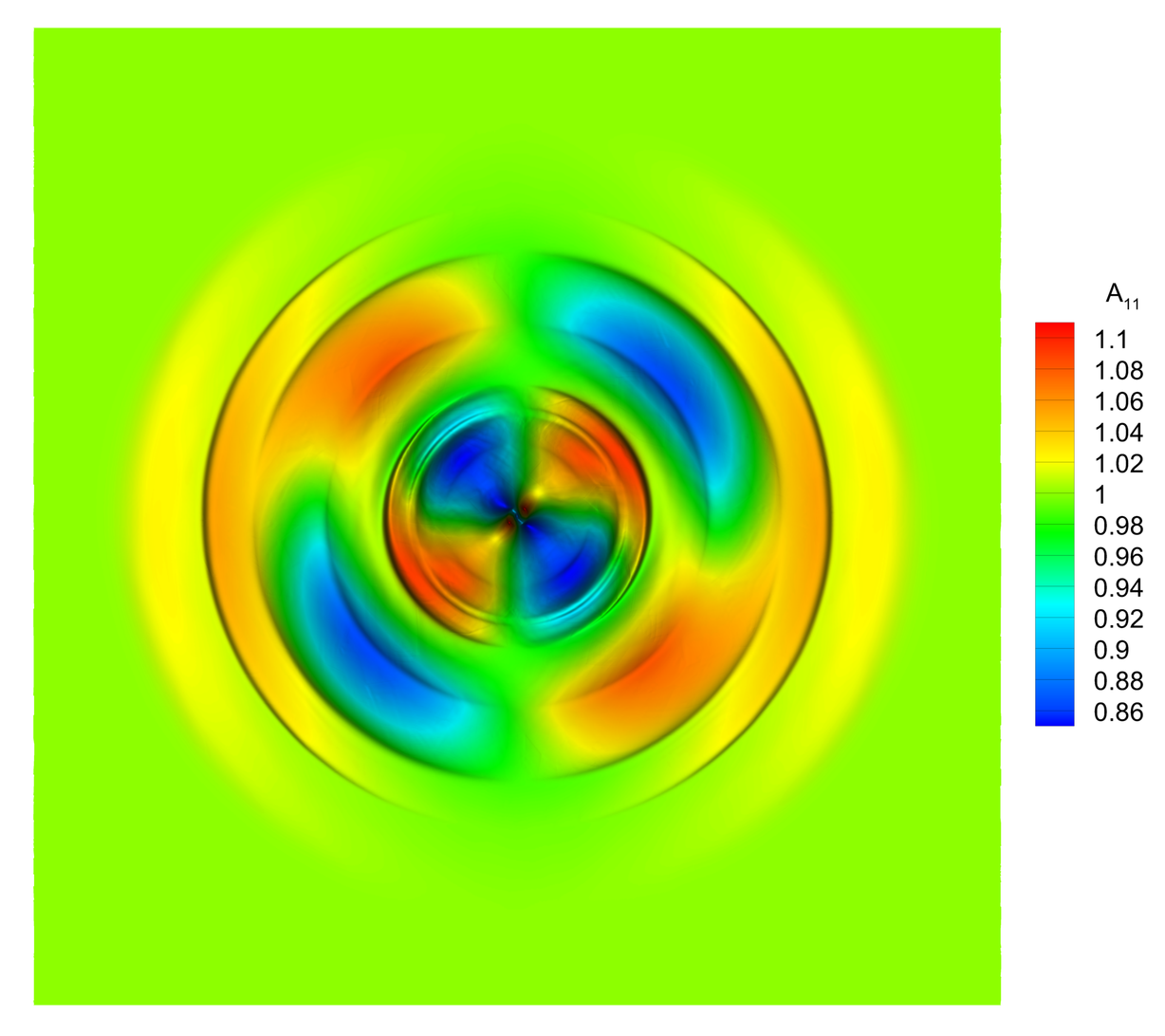}
	\includegraphics[width=0.32\linewidth]{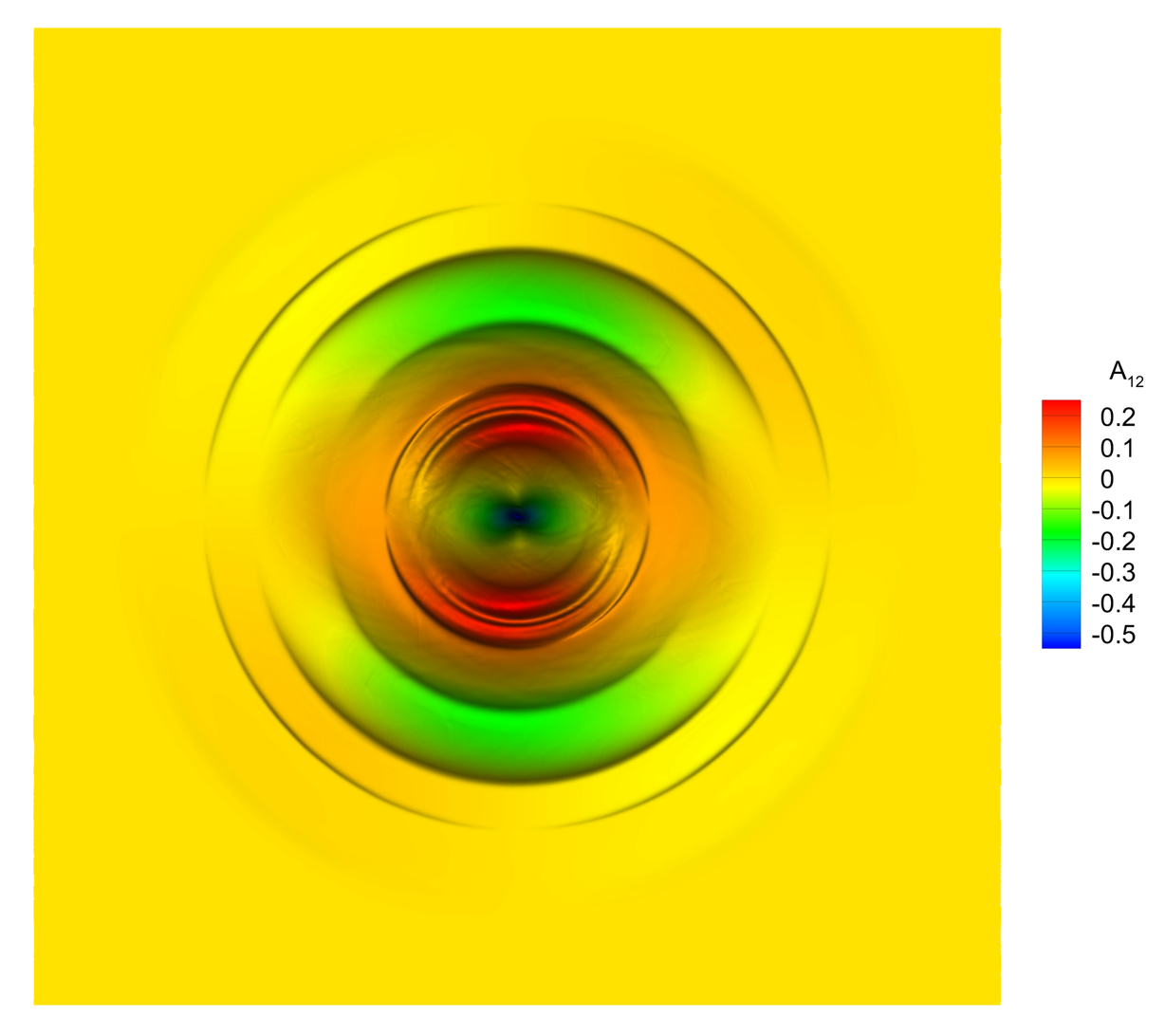}
	\includegraphics[width=0.32\linewidth]{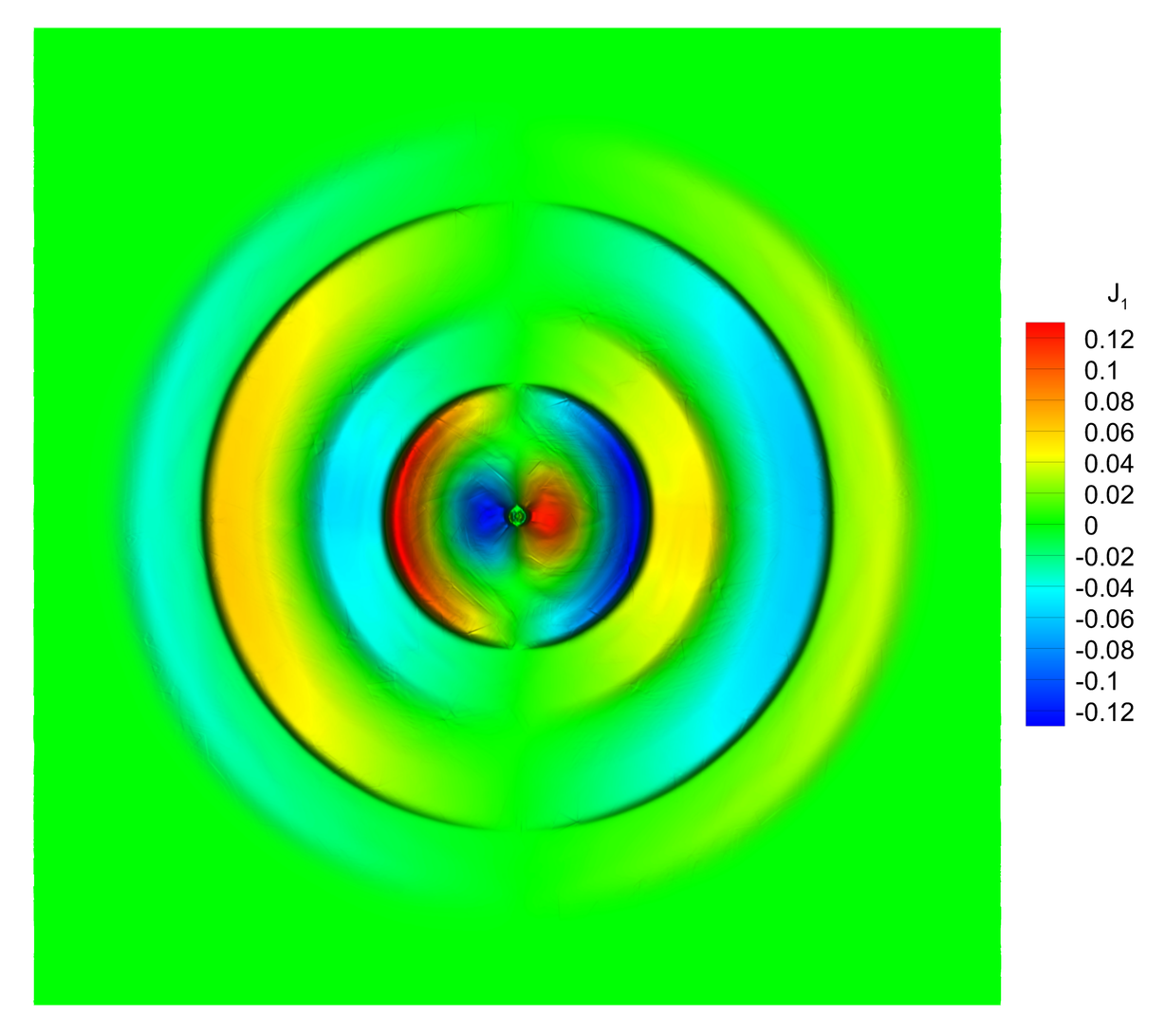}
	\includegraphics[width=0.32\linewidth]{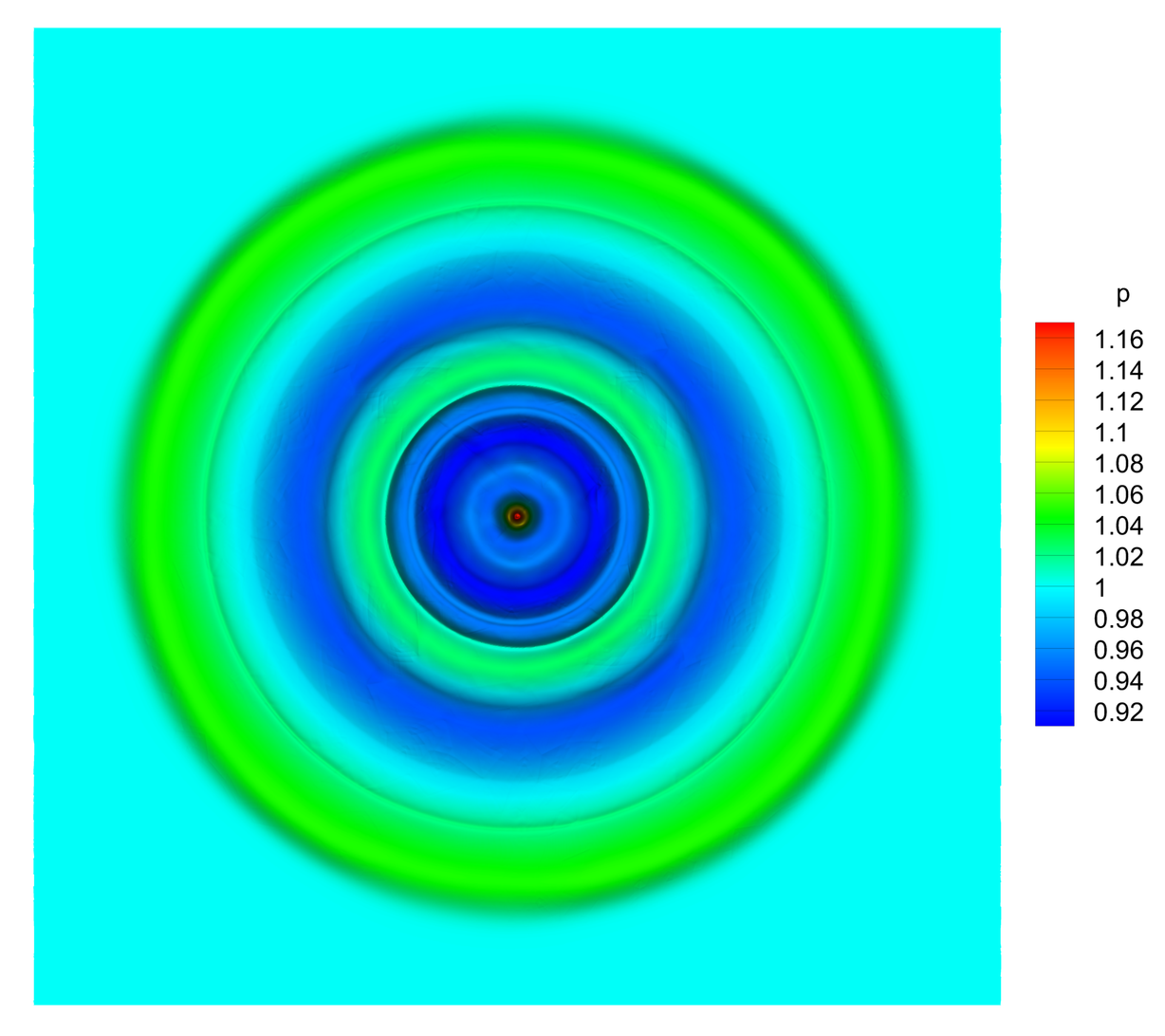}
	\caption{Solid rotor. Contour plots for the density, $\vel_{1}$, $A_{11}$, $A_{12}$, $J_{1}$ and pressure fields at $t_{e}=0.3$.}
	\label{fig:SR}
\end{figure}

\begin{figure}[H]
	\centering
	\includegraphics[width=0.32\linewidth]{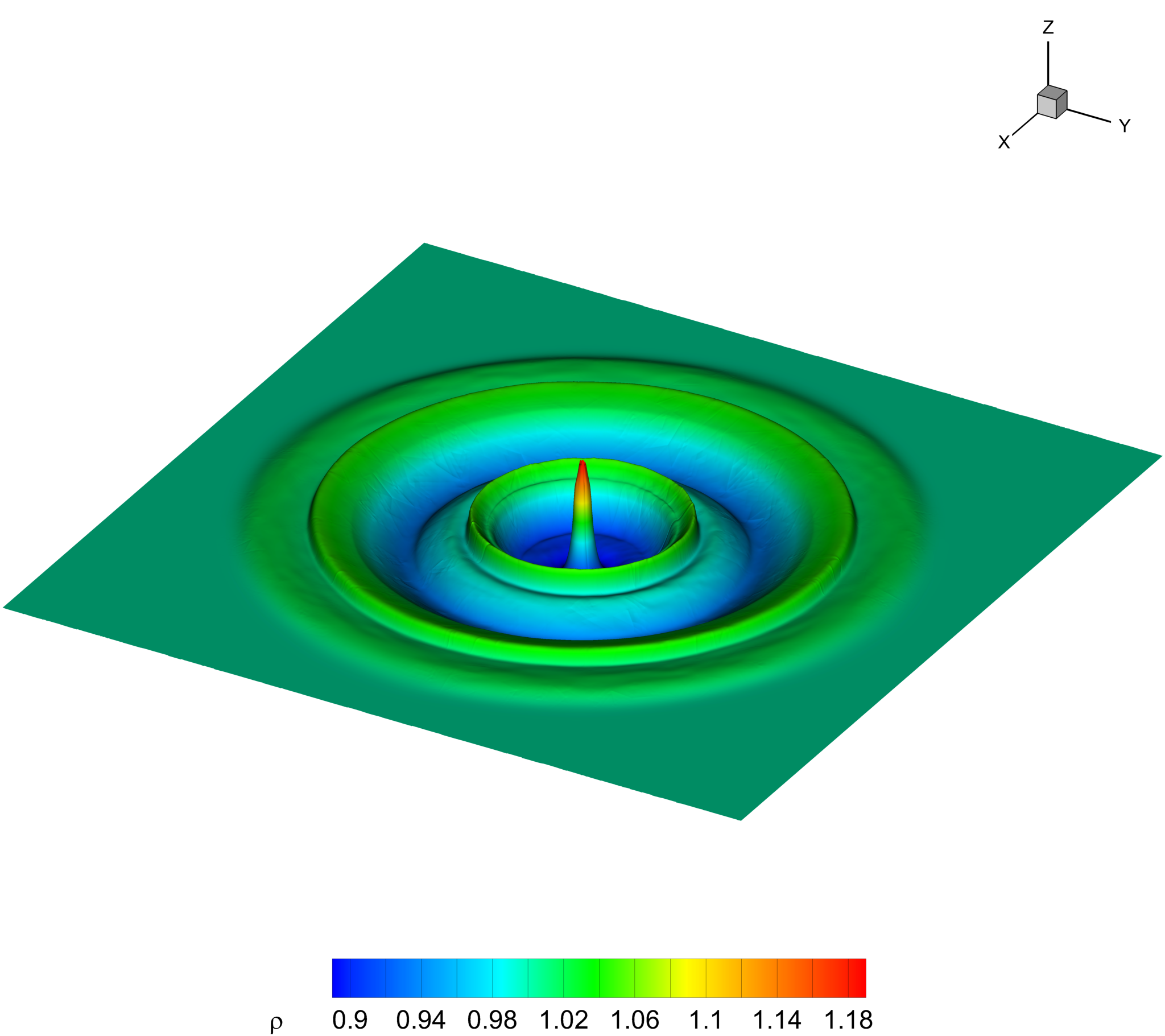}
	\includegraphics[width=0.32\linewidth]{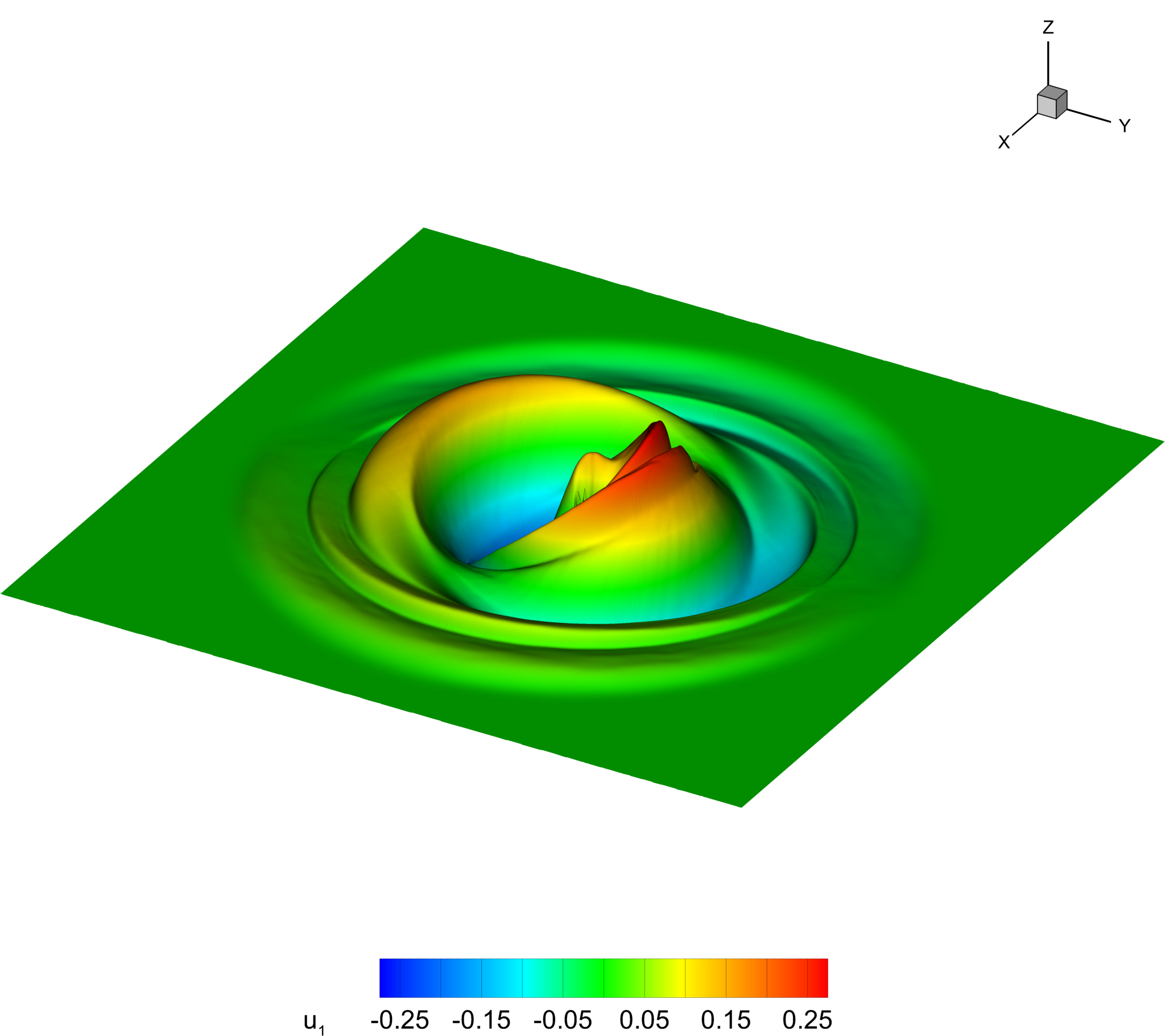}
	\includegraphics[width=0.32\linewidth]{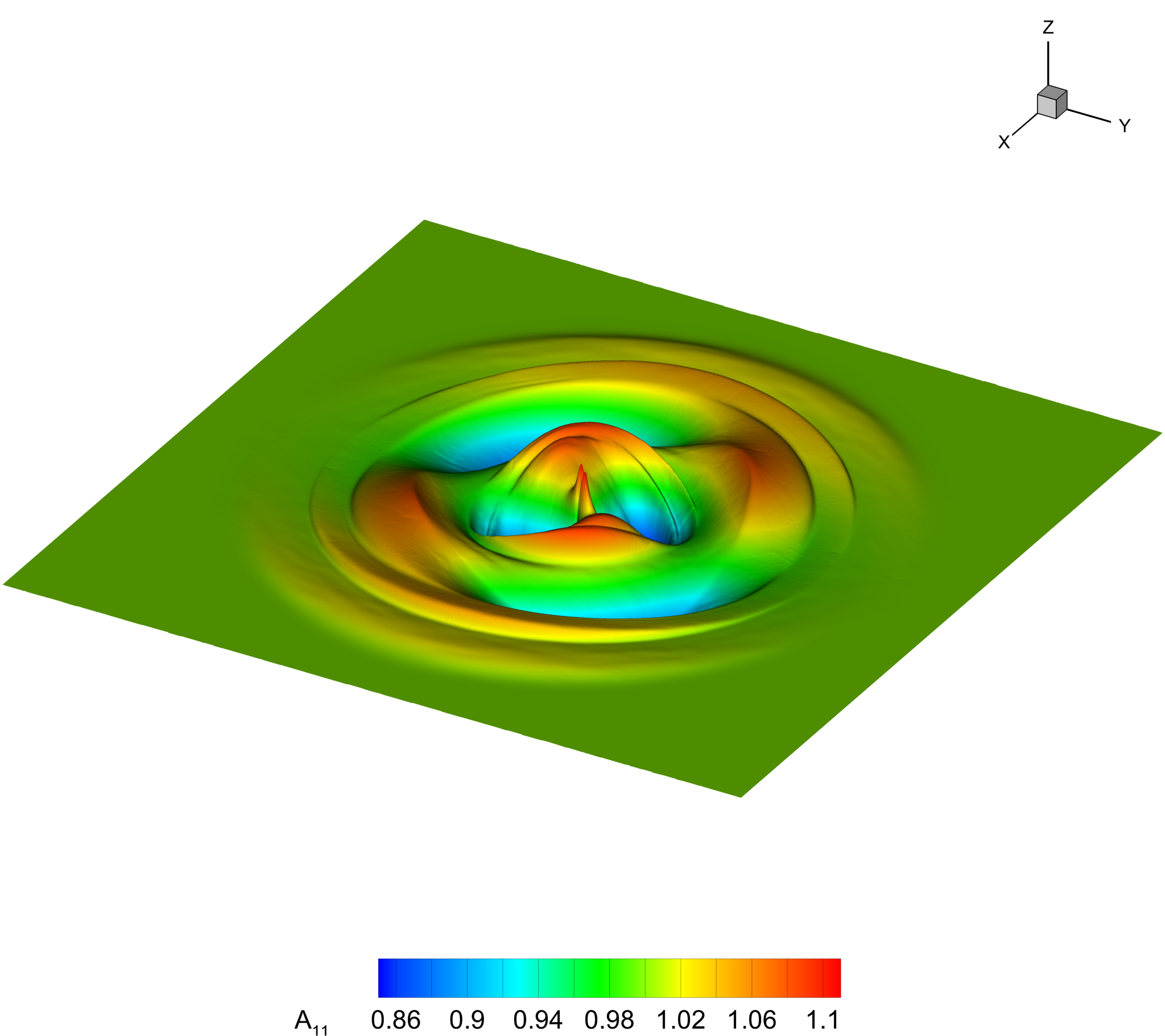}
	\includegraphics[width=0.32\linewidth]{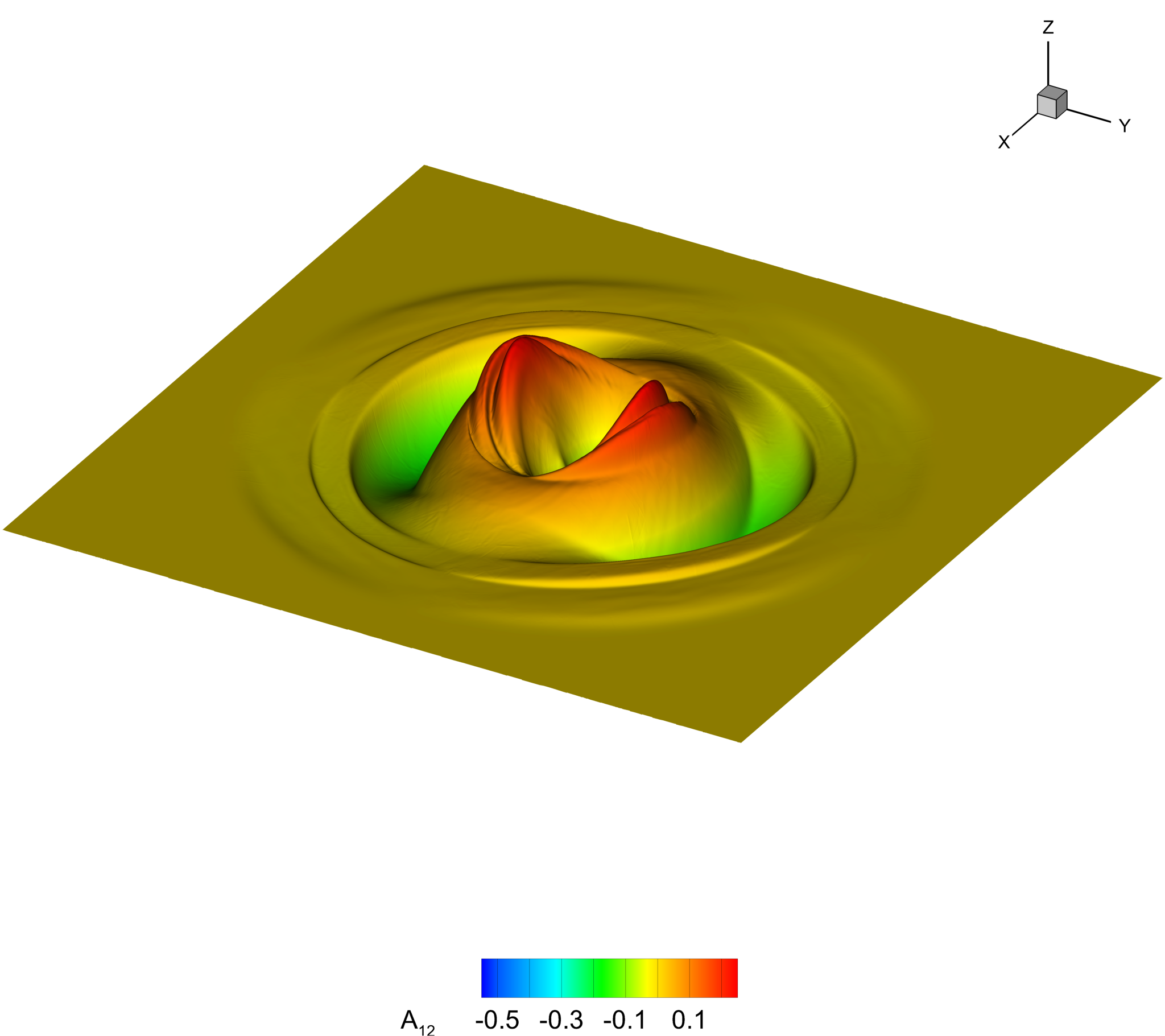}
	\includegraphics[width=0.32\linewidth]{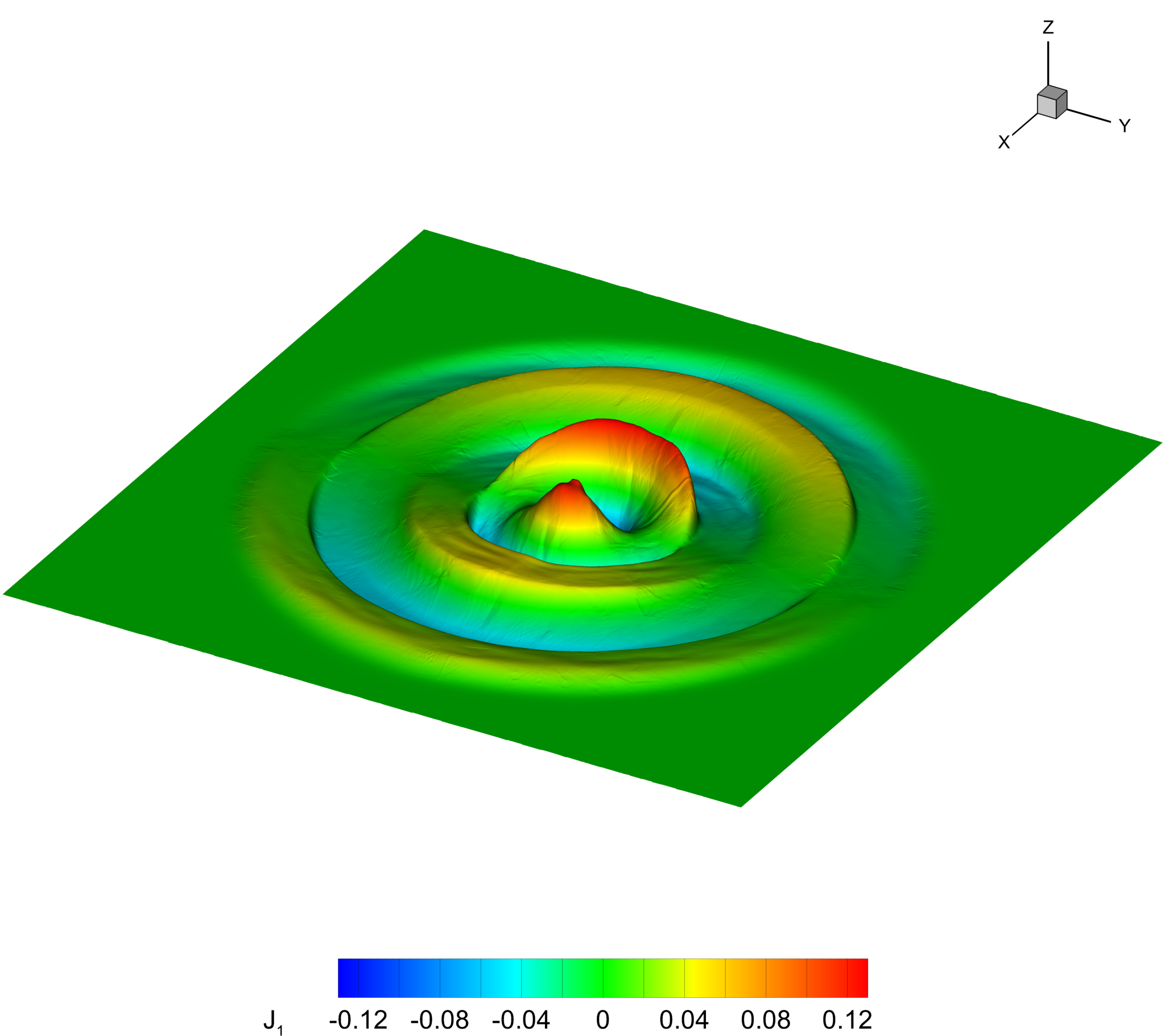}
	\includegraphics[width=0.32\linewidth]{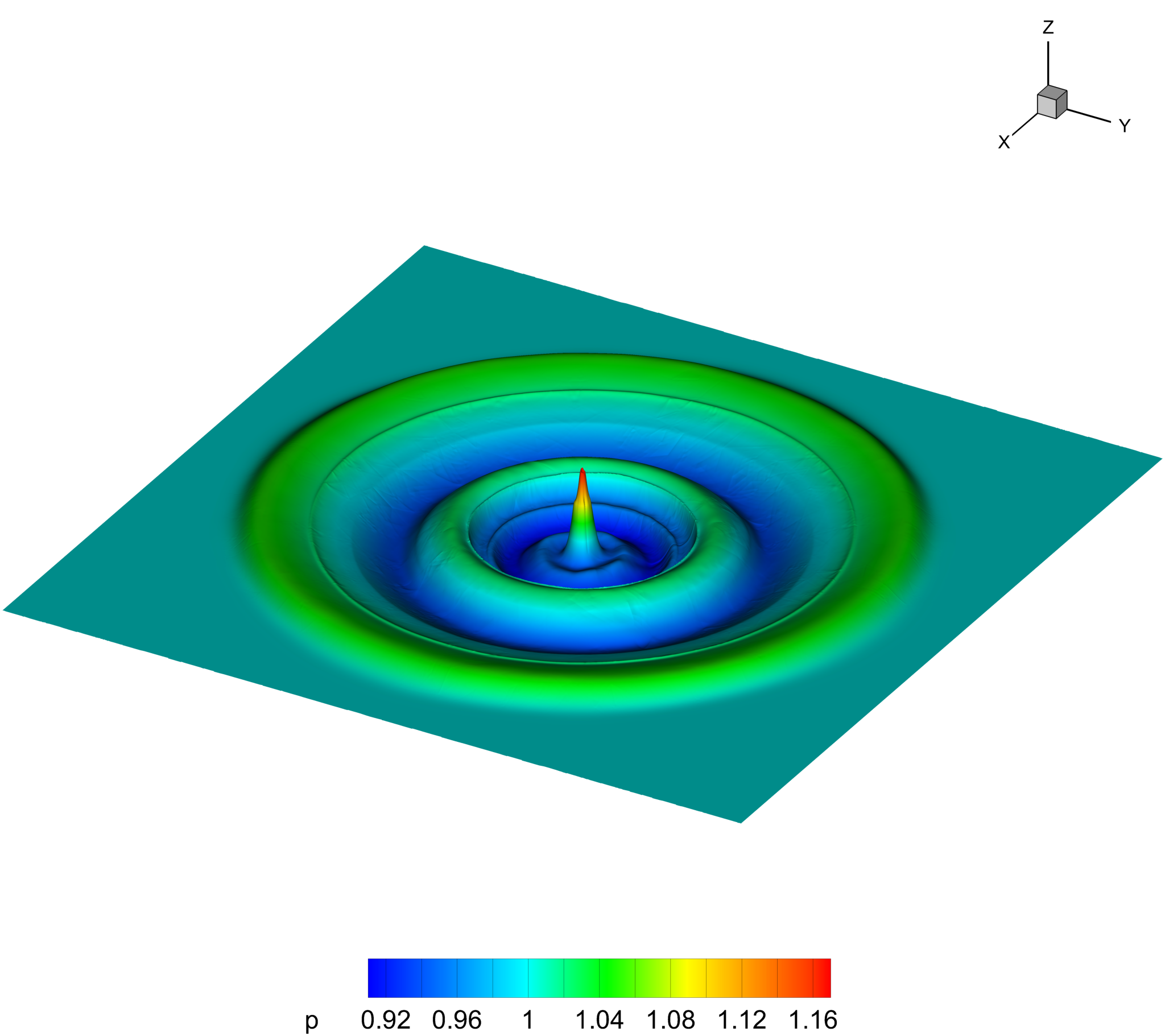}
	\caption{Solid rotor. Elevated contour plots for the density, $\vel_{1}$, $A_{11}$, $A_{12}$, $J_{1}$ and pressure fields at $t_{e}=0.3$.}
	\label{fig:SR3D}
\end{figure}

\begin{figure}[H]
	\centering
	\includegraphics[width=0.32\linewidth]{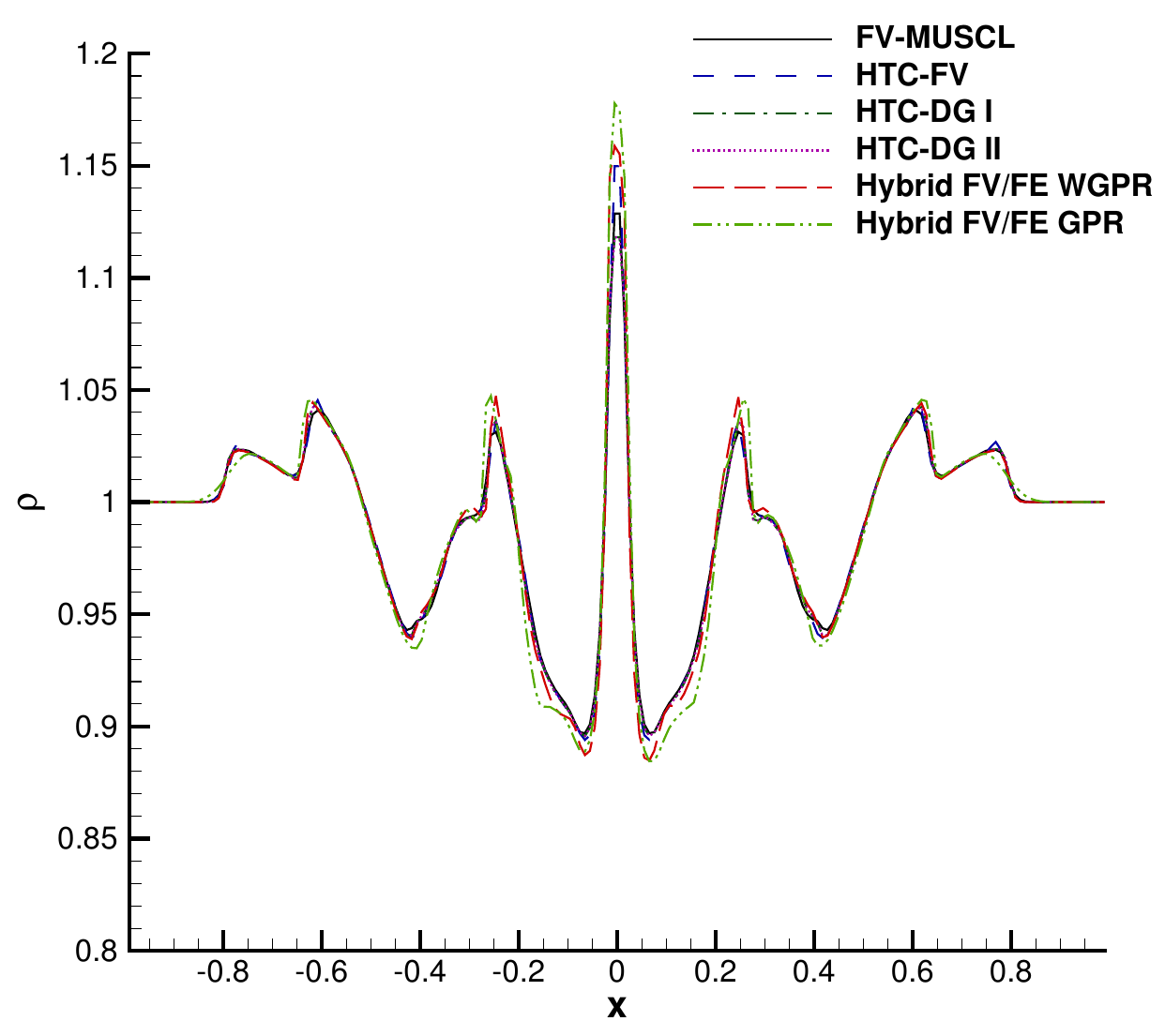}
	\includegraphics[width=0.32\linewidth]{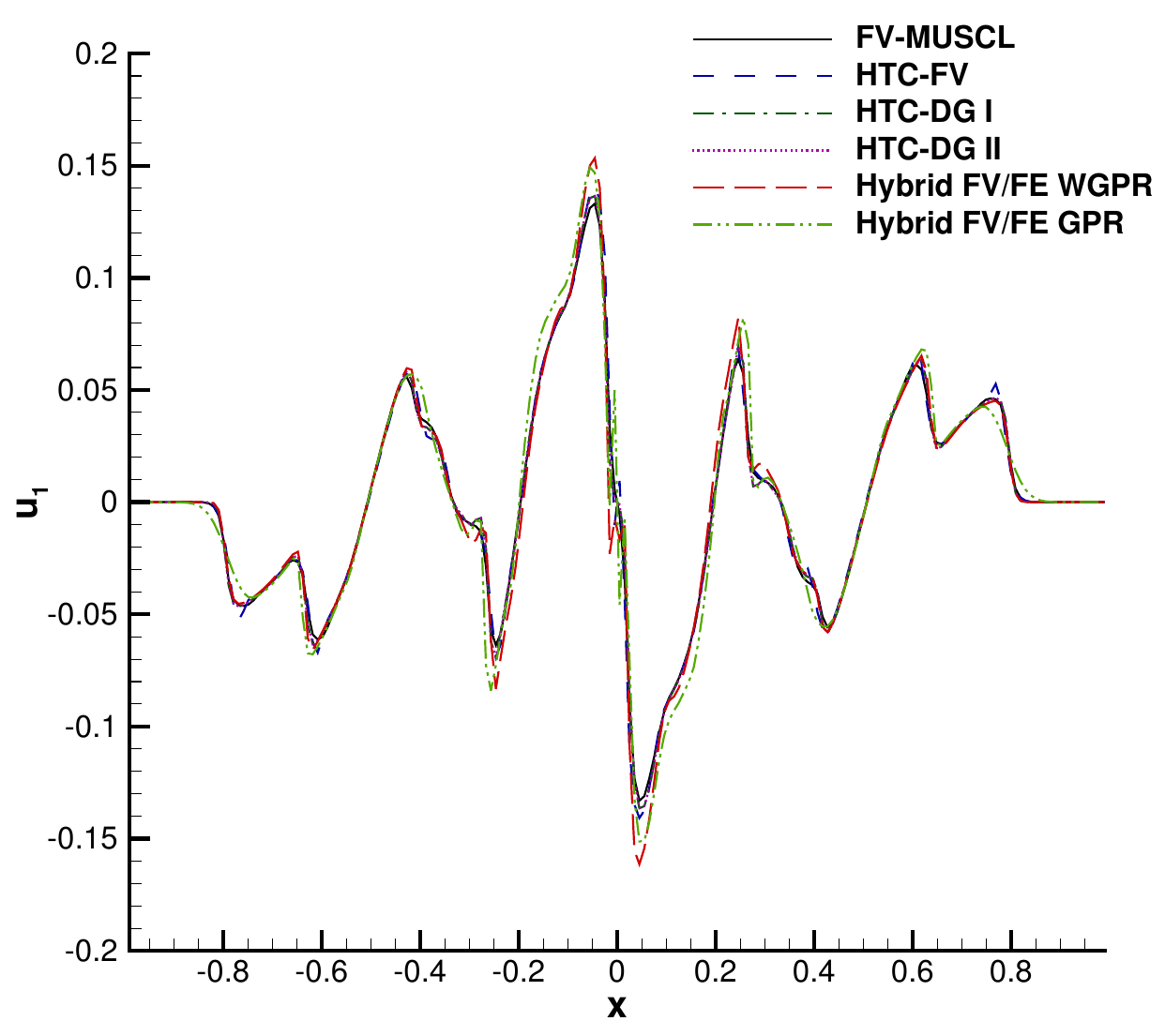}
	\includegraphics[width=0.32\linewidth]{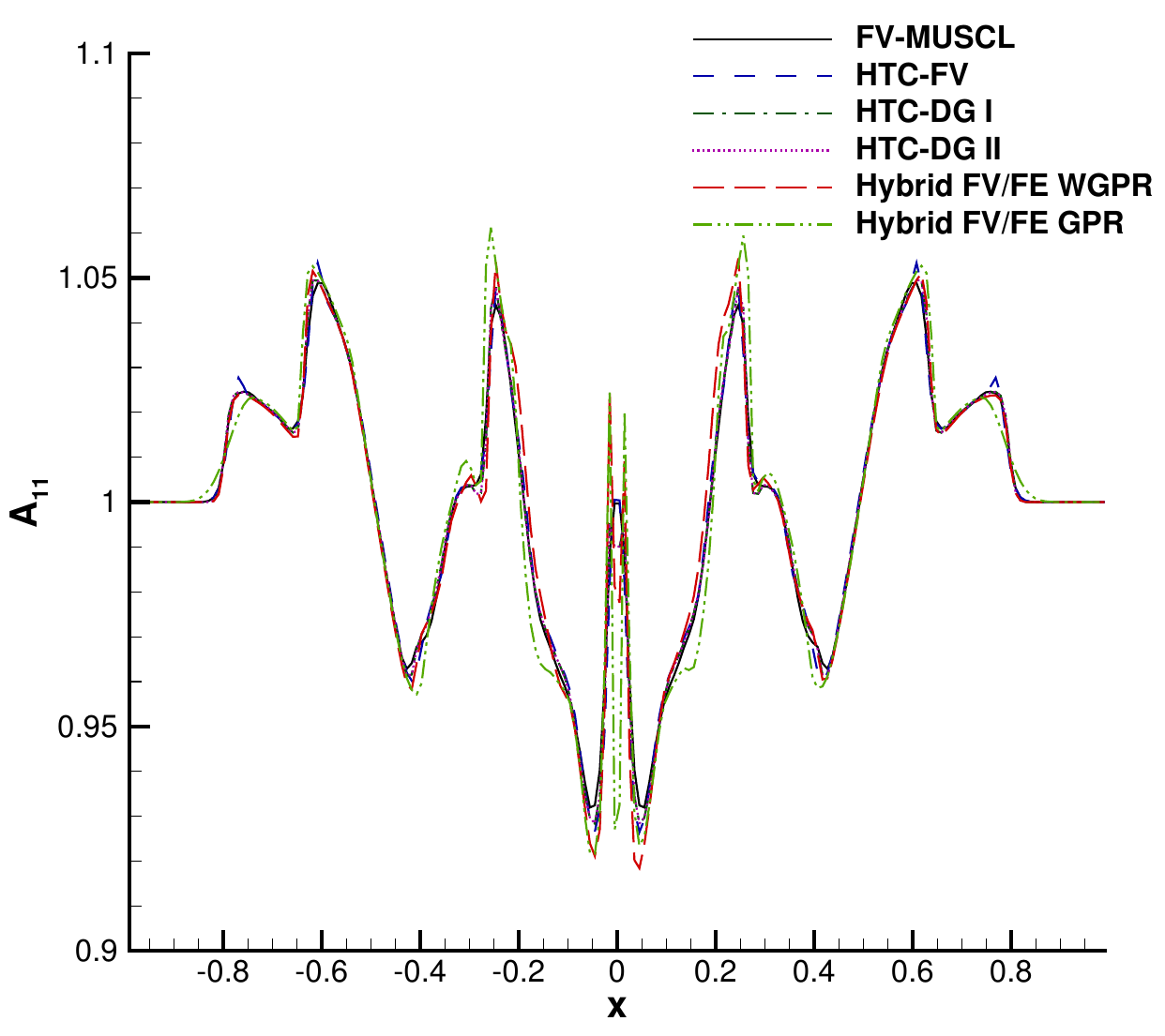}
	\includegraphics[width=0.32\linewidth]{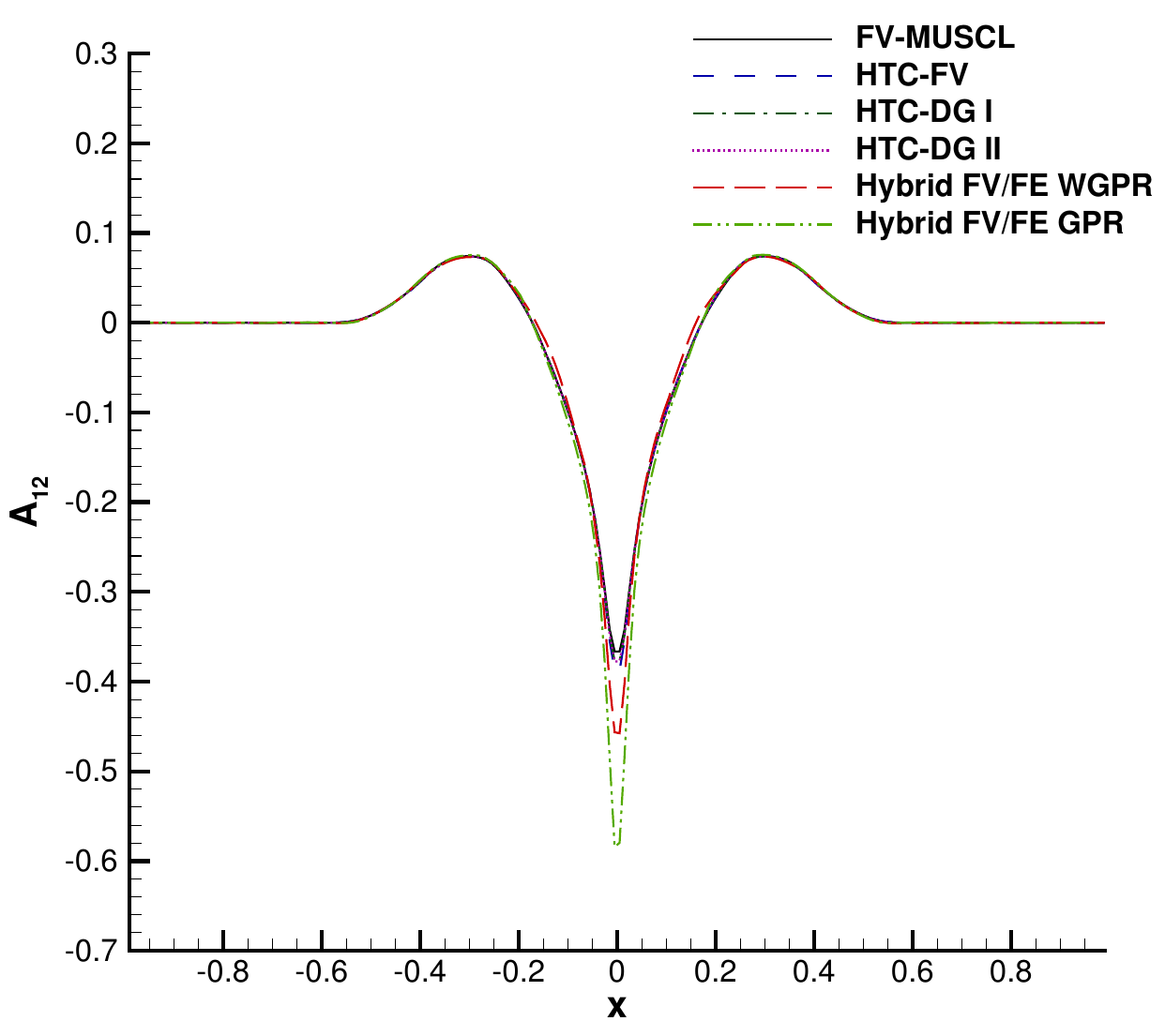}
	\includegraphics[width=0.32\linewidth]{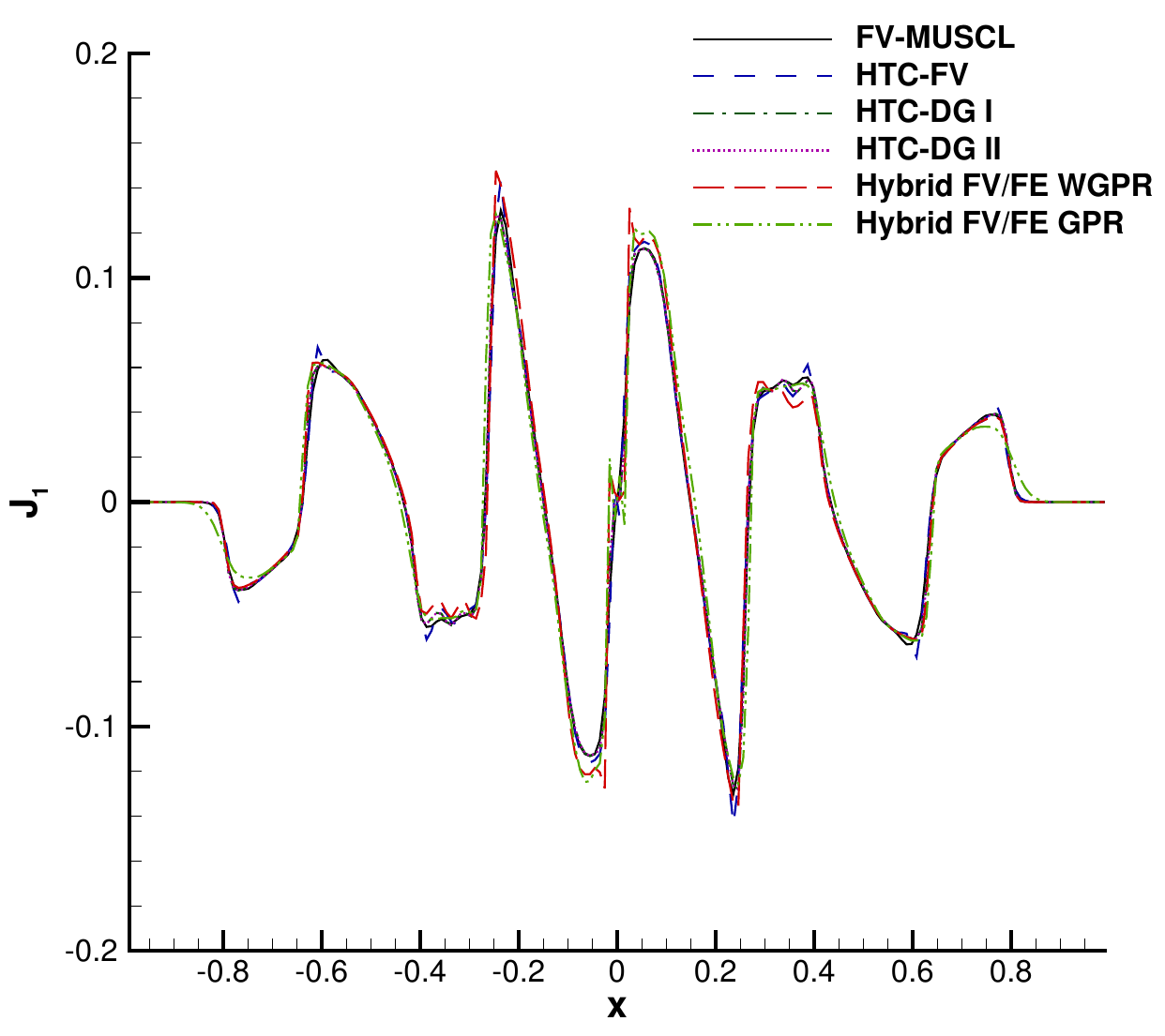}
	\includegraphics[width=0.32\linewidth]{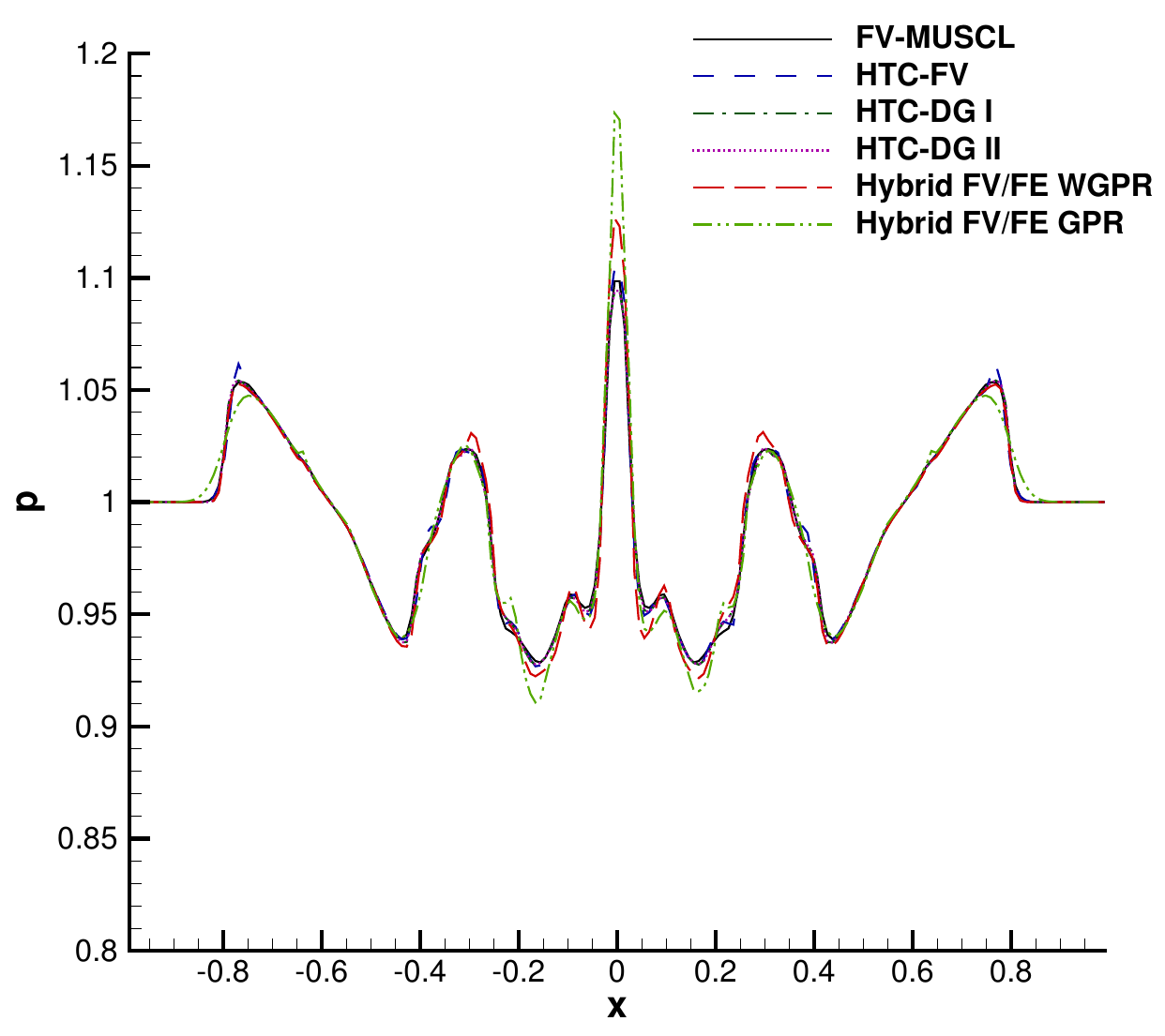}
	\caption{Solid rotor. 1D comparative of the solution along the $x$-axis obtained using the new hybrid methodology for the compressible GPR model, the hybrid approach for weakly compressible medium \cite{HybridGPR}, the FV and DG thermodynamically compatible schemes in \cite{HTCAbgrall} and the asymptotic preserving FV-MUSCL method \cite{Boscheri2021SIGPR}. From left-top to right-bottom: density, $\vel_{1}$, $A_{11}$, $A_{12}$, $J_{1}$ and pressure.}
	\label{fig:SR1D}
\end{figure}

The solid rotor test case is employed to further test the behaviour of the proposed numerical method in the solid limit of the GPR model \cite{Boscheri2021SIGPR,HTCGPR,HTCAbgrall}. The solution in the computational domain $\Omega=[-1,1]^2$ is initialized as
\begin{gather*}
	\rho\left(\mathbf{x},0\right) =\press\left(\mathbf{x},0\right) =  1,\quad \bvel\left(\mathbf{x},0\right)=\left\lbrace\begin{array}{ll}
		\left(\frac{-y}{0.2}, \frac{x}{0.2}, 0\right)^{T}  & \mathrm{if} \; \left\| \mathbf{x} \right\| \leq 0.2,\\
		\mathbf{0} & \mathrm{if} \; \left\| \mathbf{x} \right\| > 0.2.
	\end{array} \right.
\end{gather*}
while the model parameters are $c_{s}=c_{h}=1$, $\tau_{1}=6\cdot 10^{20}$, $\tau_{2}=10^{20}$, $p_0=0$, $c_{v}=1$, $c_{p}=1.4$.
A simulation is run employing the hybrid FV/FE approach with ENO limiting and an implicit discretization of the source terms. Besides, periodic boundary conditions are assumed in all boundaries and the final simulation time is $t_{\mathrm{e}}=0.3$. The obtained solution for representing variables is reported in Figures~\ref{fig:SR}--\ref{fig:SR3D}. To ease comparison with reference solutions in the bibliography, 1D cuts of $\rho$, $\vel_{1}$, $A_{11}$, $A_{12}$, $J_{1}$ and $\press$ fields are reported in Figure~\ref{fig:SR1D}. An excellent agreement is observed for all variables with the numerical approaches proposed in \cite{Boscheri2021SIGPR,HTCAbgrall,HybridGPR}. The solution obtained with the  hybrid methodology has been calculated using a very fine mesh made of $2975744$ primal triangular elements which allows a detailed definition of the steepest wavefronts.

\subsection{Sloshing}
\begin{figure}[H]
	\centering
	\includegraphics[width=0.45\linewidth]{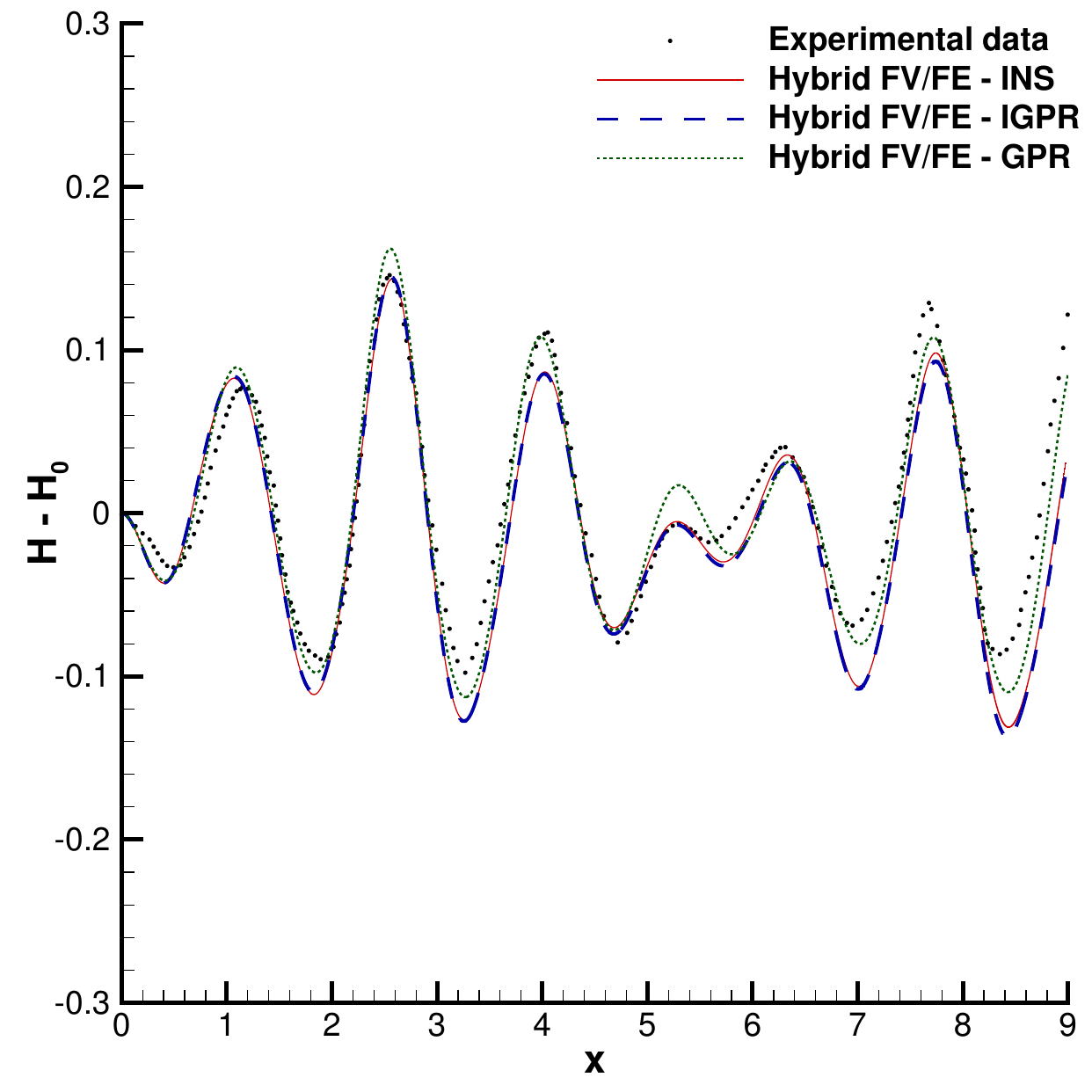}
	\caption{Sloshing. Time series of the mesh surface elevation for the ALE hybrid method solving the compressible GPR model (green dotted line), the incompressible GPR model (blue dashed line) and the incompressible Navier-Stokes equations \cite{HybridALE} (red dashed line) and experimental data from \cite{FaltinsenTimokha2000} (black circles).}
	\label{fig:Sloshing} %\cite{HybridGPR} 
\end{figure}

\begin{figure}[H]
	\centering
	\includegraphics[width=0.5\linewidth]{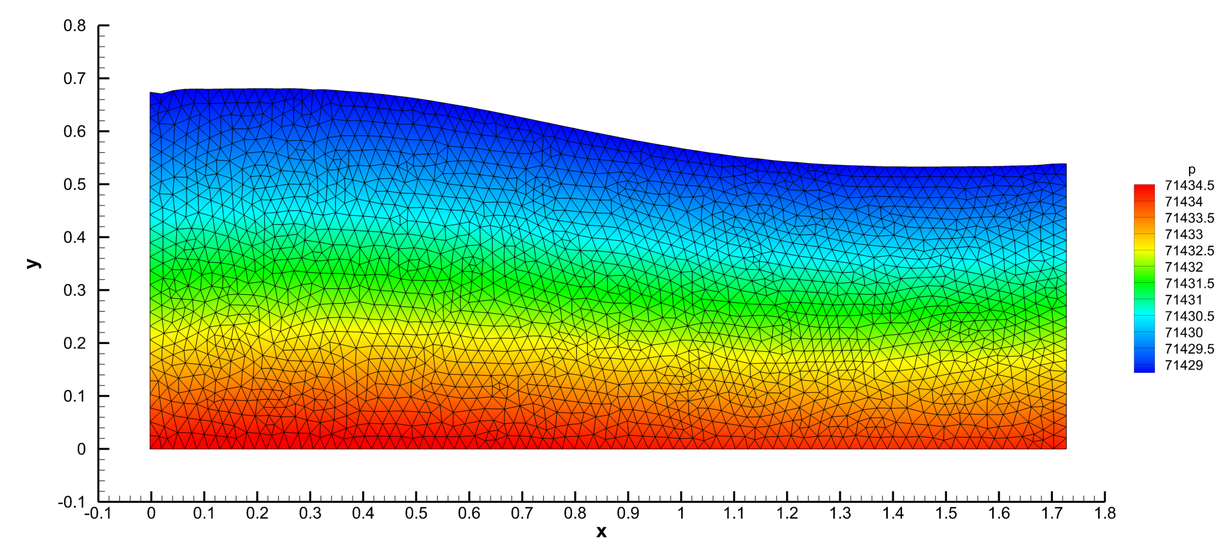}
	\includegraphics[width=0.5\linewidth]{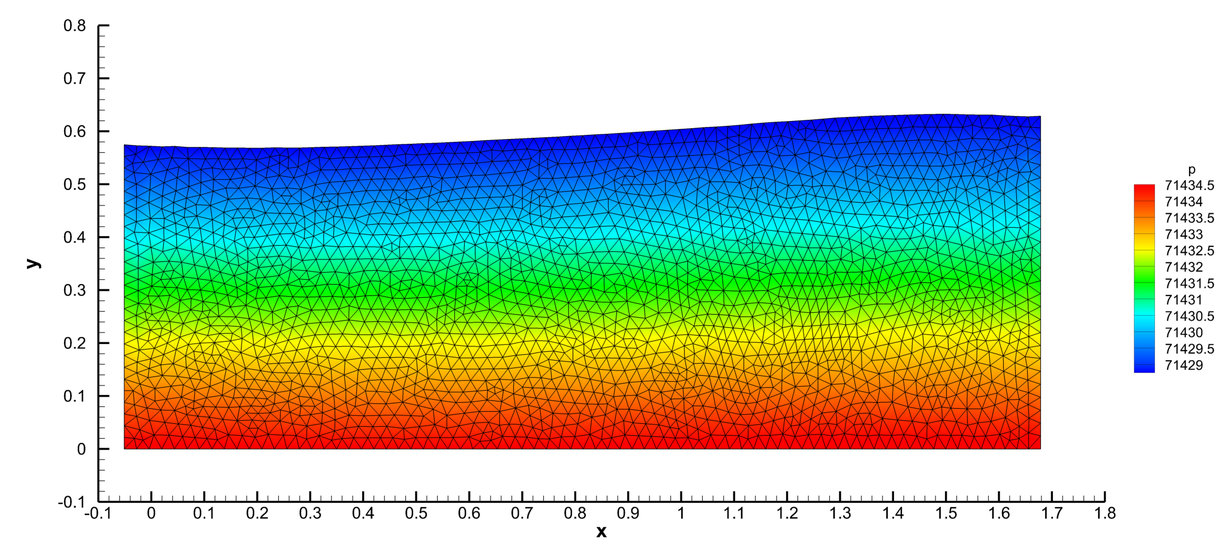}
	\includegraphics[width=0.5\linewidth]{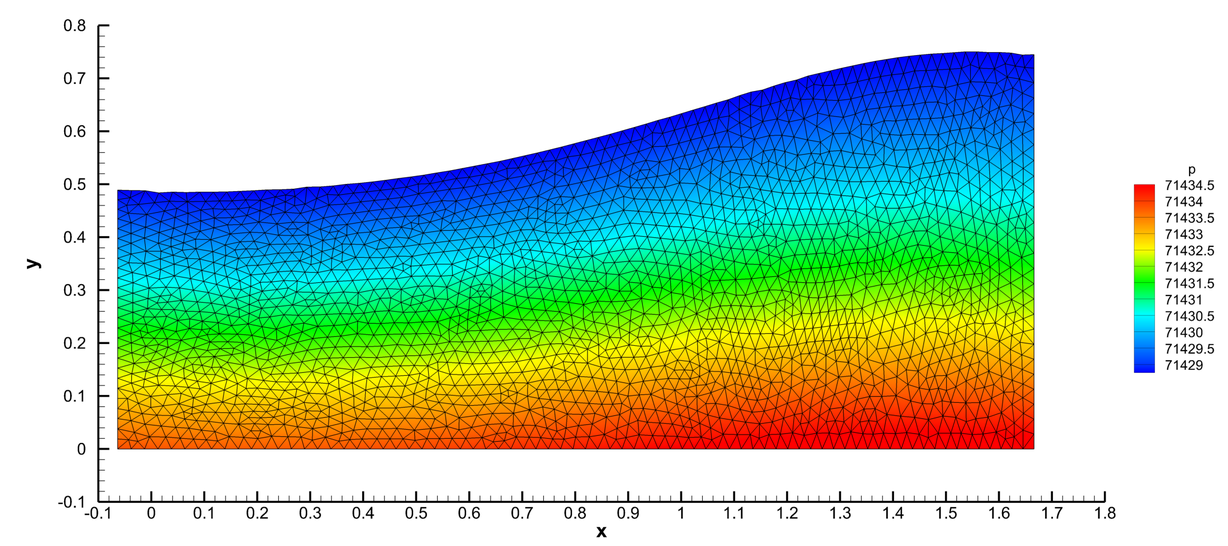}
	\caption{Sloshing. Pressure contour plot and mesh deformation at times $t\in\left\lbrace 1.2, 2.14, 3.29 \right\rbrace$ (from top to bottom). Results obtained using the first order ALE hybrid FV/FE method for the GPR model with artificial viscosity $c_{\alpha}=2$.}
	\label{fig:Sloshing2D}
\end{figure}

The last test case corresponds to a sloshing of an inviscid fluid in a moving tank \cite{FaltinsenTimokha2000,LagrangeNC,HybridALE}. The initial computational domain is $\Omega_{0}=[0,1.73]\times [0,0.6]$. The left, right and bottom boundaries of are assumed to be slip walls moving horizontally following
\begin{equation}
	\vel_1(t)= - \omega A \sin(\omega t), \qquad A=0.032, \qquad \omega = 2 \frac{\pi}{T}, \qquad T=1.3.
\end{equation}
Further, the fluid is left free in the top boundary allowing for the sloshing phenomena to occur. To this end, a pressure boundary condition with $\press_{\mathrm{top}}=0$ is imposed at the top of the domain while inside the domain the fluid moves with a local fluid velocity with speed smoothing $\varsigma=10^3$. 
As initial conditions, we consider a fluid at rest with hydrostatic pressure distribution,
\begin{equation}
		\rho(0,\x) =1,\qquad \bvel(0,\x) =\boldsymbol{0}, \qquad \press(0,\x) =\frac{ 10^{5}}{\gamma} +  g (0.6-y),  \qquad \g = (0,-g),  \qquad g=9.81.
\end{equation}
This benchmark is solved employing the second order semi-implicit ALE hybrid method with ENO limiting and auxiliary numerical viscosity coefficient $c_{\alpha}=2$. Besides, a Lax time dissipation between the staggered grids is applied over the new contributions of the transport stage at each time step. Three simulations are run addressing different models: the incompressible Navier-Stokes equations \cite{HybridALE}, the incompressible GPR model \cite{HybridGPR} and the fully compressible GPR model. All three simulations are run on an unstructured triangular mesh of $4072$ primal elements up to time $t_{\mathrm{e}}=9$. The time series of the surface elevation of a point tracer initially located at $\mathbf{x}=(0.05,0.6)$ is reported in Figure~\ref{fig:Sloshing}. A good agreement is observed among the obtained solutions and with the experimental data from \cite{FaltinsenTimokha2000}. To further illustrate the obtained solutions and ease comparison with further references, as e.g. \cite{LagrangeNC,HybridALE}, several contour plots of the pressure and the mesh deformation computed using the GPR model are depicted in Figure~\ref{fig:Sloshing2D}.

% % % % % % % % % % % % % % % % % % % % % % % % % % % % % %
% % % % % % % % % % % % % % % % % % % % % % % % % % % % % %
%                  Conclusions                            %
% % % % % % % % % % % % % % % % % % % % % % % % % % % % % %
% % % % % % % % % % % % % % % % % % % % % % % % % % % % % %
\section{Conclusions} \label{sec:conclusions}
A novel semi-implicit direct arbitrary Lagrangian-Eulerian methodology has been presented for the numerical solution of the GPR model for continuum mechanics. The proposed hybrid  methodology combines explicit finite volume methods for the transport subsystem with a continuous finite element discretization for the Poisson-type problem associated with the pressure field and for the Laplacian equation that may govern the mesh motion. This decoupling is enabled by an operator splitting approach inspired by previous all Mach number methods, allowing for a robust and efficient treatment of the underlying physical processes.

In addition, we have introduced a thermodynamically compatible augmented GLM GPR model, which has been shown to significantly reduce errors associated with curl involutions. For demanding simulations in the solid limit, employing an exactly curl free scheme may become important. Therefore, a future line of research would tackle the development of exactly involution preserving discretizations profiting form the staggered grid structure \cite{HybridHexa1}. The preservation of additional structural properties, such as thermodynamic compatibility \cite{HTCAbgrall,HTCA2}, will also take part of future research.

The numerical results presented demonstrate the %potential 
capability 
of the proposed methodology to accurately handle both solids and fluids at all Mach numbers. Future efforts will, therefore, focus on incorporating very high-order schemes, on employing alternative equations of state enabling the treatment of medium to large deformations \cite{FGN14_HTChypelast,GNH16_EOS} and on extending the methodology to solve fluid structure iteration problems.

% % % % % % % % % % % % % % % % % % % % % % % % % % % % % %
% % % % % % % % % % % % % % % % % % % % % % % % % % % % % %
%                Acknowledgment                          %
% % % % % % % % % % % % % % % % % % % % % % % % % % % % % %
% % % % % % % % % % % % % % % % % % % % % % % % % % % % % %
\section*{Acknowledgements}
The author acknowledges support from the Spanish Ministry of Science, Innovation and Universities (MCIN), the Spanish AEI (MCIN/AEI/10.13039/501100011033) and European Social Fund Plus under the project RYC2022-036355-I; from FEDER and the Spanish Ministry of Science, Innovation and Universities under project PID2021-122625OB-I00 and from the Xunta de Galicia (Spain) under project GI-1563 ED431C 2025/09. 
The author would like to acknowledge support from the CESGA, Spain, for the access to the FT3 supercomputer. 
Views and opinions expressed are however those of the author only and do not necessarily reflect those of the granting authorities. Neither the European Union nor the granting authorities can be held responsible for them.

%% % % % % % % % % % % % % % % % % % % % % % % % % % % % % %
%% % % % % % % % % % % % % % % % % % % % % % % % % % % % % %
%%              Bibliography
%% % % % % % % % % % % % % % % % % % % % % % % % % % % % % %
%% % % % % % % % % % % % % % % % % % % % % % % % % % % % % %
\bibliographystyle{plain}
\bibliography{./mibiblio}

\end{document}